\documentclass[12pt]{article}
\usepackage{amsmath,latexsym,amssymb}
\begin{document}
\headsep 0.5 true cm

\begin{center}
{\Large\bf Modularity and uniformization of a higher genus
           algebraic space curve, its distinct arithmetical
           realizations by cohomology groups and $E_6$,
           $E_7$, $E_8$-singularities}
\vskip 1.0 cm
{\large\bf Lei Yang}
\end{center}
\vskip 1.5 cm

\begin{center}
{\large\bf Abstract}
\end{center}

\vskip 0.5 cm

  We prove the modularity for an algebraic space curve $Y$ of genus $50$
in $\mathbb{P}^5$, which consists of $21$ quartic polynomials in six
variables, by means of an explicit modular parametrization by theta
constants of order $13$. This provides an example of modularity, explicit
uniformization and hyperbolic uniformization of arithmetic type for a
higher genus (arithmetic) algebraic space curve. In particular, it gives
a new example for Hilbert's 22nd problem. This gives $21$ modular equations
of order $13$, which greatly improve the result of Ramanujan and Evans on
the construction of modular equations of order $13$. We show that $Y$ is
isomorphic to the modular curve $X(13)$. The corresponding ideal $I(Y)$
is invariant under the action of $\text{SL}(2, 13)$, which leads to a
$21$-dimensional reducible representation of $\text{SL}(2, 13)$, whose
decomposition as the direct sum of $1$, $7$ and $13$-dimensional
representations gives two distinct arithmetical realizations of $X(13)$
by character fields $\mathbb{Q}(\chi)=\mathbb{Q}(\zeta_7+\zeta_7^{-1})$
or $\mathbb{Q}(\chi)=\mathbb{Q}(\sqrt{13})$ of irreducible representations
of $\text{SL}(2, 13)$ in two different ways corresponding to the decompositions
of cohomology groups of a projective or affine variety with values in a coherent
algebraic sheaf on $X(13)$ as well as the geometric construction of $Y$, the
geometric realization of the degenerate principal series and the Steinberg
representation of $\text{SL}(2, 13)$. The projection
$Y \rightarrow Y/\text{SL}(2, 13)$ (identified with $\mathbb{CP}^1$)
is a Galois covering whose generic fibre is interpreted as the Galois
resolvent of the modular equation $\Phi_{13}(\cdot, j)=0$ of level $13$,
i.e., the function field of $Y$ is the splitting field of this modular
equation over $\mathbb{C}(j)$. The ring of invariant polynomials
$(\mathbb{C}[z_1, z_2, z_3, z_4, z_5, z_6]/I(Y))^{\text{SL}(2, 13)}$
over $X(13)$ leads to a new perspective on the theory of $E_6$, $E_7$
and $E_8$-singularities.

\vskip 2.0 cm

\begin{center}
{\large\bf Contents}
\end{center}

$$\aligned
  &1. \text{\quad Introduction}\\
  &2. \text{\quad Modularity for equations of $E_6$, $E_7$ and $E_8$-singularities}\\
  &\quad \quad \text{coming from modular curves $X(3)$, $X(4)$ and $X(5)$}\\
  &3. \text{\quad Invariant theory and modular forms for $\text{SL}(2, 13)$ I:}\\
  &\quad \quad \text{the quadratic and cubic invariants}\\
  &4. \text{\quad Invariant theory and modular forms for $\text{SL}(2, 13)$ II:}\\
  &\quad \quad \text{the quartic invariants and the modular curve $X(13)$}\\
  &5. \quad \text{An invariant ideal defining the modular curve $X(13)$}\\
  &6. \quad \text{The decomposition of a $21$-dimensional reducible representation}\\
  &\quad \quad \text{of $\text{SL}(2, 13)$ and its geometric realization}\\
  &\quad 6.1. \quad \text{A $21$-dimensional reducible representation of $\text{SL}(2, 13)$}\\
  &\quad \quad \text{and its decomposition: $\mathbf{21}=\mathbf{1} \oplus \mathbf{7}
          \oplus \mathbf{13}$}\\
  &\quad 6.2. \quad \text{A geometric construction of the curve $Y$ and the ring of}\\
  &\quad \quad \text{invariant polynomials}\\
  &\quad 6.3. \quad \text{A non-standard geometric realization of the degenerate}\\
  &\quad \quad \text{principal series for $\text{SL}(2, 13)$ and a geometric realization of}\\
  &\quad \quad \text{the Steinberg representation for $\text{SL}(2, 13)$}\\
  &7. \quad \text{Modularity for an invariant ideal and an explicit uniformization}\\
  &\quad \quad \text{of algebraic space curves of higher genus}\\
  &8. \quad \text{Galois covering $Y \rightarrow Y/\text{SL}(2, 13) \cong \mathbb{CP}^1$,
                  Galois resolvent for the}\\
  &\quad \quad \text{modular equation of level $13$ and a Hauptmodul for $\Gamma_0(13)$}\\
  &9. \text{\quad Invariant theory and modular forms for $\text{SL}(2, 13)$ III: some}\\
  &\quad \quad \text{computation for invariant polynomials}\\
  &10. \text{\quad Modularity for equations of $E_6$, $E_7$ and
                  $E_8$-singularities coming}\\
  &\quad \quad \text{from $C_Y/\text{SL}(2, 13)$ and variations of $E_6$, $E_7$
                     and $E_8$-singularity}\\
  &\quad \quad \text{structures over $X(13)$}\\
  &11. \text{\quad Modularity for equations of $Q_{18}$ and $E_{20}$-singularities
                   coming}\\
  &\quad \quad \text{from $C_Y/\text{SL}(2, 13)$ and variations of $Q_{18}$ and
                     $E_{20}$-singularity}\\
  &\quad \quad \text{structures over $X(13)$}
\endaligned$$

\vskip 2.0 cm

\begin{center}
{\large\bf 1. Introduction}
\end{center}

  It is well-known that modularity provides a deep relation between
automorphic forms and algebraic varieties defined by equations. Such
modular varieties have a very rich structure that can be used to study
their arithmetic properties. For cubic equations, we have the rich theory
of modular curves of higher level and the associated modular parameterizations
of elliptic curves. In particular, there is the modularity theorem for
elliptic curves over $\mathbb{Q}$ (formerly known as Taniyama-Shimura-Weil
conjecture) (see \cite{W}, \cite{TW}, \cite{BCDT} and \cite{AC}) and its
analogue for curves of genus two (see \cite{BCGP} and \cite{BPTVY}).
Despite of this, very little is known for curves of higher genus, especially
for the explicit construction of modular parameterizations. In fact, for a
generic curve of genus $g$ (or abelian variety of dimension $g$), we expect
an automorphic representation for the (split) orthogonal group $\text{SO}_{2g+1}$,
the group whose Langlands dual is $\text{Sp}_{2g}$. This is isomorphic to
$\text{Sp}_{2g}$ if $g \leq 2$ but not other wise. More precisely, when
$g=1$ or $g=2$, there are well-known exceptional isomorphisms which allow
us to replace $\text{SO}_{2g+1}$ by the groups $\text{GL}_2$ and $\text{Sp}_4$
respectively. One may well ask whether the method of the above papers could be
used to prove (potential) modularity of curves of genus $g \geq 3$. At the moment,
this seems exceedingly unlikely without some substantial new idea (see \cite{BCGP},
page 14). In fact, the automorphic representations contributing to the coherent
cohomology of orthogonal Shimura varieties are representations of the inner form
$\text{SO}(2g-1, 2)$ of $\text{SO}_{2g+1}$ which is non-split if $g>1$, whose
infinity components $\pi_{\infty}$ are either discrete series or non-degenerate
limits of discrete series. If $g=1$, the representations are discrete series, and
if $g=2$, they are non-degenerate limits of discrete series, but if $g \geq 3$,
then neither possibility occurs, so the automorphic representations do not
contribute to the cohomology of the corresponding Shimura variety (see \cite{BCGP}).
In particular, the general modularity problem for curves of genus $g \geq 3$ seems
at least as hard as proving non-solvable cases of the Artin conjecture for totally
even representations, and even proving the modularity of a single such curve seems
completely out of reach (see also \cite{BCGP}).

  We resolve this problem in Theorem 1.12 and Corollary 1.13 below by proving
the modularity for an algebraic space curve $Y$ of genus $50$ and degree $35$
in $\mathbb{P}^5$ as well as two curves $Y_2$ and $Y_3$ lying over $Y$ which
have the modular components, i.e., some components are modular. In particular,
we give the explicit modular parameterizations for the curve $Y$ and some
components of $Y_2$ and $Y_3$. Now, let us give some background. For an algebraic
curve or a Riemann surface $X$ of genus $g$, there are two objects: one is its
fundamental group $\pi_1(X)$, the other is its homology group $H_1(X) \cong \mathbb{Z}^{2g}$
(for simplicity, we do not consider the torsion group). The first one is non-abelian,
the second one is abelian (it is the abelianization of $\pi_1(X)$). For the abelian
object $H_1(X)$, one can associate it with the following structures:

(1) Jacobian variety $\text{Jac}(X)$ of $X$;

(2) symplectic groups: $\text{Sp}(2g, \mathbb{Z})$, $\text{Sp}(2g, \mathbb{F}_q)$
    and symplectic congruence subgroups $\text{Sp}(2g, \mathbb{Z})(N)$;

(3) Siegel modular forms of genus $g$.

  In particular, we have the global Langlands correspondence for discrete
symplectic motives of rank $2g$ over $\mathbb{Q}$. To such a motive Langlands
conjecturally associates a generic automorphic representation $\pi$ of the
split orthogonal group $\text{SO}_{2g+1}$ over $\mathbb{Q}$, which appears
with multiplicity one in the cuspidal spectrum (see \cite{Gro}). More
precisely, let $M$ be a pure motive of weight $-1$ and rank $2g$ over
$\mathbb{Q}$ with a non-degenerate symplectic polarization
$$\psi: \Lambda^2 M \rightarrow \mathbb{Q}(1).$$
The simplest examples come from complete, non-singular curves $X$ of genus
$g$ over $\mathbb{Q}$, where $M=H^1(X)(1)=H_1(X)$. Let $(A, \psi)$ be
the corresponding Jacobian variety (a polarized abelian variety) of
dimension $g$ over $\mathbb{Q}$ and let $(H_1(A), \psi)$ be the associated
polarized symplectic motive $M$ of rank $2g$.

  The corresponding real Langlands parameter $W_{\mathbb{R}} \rightarrow
\text{Sp}_{2g}(\mathbb{C})$ is given by $g$ copies of the $2$-dimensional
representation $\text{Ind}(z/\overline{z})^{1/2}$. The centralizer of the
image in $\text{Sp}_{2g}(\mathbb{C})$ is isomorphic to the orthogonal group
$\text{O}_g(\mathbb{C})$ and has component group of order $2$. Hence there
are two representations in the corresponding Vogan $L$-packet of limit
discrete series. The first is the genetic limit discrete series
$\pi_{\infty}$ for $\text{SO}(g+1, g)$, which corresponds to the trivial
character of the component group. The second representation in the Vogan
$L$-packet is also a limit discrete series, although it is not generic.
It gives a representation of the split group $\text{SO}(g+1, g)$ when
$g=2m$ is even. When $g=2m+1$ is odd the element $-1$ in the center of
$\text{Sp}_{2g}(\mathbb{C})$ lies in the non-trivial coset of the
connected component of the centralizer and the non-generic limit discrete
series is a representation of the non-split group $\text{SO}(g-1, g+2)$
(see \cite{Gro}).

  In fact, when the genus of $X$ increases, $\pi_1(X)$ is far away from
abelian. Hence, $H_1(X)$ is far away from $\pi_1(X)$. Since the fundamental
group $\pi_1(X)$ captures the non-abelian aspect of $X$, while the homology
group $H_1(X)$ can only give the abelian aspect of $X$. It is better to study
$\pi_1(X)$ instead of $H_1(X)$, i.e., study $X$ instead of its Jacobian
variety $\text{Jac}(X)$ and the symplectic motives. However, for the non-abelian
object $\pi_1(X)$, very little structure can be associated with it except for
$X$ itself. By the celebrated uniformization theorem, any compact complex
Riemann surface $X$ of genus $g \geq 2$ can be represented by
$$X=\mathbb{H}/\pi_1(X),$$
where $\pi_1(X)$ is the fundamental group of $M$. In particular, the modular
curves $X(N)$ of level $N$ is given by
$$X(N)=\overline{\mathbb{H}}/\Gamma(N)=\mathbb{H}/\pi_1(X(N)),$$
where
$$\Gamma(N)=\left\{ \left(\begin{matrix} a & b\\
            c & d \end{matrix}\right) \in \Gamma=\text{SL}(2, \mathbb{Z}):
            \left(\begin{matrix} a & b\\
            c & d \end{matrix}\right) \equiv \left(\begin{matrix}
            1 & 0\\ 0 & 1 \end{matrix}\right) (\text{mod $N$}) \right\}$$
is the principal congruence subgroup of level $N$. Hence,
$\Gamma(N)$ is closely related to the fundamental group $\pi_1(X(N))$ and
give the geometry and topology of $X(N)$ due to Riemann and Klein. This
leads to the uniformization theorem which was at the centre of the evolution
of mathematics in the 19th century (see \cite{dSG}). In the diversity of its
algebraic, geometric, analytic, topological, and even number-theoretic aspects,
the uniformization theorem is in some sense symbolic of the mathematics of
that century. It took a whole century to get to the point of stating this
theorem and providing a convincing proof of it, relying as it did on prior
work of Gauss, Riemann, Schwarz, Klein, Poincar\'{e}, and Koebe, among others.
By 1882 Klein and Poincar\'{e} had become fully convinced of the truth of the
following uniformization theorem (see \cite{dSG}).

\textbf{Theorem 1.1.} {\it Let $X$ be any compact Riemann surface of genus
at least two. There exists a discrete subgroup $\Gamma$ of $\text{PSL}(2, \mathbb{R})$
acting freely and properly on $\mathbb{H}$ such that $X$ is isomorphic to the
quotient $\mathbb{H}/\Gamma$. In other words, the universal cover of $X$ is
holomorphically isomorphic to $\mathbb{H}$.}

  In summary, Klein and Poincar\'{e} had effectively solved one of the
main problem handed down by the founders of algebraic geometry, that is,
to parametrize an algebraic curve
$$F(x, y)=0\eqno{(1.1)}$$
of genus at least two by single-valued meromorphic functions $x, y:
\mathbb{H} \rightarrow \mathbb{C}$.

  The uniformization of the plane algebraic curves (see \cite{Br}), defined
by irreducible polynomials (1.1) over $\mathbb{C}$ is a union of analytic
theory and theory of infinite discontinuous groups of linear transformations
of the form
$$\tau \mapsto \frac{a \tau+b}{c \tau+d}.$$
Here, the algebraic and geometric aspects include theory of Fuchsian and
Kleinian groups, their subgroups, hyperbolic geometry, and fundamental
polygons. Fuchsian groups of higher genera have received less attention for
the reason of computation complexity. However, pure algebraic and geometric
approaches do not touch the question of explicit parametrization of the
implicit relation $F(x, y)=0$. This point is the central problem in the
uniformization and has essentially analytic characterization. In particular,
in his pioneering work \cite{We}, Hermann Weyl stated that ``To crown the
whole development of this part of Riemannian function theory, the solution
of the following problem would be desired. Given a Riemann surface in its
normal form (that is, the associated group of motions of the non-Euclidean
plane), to express each of the uniform or many-valued un-branched functions
on the surface by means of a closed analytic formula in terms of the
coordinate $t$ of the points in the Lobachevsky plane.'' More precisely,
the analytical description implies the presence of two analytic functions,
single-valued in a domain of their existence, of the global uniformizing
parameter $\tau$:
$$x=x(\tau), \quad y=y(\tau).$$

  First examples of explicit parametrizations of higher genera curves were
obtained by Jacobi in his celebrated Fundamenta Nova (see \cite{Ja}). After
Jacobi, his followers (Schl\"{a}fli, Sohnke et all) expanded this list and,
nowadays, it is known as classical modular equations. In particular, Jacobi
established modular equations of degrees $3$ and $5$ in his Fundamenta Nova
(see \cite{Ja}). Modern achievements in the explicit construction of the
uniformizing map are also related to particular examples of modular curves.
There are two key points: one is in the spirit of Fricke and Klein, to consider
subgroups of the full modular group $\text{SL}(2, \mathbb{Z})$, the other is
Jacobi's theta constants and representations in terms of theta constants. In
particular, the parametrization of the Klein quartic curve represents an
important stage in the explicit uniformization of algebraic plane curves along
these lines.

  However, algebraic plane curves do not exhaust all of the algebraic curves.
When the genus $g \geq 4$, the generic curves are algebraic space curves. In
fact, according to Weyl's comments mentioned as above, given a Riemann surface
in its normal form $\Gamma \backslash \mathbb{H}$, one will obtain, in general,
an algebraic space curve, not necessarily an algebraic plane curve. Therefore,
the uniformization of algebraic space curves in $n$-dimensional affine or
projective space presents a major challenge that we do not need to deal with
for algebraic plane curves (1.1) which we mention as above. Instead of a single
equation (1.1), a system of $n-1$ or more equations is needed. This system
is far from unique and, in many cases, may be over-determined. General
algebraic space curves is a topic with various unresolved issues of mathematical
and computational interest and an area of important future research. In particular,
representing algebraic space curves is a fundamental problem of algebraic geometry.
These curves are defined as the intersection curves of algebraic hyper-surfaces.
Representation issues include problems such as finding the minimum number of
equations needed to define an algebraic space curve in affine and projective
spaces as well as the degrees of these defining equations. However, irreducible
space curves in general, defined by more than two hyper-surfaces are difficult
to handle equationally and one needs to resort to ideal-theoretic methods.
Therefore, finding a small set of generators for a polynomial ideal defining
an algebraic space curve and constructing an explicit uniformization is a
challenge problem in algebraic geometry.

  In fact, the uniformization of algebraic space curves by automorphic functions
also appears in the celebrated Hilbert's problems, i.e., uniformization of
analytic relations by means of automorphic functions. Hilbert posed this question
as the twenty-second in his list of $23$ problems. The following quotes are excerpts
from a translation of Hilbert's ICM article \cite{Hil}.

  ``As Poincar\'{e} was the first to prove, it is always possible to reduce any
algebraic relation between two variables to uniformity by the use of automorphic
functions of one variable. That is, if any algebraic equation in two variables
be given, there can always be found for these variables two such single valued
automorphic functions of a single variable that their substitution renders the
given algebraic equation an identity.''

  ``In conjunction with this problem comes up the problem of reducing to uniformity
an algebraic or any other analytic relation among three or more complex variables--
a problem which is known to be solvable in many particular cases. Toward the
solution of this the recent investigations of Picard on algebraic functions of two
variables are to be regarded as welcome and important preliminary studies.''

  In summary, Hilbert's 22nd problem asks whether every algebraic or analytic
curve, i.e., solutions to polynomial equations, can be written in terms of
single-valued functions. The problem has been resolved in the one-dimensional
case, i.e., a single polynomial equation in two variables given by (1.1), and
continues to be studied in other cases (see \cite{Ber1} and \cite{Ber2}). More
precisely, Poincar\'{e} studied $F(x, y)=0$ with $F$ algebraic or analytic,
while Picard studied $F(x, y, z)=0$ with $F$ algebraic, by means of automorphic
functions.

  Despite of this, very little is known for reducing to uniformity by means
of automorphic functions on algebraic functions of three or more variables.
In particular, one does not know whether solutions to a system of polynomial
equations can be written uniformly (not locally) in terms of single-valued
automorphic functions, which plays a central role in Hilbert's 22nd problem.
The main difficulty is that the uniformization for a system of polynomial
equations (e.g. algebraic space curves) leads to the simultaneous
uniformization for several algebraic relations among three or more variables
by means of automorphic functions (e.g. the space curves are represented as
the intersection of higher-dimensional algebraic varieties). Thus, we have
to deal with the uniformization for algebraic functions of three or more
variables by means of automorphic functions. In the present paper, we give
an explicit construction of uniformity for an over-determined system of
algebraic relations among five variables by theta constants (automorphic
functions), which leads to a new solution for Hilbert's 22nd problem.

  Moreover, the uniformization theorem is closely related to the modularity by
an equivalent statement of the Belyi's theorem (see \cite{Groth}):

\textbf{Theorem 1.2.} {\it For every algebraic curve defined over
$\overline{\mathbb{Q}}$, there is a finite index subgroup $\Gamma^{\prime}$
of $\Gamma=\text{SL}(2, \mathbb{Z})$, such that
$$\overline{\mathbb{H}}/\Gamma^{\prime} \cong X(\mathbb{C}),$$
i.e., the corresponding complex algebraic curve is uniformized by a finite
index subgroup of $\Gamma=\text{SL}(2, \mathbb{Z})$.}

  Despite that the majority of these groups $\Gamma^{\prime}$ are non-congruence,
the congruence subgroups give rise to the key.

{\it Definition 1.3.} (see \cite{Ma}). Let $C$ be an algebraic curve. A hyperbolic
uniformization (of $C$) of arithmetic type is a hyperbolic uniformization of
the algebraic curve $C$ which is periodic with respect to a congruence
subgroup $\Gamma^{\prime} \subset \Gamma=\text{SL}(2, \mathbb{Z})$.

  In particular, the Shimura-Taniyama-Weil conjecture asserts that any
arithmetic elliptic curve (i.e., any elliptic curve whose defining equation
can be taken with coefficients in $\mathbb{Q}$) admits a hyperbolic
uniformization of arithmetic type. In general, let $C$ be an arithmetic
algebraic curve (i.e., an algebraic curve whose defining ideal or defining
equations can be taken with coefficients in $\mathbb{Q}$). A hyperbolic
uniformization of arithmetic type will put an exceedingly rich geometric
structure on this arithmetic algebraic curve $C$ and carry deep implications
for arithmetic questions. Thus, in order to generalize the Shimura-Taniyama-Weil
conjecture for elliptic curves to curves with higher genus ($g \geq 4$), a
challenge problem is the following:

\textbf{Problem 1.4.} Find an arithmetic algebraic space curve $C$ with a
hyperbolic uniformization of arithmetic type, where the genus $g(C) \geq 4$.
In particular, for such an arithmetic algebraic curve $C$, find the level
$N$ such that $\Gamma(N)$ is involved in a hyperbolic uniformization of
arithmetic type for $C$.

  In the present paper, we prove the modularity for an algebraic space
curve $Y$ of genus $50$ and degree $35$ in $\mathbb{P}^5$ and construct
the associated modular parameterizations. We show that this curve $Y$ is
isomorphic to the modular curve $X(13)$. We can associate $\Gamma(13)$
and hence $\pi_1(X(13))$ with an invariant ideal $I(Y)$ which is invariant
under the action of $\text{SL}(2, 13)$ and establish its modularity by
theta constants of order $13$. In particular, given an algebraic curve
or a compact Riemann surface as simple as the modular curve $X(13)$,
there are two distinct arithmetical realizations of $X(13)$ by character
fields $\mathbb{Q}(\chi)=\mathbb{Q}(\zeta_7+\zeta_7^{-1})$ or
$\mathbb{Q}(\chi)=\mathbb{Q}(\sqrt{13})$ (i.e. two distinct algebraic
number fields) of irreducible representations of $\text{SL}(2, 13)$ in two
different ways corresponding to the decompositions of cohomology groups of
a projective or affine variety with values in a coherent algebraic sheaf
on $X(13)$. In fact, after we have proved the modularity of $Y$, obtained
an explicit uniformization of $Y$ as well as the hyperbolic uniformization
of arithmetic type for this higher genus arithmetic algebraic curve $Y$, a
major problem as well as a major application is to study its arithmetic.
That is, there are two distinct arithmetical realizations by cohomology
groups whose fields of definition are distinct algebraic number fields.
Hence, we have the following picture:
$$\begin{array}{rcl}
  \text{algebraic curves} & \underleftrightarrow{\text{modularity/uniformization}}
                          & \text{modular forms}\\
  \searrow &              & \swarrow\\
           & \text{modular curves} &\\
           & \downarrow &\\
           & \text{arithmetical realizations} &\\
           & \text{(fields of definition/character fields)} &
\end{array}$$
More precisely, for modular curves $X(p)$, we have the following two
different objects associated with it:
$$X(p): \left\{\aligned
  &H^0(X(p), \Omega_{X(p)}^1)\\
  &\mathcal{L}(X(p))
  \endaligned\right. \quad
  \text{and their decompositions into irreducible parts}.$$
The first one $H^0(X(p), \Omega_{X(p)}^1)$ coming from geometry, which is
a linear object, gives a finite dimensional (reducible) representations of
the finite group $\text{PSL}(2, \mathbb{F}_p)$ whose dimension is equal to
the genus of $X(p)$, which was studied by Hecke in his remarkable papers
\cite{Hec1}, \cite{Hec2}, \cite{Hec3}, \cite{Hec4}, and \cite{Hec5}, where
he combined geometry, topology and complex analysis, i.e, Riemann-Weierstrass
function theory with Frobenius's character theory (see \cite{Fr1}) to give
its decomposition into the direct sum of irreducible representations. The
second one $\mathcal{L}(X(p))$ coming from algebra is the locus of $X(p)$
in projective spaces, i.e., a highly nonlinear object, which is given by
Klein's $z$-model of the modular curve $X(p)$ in $\mathbb{CP}^{\frac{p-3}{2}}$
(see \cite{K4}).

  In particular, when $p=13$, both $H^0(X(13), \Omega_{X(13)}^1)$ and
$\mathcal{L}(X(13))$, i.e., a linear object and a highly nonlinear object
give rise to finite dimensional reducible representations of $\text{SL}(2, 13)$
or $\text{PSL}(2, 13)$ with dimensions $50$ and $21$, respectively. The first
one is the genus of $X(13)$, i.e., the number of linearly independent integrals
of the first kind. The second one is the dimension of the invariant ideal,
i.e., the number of quartic equations which define $\mathcal{L}(X(13))$,
where representation theory, algebra and arithmetic are put together in
the study of its decomposition into irreducible parts. In fact, both
$H^0(X(13), \Omega_{X(13)}^1)$ and $\mathcal{L}(X(13))$ are unified in the
context of cohomology groups of a projective and affine variety with values
in a coherent algebraic sheaf on $X(13)$, respectively. This leads to the
following correspondence:
$$\left\{\aligned
  \text{The homology group (a motive)}: H_1(X(13)) &\longleftrightarrow
                             H^0(X(13), \Omega_{X(13)}^1),\\
  \text{an abelian object} &\longleftrightarrow \text{a linear object}.\\
  \text{The fundamental group}: \pi_1(X(13)) \leftrightarrow \Gamma(13)
  &\longleftrightarrow \mathcal{L}(X(13)),\\
  \text{a non-abelian object} &\longleftrightarrow \text{a nonlinear object}.
\endaligned\right.$$

  According to the above correspondence, we can associate $\Gamma(13)$ with
the following thirteen kinds of structures:

(1) defining invariant ideal $I(Y)$ for $X(13)$ in $\mathbb{P}^5$, i.e.,
    the locus $\mathcal{L}(X(13))$, which is generated by $21$ quartic
    polynomials in six variables and the corresponding invariant theory;

(2) invariant curve $Y$ of genus $50$ and degree $35$ in $\mathbb{P}^5$ as
    well as two invariant curves $Y_2$, $Y_3$ lying over $Y$;

(3) the ring of invariant polynomials in six variables
$$\left(\mathbb{C}[z_1, z_2, z_3, z_4, z_5, z_6]/I(Y)\right)^{\text{SL}(2, 13)};$$

(4) invariant cones $C_Y/\text{SL}(2, 13)$, $C_{Y_2}/\text{SL}(2, 13)$
    and $C_{Y_3}/\text{SL}(2, 13)$ over the curve $Y$;

(5) invariant quartic Fano four-fold $Y_1$, invariant sextic Calabi-Yau
    four-fold $W_2$ and invariant octic general type four-fold $W_3$;

(6) a $21$-dimensional reducible representation of $\text{SL}(2, 13)$,
    its decomposition as the direct sum of $1$, $7$ and $13$-dimensional
    representations, which leads to two distinct arithmetical realizations
    of $X(13)$ by character fields $\mathbb{Q}(\chi)=\mathbb{Q}(\zeta_7+\zeta_7^{-1})$
    or $\mathbb{Q}(\chi)=\mathbb{Q}(\sqrt{13})$ (i.e. two distinct algebraic
    number fields) of irreducible representations of $\text{SL}(2, 13)$ in
    two different ways corresponding to the decompositions of cohomology
    groups of a projective or affine variety with values in a coherent
    algebraic sheaf on $X(13)$ and the corresponding geometric construction
    of $Y$;

(7) geometric realization of the degenerate principal series representation
    and the Steinberg representation for $\text{SL}(2, 13)$;

(8) modularity of $Y$ as well as some components of $Y_2$ and $Y_3$ in
    $\mathbb{P}^5$ by theta constants of order $13$;

(9) an explicit uniformization of algebraic space curves of higher genus;

(10) a hyperbolic uniformization of arithmetic type for a higher genus
     arithmetic space curve;

(11) Galois covering $Y \rightarrow Y/\text{SL}(2, 13) \cong \mathbb{P}^1$;

(12) Galois resolvent for the modular equation of order $13$ and a Hauptmodul
    for $\Gamma_0(13)$;

(13) $E_6$, $E_7$, $E_8$-singularities and their variation structures.

  From the viewpoint of automorphic representations, for such a curve
of genus $g=50$, one should study the following objects:

(1) $H_1(X(13)) \cong \mathbb{Z}^{100}$ (an abelian object);

(2) the symplectic motives;

(3) orthogonal Shimura varieties for $\text{SO}(99, 2)$

(4) generic and non-generic limit discrete series with multiplicity
    one for $\text{SO}(51, 50)$, whose rank is $50$;

(5) modularity of orthogonal Galois representations:
$$\rho: \text{Gal}(\overline{\mathbb{Q}}/\mathbb{Q})
  \rightarrow \text{GSO}_{101}(\overline{\mathbb{Q}_{\ell}}).$$

(6) symplectic Galois representations:
$$\rho: \text{Gal}(\overline{\mathbb{Q}}/\mathbb{Q})
  \rightarrow \text{Sp}_{100}(\overline{\mathbb{Q}_{\ell}}).$$
Here, the symplectic Galois representations coming from the \'{e}tale
cohomology group $H^1$ of a curve have weights $0$, $1$ each occurring
with multiplicity $50$.

  Comparing with the above approach of automorphic representations, we
study the following objects:

(1) $\pi_1(X(13)) \leftrightarrow \mathcal{L}(X(13))$, that is, a non-abelian
    object corresponds to a highly nonlinear object;

(2) an invariant ideal $I(Y)$ and the ring of invariant polynomials;

(3) a $21$-dimensional reducible representation of $\text{SL}(2, 13)$,
    its decomposition as the direct sum of $1$, $7$ and $13$-dimensional
    representations, which leads to two distinct arithmetical realizations
    of $X(13)$ by character fields $\mathbb{Q}(\chi)=\mathbb{Q}(\zeta_7+\zeta_7^{-1})$
    or $\mathbb{Q}(\chi)=\mathbb{Q}(\sqrt{13})$ (i.e. two distinct algebraic
    number fields) of irreducible representations of $\text{SL}(2, 13)$ in
    two different ways corresponding to the decompositions of cohomology
    groups of a projective or affine variety with values in a coherent
    algebraic sheaf on $X(13)$ and the corresponding geometric construction
    of $Y$;

(4) modular parametrization of $Y$ and some components of $Y_2$ and $Y_3$
    by theta constants of order $13$;

(5) explicit uniformization of $Y$;

(6) hyperbolic uniformization of arithmetic type for $Y$;

(7) Galois covering $Y \rightarrow Y/\text{SL}(2, 13) \cong \mathbb{P}^1$.

  Arithmetic properties of Shimura varieties are an exciting topic which
has roots in classical topics of algebraic geometry and of number theory
such as modular curves and modular forms. In particular, given Shimura
varieties, a fundamental problem is to determine their defining ideals or
defining equations as algebraic varieties and establish their modularity.
However, even for the special case of modular curves, this is far away
from being solved. In this paper, at least for the modular curve $X(13)$,
we can solve this problem and obtain the following pictures for $X(13)$:

(1) Riemann surfaces $X(13)$ (geometry and topology);

(2) curves given by a system of quartic equations (algebraic curves);

(3) representing algebraic space curves with the minimum number of defining
    equations and finding a small set of generators for a polynomial defining
    ideal (algebraic geometry);

(4) an invariant quartic Fano four-fold, an invariant sextic Calabi-Yau four-fold
    and an invariant octic general type four-fold (higher-dimensional algebraic geometry);

(5) an invariant ideal and the ring of invariant polynomials
   (commutative algebra and invariant theory);

(6) three invariant cones $C_Y$, $C_{Y_2}$ and $C_{Y_3}$ over $Y$
    (algebraic geometry);

(7) a $21$-dimensional reducible representation of $\text{SL}(2, 13)$, its
    decomposition as the direct sum of $1$, $7$ and $13$-dimensional
    representations, which leads to two distinct arithmetical realizations
    of $X(13)$ by character fields $\mathbb{Q}(\chi)=\mathbb{Q}(\zeta_7+\zeta_7^{-1})$
    or $\mathbb{Q}(\chi)=\mathbb{Q}(\sqrt{13})$ (i.e. two distinct algebraic
    number fields) of irreducible representations of $\text{SL}(2, 13)$ in
    two different ways corresponding to the decompositions of cohomology
    groups of a projective or affine variety with values in a coherent
    algebraic sheaf on $X(13)$ and the corresponding geometric construction
    of $Y$ (representation theory and arithmetic geometry);

(8) geometric realization of the degenerate principal series representation
    and the Steinberg representation for $\text{SL}(2, 13)$ (geometric
    representation theory);

(9) modular parametrization of the curves $Y$ and some components of $Y_2$
    and $Y_3$ (number theory);

(10) an explicit uniformization of algebraic space curves of higher genus
   (complex analysis);

(11) a hyperbolic uniformization of arithmetic type for a higher genus
     arithmetic space curve (arithmetic geometry);

(12) Galois covering $Y \rightarrow Y/\text{SL}(2, 13) \cong \mathbb{P}^1$
   (arithmetic geometry);

(13) Galois resolvent for the modular equation of order $13$ and a Hauptmodul
     for $\Gamma_0(13)$ (arithmetic geometry).

\noindent Here, arithmetic, algebraic, geometric and analytic facets are
tightly imbricated in the study of the modular curve $X(13)$ and its
projective model $Y$ in $\mathbb{P}^5$.

\textbf{Main results.}

  Let us begin with the invariant theory for $\text{SL}(2, 13)$, which
we developed in \cite{Y1}, \cite{Y2}, \cite{Y3}, \cite{Y4} and \cite{Y5}.
The representation of $\text{SL}(2, 13)$ we will consider is the unique
six-dimensional irreducible complex representation for which the
eigenvalues of $\left(\begin{matrix} 1 & 1\\ 0 & 1 \end{matrix}\right)$
are the $\exp (\underline{a} . 2 \pi i/13)$ for $\underline{a}$ a
non-square mod $13$. We will give an explicit realization of this
representation. This explicit realization will play a major role for
giving a complete system of invariants associated to $\text{SL}(2, 13)$.
Recall that the six-dimensional representation of the finite group
$\text{SL}(2, 13)$ of order $2184$, which acts on the five-dimensional
projective space
$$\mathbb{CP}^5=\{ (z_1, z_2, z_3, z_4, z_5, z_6): z_i \in \mathbb{C}
  \quad (i=1, 2, 3, 4, 5, 6) \}.$$
This representation is defined over the cyclotomic field
$\mathbb{Q}(e^{\frac{2 \pi i}{13}})$. Put
$$S=-\frac{1}{\sqrt{13}} \begin{pmatrix}
  \zeta^{12}-\zeta & \zeta^{10}-\zeta^3 & \zeta^4-\zeta^9
& \zeta^5-\zeta^8 & \zeta^2-\zeta^{11} & \zeta^6-\zeta^7\\
  \zeta^{10}-\zeta^3 & \zeta^4-\zeta^9 & \zeta^{12}-\zeta
& \zeta^2-\zeta^{11} & \zeta^6-\zeta^7 & \zeta^5-\zeta^8\\
  \zeta^4-\zeta^9 & \zeta^{12}-\zeta & \zeta^{10}-\zeta^3
& \zeta^6-\zeta^7 & \zeta^5-\zeta^8 & \zeta^2-\zeta^{11}\\
  \zeta^5-\zeta^8 & \zeta^2-\zeta^{11} & \zeta^6-\zeta^7
& \zeta-\zeta^{12} & \zeta^3-\zeta^{10} & \zeta^9-\zeta^4\\
  \zeta^2-\zeta^{11} & \zeta^6-\zeta^7 & \zeta^5-\zeta^8
& \zeta^3-\zeta^{10} & \zeta^9-\zeta^4 & \zeta-\zeta^{12}\\
  \zeta^6-\zeta^7 & \zeta^5-\zeta^8 & \zeta^2-\zeta^{11}
& \zeta^9-\zeta^4 & \zeta-\zeta^{12} & \zeta^3-\zeta^{10}
\end{pmatrix}$$
and
$$T=\text{diag}(\zeta^7, \zeta^{11}, \zeta^8, \zeta^6, \zeta^2, \zeta^5)$$
where $\zeta=\exp(2 \pi i/13)$. We have
$$S^2=-I, \quad T^{13}=(ST)^3=I.\eqno{(1.2)}$$
Let
$G=\langle S, T \rangle$, then
$$G \cong \text{SL}(2, 13).$$
We construct some $G$-invariant polynomials in six variables
$z_1, \ldots, z_6$. Let
$$\left\{\aligned
  w_{\infty} &=13 \mathbf{A}_0^2,\\
  w_{\nu} &=(\mathbf{A}_0+\zeta^{\nu} \mathbf{A}_1+\zeta^{4 \nu} \mathbf{A}_2
  +\zeta^{9 \nu} \mathbf{A}_3+\zeta^{3 \nu} \mathbf{A}_4+\zeta^{12 \nu}
  \mathbf{A}_5+\zeta^{10 \nu} \mathbf{A}_6)^2
\endaligned\right.\eqno{(1.3)}$$
for $\nu=0, 1, \ldots, 12$, where the senary quadratic forms (quadratic
forms in six variables) $\mathbf{A}_j$ $(j=0, 1, \ldots, 6)$ are given by
$$\left\{\aligned
  \mathbf{A}_0 &=z_1 z_4+z_2 z_5+z_3 z_6,\\
  \mathbf{A}_1 &=z_1^2-2 z_3 z_4,\\
  \mathbf{A}_2 &=-z_5^2-2 z_2 z_4,\\
  \mathbf{A}_3 &=z_2^2-2 z_1 z_5,\\
  \mathbf{A}_4 &=z_3^2-2 z_2 z_6,\\
  \mathbf{A}_5 &=-z_4^2-2 z_1 z_6,\\
  \mathbf{A}_6 &=-z_6^2-2 z_3 z_5.
\endaligned\right.\eqno{(1.4)}$$
Then $w_{\infty}$, $w_{\nu}$ for $\nu=0, \ldots, 12$ are the roots of a
polynomial of degree fourteen. The corresponding equation is just the
Jacobian equation of degree fourteen (see \cite{K}, pp.161-162). On the
other hand, set
$$\left\{\aligned
  \delta_{\infty} &=13^2 \mathbf{G}_0,\\
  \delta_{\nu} &=-13 \mathbf{G}_0+\zeta^{\nu} \mathbf{G}_1+\zeta^{2 \nu}
  \mathbf{G}_2+\cdots+\zeta^{12 \nu} \mathbf{G}_{12}
\endaligned\right.\eqno{(1.5)}$$
for $\nu=0, 1, \ldots, 12$, where the senary sextic forms (i.e., sextic
forms in six variables) $\mathbf{G}_j$ $(j=0, 1, \ldots, 12)$ are given
by
$$\left\{\aligned
  \mathbf{G}_0 =&\mathbf{D}_0^2+\mathbf{D}_{\infty}^2,\\
  \mathbf{G}_1 =&-\mathbf{D}_7^2+2 \mathbf{D}_0 \mathbf{D}_1+10 \mathbf{D}_{\infty}
                 \mathbf{D}_1+2 \mathbf{D}_2 \mathbf{D}_{12}+\\
                &-2 \mathbf{D}_3 \mathbf{D}_{11}-4 \mathbf{D}_4 \mathbf{D}_{10}
                 -2 \mathbf{D}_9 \mathbf{D}_5,\\
  \mathbf{G}_2 =&-2 \mathbf{D}_1^2-4 \mathbf{D}_0 \mathbf{D}_2+6 \mathbf{D}_{\infty}
                 \mathbf{D}_2-2 \mathbf{D}_4 \mathbf{D}_{11}+\\
                &+2 \mathbf{D}_5 \mathbf{D}_{10}-2 \mathbf{D}_6 \mathbf{D}_9-2
                 \mathbf{D}_7 \mathbf{D}_8,\\
  \mathbf{G}_3 =&-\mathbf{D}_8^2+2 \mathbf{D}_0 \mathbf{D}_3+10 \mathbf{D}_{\infty}
                 \mathbf{D}_3+2 \mathbf{D}_6 \mathbf{D}_{10}+\\
                &-2 \mathbf{D}_9 \mathbf{D}_7-4 \mathbf{D}_{12} \mathbf{D}_4
                 -2 \mathbf{D}_1 \mathbf{D}_2,\\
  \mathbf{G}_4 =&-\mathbf{D}_2^2+10 \mathbf{D}_0 \mathbf{D}_4-2 \mathbf{D}_{\infty}
                 \mathbf{D}_4+2 \mathbf{D}_5 \mathbf{D}_{12}+\\
                &-2 \mathbf{D}_9 \mathbf{D}_8-4 \mathbf{D}_1 \mathbf{D}_3-2
                 \mathbf{D}_{10} \mathbf{D}_7,\\
  \mathbf{G}_5 =&-2 \mathbf{D}_9^2-4 \mathbf{D}_0 \mathbf{D}_5+6 \mathbf{D}_{\infty}
                 \mathbf{D}_5-2 \mathbf{D}_{10} \mathbf{D}_8+\\
                &+2 \mathbf{D}_6 \mathbf{D}_{12}-2 \mathbf{D}_2 \mathbf{D}_3
                 -2 \mathbf{D}_{11} \mathbf{D}_7,\\
  \mathbf{G}_6 =&-2 \mathbf{D}_3^2-4 \mathbf{D}_0 \mathbf{D}_6+6 \mathbf{D}_{\infty}
                 \mathbf{D}_6-2 \mathbf{D}_{12} \mathbf{D}_7+\\
                &+2 \mathbf{D}_2 \mathbf{D}_4-2 \mathbf{D}_5 \mathbf{D}_1-2
                 \mathbf{D}_8 \mathbf{D}_{11},\\
  \mathbf{G}_7 =&-2 \mathbf{D}_{10}^2+6 \mathbf{D}_0 \mathbf{D}_7+4 \mathbf{D}_{\infty}
                 \mathbf{D}_7-2 \mathbf{D}_1 \mathbf{D}_6+\\
                &-2 \mathbf{D}_2 \mathbf{D}_5-2 \mathbf{D}_8 \mathbf{D}_{12}-2
                 \mathbf{D}_9 \mathbf{D}_{11},\\
  \mathbf{G}_8 =&-2 \mathbf{D}_4^2+6 \mathbf{D}_0 \mathbf{D}_8+4 \mathbf{D}_{\infty}
                 \mathbf{D}_8-2 \mathbf{D}_3 \mathbf{D}_5+\\
                &-2 \mathbf{D}_6 \mathbf{D}_2-2 \mathbf{D}_{11} \mathbf{D}_{10}-2
                 \mathbf{D}_1 \mathbf{D}_7,\\
  \mathbf{G}_9 =&-\mathbf{D}_{11}^2+2 \mathbf{D}_0 \mathbf{D}_9+10 \mathbf{D}_{\infty}
                 \mathbf{D}_9+2 \mathbf{D}_5 \mathbf{D}_4+\\
                &-2 \mathbf{D}_1 \mathbf{D}_8-4 \mathbf{D}_{10} \mathbf{D}_{12}-2
                 \mathbf{D}_3 \mathbf{D}_6,\\
  \mathbf{G}_{10} =&-\mathbf{D}_5^2+10 \mathbf{D}_0 \mathbf{D}_{10}-2 \mathbf{D}_{\infty}
                    \mathbf{D}_{10}+2 \mathbf{D}_6 \mathbf{D}_4+\\
                   &-2 \mathbf{D}_3 \mathbf{D}_7-4 \mathbf{D}_9 \mathbf{D}_1-2
                    \mathbf{D}_{12} \mathbf{D}_{11},
\endaligned\right.\eqno{(1.6)}$$
$$\left\{\aligned
  \mathbf{G}_{11} =&-2 \mathbf{D}_{12}^2+6 \mathbf{D}_0 \mathbf{D}_{11}+4 \mathbf{D}_{\infty}
                    \mathbf{D}_{11}-2 \mathbf{D}_9 \mathbf{D}_2+\\
                   &-2 \mathbf{D}_5 \mathbf{D}_6-2 \mathbf{D}_7 \mathbf{D}_4-2
                    \mathbf{D}_3 \mathbf{D}_8,\\
  \mathbf{G}_{12} =&-\mathbf{D}_6^2+10 \mathbf{D}_0 \mathbf{D}_{12}-2 \mathbf{D}_{\infty}
                    \mathbf{D}_{12}+2 \mathbf{D}_2 \mathbf{D}_{10}+\\
                   &-2 \mathbf{D}_1 \mathbf{D}_{11}-4 \mathbf{D}_3 \mathbf{D}_9-2
                    \mathbf{D}_4 \mathbf{D}_8.
\endaligned\right.\eqno{(1.6)}$$
Here, the senary cubic forms (cubic forms in six variables)
$\mathbf{D}_j$ $(j=0$, $1$, $\ldots$, $12$, $\infty)$ are given as follows:
$$\left\{\aligned
  \mathbf{D}_0 &=z_1 z_2 z_3,\\
  \mathbf{D}_1 &=2 z_2 z_3^2+z_2^2 z_6-z_4^2 z_5+z_1 z_5 z_6,\\
  \mathbf{D}_2 &=-z_6^3+z_2^2 z_4-2 z_2 z_5^2+z_1 z_4 z_5+3 z_3 z_5 z_6,\\
  \mathbf{D}_3 &=2 z_1 z_2^2+z_1^2 z_5-z_4 z_6^2+z_3 z_4 z_5,\\
  \mathbf{D}_4 &=-z_2^2 z_3+z_1 z_6^2-2 z_4^2 z_6-z_1 z_3 z_5,\\
  \mathbf{D}_5 &=-z_4^3+z_3^2 z_5-2 z_3 z_6^2+z_2 z_5 z_6+3 z_1 z_4 z_6,\\
  \mathbf{D}_6 &=-z_5^3+z_1^2 z_6-2 z_1 z_4^2+z_3 z_4 z_6+3 z_2 z_4 z_5,\\
  \mathbf{D}_7 &=-z_2^3+z_3 z_4^2-z_1 z_3 z_6-3 z_1 z_2 z_5+2 z_1^2 z_4,\\
  \mathbf{D}_8 &=-z_1^3+z_2 z_6^2-z_2 z_3 z_5-3 z_1 z_3 z_4+2 z_3^2 z_6,\\
  \mathbf{D}_9 &=2 z_1^2 z_3+z_3^2 z_4-z_5^2 z_6+z_2 z_4 z_6,\\
  \mathbf{D}_{10} &=-z_1 z_3^2+z_2 z_4^2-2 z_4 z_5^2-z_1 z_2 z_6,\\
  \mathbf{D}_{11} &=-z_3^3+z_1 z_5^2-z_1 z_2 z_4-3 z_2 z_3 z_6+2 z_2^2 z_5,\\
  \mathbf{D}_{12} &=-z_1^2 z_2+z_3 z_5^2-2 z_5 z_6^2-z_2 z_3 z_4,\\
  \mathbf{D}_{\infty}&=z_4 z_5 z_6.
\endaligned\right.\eqno{(1.7)}$$
Then $\delta_{\infty}$, $\delta_{\nu}$ for $\nu=0, \ldots, 12$ are the
roots of a polynomial of degree fourteen. The corresponding equation is
not the Jacobian equation. Let $S_d=\mathbb{C}[\Phi_{m, n}]$ be the
subalgebra of $\mathbb{C}[z_1, \ldots, z_6]$ generated by the invariant
homogeneous polynomials $\Phi_{m, n}$ given by
$$\Phi_{m, n}=w_0^m \delta_0^n+w_1^m \delta_1^n+\cdots+w_{12}^m \delta_{12}^n+
              w_{\infty}^m \delta_{\infty}^n,\eqno{(1.8)}$$
with degree $d=4m+6n$. Let
$$S=\bigoplus_{d \geq 0} S_d.$$
Then
$$S \subseteq \mathbb{C}[\mathbf{A}_0, \mathbf{A}_1, \ldots, \mathbf{A}_6,
  \mathbf{G}_0, \mathbf{G}_1, \ldots, \mathbf{G}_{12}]^G \subseteq
  \mathbb{C}[z_1, \ldots, z_6]^G.$$

  Let $x_i(z)=\eta(z) a_i(z)$ $(1 \leq i \leq 6)$, where
$$\left\{\aligned
  a_1(z) &:=e^{-\frac{11 \pi i}{26}} \theta \begin{bmatrix}
            \frac{11}{13}\\ 1 \end{bmatrix}(0, 13z),\\
  a_2(z) &:=e^{-\frac{7 \pi i}{26}} \theta \begin{bmatrix}
            \frac{7}{13}\\ 1 \end{bmatrix}(0, 13z),\\
  a_3(z) &:=e^{-\frac{5 \pi i}{26}} \theta \begin{bmatrix}
            \frac{5}{13}\\ 1 \end{bmatrix}(0, 13z),\\
  a_4(z) &:=-e^{-\frac{3 \pi i}{26}} \theta \begin{bmatrix}
            \frac{3}{13}\\ 1 \end{bmatrix}(0, 13z),\\
  a_5(z) &:=e^{-\frac{9 \pi i}{26}} \theta \begin{bmatrix}
            \frac{9}{13}\\ 1 \end{bmatrix}(0, 13z),\\
  a_6(z) &:=e^{-\frac{\pi i}{26}} \theta \begin{bmatrix}
            \frac{1}{13}\\ 1 \end{bmatrix}(0, 13z)
\endaligned\right.\eqno{(1.9)}$$
are theta constants of order $13$ and
$\eta(z):=q^{\frac{1}{24}} \prod_{n=1}^{\infty} (1-q^n)$ with
$q=e^{2 \pi i z}$ is the Dedekind eta function which are all
defined in the upper-half plane
$\mathbb{H}=\{ z \in \mathbb{C}: \text{Im}(z)>0 \}$.
In fact, the weight of $x_i(z)$ is $1$ and the parabolic modular
forms $a_i(z)$ of weight $\frac{1}{2}$ given by (1.9) form a
multiplier-system in the sense of the following (see (7.2) in
Proposition 7.2):
$$\mathbf{A}(z+1)=e^{-\frac{3 \pi i}{4}} T \mathbf{A}(z), \quad
  \mathbf{A}\left(-\frac{1}{z}\right)=e^{\frac{\pi i}{4}} \sqrt{z}
  S \mathbf{A}(z),\eqno{(1.10)}$$
where $S$ and $T$ are given as above, $0<\text{arg} \sqrt{z} \leq \pi/2$
and
$$\mathbf{A}(z):=(a_1(z), a_2(z), a_3(z), a_4(z), a_5(z), a_6(z))^{T}.$$
We will show that there is a morphism
$$\Phi: X \to Z \subset \mathbb{CP}^5$$
with $\Phi(z)=(x_1(z), \ldots, x_6(z))$, where $X=X(13)$ is the modular
curve $\overline{\Gamma(13) \backslash \mathbb{H}}$ and $Z$ is a complete
intersection algebraic curve with multi-degree $(4, 8, 10, 14)$ corresponding
to the ideal
$$I=I(Z)=(\Phi_4, \Phi_8, \Phi_{10}, \Phi_{14}),$$
where
$$\Phi_4=\Phi_{1, 0}, \quad \Phi_8=\Phi_{2, 0}, \quad
  \Phi_{10}=\Phi_{1, 1}, \quad \Phi_{14}=\Phi_{2, 1}.$$
Each $\Phi_i$ $(i=4, 8, 10, 14)$ corresponds to a unique $\Phi_{m, n}$ with
degree $i=4m+6n$. The significance of the algebraic curve $Z$ is that the
finite group $G$ acts linearly on $\mathbb{C}^6$ and on $\mathbb{CP}^5$
leaving invariant $Z \subset \mathbb{CP}^5$. Moreover, it is (1.10) that
gives an explicit realization of the isomorphism between the unique
sub-representation of parabolic modular forms of weight $\frac{1}{2}$
on $X(13)$ and the above six-dimensional complex representation of
$\text{SL}(2, 13)$ generated by $S$ and $T$.

  The beautiful interplay among the theory of automorphic forms,
invariant theory and classical algebraic geometry goes back to Felix
Klein, who studied the following:

\textbf{Problem 1.5.} Let $p$ be a prime. Give an explicit construction
of the modular curve of level $p$ from the invariant theory associated
with $\text{SL}(2, \mathbb{F}_p)$ and its representations using projective
algebraic geometry.

  In particular, when $p=5$, $7$ and $11$, Klein (see \cite{K}, \cite{K1},
\cite{K2}, \cite{K3} and \cite{K4}) created a legacy of beautiful geometry
on modular curves $X(5)$, $X(7)$ and $X(11)$, which even now is a source
of considerable inspiration. However, when $p \geq 13$, especially for
$p=13$, very little is known (see \cite{AR} and \cite{Ad2}). In fact, even
under some additional assumption, the case for $p=13$ is the most difficult
and intractable (see section 4, Theorem 4.4 and Theorem 4.5) because that it
is related to the exceptional Lie group $G_2$.

  In this paper, we answer Problem 1.5 for $p=13$. Our main theorems are
the following: In order to reduce a non-linear situation, the study of
quartic polynomials in six variables, to a linear situation, we produce
a system of linear relations, i.e., the $\mathbf{B}$-terms with the action
of $\text{SL}(2, 13)$ and hence a reducible representation of $\text{SL}(2, 13)$
and its decomposition into irreducible parts, that encapsulates some of
the information about non-linear objects, the quartic polynomials. More
precisely, let $I=I(Y)$ be an ideal generated by the twenty-one quartic
polynomials ($\mathbf{B}$-terms):
$$I=\langle \mathbf{B}_0^{(i)}, \mathbf{B}_1^{(j)}, \mathbf{B}_3^{(j)},
    \mathbf{B}_9^{(j)}, \mathbf{B}_{12}^{(j)}, \mathbf{B}_{10}^{(j)},
    \mathbf{B}_4^{(j)}, \mathbf{B}_5, \mathbf{B}_2, \mathbf{B}_6,
    \mathbf{B}_8, \mathbf{B}_{11}, \mathbf{B}_7 \rangle\eqno{(1.11)}$$
with $i=0. 1, 2$ and $j=1, 2$, where
$$\left\{\aligned
  \mathbf{B}_0^{(0)} &=z_1 z_2 z_4 z_5+z_2 z_3 z_5 z_6+z_3 z_1 z_6 z_4,\\
  \mathbf{B}_0^{(1)} &=z_1 z_5^3+z_2 z_6^3+z_3 z_4^3,\\
  \mathbf{B}_0^{(2)} &=z_1^3 z_6+z_2^3 z_4+z_3^3 z_5,
\endaligned\right.\eqno{(1.12)}$$
$$\left\{\aligned
  \mathbf{B}_1^{(1)} &=z_3 z_5^3+z_1^3 z_4-z_1 z_2^3,\\
  \mathbf{B}_1^{(2)} &=z_2 z_4 z_6^2-z_3^2 z_6 z_4-z_1^2 z_2 z_5,\\
  \mathbf{B}_3^{(1)} &=z_2 z_4^3+z_3^3 z_6-z_3 z_1^3,\\
  \mathbf{B}_3^{(2)} &=z_1 z_6 z_5^2-z_2^2 z_5 z_6-z_3^2 z_1 z_4,\\
  \mathbf{B}_9^{(1)} &=z_1 z_6^3+z_2^3 z_5-z_2 z_3^3,\\
  \mathbf{B}_9^{(2)} &=z_3 z_5 z_4^2-z_1^2 z_4 z_5-z_2^2 z_3 z_6,\\
  \mathbf{B}_{12}^{(1)} &=z_1 z_4^3+z_2^3 z_6+z_4 z_5^3,\\
  \mathbf{B}_{12}^{(2)} &=z_2 z_5 z_4^2-z_3^2 z_1 z_5-z_6^2 z_3 z_1,\\
  \mathbf{B}_{10}^{(1)} &=z_3 z_6^3+z_1^3 z_5+z_6 z_4^3,\\
  \mathbf{B}_{10}^{(2)} &=z_1 z_4 z_6^2-z_2^2 z_3 z_4-z_5^2 z_2 z_3,\\
  \mathbf{B}_4^{(1)} &=z_2 z_5^3+z_3^3 z_4+z_5 z_6^3,\\
  \mathbf{B}_4^{(2)} &=z_3 z_6 z_5^2-z_1^2 z_2 z_6-z_4^2 z_1 z_2,
\endaligned\right.\eqno{(1.13)}$$
and
$$\left\{\aligned
  \mathbf{B}_5 &=-z_2^2 z_1 z_5+z_4 z_5 z_6^2+z_2 z_3 z_4^2,\\
  \mathbf{B}_2 &=-z_1^2 z_3 z_4+z_6 z_4 z_5^2+z_1 z_2 z_6^2,\\
  \mathbf{B}_6 &=-z_3^2 z_2 z_6+z_5 z_6 z_4^2+z_3 z_1 z_5^2,\\
  \mathbf{B}_8 &=z_2 z_4 z_5^2+z_1 z_2 z_3^2+z_1^2 z_5 z_6,\\
  \mathbf{B}_{11} &=z_1 z_6 z_4^2+z_3 z_1 z_2^2+z_3^2 z_4 z_5,\\
  \mathbf{B}_7 &=z_3 z_5 z_6^2+z_2 z_3 z_1^2+z_2^2 z_6 z_4.
\endaligned\right.\eqno{(1.14)}$$
The corresponding curve associated to the ideal $I(Y)$ is denoted
by $Y$.

\textbf{Theorem 1.6.} (Main Theorem 1) {\it The modular curve $X=X(13)$
is isomorphic to the curve $Y$ in $\mathbb{CP}^5$.}

  Hence, Theorem 1.6 gives a projective model of the modular curve $X(13)$
in $\mathbb{P}^5$. Moreover, the following theorem shows that the modular
curve $X(13)$ is connected with a Fano four-fold (higher dimensional
algebraic variety) by the invariant theory. Consequently, it shows that
where the above $\mathbf{B}$-terms come from.

\textbf{Theorem 1.7.} (Main Theorem 2) {\it The curve $Y$ can be constructed
from the invariant quartic Fano four-fold $\Phi_4=0$ in $\mathbb{CP}^5$, i.e.,
$$(z_3 z_4^3+z_1 z_5^3+z_2 z_6^3)-(z_6 z_1^3+z_4 z_2^3+z_5 z_3^3)
  +3(z_1 z_2 z_4 z_5+z_2 z_3 z_5 z_6+z_3 z_1 z_6 z_4)=0.\eqno{(1.15)}$$}

  Note that there is a close relationship between ideals $I(Y)$ and varieties
$Y$ which reveals the intimate link between algebra and geometry. The solutions
of a system of polynomial equations form a geometric object called a variety,
the corresponding algebraic object is an ideal. Moreover, invariant theory play
an important role, which has had a profound effect on the development of
algebraic geometry. Theorem 1.6 and Theorem 1.7 give the geometric aspect of
the algebraic space curve $Y$. For its algebraic aspect, the corresponding ideal
$I(Y)$ admits an action of the group $G \cong \text{SL}(2, 13)$ by the theory
of invariants.

\textbf{Theorem 1.8.} (Main Theorem 3) {\it The ideal $I(Y)$ is invariant
under the action of $G$, which gives a twenty-one dimensional representation
of $G$.}

  Theorem 1.8 leads to a $21$-dimensional representation of $\text{SL}(2, 13)$.
Here, the irreducible representations of $\text{SL}(2, q)$ are well known for
a very long time (see \cite{Fr1}, \cite{J} and \cite{S2}). In fact, Schur (see
\cite{S2}) derived the character tables of the groups $\text{SL}(2, q)$ for all
values of $q$. For $G=\text{SL}(2, q)$, there are at most
$$1+1+\frac{q-3}{2}+2+\frac{q-1}{2}+2=q+4$$
conjugacy classes which can appear in the following formula on the
dimensions of ordinary representations:
$$|G|=1^2+q^2+\frac{q-3}{2} (q+1)^2+2 \cdot \left(\frac{q+1}{2}\right)^2
     +\frac{q-1}{2} (q-1)^2+2 \cdot \left(\frac{q-1}{2}\right)^2.$$
In particular, when $q=13$, $|G|=2184$. There are at most $q+4=17$
conjugacy classes and
$$1^2+13^2+5 \times 14^2+2 \times 7^2+6 \times 12^2+2 \times 6^2=2184$$
corresponding to the decomposition of conjugacy classes:
$$1+1+5+2+6+2=17.$$

  It is known that there are the following representations $\rho$
for $\text{SL}(2, 13)$:

(1) discrete series
$$\text{dim}(\rho)=\left\{\aligned
  &12 \quad \text{(super-cuspidal representations) or}\\
  &6 \quad \text{(two degenerate discrete series representations)}
\endaligned\right.$$

(2) principal series
$$\text{dim}(\rho)=\left\{\aligned
  &14 \quad \text{(principal series representations) or}\\
  &13 \quad \text{(Steinberg representation) or}\\
  &7 \quad \text{(two degenerate principal series representations)}
\endaligned\right.$$

\textbf{Theorem 1.9.} (Main Theorem 4) {\it The $21$-dimensional representation
is reducible, which can be decomposed as the direct sum of a $1$-dimensional
representation $($the trivial representation$)$, a $7$-dimensional representation
$($the degenerate principal series representation$)$ and a $13$-dimensional
representation $($the Steinberg representation$)$$:$
$$\mathbf{21}=\mathbf{1} \oplus \mathbf{7} \oplus \mathbf{13}.\eqno{(1.16)}$$
More precisely, let $V$ be a complex vector space generated by the $21$ quartic
polynomials $\mathbf{B}_0^{(0)}$, $\mathbf{B}_0^{(1)}$, $\ldots$, $\mathbf{B}_7$.
Then $V$ has the following decomposition:
$$V=V_1 \oplus V_7 \oplus V_{13},\eqno{(1.17)}$$
where the $1$-dimensional subspace
$$V_1=\langle 3 \mathbf{B}_0^{(0)}+\mathbf{B}_0^{(1)}-\mathbf{B}_0^{(2)} \rangle,
      \eqno{(1.18)}$$
the $7$-dimensional subspace
$$\aligned
  V_7=\langle &\mathbf{B}_0^{(1)}+\mathbf{B}_0^{(2)},
               \mathbf{B}_1^{(1)}-3 \mathbf{B}_1^{(2)},
               \mathbf{B}_3^{(1)}-3 \mathbf{B}_3^{(2)},
               \mathbf{B}_9^{(1)}-3 \mathbf{B}_9^{(2)},\\
              &\mathbf{B}_{12}^{(1)}+3 \mathbf{B}_{12}^{(2)},
               \mathbf{B}_{10}^{(1)}+3 \mathbf{B}_{10}^{(2)},
               \mathbf{B}_4^{(1)}+3 \mathbf{B}_4^{(2)}
\rangle\endaligned\eqno{(1.19)}$$
and the $13$-dimensional subspace
$$\aligned
  V_{13}=\langle &4 \mathbf{B}_0^{(0)}-\mathbf{B}_0^{(1)}+\mathbf{B}_0^{(2)},
               \mathbf{B}_5, \mathbf{B}_2, \mathbf{B}_6, \mathbf{B}_8,
               \mathbf{B}_{11}, \mathbf{B}_7,\\
              &\mathbf{B}_1^{(1)}+\mathbf{B}_1^{(2)},
               \mathbf{B}_3^{(1)}+\mathbf{B}_3^{(2)},
               \mathbf{B}_9^{(1)}+\mathbf{B}_9^{(2)},\\
              &-\mathbf{B}_{12}^{(1)}+\mathbf{B}_{12}^{(2)},
               -\mathbf{B}_{10}^{(1)}+\mathbf{B}_{10}^{(2)},
               -\mathbf{B}_4^{(1)}+\mathbf{B}_4^{(2)} \rangle
\endaligned\eqno{(1.20)}$$
are stable under the action of $G \cong \text{SL}(2, 13)$.}

  Note that the basis of $V_1$ is just the invariant quartic four-fold
$\Phi_4=0$. This gives a natural viewpoint from representation theory
that it is an $\text{SL}(2, 13)$-invariant quartic four-fold.

  In order to establish an algebraic theory of algebraic functions,
in a series of remarkable papers \cite{Hec1}, \cite{Hec2}, \cite{Hec3},
\cite{Hec4} and \cite{Hec5}, Hecke studied the arithmetic theory of
elliptic modular functions by the methods of function theory
(Riemann-Weierstrass), the methods of arithmetic (Dedekind-Weber),
Galois theory and Frobenius's character theory.

  In the introduction of his work \cite{Hec5}, Hecke stated that
``If an algebraic structure in a variable of genus $p$ has a group $\mathfrak{G}$
of one-to-one analytical mappings (automorphisms), then the $p$ linearly
independent differentials of the first kind experience linear homogeneous
substitutions in these mappings, which together form a representation
$\mathfrak{S}$ of the abstract group $\mathfrak{G}$. A general method of
determining this representation, i.e. determining its character, was first
given by Hurwitz. For the special structure, which is defined by the elliptic
modular functions of a prime number level $q$, I determined the decomposition
of the above representation into its irreducible components a few years ago
in a somewhat different way, whereby a strange connection with the class number
of the quadratic number field $K(\sqrt{-q})$ yielded. Since then, the question
has been taken up again by C. Chevalley and A. Weil, and a very general theorem
about these and related groups has been proved by algebraic methods.''

  In the other work \cite{Hec3}, page 551, Hecke stated that ``We now consider
the algebraic structure, which is defined by the elliptic modular functions of
the prime level $q$ ($q \geq 3$). The structure has a group of mappings into itself,
which are explained by the totality of the modular substitutions applied to the
argument $\tau$; the mappings are thus reduced to the system of modular substitutions,
mod $q$, and isomorphic with the modular group $\mathfrak{M}(q)$ of order
$\frac{q(q^2-1)}{2}$. The group of integrals of the first kind induced in
this way is decomposed into its irreducible components
$$\mathfrak{S}=x \mathfrak{G}_q+y_1 \mathfrak{G}_{\frac{q+\varepsilon}{2}}+
  y_2 \mathfrak{G}_{\frac{q+\varepsilon}{2}}^{\prime}+\sum_{i} u_i
  \mathfrak{G}_{q+1}^{(i)}+\sum_{i} v_i \mathfrak{G}_{q-1}^{(i)}.$$
Here the $\mathfrak{G}_n$, $\mathfrak{G}_n^{(i)}$ are the various possible
irreducible representations of $\mathfrak{M}(q)$, $n$ being the degree,
$\varepsilon=(-1)^{\frac{q-1}{2}}$; the multiplicities $x$, $u_i$, $v_i$
and also $y_1$, $y_2$ have already been calculated by me.''

  In modern terminology, using modular curves $X(p)$ with $p$ a prime,
Hecke constructed a finite dimensional (reducible) representation
$H^0(X(p), \Omega_{X(p)}^{1})$ of the finite group $\text{PSL}(2, \mathbb{F}_p)$
whose dimension is equal to the genus of $X(p)$, i.e. the number of
linearly independent differentials of the first kind:
$$\dim H^0(X(p), \Omega_{X(p)}^{1})=g(X(p)).$$
This construction has the rich structure: the finite dimensional representation
of $\text{PSL}(2, \mathbb{F}_p)$ can be obtained by the methods from geometry,
topology and complex analysis, i.e. Riemann-Weierstrass function theory. In
particular, with the help of the character theory established by Frobenius
(see \cite{Fr1}), this leads Hecke to study the decomposition of the vector
space of regular differentials on $X(p)$: $H^0(X(p), \Omega_{X(p)}^{1})$ under
the action of $\text{PSL}(2, \mathbb{F}_p)$ as the direct sum of spaces of
irreducible representations for $\text{PSL}(2, \mathbb{F}_p)$ defined over
some number fields (see \cite{Hec1}, \cite{Hec2}, \cite{Hec3}, \cite{Hec4}, and
\cite{Hec5}). In particular, he proved that (see \cite{Hec5}, p.761, Satz 14)
that every irreducible representation of $\text{PSL}(2, \mathbb{F}_p)$ is
equivalent to one whose coefficients lie in the field of its character
$\mathbb{Q}(\chi)$.

  According to \cite{Se2}, Chap. 12, let $K$ denote a field of characteristic
zero and $C$ an algebraic closure of $K$. If $V$ is a $K$-vector space, we let
$V_C$ denote the $C$-vector space $C \otimes_K V$ obtained from $V$ by extending
scalars from $K$ to $C$. If $G$ is a finite group, each linear representation
$\rho: G \rightarrow \text{GL}(V)$ over the field $K$ defines a representation
$$\rho_C: G \rightarrow \text{GL}(V) \rightarrow \text{GL}(V_C)$$
over the field $C$. In terms of modules, we have
$$V_C=C[G] \otimes_{K[G]} V.$$
The character $\chi_{\rho}=\text{Tr}(\rho)$ of $\rho$ is the same as for $\rho_C$.
It is a class function on $G$ with values in $K$. We denote by $R_K(G)$ the group
generated by the characters of the representations of $G$ over $K$. It is a sub-ring
of the ring $R(G)=R_C(G)$.

  A linear representation of a finite group $G$ over $C$ is said to be realizable
over $K$ (or rational over $K$) if it is isomorphic to a representation of the form
$\rho_C$, where $\rho$ is a linear representation of $G$ over $K$; this amounts to
saying that it can be realized by matrices having coefficients in $K$.

  Now, let us consider the realizability over cyclotomic fields. Denote by $m$ the
least common multiple of the orders of the elements of $G$. It is a divisor of $g$,
the order of $G$. We have Brauer's theorem (see \cite{Se2}, 12.3, Theorem 24):
If $K$ contains the $m$-th roots of unity, then $R_K(G)=R(G)$. Consequently,
each linear representation of $G$ can be realized over $K$.

  In his arithmetic study of complex representations of finite groups, Schur
(see \cite{S1}) introduced the notion of the index. Roughly speaking, suppose
$G$ is a finite group, $K$ is a splitting field for $G$, and $\chi$ is the
character of an irreducible linear representation $\rho$ of $G$ over $K$.
Suppose $k$ is the subfield of $K$ generated by the character values $\chi(g)$,
$g \in G$. The Schur index of $\chi$ or the Schur index of $\rho$ is defined as
the smallest positive integer $m$ such that there exists a degree $m$ extension
$L$ of $k$ such that $\rho$ can be realized over $L$, i.e., we can change basis
so that all the matrix entries are from $L$. The Schur index of a character
$\chi$ is denoted $m(\chi)$ (see also \cite{Se2} for the definition using the
theory of semi-simple algebras). Many important results about Schur indices appear
to depend on deep facts about division algebras and number theory (see \cite{Fei}
and \cite{I}, Chap. 10). Note that if the representation can be realized over the field
generated by the character values for that representation, the Schur index is one.
In particular, all irreducible complex characters of $\text{PSL}(2, \mathbb{F}_q)$
have Schur index $1$ (see \cite{Jan}, also \cite{S2}, \cite{Dor}, \S 38 and
\cite{Se2}, \S 12.6).

  In fact, let $G$ be a finite group and $\rho: G \rightarrow \text{GL}(n, \mathbb{C})$
be an irreducible complex representation. Then $\rho$ descends to a representation
over $\overline{\mathbb{Q}}$ uniquely up to isomorphism. Hence, we can consider
$\rho$ as a representation with values in $\text{GL}(n, \overline{\mathbb{Q}})$.
Under $\overline{\mathbb{Q}}$, any field over which $\rho$ is definable must
contain the trace field of $\rho$, defined as the extension of $\mathbb{Q}$
generated by the traces of all the matrices $\rho(g)$, $g \in G$, that is,
$\mathbb{Q}(\text{tr}(\rho(g)), \rho \in G)$, which is just the character field
$\mathbb{Q}(\chi)$. The trace field of $\rho$ is characterized as the fixed field
of Galois automorphisms $\sigma \in \text{Gal}(\overline{\mathbb{Q}}/\mathbb{Q})$
where $\rho \cong \rho^{\sigma}$, the Galois twist of $\rho$ by $\sigma$. Thus,
we can denote the trace field of $\rho$ as $\mathbb{Q}(\rho)$. The obstruction to
the representation descending to its trace field can be described by a Galois
$2$-cocycle, which represents an element of the Brauer group of the trace field,
i.e., the Brauer obstruction of $\rho$. More precisely, the Brauer obstruction
$[\psi_{\rho}] \in H^2(\text{Gal}(\overline{\mathbb{Q}}/\mathbb{Q}(\rho)),
\overline{Q}^{\times})=\text{Br}(\mathbb{Q}(\rho))$ is defined as follows: for
$\sigma \in \text{Gal}(\overline{\mathbb{Q}}/\mathbb{Q}(\rho))$, the isomorphism
$\rho \cong \rho^{\sigma}$ is unique up to $\overline{\mathbb{Q}}^{\times}$ by
Schur lemma. When we choose isomorphisms $\rho \cong \rho^{\sigma}$, we obtain
the Brauer obstruction in $H^2(\text{Gal}(\overline{\mathbb{Q}}/\mathbb{Q}(\rho)),
\overline{Q}^{\times})$, which is well-defined as the cohomology class is independent
of the choice of isomorphisms.

  Despite of this, we still need the realization in the field $\mathbb{Q}(\chi)$
and the exact knowledge of it. Thanks to Deligne (see \cite{De2} and \cite{De4}),
who gives an elegant and different proof using the Kirillov model in representation
theory.

\textbf{Deligne's realization theorem.} (see \cite{De4}). {\it Any representation
of $\text{GL}(2, \mathbb{F}_q)$ can be realized over the field of its character.}

  In fact, Deligne proved the more general result: The split (resp. non-split)
torus $T_0$ of $\text{SL}(2, \mathbb{F}_q)$ are, up to $x \mapsto x^{-1}$,
naturally isomorphic to $\mathbb{F}_q^{\times}$ (diagonal matrices $\left(\begin{matrix}
x & 0\\ 0 & x^{-1} \end{matrix}\right)$), resp. Ker(N: $\mathbb{F}_{q^2}^{\times}
\rightarrow \mathbb{F}_q^{\times})$. Let $\chi_0$ be a nontrivial character
of $T_0$, $t_0$ be its order, and $X_{\chi_0}$ be the corresponding
representation of $\text{SL}(2, \mathbb{F}_q)$, well defined up to
isomorphism. For $\chi_0$ of order $2$, we mean the sum of the two
representations of dimension $\frac{q+1}{2}$ (resp. $\frac{q-1}{2}$).
Assume that $\chi_0(-1)=1$. This is the case for $q$ even. For $q$ odd,
it is equivalent to $t_0| \frac{q-1}{2}$ (resp. $t_0 | \frac{q+1}{2}$).
It is also equivalent to $X_{\chi_0}$ factoring through a representation
of $\text{PSL}(2, \mathbb{F}_q)$.

  Under this assumption, Hecke states that $X_{\chi_0}$ can be realized
over the field (generated by values) of his character, which is the real
subfield $F_0$ of the $t_0$-cyclotomic field.

   Let $H$ be the following quaternion algebra over $F_0$: it is generated
by the quadratic extension $F$ of $F_0$, and by an element $u$, such that
conjugation by $u$ induces on $F$ the generator of $\text{Gal}(F/F_0)$, and
that
$$\left\{\aligned
 &\text{in the split case (principal series)} \quad : u^2=q,\\
 &\text{in the non-split case (discrete series)} : u^2=-q.
\endaligned\right.$$
This algebra is a matrix algebra if and only if $q$ (resp. $-q$) is
in the image of the norm $\text{N}_{F/F_0}$.

  In particular, Deligne proved that the obstruction to realize the
representation of $\text{PSL}(2, \mathbb{F}_q)$ corresponding to $\chi_0$
over the field of its character, that is over $F_0$, is the class of the
quaternion algebra $H$ in the Brauer group of $F_0$.

  Hence, in order to prove (at least for $t_0 \neq 2$) Hecke's statement
that representations of $\text{PSL}(2, \mathbb{F}_q)$ are realizable over
the field of their characters, it remains to prove that the obstruction made
explicit as above vanishes. In particular, for $X$ a discrete series representation
of $\text{GL}(2, \mathbb{F}_q)$, it essentially leads to the Kirillov model of
the representation $X$.

  Recall that the main result of \cite{Hec5} is expressed in the following
simple formula: The multiplicity $r$ of an irreducible representation $\rho$ of
$\text{PSL}(2, \mathbb{F}_p)$ with the character $\chi$ within the group
$H^0(X(p), \Omega_{X(p)}^{1})+\overline{H^0(X(p), \Omega_{X(p)}^{1})}$
(where $\rho$ is not the identical representation) is
$$\kappa(\rho)+\kappa(\overline{\rho})=r=
 f-\frac{1}{p} \sum_{\text{$n$ mod $p$}} \chi(T^n)-\frac{1}{2} \sum_{\text{$n$ mod $2$}}
  \chi(S^n)-\frac{1}{3} \sum_{\text{$n$ mod $3$}} \chi((ST)^n).$$
Here $T$, $S$ are the two modular substitutions $\tau^{\prime}=\tau+1$,
$\tau^{\prime}=-\frac{1}{\tau}$, $f$ is the degree of $\rho$. Therefore,
Hecke obtained the following beautiful decomposition formula:
$$H^0(X(p), \Omega_{X(p)}^{1})=m V_p \oplus n_1 V_{\frac{p+\varepsilon}{2}}
  \oplus n_2 V_{\frac{p+\varepsilon}{2}}^{\prime} \oplus \bigoplus_i u_i
  V_{p+1}^{(i)} \oplus \bigoplus_j v_j V_{p-1}^{(j)},$$
where $\varepsilon=(-1)^{\frac{p-1}{2}}$, $i=1, 2, \ldots, \frac{q-\varepsilon}{4}-1$
and $j=1, 2, \ldots, \frac{q+\varepsilon}{4}-\frac{1}{2}$. The multiplicities
$m$, $n_1$, $n_2$, $u_i$ and $v_j$ can be calculated by the above formula on $r$.
Note that for small primes $p$, the irreducible representations appearing in the
decomposition formula  of $H^0(X(p), \Omega_{X(p)}^{1})$ can not exhaust all
of the irreducible representations of $\text{PSL}(2, \mathbb{F}_p)$, especially
for $p=7$, $11$ and $13$.

  Specifically, for every $p>3$, $p \equiv 3$ (mod $4$), Hecke (see \cite{Hec5},
p.768, Satz 16) described a $\text{SL}(2, \mathbb{F}_p)$-invariant subspace of
dimension $h \cdot \frac{p-1}{2}$ in $H^0(X(p), \Omega_{X(p)}^{1})$, where $h$
is the class number of the imaginary quadratic field $K=\mathbb{Q}(\sqrt{-p})$,
the group $\text{SL}(2, \mathbb{F}_p)$ acts on this subspace by $h$ copies of $W$,
one of the two irreducible degenerate discrete series representations of
$\text{SL}(2, \mathbb{F}_p)$ of dimension $\frac{p-1}{2}$ (see also \cite{Sh1},
\cite{Sh2} and \cite{Gro2} for the later development). This subspace is spanned
by weighted binary theta series, and Hecke identified certain periods of these
differentials as the periods of elliptic curves with complex multiplication
by $K$. Namely, Hecke showed that (see \cite{Hec5}, p.768, Satz 16) among the
$h \cdot \frac{p-1}{2}$ integrals of the first kind of this family, the independent
ones can be selected in such a way that their periods at $\Gamma(p)$ are always
integers from $\mathbb{Q}(\sqrt{-p})$. These elliptic curves appear as factors
of the Jacobian of the modular curve $X(p)$.

  In particular, for the smallest such primes $p=7$ and $p=11$ with class number
one, Hecke gave the following remarkable decompositions:
$$H^0(X(7), \Omega_{X(7)}^1)=V_3,$$
$$H^0(X(11), \Omega_{X(11)}^1)=V_5 \oplus V_{11} \oplus V_{10}.$$

  Here, when $p=7$, the curve $X(7)$ has genus three and is isomorphic to the
Klein quartic curve. In this case, Hecke showed that the multiplicity $m(W)$
for the irreducible representations $W$ of $\text{PSL}(2, 7)$ is equal to one,
for $W$ is one of the degenerate discrete series representations of dimension
$3$, and equal to zero for all other irreducible representations of $\text{PSL}(2, 7)$.
Moreover, the integrals of $V_3$ only have periods that are integers from
$\mathbb{Q}(\sqrt{-7})$.

  When $p=11$ the curve $X(11)$ has genus $26$. In this case, Hecke showed
that the multiplicity $m(W)$ is equal to one for $W$ is the degenerate discrete
series representation of dimension $5$, the Steinberg representation of
dimension $11$, and the discrete series representation of dimension $10$,
and equal to zero for all other irreducible representations of
$\text{PSL}(2, 11)$. Hecke also pointed out that (see \cite{Hec5}, p.769,
Satz 18) for the modular curve $X(11)$ of level $11$ and genus $g=26$, the
$26$ integrals of the first kind can be selected as elliptic integrals.
Moreover, the integrals of $V_5$ only have periods that are integers from
$\mathbb{Q}(\sqrt{-11})$, and $V_{11}$, $V_{10}$ have a rational character
and the multiplicity $1$.

  However, for every $p>5$, $p \equiv 1$ (mod $4$), it becomes completely
different. All of the conjugacy classes in $\text{SL}(2, \mathbb{F}_p)$ are
real and all of the irreducible characters take real values. More precisely,
Hecke (see \cite{Hec5}, p.762, p.763 and p.766) showed that for the Steinberg
representation of dimension $p$, the only such representation is of a rational
character. For the degenerate principal series representation of dimension
$\frac{p+1}{2}$, the two representations of this degree are first known in a
form from $\mathbb{Q}(\zeta_p)$, where $\zeta_p=e^{\frac{2 \pi i}{p}}$, while
the field $\mathbb{Q}(\chi)$ is $\mathbb{Q}(\sqrt{p})$. For the principal series
representation of dimension $p+1$, the characters of such a representation are
certain numbers from the real subfield of the field of the primitive $t$-th roots
of unity $\varrho$, where $t$ passes through all divisors of $\frac{p-1}{2}$
except $1$ and $2$. In fact, there are exactly $\frac{1}{2} \varphi(t)$ different
simple characters of degree $p+1$ for each of these $t$. The field
$\mathbb{Q}(\chi)$ is the real subfield $\mathbb{Q}(\varrho+\varrho^{-1})$.
For the discrete series representation of dimension $p-1$, the characters of
such a representation is again a primitive $t$-th roots of unity $\sigma$, where
now $t$ must be a divisor of $\frac{p+1}{2}$ and $t>2$, such that a total of
$\frac{1}{2} \varphi(t)$ different representations belong to the different
$\sigma$ with the same $t$, where $\sigma$ and $\sigma^{-1}$ give the same
representation. The character $\chi$ generates the real field
$\mathbb{Q}(\sigma+\sigma^{-1})$.

  According to \cite{Hec5}, p.770), the general case of the representations
$V_f^{(l)}(t)$ is that the field $\mathbb{Q}(\chi)$ of their character is a
totally real (abelian) number field of degree $n=\frac{1}{2} \varphi(t)$.
First let $\kappa=1$, such as $p=13$. Each of the three representations
$V_{12}^{(l)}(7)$ ($l=1, 2, 3$) with conjugate coefficients from the real
cubic subfield of the $7$th roots of unity $\mathbb{Q}(\sigma+\sigma^{-1})$
is realized by an integral vector of the first kind with $12$ components. In
general, Hecke (see \cite{Hec5}, p.771, Satz 19) prove that the period problem
for the $n=\frac{1}{2} \varphi(t)$ algebraic conjugate representations $V_f^{(l)}(t)$
$(l=1, \ldots, n)$ is, if the multiplicity $\kappa$ of them in the integral group
is $1$, equivalent to determining the values of Hilbert modular functions of $n$
variables over the totally real field of the $n$-th degree which is determined
by the characters of the representations $V_f^{(l)}(t)$.

  In particular, for the smallest such prime $p=13$, Hecke showed that the
following deep result:
$$H^0(X(13), \Omega_{X(13)}^1)=V_{14} \oplus V_{12}^{(1)} \oplus V_{12}^{(2)}
                        \oplus V_{12}^{(3)},$$
where $q=13$, $\frac{q-1}{2}=6$.  Hence, $t=3$ or $6$. This implies that
the field $\mathbb{Q}(\chi)$ is $\mathbb{Q}$ and $V_{14}$ is defined over
$\mathbb{Q}$. On the other hand, $V_{12}^{(i)}$ $(i=1, 2, 3)$ are conjugate
over the totally real field $\mathbb{Q}(\zeta_7+\zeta_7^{-1})$ with $\zeta_7$
the seventh root of unity. In this case, Hecke showed that the multiplicity
$m(W)$ is equal to one for $W$ is the principal series representation of
dimension $14$, and the three conjugate discrete series representations of
dimension $12$, and equal to zero for all other irreducible representations
of $\text{SL}(2, 13)$. In particular, the trivial representation of dimension
one, the degenerate principal series of dimension $7$ and the Steinberg
representation of dimension $13$ do not appear in the above decomposition
corresponding to the modular curve $X(13)$. A natural problem is {\it why
these three irreducible representations do not appear and where do they appear?}

  What makes Hecke's decomposition so appealing is the wide range of mathematics
it touches: it can be regarded as a result in modular forms and modular curves
(for $X(p)$), number theory (for the field of definition), representation theory
(for $\text{SL}(2, \mathbb{F}_p)$), algebraic geometry (for the genus
$g=\dim H^0(X(p), \Omega_{X(p)}^1)$), and homological algebra (for the direct
sum decomposition for the cohomology groups $H^0(X(p), \Omega_{X(p)}^1)$).

  Theorem 1.9 gives a distinct decomposition (1.17) defined over the cyclotomic
field $\mathbb{Q}(\zeta)$, where the modular forms and modular curves (for
$X(13)$), number theory (for the field $\mathbb{Q}(\zeta)$ of definition),
representation theory (for $\text{SL}(2, 13)$), algebraic geometry ($21$
quartic equations give the locus for $X(13)$ in $\mathbb{P}^5$), and
homological algebra (for the direct sum decomposition for the linear space
$V$ associated with the ideal $I(Y)$ for $X(13)$ corresponding to
$\Gamma(Y, \mathcal{O}_Y)$) are put together. In fact, let
$A(Y)=\mathbb{C}[z_1, z_2, z_3, z_4, z_5, z_6]/I(Y)$ be the affine coordinate
ring of $Y$ where we consider $Y \subset \mathbb{A}^6(\mathbb{C})$ as an affine
variety, then for coherent algebraic sheaves on affine varieties we have (see
\cite{Se1}, p.237, Corollaire 3) $A(Y) \cong \Gamma(Y, \mathcal{O}_Y)$, hence
it is an invariant up to isomorphism. In particular, the three missing
irreducible representations of dimension $1$, $7$ and $13$ in Hecke's
decomposition corresponding to the modular curve $X(13)$ do appear in the
decomposition (1.17). Hence, we have the following theorem which shows
that $X(13)$ is very exceptional in the sense that even for $p \equiv 1$
(mod $4$), such as $p=13$, there exist two distinct decompositions
corresponding to the same modular curve $X(13)$ with different fields
of definition: one is associated with the $7$-th cyclotomic field, the
other is associated with the $13$-th cyclotomic field.

\textbf{Theorem 1.10.} (Main Theorem 5) (Distinct arithmetical realizations
by cohomology groups of a projective or affine variety with values in a
coherent algebraic sheaf corresponding to the same modular curve $X(13)$)
{\it There are two distinct arithmetical realizations of $X(13)$ by
character fields $\mathbb{Q}(\chi)=\mathbb{Q}(\zeta_7+\zeta_7^{-1})$ or
$\mathbb{Q}(\chi)=\mathbb{Q}(\sqrt{13})$ $($i.e., two distinct algebraic
number fields$)$ of irreducible representations of $\text{SL}(2, 13)$ in
two different ways corresponding to the decompositions of some spaces on
$X(13)$:

(1) Hecke's arithmetical realization $($the standard realization$)$ over
$\mathbb{Q}(\chi)=\mathbb{Q}(\zeta_7+\zeta_7^{-1})$ by cohomology groups
of a projective variety with values in a coherent algebraic sheaf
$\Omega_{X(13)}^1$ corresponding to Hecke's decomposition of the space of
regular differentials on $X(13)$:
$$H^0(X(13), \Omega_{X(13)}^1)=V_{14} \oplus V_{12}^{(1)} \oplus V_{12}^{(2)}
  \oplus V_{12}^{(3)}, \quad (\text{genus: $g(X(13))=50$}),\eqno{(1.21)}$$
which involves the principal series representation and the discrete series
representations, the former one is defined over $\mathbb{Q}$, the latter
three are defined over the totally real subfield
$\mathbb{Q}(\zeta_7+\zeta_7^{-1})$ of the $7$-th cyclotomic field
$\mathbb{Q}(\zeta_7)$. In particular, the character field of $V_{14}$ is
$\mathbb{Q}$, while the character fields of $V_{12}^{(i)}$ $(i=1, 2, 3)$
are $\mathbb{Q}(\zeta_7+\zeta_7^{-1})$.

(2) Our arithmetical realization $($the nonstandard arithmetical realization$)$
over $\mathbb{Q}(\chi)=\mathbb{Q}(\sqrt{13})$ by cohomology groups of an
affine variety with values in a coherent algebraic sheaf $\mathcal{O}_Y$
corresponding to our decomposition of the space of $21$ quartic equations
$($i.e. the locus $\mathcal{L}$$)$ on $X(13)$:
$$\Gamma(Y, \mathcal{O}_Y) \cong A(Y):=
  \mathbb{C}[z_1, z_2, z_3, z_4, z_5, z_6]/I(Y), \eqno{(1.22)}$$
where the vector space $V_{I(Y)}$ generated by the generators of the ideal
$I(Y)$ has the following decomposition:
$$\text{$V_{I(Y)}=V_1 \oplus V_7 \oplus V_{13}$, i.e.
  $\mathbf{21}=\mathbf{1} \oplus \mathbf{7} \oplus \mathbf{13}$,
  $($dimension of the ideal: $21$$)$}\eqno{(1.23)}$$
which involves the trivial representation, the degenerate principal series
representation and the Steinberg representation, all of them are defined
over the $13$-th cyclotomic field $\mathbb{Q}(\zeta)$. In particular, the
character fields of $V_1$ and $V_{13}$ are $\mathbb{Q}$, while the character
field of $V_7$ is $\mathbb{Q}(\sqrt{13})$.}

  In fact, Hecke's arithmetical realization is in the spirit of the topological
ideas of Riemann and Poincar\'{e}, i.e., the concept of genus. Our arithmetical
realization is in the spirit of the algebraic ideas of Hilbert, i.e., the
concept of ideals. Both of these two approaches can be unified in the context
of cohomology groups of a projective or affine variety with values in a
coherent algebraic sheaf $\mathcal{O}_{X(13)}$. Namely, by the Serre duality
theorem for projective varieties,
$$H^0(X(13), \Omega_{X(13)}^1) \cong H^1(X(13), \mathcal{O}_{X(13)}).$$
Therefore, (1.21) and (1.23) in Theorem 1.10 can be rewritten in a unified
way:
$$\left\{\aligned
  &H^0(Y, \mathcal{O}_Y)=\mathbb{C}[z_1, z_2, z_3, z_4, z_5, z_6]/I(Y),\\
  &\text{(as coherent sheaves on affine varieties)},\\
  &\text{where $V_{I(Y)}=V_1 \oplus V_7 \oplus V_{13}$}.\\
  &H^0(X(13), \mathcal{O}_{X(13)})=\mathbb{C},\\
  &H^1(X(13), \mathcal{O}_{X(13)})=V_{14} \oplus V_{12}^{(1)} \oplus
                                   V_{12}^{(2)} \oplus V_{12}^{(3)},\\
  &\text{(as coherent sheaves on projective varieties)}\\
\endaligned\right.\eqno{(1.24)}$$
with
$$\dim V_{I(Y)}=21, \dim H^0(X(13), \mathcal{O}_{X(13)})=1,
  \dim H^1(X(13), \mathcal{O}_{X(13)})=50.$$
Furthermore, (1.24) exhausts all of the cohomology groups of a projective or
affine variety with values in a coherent algebraic sheaf $\mathcal{O}_{X(13)}$
corresponding to the same modular curve $X(13)$. In particular, Theorem 1.10
shows that there are two distinct types of arithmetic corresponding to the same
modular curve $X(13)$: $7$ and $13$, i.e., one is the $7$-th cyclotomic field,
the other is the $13$-th cyclotomic field. They are distinguished by their
character fields $\mathbb{Q}(\zeta_7+\zeta_7^{-1})$ and $\mathbb{Q}(\sqrt{13})$,
i.e. two distinct algebraic number fields. Now, we give some of the character
table of $\text{PSL}(2, 13)$ which are associated with theses two distinct
algebraic number fields (see \cite{CC}).

$$\text{Table $1$. The character table of $\text{PSL}(2, 13)$
        associated with $\mathbb{Q}(\sqrt{13})$}$$
$$\begin{matrix}
            & 1A & 13A & 13B\\
     \chi_2 &  7 & \frac{1-\sqrt{13}}{2} & \frac{1+\sqrt{13}}{2}\\
     \chi_3 &  7 & \frac{1+\sqrt{13}}{2} & \frac{1-\sqrt{13}}{2}\\
  \chi_{10} &  6 & \frac{-1+\sqrt{13}}{2} & \frac{-1-\sqrt{13}}{2}\\
  \chi_{11} &  6 & \frac{-1-\sqrt{13}}{2} & \frac{-1+\sqrt{13}}{2}
\end{matrix}$$

$$\text{Table $2$. The character table of $\text{PSL}(2, 13)$
        associated with $\mathbb{Q}(\zeta_7+\zeta_7^{-1})$}$$
$$\begin{matrix}
            & 1A & 7A & 7B & 7C\\
     \chi_4 & 12 & -(\zeta_7+\zeta_7^{-1}) & -(\zeta_7^2+\zeta_7^{-2}) &
                   -(\zeta_7^4+\zeta_7^{-4})\\
     \chi_5 & 12 & -(\zeta_7^4+\zeta_7^{-4}) & -(\zeta_7+\zeta_7^{-1}) &
                   -(\zeta_7^2+\zeta_7^{-2})\\
     \chi_6 & 12 & -(\zeta_7^2+\zeta_7^{-2}) & -(\zeta_7^4+\zeta_7^{-4}) &
                   -(\zeta_7+\zeta_7^{-1})\\
  \chi_{12} & 12 &  -(\zeta_7+\zeta_7^{-1}) & -(\zeta_7^2+\zeta_7^{-2}) &
                   -(\zeta_7^4+\zeta_7^{-4})\\
  \chi_{13} & 12 &  -(\zeta_7^4+\zeta_7^{-4}) & -(\zeta_7+\zeta_7^{-1}) &
                   -(\zeta_7^2+\zeta_7^{-2})\\
  \chi_{14} & 12 & -(\zeta_7^2+\zeta_7^{-2}) & -(\zeta_7^4+\zeta_7^{-4}) &
                   -(\zeta_7+\zeta_7^{-1})
\end{matrix}$$

  The significance of Theorem 1.10 also comes from that the irreducible
representations appearing in Theorem 1.10 exhaust all of the types of the
irreducible representations of $\text{SL}(2, 13)$. In fact, the trivial
representation of dimension $1$, the degenerate principal series representation
of dimension $7$ and the Steinberg representation of dimension $13$ are
induced from the degenerate discrete series of dimension $6$ by the action
on the $21$ quartic equations, i.e. the locus $\mathcal{L}$ for the modular
curve $X(13)$.

  Let us recall some basic facts about the modular curve $X(N)$ which classifies
elliptic curves with full level $N$ structure. Here, a level $N$ structure of an
elliptic curve $E$ is an isomorphism from $\mathbb{Z}/N \mathbb{Z} \times
\mathbb{Z}/N \mathbb{Z}$ to the group of $N$-torsion points on $E$. Moreover,
it satisfies the additional condition that the two base points we choose map
to a certain $N$-th root of unity under the Weil pairing. Hence, this curve
(regarding as a scheme) is defined over $\mathbb{Q}(\zeta_N)$. More precisely,
we consider the moduli problem associated to the congruence subgroup $\Gamma(N)$.
If $\gamma \in \Gamma(N)$ and $\tau \in \mathbb{H}/\Gamma(N)$, then one checks
that the points $1/N$ and $\tau/N$ in $\mathbb{C}/(\mathbb{Z}+\mathbb{Z} \tau)$
remain invariant under the action of $\gamma$. Thus associated to a point of
$\mathbb{H}/\Gamma(N)$ is an elliptic curve $\mathbb{C}/(\mathbb{Z}+\mathbb{Z}
\tau)$, together with a basis $\{ 1/N, \tau/N \}$ for the group of $N$-torsion
points. However, a point of $\mathbb{H}/\Gamma(N)$ contains one further piece
of information. Recall that there is a pairing $e_N$ on the group of $N$-torsion
points of an elliptic curve. Then one can check that
$$e_N \left(\frac{1}{N}, \frac{\tau}{N}\right)=e^{\frac{2 \pi i}{N}}.$$
Thus not only do we get a basis for the $N$-torsion, but also the two points
making up that basis pair, via the Weil pairing, to a specific primitive $N$-th
root of unity.

  In conclusion, there is a smooth projective curve $X(N)/\mathbb{Q}(\zeta_N)$
and a complex analytic isomorphism
$$j_N: \overline{\mathbb{H}}/\Gamma(N) \longrightarrow X(N)(\mathbb{C})$$
such that the following holds: Let $\tau \in \mathbb{H}/\Gamma(N)$, and
let $K=\mathbb{Q}(\zeta_N, j_N(\tau))$. $\tau$ corresponds to an equivalence
class of triples $(E, T_1, T_2)$, where $E$ is an elliptic curve, and $\{ T_1,
T_2 \}$ are generators for $E(N)$ satisfying $e_N(T_1, T_2)=\zeta_N$. Here
$e_N$ is the Weil pairing. Then this equivalence class contains a triple such
that $E$ is defined over $K$ and $T_1, T_2 \in E(K)$.

  The curve $Y(N)$ as a moduli variety: the map
$$\aligned
  \mathbb{H} &\longrightarrow \mathcal{E}_N(\mathbb{C})\\
  z &\mapsto (\mathbb{C}/\Lambda(z, 1), \quad \text{$(z/N, 1/N)$ mod $\Lambda(z, 1)$})
\endaligned$$
induces a bijection $\Gamma(N) \backslash \mathbb{H} \longrightarrow
\mathcal{E}_N(\mathbb{C})$. Let $k$ be a field containing $\mathbb{Q}(\zeta_N)$,
where $\zeta_N$ is a primitive $N$-th root of unity. The moduli problem
$\mathcal{E}_N$ has a solution over $k$. When $k=\mathbb{C}$, $M$ is canonically
isomorphic to $Y(N)_{\mathbb{C}}$($=X(N)_{\mathbb{C}}$ with the cusps removed).
Let $M$ be the solution to the moduli problem $\mathcal{E}_N$ over
$\mathbb{Q}(\zeta_N)$, then $M$ has good reduction at the prime ideals not
dividing $N$.

  Let $F_N$ be the modular function field of level $N$, i.e. the function
field of the modular curve $X(N)$ over the field of $N$-th root of unity
over $\mathbb{Q}$. The algebraic closure of $\mathbb{Q}$ in $F_N$ is
$\mathbb{Q}(\zeta_N)$.

  In fact, as Deligne pointed out to the author (see \cite{De1} and \cite{De3}),
$X(N)$ can always be defined over $\mathbb{Q}$. His proof is as follows (see
\cite{De3}): using the map
$$\aligned
  &\alpha: \mathbb{Z}/N \mathbb{Z} \times \mu_N \stackrel{\sim}{\longrightarrow} E_N,
   \quad \text{with for the Weil pairing}\\
  &(\alpha(1 \in \mathbb{Z}/N \mathbb{Z}), \alpha(\zeta \in \mu_N))=\zeta,
   \quad \text{for any $\zeta \in \mu_N$}
\endaligned$$
which does make sense over $\mathbb{Q}$. There is the other description:
we consider first
$$\alpha: \mathbb{Z}/N \mathbb{Z} \times \mathbb{Z}/N \mathbb{Z}
  \stackrel{\sim}{\longrightarrow} E_N \quad \text{(no restriction)}.$$
Over $\mathbb{C}$, we get a disconnected moduli scheme, with connected
component indexed by the primitive $N$-th roots of unity. This disconnected
scheme $\mathcal{M}_N$ is defined over $\mathbb{Q}$, and one has
$$\mathcal{M}_N \longrightarrow \text{Spec($N$-cyclotomic)} \longrightarrow
  \text{Spec}(\mathbb{Q}).$$
The scheme $\mathcal{M}_N$ is acted by the group $\text{GL}(2, \mathbb{Z}/N \mathbb{Z})$,
with an action of $\gamma \in \text{GL}(2, \mathbb{Z}/N \mathbb{Z})$ on
Spec($N$-cyclotomic) given by $\det(\gamma)$. The $\mathbb{Q}$-form
of $X(N)$ that we consider can also be described as
$$\text{$\mathcal{M}_N/$ subgroup $\left(\begin{matrix} 1 & 0\\
  0 & * \end{matrix}\right)$}.$$
In other words, the level structure considered is
$$\alpha: \mathbb{Z}/N \mathbb{Z} \times \mathbb{Z}/N \mathbb{Z}
  \stackrel{\sim}{\longrightarrow} E_N$$
given up to $\left(\begin{matrix} 1 & 0\\ 0 & *
\end{matrix}\right)$. This form is a scheme over $\mathbb{Q}$ which over
$\mathbb{C}$ becomes $X(N)$. However, the action of $\text{SL}(2, \mathbb{Z}/N \mathbb{Z})$
which we have over $\mathbb{C}$ is not defined over $\mathbb{Q}$, but only over
the cyclotomic field. In our case, $N=13$, the $21$ quartic equations given by
(1.12), (1.13) and (1.14) have rational coefficients, but the action of some
elements of $\text{SL}(2, 13)$ is given by matrices with the $13$-th cyclotomic
entries in the decomposition formula (1.23).

  Note that the simple group $\text{PSL}(2, 13)$ of order $1092$ corresponds to
four different curves:

\noindent (1) the modular curve $X(13)$ or its algebraic model $Y \subset \mathbb{P}^5$
with genus $50$ defined over $\mathbb{Q}$, where $\text{PSL}(2, 13) \cong (2, 3, 13)$;

\noindent (2) the first Hurwitz triplet $X_1$, $X_2$ and $X_3$, which can
be realized as the Shimura curves whose levels are with norm $13$ or realized
as three non-congruence modular curves of level $7$ (see \cite{Y1}) or their
algebraic models $Y_i \subset \mathbb{P}^{13}$ $(i=1, 2, 3)$ defined over
$\mathbb{Q}(\zeta_7+\zeta_7^{-1})$ with genus $14$, where $\text{PSL}(2, 13)
\cong (2, 3, 7; 6)$, $(2, 3, 7; 7)$ or $(2, 3, 7; 13)$, respectively.

  For $X(13)$ or $Y$, we have the decomposition formulas (1.21) and (1.23).
However, for $X_i$ or $Y_i$ $(i=1, 2, 3)$, one only has
$$H^0(X_i, \Omega_{X_i}^1)=V_{14}, \quad (i=1, 2, 3),$$
which only involves the principal series representation $V_{14}$ defined over
$\mathbb{Q}$. In particular, the character field of $V_{14}$ is $\mathbb{Q}$.
Moreover, the three explicit algebraic models for $Y_i$ (i.e., the three
defining ideals $I(Y_i)$) $(i=1, 2, 3)$ defined over $\mathbb{Q}(\zeta_7+\zeta_7^{-1})$
are still unknown. In this sense, the modular curve $X(13)$ has much richer
structure than the first Hurwitz triplet $X_i$ $(i=1, 2, 3)$.

  In fact, we have constructed four kinds of representations of $\text{SL}(2, 13)$
corresponding to the following bases:

(1) $(z_1, \ldots, z_6)$: discrete series of dimension six

(2) $(\mathbf{A}_0, \mathbf{A}_1, \ldots, \mathbf{A}_6)$: principal series
    of dimension seven (see \cite{Y2})

(3) $(\mathbf{D}_0, \mathbf{D}_1, \ldots, \mathbf{D}_{12}, \mathbf{D}_{\infty})$:
    principal series of dimension fourteen (see \cite{Y2})

(4) $(\mathbf{B}_0^{(0)}, \mathbf{B}_0^{(1)}, \mathbf{B}_0^{(2)}, \ldots, \mathbf{B}_7)$:
    reducible representation of dimension $21$

  The significance of our $21$-dimensional reducible representation comes
from that it is closely related to the construction of the curve $Y$
(invariant theory and algebraic geometry) and modularity as well as
uniformization of this curve (number theory and arithmetic geometry).
Consequently, the $21$-dimensional reducible projective representation of
$\text{PSL}(2, 13)$ is obtained from $21$-dimensional ordinary reducible
representation of $\text{SL}(2, 13)$.

  Let $\mathfrak{a}_i$ be the ideals generated by the basis of $V_{6i-5}$
$(i=1, 2, 3)$ and $Y_i$ be the varieties corresponding to the ideals
$\mathfrak{a}_i$. Theorem 1.9 implies that $\mathfrak{a}_i$ $(i=1, 2, 3)$
are $\text{SL}(2, 13)$-invariant ideals and $Y_i$ $(i=1, 2, 3)$ are
$\text{SL}(2, 13)$-invariant varieties. In particular, the decomposition
(1.20) and (1.21) leads to the corresponding geometric construction of
the curve $Y$:

\textbf{Theorem 1.11.} (Main Theorem 6) {\it The ideal $I(Y)$ can be decomposed
as the sum of ideals $\mathfrak{a}_i$ $(i=1, 2, 3)$. The curve $Y$ can be
expressed as the intersection of $Y_1$, $Y_2$ and $Y_3$:
$$I(Y)=\mathfrak{a}_1+\mathfrak{a}_2+\mathfrak{a}_3.\eqno{(1.25)}$$
$$Y=Y_1 \cap Y_2 \cap Y_3.\eqno{(1.26)}$$
In particular, the curve $Y_2$ is in a $\text{SL}(2, 13)$-invariant sextic
four-fold $W_2 \subset \mathbb{P}^5$, i.e., a Calabi-Yau four-fold, given by an
invariant homogeneous polynomial $\Psi_{1, 1}=0$, where
$$\aligned
  \Psi_{1, 1} =&2 \mathbf{A}_0 (\mathbf{B}_0^{(1)}+\mathbf{B}_0^{(2)})
                +\mathbf{A}_1 (\mathbf{B}_{12}^{(1)}+3 \mathbf{B}_{12}^{(2)})+\\
               &+\mathbf{A}_4 (\mathbf{B}_{10}^{(1)}+3 \mathbf{B}_{10}^{(2)})
                +\mathbf{A}_3 (\mathbf{B}_4^{(1)}+3 \mathbf{B}_4^{(2)})+\\
               &+\mathbf{A}_5 (\mathbf{B}_1^{(1)}-3 \mathbf{B}_1^{(2)})
                +\mathbf{A}_6 (\mathbf{B}_3^{(1)}-3 \mathbf{B}_3^{(2)})+\\
               &+\mathbf{A}_2 (\mathbf{B}_9^{(1)}-3 \mathbf{B}_9^{(2)})
                \in \mathfrak{a}_2.
\endaligned\eqno{(1.27)}$$
The curve $Y_3$ is in a $\text{SL}(2, 13)$-invariant octic four-fold $W_3
\subset \mathbb{P}^5$, i.e., a general type four-fold, given by an invariant
homogeneous polynomial $\Omega_{1, 1}=0$, where
$$\aligned
  &\Omega_{1, 1}=(4 \mathbf{B}_0^{(0)}-\mathbf{B}_0^{(1)}+\mathbf{B}_0^{(2)})
                 (6 \mathbf{A}_0^2-\mathbf{A}_1 \mathbf{A}_5-\mathbf{A}_2 \mathbf{A}_3
                  -\mathbf{A}_4 \mathbf{A}_6)+\\
                &+7 \mathbf{B}_5 (\mathbf{A}_2^2+2 \mathbf{A}_3 \mathbf{A}_5)
                 +7 \mathbf{B}_2 (\mathbf{A}_5^2+2 \mathbf{A}_1 \mathbf{A}_6)
                 +7 \mathbf{B}_6 (\mathbf{A}_6^2+2 \mathbf{A}_4 \mathbf{A}_2)+\\
                &+7 \mathbf{B}_8 (\mathbf{A}_3^2+2 \mathbf{A}_1 \mathbf{A}_2)+
                 +7 \mathbf{B}_{11} (\mathbf{A}_1^2+2 \mathbf{A}_4 \mathbf{A}_5)
                 +7 \mathbf{B}_7 (\mathbf{A}_4^2+2 \mathbf{A}_3 \mathbf{A}_6)+\\
                &+7 (\mathbf{B}_1^{(1)}+\mathbf{B}_1^{(2)}) (\mathbf{A}_0
                 \mathbf{A}_5+\mathbf{A}_3 \mathbf{A}_4)
                 +7 (\mathbf{B}_3^{(1)}+\mathbf{B}_3^{(2)}) (\mathbf{A}_0
                 \mathbf{A}_6+\mathbf{A}_1 \mathbf{A}_3)+\\
                &+7 (\mathbf{B}_9^{(1)}+\mathbf{B}_9^{(2)}) (\mathbf{A}_0
                 \mathbf{A}_2+\mathbf{A}_1 \mathbf{A}_4)
                 +7 (-\mathbf{B}_{12}^{(1)}+\mathbf{B}_{12}^{(2)}) (\mathbf{A}_0
                 \mathbf{A}_1+\mathbf{A}_2 \mathbf{A}_6)+\\
                &+7 (-\mathbf{B}_{10}^{(1)}+\mathbf{B}_{10}^{(2)}) (\mathbf{A}_0
                 \mathbf{A}_4+\mathbf{A}_2 \mathbf{A}_5)
                 +7 (-\mathbf{B}_4^{(1)}+\mathbf{B}_4^{(2)}) (\mathbf{A}_0
                 \mathbf{A}_3+\mathbf{A}_5 \mathbf{A}_6)\\
                &\in \mathfrak{a}_3.
\endaligned\eqno{(1.28)}$$}

  Theorem 1.6, Theorem 1.7, Theorem 1.8, Theorem 1.9, Theorem 1.10 and
Theorem 1.11 give an answer on the construction of the curve $Y$, i.e.
Problem 1.5. In particular, two invariant curves $Y_2$ and $Y_3$ are
lying over the curve $Y$.

  The structure of the homogeneous coordinate ring of a canonical curve is
known for genus $g \leq 5$. For example, when $g=4$, $X_6=(2) \cap (3)
\subset \mathbb{P}^4$. When $g=5$, $X_8=(2) \cap (2) \cap (2) \subset
\mathbb{P}^5$. When the genus is higher some difficulties arise. In \cite{Mu1},
\cite{Mu2} and \cite{Mu3}, the cases $g=7$, $8$ and $9$ are treated by means
of the consideration of the canonical curve as a transversal linear section
of a symmetric space. Despite of this, very little is known for curves with
much higher genus. Theorem 1.11 provides an example of the construction of
the curve with genus $50$. Moreover, the decomposition (1.26) has the
invariance under the action of groups.

  In particular, Theorem 1.11 provides three distinct realizations of the
modular curve $X(13)$:

(1) in the space $V$ of the $21$-dimensional reducible representation;

(2) in the space $V_7$ of the degenerate principal series;

(3) in the space $V_{13}$ of the Steinberg representation.

  Theorem 1.11 also provides an example of Fano four-folds $Y_1$,
an example of Calabi-Yau four-folds $W_2$ as well as an example of
general type four-folds $W_3$ which admit the same (large) automorphism
group $\text{SL}(2, 13)$.

  It is well-known that the geometric realization of representations often
reveals deeper layers of the structure in the form of representation theory.
Klein showed that the modular curve $X(p)$ could be $\text{PSL}(2, \mathbb{F}_p)$
equivariantly mapped to the projective space $\mathbb{P}(V^{+})$ associated to
the even part $V^{+}$ of the Weil representation of $\text{SL}(2, \mathbb{F}_p)$.
The resulting curve, which is called the $\mathbf{A}$-curve of level $p$ and
denote $\mathbf{A}(p)$, has the same genus as $X(p)$. Denote by $\mathbf{A}_0$,
$\mathbf{A}_1$, $\ldots$, $\mathbf{A}_{\frac{p-1}{2}}$ the coordinates of the
$\mathbf{A}(p)$ and by $z_1$, $\ldots$, $z_{\frac{p-1}{2}}$ those of the
$z$-curve in $\mathbb{P}(V^{-})$. In particular, the basis $\mathbf{A}_0$,
$\mathbf{A}_1$, $\ldots$, $\mathbf{A}_{\frac{p-1}{2}}$ provides the standard
geometric realization of the degenerate principal series representations of
$\text{SL}(2, \mathbb{F}_p)$ with dimension $\frac{p+1}{2}$ by a concrete
vector space with a concrete $\text{SL}(2, \mathbb{F}_p)$-action. Now, we
give a non-standard geometric realization of degenerate principal series
for $\text{SL}(2, 13)$ by a concrete vector space corresponding to the curve
$Y_2$ lying over $Y$ with a concrete $\text{SL}(2, 13)$-action on this curve.
This provides a new viewpoint to the geometric realization of degenerate
principal series for $\text{SL}(2, 13)$ which is distinct from the Klein
$\mathbf{A}$-basis realization. Note that our realization has a rich
geometric structure of invariant algebraic curves, while the standard
realization only has the structure of vector spaces.

\textbf{Theorem 1.12.} (Main Theorem 7) {\it The degenerate principal series
of $\text{SL}(2, 13)$ has two distinct geometric realizations: The standard
realization is given by the basis $(\mathbf{A}_0, \mathbf{A}_1$, $\mathbf{A}_2$,
$\mathbf{A}_3$, $\mathbf{A}_4$, $\mathbf{A}_5$, $\mathbf{A}_6)$ corresponding
to the Klein $\mathbf{A}$-curve. The non-standard realization is given by
the basis $(\mathbf{B}_0^{(1)}+\mathbf{B}_0^{(2)}$, $\mathbf{B}_1^{(1)}-3
\mathbf{B}_1^{(2)}$, $\mathbf{B}_3^{(1)}-3 \mathbf{B}_3^{(2)}$,
$\mathbf{B}_9^{(1)}-3 \mathbf{B}_9^{(2)}$, $\mathbf{B}_{12}^{(1)}+3
\mathbf{B}_{12}^{(2)}$, $\mathbf{B}_{10}^{(1)}+3 \mathbf{B}_{10}^{(2)}$,
$\mathbf{B}_4^{(1)}+3 \mathbf{B}_4^{(2)})$ corresponding to the curve $Y_2$
lying over $Y$. In particular,
$$\left\{\aligned
   \sqrt{13} ST^{\nu}(\mathbf{A}_0)&=\mathbf{A}_0+\zeta^{\nu} \mathbf{A}_1
   +\zeta^{3 \nu} \mathbf{A}_4+\zeta^{9 \nu} \mathbf{A}_3\\
  &+\zeta^{12 \nu} \mathbf{A}_5+\zeta^{10 \nu} \mathbf{A}_6+\zeta^{4 \nu}
   \mathbf{A}_2,\\
   \sqrt{13} ST^{\nu}(\mathbf{B}_0^{(1)}+\mathbf{B}_0^{(2)})
  &=(\mathbf{B}_0^{(1)}+\mathbf{B}_0^{(2)})
  +\zeta^{\nu} (\mathbf{B}_1^{(1)}-3 \mathbf{B}_1^{(2)})\\
  &+\zeta^{3 \nu} (\mathbf{B}_3^{(1)}-3 \mathbf{B}_3^{(2)})
  +\zeta^{9 \nu} (\mathbf{B}_9^{(1)}-3 \mathbf{B}_9^{(2)})\\
  &+\zeta^{12 \nu} (\mathbf{B}_{12}^{(1)}+3 \mathbf{B}_{12}^{(2)})
  +\zeta^{10 \nu} (\mathbf{B}_{10}^{(1)}+3 \mathbf{B}_{10}^{(2)})\\
  &+\zeta^{4 \nu} (\mathbf{B}_4^{(1)}+3 \mathbf{B}_4^{(2)}).
\endaligned\right.\eqno{(1.29)}$$}

  Similarly, the curve $Y_3$ over $Y$ gives a geometric realization of the
Steinberg representation for $\text{SL}(2, 13)$. In particular, comparing
with (1.29), we have
$$\aligned
  &13S T^{\nu}(4 \mathbf{B}_0^{(0)}-\mathbf{B}_0^{(1)}+\mathbf{B}_0^{(2)})\\
 =&-(4 \mathbf{B}_0^{(0)}-\mathbf{B}_0^{(1)}+\mathbf{B}_0^{(2)})\\
  &+14 \zeta^{5 \nu} \mathbf{B}_5
   +14 \zeta^{2 \nu} \mathbf{B}_2
   +14 \zeta^{6 \nu} \mathbf{B}_6
   +14 \zeta^{8 \nu} \mathbf{B}_8
   +14 \zeta^{11 \nu} \mathbf{B}_{11}
   +14 \zeta^{7 \nu} \mathbf{B}_7\\
  &+7 \zeta^{\nu} (\mathbf{B}_1^{(1)}+\mathbf{B}_1^{(2)})
   +7 \zeta^{3 \nu} (\mathbf{B}_3^{(1)}+\mathbf{B}_3^{(2)})
   +7 \zeta^{9 \nu} (\mathbf{B}_9^{(1)}+\mathbf{B}_9^{(2)})\\
  &+7 \zeta^{12 \nu} (-\mathbf{B}_{12}^{(1)}+\mathbf{B}_{12}^{(2)})
   +7 \zeta^{10 \nu} (-\mathbf{B}_{10}^{(1)}+\mathbf{B}_{10}^{(2)})
   +7 \zeta^{4 \nu} (-\mathbf{B}_4^{(1)}+\mathbf{B}_4^{(2)}).
\endaligned\eqno{(1.30)}$$

\textbf{Theorem 1.13.} (Main Theorem 8) (Modularity for the ideal $I(Y)$)
{\it The invariant ideal $I(Y)$ can be parameterized by theta constants
of order $13$, i.e., there are the following twenty-one modular equations
of order $13$:
$$\left\{\aligned
  \mathbf{B}_0^{(i)}(a_1(z), \ldots, a_6(z)) &=0,\\
  \mathbf{B}_1^{(j)}(a_1(z), \ldots, a_6(z)) &=0,\\
  \mathbf{B}_3^{(j)}(a_1(z), \ldots, a_6(z)) &=0,\\
  \mathbf{B}_9^{(j)}(a_1(z), \ldots, a_6(z)) &=0,\\
  \mathbf{B}_{12}^{(j)}(a_1(z), \ldots, a_6(z)) &=0,\\
  \mathbf{B}_{10}^{(j)}(a_1(z), \ldots, a_6(z)) &=0,\\
  \mathbf{B}_4^{(j)}(a_1(z), \ldots, a_6(z)) &=0,
\endaligned\right.\eqno{(1.31)}$$
$$\left\{\aligned
  \mathbf{B}_5(a_1(z), \ldots, a_6(z)) &=0,\\
  \mathbf{B}_2(a_1(z), \ldots, a_6(z)) &=0,\\
  \mathbf{B}_6(a_1(z), \ldots, a_6(z)) &=0,\\
  \mathbf{B}_8(a_1(z), \ldots, a_6(z)) &=0,\\
  \mathbf{B}_{11}(a_1(z), \ldots, a_6(z)) &=0,\\
  \mathbf{B}_7(a_1(z), \ldots, a_6(z)) &=0,
\endaligned\right.\eqno{(1.31)}$$
where $i=0, 1, 2$ and $j=1, 2$.}

\textbf{Corollary 1.14.} {\it The space curve $Y$ is modular, i.e., it
is parameterized by theta constants of order $13$. Two curves $Y_2$ and
$Y_3$ lying over $Y$ have modular components given by (1.31).}

  Corollary 1.14 provides an example of modularity for higher genus space
curve $(g \geq 4)$ by means of the explicit modular parametrization.

\textbf{Corollary 1.15.} {\it The space curve $Y$ admits a hyperbolic
uniformization of arithmetic type with respect to the congruence
subgroup $\Gamma(13)$.}

  Corollary 1.15 gives an example of hyperbolic uniformization of arithmetic
type with respect to a congruence subgroup of $\text{SL}(2, \mathbb{Z})$
as well as an example of the explicit uniformization of algebraic space
curves of higher genus. This gives an answer to Problem 1.4.

\textbf{Corollary 1.16.} {\it The $21$ quartic polynomials in six variables
given by (1.32) form a system of over-determined algebraic equations which
can be uniformized by theta constants of order $13$.}

  Corollary 1.16 gives a new example for Hilbert's 22nd problem, i.e., an
explicit uniformization for algebraic relations among five variables by means
of automorphic functions.

  It is known that Ramanujan (see \cite{R1}, p.326, \cite{R2}, p.244 and
\cite{B}, p.372) constructed some modular equations of order $13$, most of
these results show Ramanujan at his very best. Evans (see \cite{Ev}) gave
beautiful proofs of some of Ramanujan's modular equations of order $13$.
Moreover, inspired by Ramanujan's work, Evans proved a beautiful formula
in the spirit of Ramanujan (see \cite{Ev} or \cite{B}, p.375-376):
$$\frac{1}{a_1(z) a_4(z)}+\frac{1}{a_2(z) a_5(z)}+\frac{1}{a_3(z) a_6(z)}=0.
  \eqno{(1.32)}$$
We show that, in fact, there exist $21$ such kinds of formulas, all of
which have a natural geometric interpretation in terms of the projective
model of the modular curve $X(13)$ in $\mathbb{P}^5$ and as a whole, they
are invariant under the action of $G \cong \text{SL}(2, 13)$.

\textbf{Theorem 1.17.} (Main Theorem 9) {\it Besides $(1.32)$, there are
twenty formulas given by $(1.31)$ in the spirits of $(1.32)$ as follows:

\noindent (1) Quadratic fractions:
$$\mathbf{B}_0^{(0)}: \quad
  \frac{1}{a_1(z) a_4(z)}+\frac{1}{a_2(z) a_5(z)}+\frac{1}{a_3(z) a_6(z)}=0.$$
(2) Cubic fractions:
$$\left\{\aligned
  \mathbf{B}_1^{(2)}: \quad
 &\frac{a_6(z)}{a_1(z) a_3(z) a_5(z)}-\frac{a_3(z)}{a_1(z) a_2(z) a_5(z)}-
  \frac{a_1(z)}{a_3(z) a_4(z) a_6(z)}=0,\\
  \mathbf{B}_3^{(2)}: \quad
 &\frac{a_5(z)}{a_2(z) a_3(z) a_4(z)}-\frac{a_2(z)}{a_1(z) a_3(z) a_4(z)}-
  \frac{a_3(z)}{a_2(z) a_5(z) a_6(z)}=0,\\
  \mathbf{B}_9^{(2)}: \quad
 &\frac{a_4(z)}{a_1(z) a_2(z) a_6(z)}-\frac{a_1(z)}{a_2(z) a_3(z) a_6(z)}-
  \frac{a_2(z)}{a_1(z) a_4(z) a_5(z)}=0.
\endaligned\right.\eqno{(1.33)}$$
$$\left\{\aligned
  \mathbf{B}_{12}^{(2)}: \quad
 &\frac{a_4(z)}{a_1(z) a_3(z) a_6(z)}-\frac{a_3(z)}{a_2(z) a_4(z) a_6(z)}-
  \frac{a_6(z)}{a_2(z) a_4(z) a_5(z)}=0,\\
  \mathbf{B}_{10}^{(2)}: \quad
 &\frac{a_6(z)}{a_2(z) a_3(z) a_5(z)}-\frac{a_2(z)}{a_1(z) a_5(z) a_6(z)}-
  \frac{a_5(z)}{a_1(z) a_4(z) a_6(z)}=0,\\
  \mathbf{B}_4^{(2)}: \quad
 &\frac{a_5(z)}{a_1(z) a_2(z) a_4(z)}-\frac{a_1(z)}{a_3(z) a_4(z) a_5(z)}-
  \frac{a_4(z)}{a_3(z) a_5(z) a_6(z)}=0.
\endaligned\right.\eqno{(1.34)}$$
$$\left\{\aligned
  \mathbf{B}_5: \quad
 &-\frac{a_2(z)}{a_3(z) a_4(z) a_6(z)}+\frac{a_6(z)}{a_1(z) a_2(z) a_3(z)}+
  \frac{a_4(z)}{a_1(z) a_5(z) a_6(z)}=0,\\
  \mathbf{B}_2: \quad
 &-\frac{a_1(z)}{a_2(z) a_5(z) a_6(z)}+\frac{a_5(z)}{a_1(z) a_2(z) a_3(z)}+
  \frac{a_6(z)}{a_3(z) a_4(z) a_5(z)}=0,\\
  \mathbf{B}_6: \quad
 &-\frac{a_3(z)}{a_1(z) a_4(z) a_5(z)}+\frac{a_4(z)}{a_1(z) a_2(z) a_3(z)}+
  \frac{a_5(z)}{a_2(z) a_4(z) a_6(z)}=0.
\endaligned\right.\eqno{(1.35)}$$
$$\left\{\aligned
  \mathbf{B}_8: \quad
 &\frac{a_5(z)}{a_1(z) a_3(z) a_6(z)}+\frac{a_3(z)}{a_4(z) a_5(z) a_6(z)}+
  \frac{a_1(z)}{a_2(z) a_3(z) a_4(z)}=0,\\
  \mathbf{B}_{11}: \quad
 &\frac{a_4(z)}{a_2(z) a_3(z) a_5(z)}+\frac{a_2(z)}{a_4(z) a_5(z) a_6(z)}+
  \frac{a_3(z)}{a_1(z) a_2(z) a_6(z)}=0,\\
  \mathbf{B}_7: \quad
 &\frac{a_6(z)}{a_1(z) a_2(z) a_4(z)}+\frac{a_1(z)}{a_4(z) a_5(z) a_6(z)}+
  \frac{a_2(z)}{a_1(z) a_3(z) a_5(z)}=0.
\endaligned\right.\eqno{(1.36)}$$
(3) Quartic fractions:
$$\left\{\aligned
  \mathbf{B}_0^{(1)}: \quad
 &\frac{a_5(z)^2}{a_2(z) a_3(z) a_4(z) a_6(z)}+
  \frac{a_6(z)^2}{a_1(z) a_3(z) a_4(z) a_5(z)}+\\
 &+\frac{a_4(z)^2}{a_1(z) a_2(z) a_5(z) a_6(z)}=0,\\
  \mathbf{B}_0^{(2)}: \quad
 &\frac{a_1(z)^2}{a_2(z) a_3(z) a_4(z) a_5(z)}+
  \frac{a_2(z)^2}{a_1(z) a_3(z) a_5(z) a_6(z)}+\\
 &+\frac{a_3(z)^2}{a_1(z) a_2(z) a_4(z) a_6(z)}=0.
\endaligned\right.\eqno{(1.37)}$$
$$\left\{\aligned
  \mathbf{B}_1^{(1)}: \quad
 &\frac{a_5(z)^2}{a_1(z) a_2(z) a_4(z) a_6(z)}+
  \frac{a_1(z)^2}{a_2(z) a_3(z) a_5(z) a_6(z)}+\\
 &-\frac{a_2(z)^2}{a_3(z) a_4(z) a_5(z) a_6(z)}=0,\\
  \mathbf{B}_3^{(1)}: \quad
 &\frac{a_4(z)^2}{a_1(z) a_3(z) a_5(z) a_6(z)}+
  \frac{a_3(z)^2}{a_1(z) a_2(z) a_4(z) a_5(z)}+\\
 &-\frac{a_1(z)^2}{a_2(z) a_4(z) a_5(z) a_6(z)}=0,\\
  \mathbf{B}_9^{(1)}: \quad
 &\frac{a_6(z)^2}{a_2(z) a_3(z) a_4(z) a_5(z)}+
  \frac{a_2(z)^2}{a_1(z) a_3(z) a_4(z) a_6(z)}+\\
 &-\frac{a_3(z)^2}{a_1(z) a_4(z) a_5(z) a_6(z)}=0.
\endaligned\right.\eqno{(1.38)}$$
$$\left\{\aligned
  \mathbf{B}_{12}^{(1)}: \quad
 &\frac{a_4(z)^2}{a_2(z) a_3(z) a_5(z) a_6(z)}+
  \frac{a_2(z)^2}{a_1(z) a_3(z) a_4(z) a_5(z)}+\\
 &+\frac{a_5(z)^2}{a_1(z) a_2(z) a_3(z) a_6(z)}=0,\\
  \mathbf{B}_{10}^{(1)}: \quad
 &\frac{a_6(z)^2}{a_1(z) a_2(z) a_4(z) a_5(z)}+
  \frac{a_1(z)^2}{a_2(z) a_3(z) a_4(z) a_6(z)}+\\
 &+\frac{a_4(z)^2}{a_1(z) a_2(z) a_3(z) a_5(z)}=0,\\
  \mathbf{B}_4^{(1)}: \quad
 &\frac{a_5(z)^2}{a_1(z) a_3(z) a_4(z) a_6(z)}+
  \frac{a_3(z)^2}{a_1(z) a_2(z) a_5(z) a_6(z)}+\\
 &+\frac{a_6(z)^2}{a_1(z) a_2(z) a_3(z) a_4(z)}=0.
\endaligned\right.\eqno{(1.39)}$$}

  Theorem 1.17 greatly improves the result of Ramanujan and Evans on
modular equations of order $13$. Moreover, the significance of these
quartic fractions come from that they can define the locus $\mathcal{L}$
of the modular curve $X(13)$ in $\mathbb{P}^5$ (see Theorem 1.19).

  In fact, combining Theorem 1.10 with Theorem 1.13, Corollary 1.14,
Corollary 1.15, Corollary 1.16 and Theorem 1.17, a major application
of the modularity of $Y$, an explicit uniformization of $Y$ as well
as the hyperbolic uniformization of arithmetic type for this higher
genus arithmetic algebraic curve $Y$, is the arithmetic of $Y$. That
is, there are two distinct arithmetical realizations whose fields of
definition are distinct algebraic number fields in the sense of
Theorem 1.10.

\textbf{Problem 1.18.} Determine minimal generators for the ideal
of identities and whether the identities determine the modular curve
$X(p)$ in $\mathbb{CP}^{\frac{p-3}{2}}$ for prime number $p \geq 5$.

  In fact, when $p=5$, $7$ and $11$, the above problem was studied
by Klein in \cite{K1}, \cite{K2} and \cite{K3}. In particular, when $p=5$,
this is trivial since the modular $X(5) \cong \mathbb{CP}^1$. When $p=7$,
the number of minimal generators for the ideal of identities is one,
this identity is just the Klein quartic curve $x^3 y+y^3 z+z^3 x=0$,
which can determine the modular curve $X(7)$. When $p=11$, the number
of minimal generators for the ideal of identities is ten, all of which
can be constructed from the Klein cubic threefold
$v^2 w+w^2 x+x^2 y+y^2 z+z^2 v=0$.
Despite of this, when $p \geq 13$, little is known, especially for the
explicit construction of the minimal generators for the ideal of identities
which can determine the modular curve $X(p)$ (see \cite{AR} and \cite{FK}).
In the present paper, we solve Problem 1.18 for $p=13$. In fact, Ramanan
(see \cite{AR}, p.59) remarked that the locus $\mathcal{L}$ of the modular
curve $X(p)$ can be defined by $(p-1)/2$ quartics. For $p=13$, this number
is $6$. However, this does not lead to explicit equations. Now, we give an
explicit construction for $p=13$:

\textbf{Theorem 1.19.} {\it The locus $\mathcal{L}$ of the modular curve
$X(13)$ can be defined by six quartics which can be given as follows:
$$(\mathbf{B}_0^{(1)}, \mathbf{B}_0^{(2)}, \mathbf{B}_3^{(1)},
   \mathbf{B}_{10}^{(1)}, \mathbf{B}_9^{(1)}, \mathbf{B}_4^{(1)})$$
for $z_6 \neq 0$ or $z_3 \neq 0$.
$$(\mathbf{B}_0^{(1)}, \mathbf{B}_0^{(2)}, \mathbf{B}_9^{(1)},
   \mathbf{B}_4^{(1)}, \mathbf{B}_1^{(1)}, \mathbf{B}_{12}^{(1)})$$
for $z_5 \neq 0$ or $z_2 \neq 0$.
$$(\mathbf{B}_0^{(1)}, \mathbf{B}_0^{(2)}, \mathbf{B}_1^{(1)},
   \mathbf{B}_{12}^{(1)}, \mathbf{B}_3^{(1)}, \mathbf{B}_{10}^{(1)})$$
for $z_4 \neq 0$ or $z_1 \neq 0$.}

  Following on foundational work, notably of Gauss, Abel and Jacobi, a basic
problem became that of linking $j(\tau)$ and $j^{\prime}(\tau)=j(N \tau)$
for $\tau \in \mathbb{H}$ and $N$ an integer $\geq 2$. It can be shown that
there exists a polynomial $\Phi_N \in \mathbb{C}[X, Y]$ such that
$$\Phi_N(j^{\prime}, j)=0.$$
When $\Phi_N$ is minimal, this is called the modular equation associated
with transformations of order $N$. The field of meromorphic
$\text{SL}(2, \mathbb{Z})$-invariant functions coincides with
$\mathbb{C}(j)$, which is isomorphic to the field of rational functions
over $\mathbb{C}$ in a single variable. Hence for every finite-index
subgroup $\Gamma$ of $\text{SL}(2, \mathbb{Z})$, the field of functions
$K(\Gamma)$ is a finite extension of $\mathbb{C}(j)$, Galois if and only if
$\Gamma$ is a normal subgroup of $\Gamma(1)$, which is the case for the
principal congruence subgroup $\Gamma(N)$, whose degree is equal to that
of the branched covering $X_{\Gamma} \rightarrow X(1)$. The splitting
field of $\Phi_N \in \mathbb{C}[j][X]$ is the function field $K(\Gamma(N))$.
In a series of work \cite{K1}, \cite{K2} and \cite{K}, Klein studied the
geometric realization of $K(\Gamma(N))$ for $N=2$, $3$, $4$, $5$, $7$ and
$13$ for which $g(X_0(N))=0$. In fact, the modular polynomial $\Phi_N(X, Y)$
is a model for the modular curve $X_0(N)$ parametrizing cyclic $N$-isogenies
between elliptic curves. However, the modular polynomial can quickly get very
complicated (see section 8 for details). It is better to parametrize the curve
$\Phi_N(X, Y)=0$ by a different modular function, i.e., the Hauptmodul $\tau$,
and then write $X$ and $Y$ in terms of this parameter:
$$X=j(\tau), \quad Y=j(N \tau).$$
Here, the function $j(Nz)$ is algebraic over $\mathbb{C}(j)$, and its minimal
polynomial $\Phi_N(X)$ has coefficients in $\mathbb{Z}[j]$. In fact, $\Phi_N(X, Y)$
is symmetric, i.e. $\Phi_N(X, Y)=\Phi_N(Y, X)$, and $\deg_X \Phi_N=N+1$. In
particular, when $N=5$, this leads to the quintic equation and the icosahedron
(see \cite{K}). When $N=7$, this leads to the Klein quartic curve (see \cite{K2}),
where arithmetic, algebraic, geometric and combinatorial facets are tightly
imbricated in this work of Klein, revealing the Klein quartic to be a central
and fascinating mathematical object. However, Klein left the last case $N=13$
open, which is a long-standing problem since his paper \cite{K1}. In the present
paper, we solve this problem by the following:

\textbf{Theorem 1.20.} (Main Theorem 10) (Galois covering, Galois resolvent
and their geometric realizations). {\it In the projective model, $X(13)$ is
isomorphic to $Y$ in $\mathbb{CP}^5$. The natural morphism from $X(13)$ onto
$X(1) \simeq \mathbb{CP}^1$ is realized as the projection of $Y$ onto $Y/G$
$($identified with $\mathbb{CP}^1$$)$. This is a Galois covering whose generic
fibre is interpreted as the Galois resolvent of the modular equation
$\Phi_{13}(\cdot, j)=0$ of level $13$, i.e., the function field of $Y$ is the
splitting field of this modular equation over $\mathbb{C}(j)$.}

  The central relation between automorphic forms and algebraic varieties
provides a passage from a context where a given critical assertion is
difficult, even impossible, to one in which, it is almost transparent.
In our case, because of the complicated expression of $Y$ by $21$ quartic
polynomials in six variables, it may be not possible to prove directly
that $Y/\text{SL}(2, 13) \cong \mathbb{CP}^1$ with the desired ramification
property, and even impossible that this is a Galois covering whose generic
fibre is the Galois resolvent of the modular equation of level $13$, i.e.,
the function field of $Y$ is the splitting field of this equation over
$\mathbb{C}(j)$. But once we know that $Y \cong X(13)$, the same conclusion
can be immediate.

  It is known that if $G$ is simple and $|G|$ is composite but a product
of at most $5$ prime numbers, then $|G|$ is $60$, $168$, $660$ or $1092$
and in fact $G$ is $\text{PSL}(2, p)$ for $p=5$, $7$, $11$, $13$, respectively.
Amongst all groups, whose order is a product of $5$ prime numbers, there
are only $3$ simple groups, consisting of the non-singular fractional
linear substitutions with respect to a prime modulus $p$, of order
$\frac{1}{2} p (p^2-1)$ for $p=7$, $11$ and $13$, hence of order
$$168=2^3 \cdot 3 \cdot 7, \quad
  660=2^2 \cdot 3 \cdot 5 \cdot 11, \quad
  1092=2^2 \cdot 3 \cdot 7 \cdot 13.$$
Now, we have a complete picture on the trilogy ($X_0(5)$, $X_0(7)$,
$X_0(13)$):

(1) $\text{SL}(2, 5)$: $g(X_0(5))=0$.

  In the modular equation $\Phi_5(X, Y)=0$ of level five:
$$X=j(z)=-\frac{(\tau^2-250 \tau+3125)^3}{\tau^5}
        =-\frac{({\tau^{\prime}}^2-10 \tau^{\prime}+5)^3}{\tau^{\prime}},$$
$$Y=j(5z)=-\frac{(\tau^2-10 \tau+5)^3}{\tau},$$
where the Hauptmodul is given by
$$\tau=\left(\frac{\eta(z)}{\eta(5 z)}\right)^6,$$
and the Fricke involution $\tau^{\prime}=\tau(-\frac{1}{5z})$ satisfies
that $\tau \tau^{\prime}=125$.
$$\mathbb{C}^2/\text{SL}(2, 5) \cong \text{Spec}\left(
  \mathbb{C}[z_1, z_2]^{\text{SL}(2, 5)}\right).\eqno{(1.40)}$$
$$\mathbb{P}^1 \rightarrow \mathbb{P}^1/\text{SL}(2, 5) \cong \mathbb{P}^1.
  \eqno{(1.41)}$$
Galois resolvent of the modular equation $\Phi_5(\cdot, j)=0$ of level five.

(2) $\text{PSL}(2, 7)$: $g(X_0(7))=0$.

  In the modular equation $\Phi_7(X, Y)=0$ of level seven:
$$\aligned
  X=j(z) &=\frac{(\tau^2+13 \tau+49) (\tau^2+245 \tau+2401)^3}{\tau^7}\\
         &=\frac{({\tau^{\prime}}^2+13 \tau^{\prime}+49)({\tau^{\prime}}^2
           +5 \tau^{\prime}+1)^3}{\tau^{\prime}},
\endaligned$$
$$Y=j(7z)=\frac{(\tau^2+13 \tau+49)(\tau^2+5 \tau+1)^3}{\tau},$$
where the Hauptmodul is given by
$$\tau=\left(\frac{\eta(z)}{\eta(7 z)}\right)^4,$$
and the Fricke involution $\tau^{\prime}=\tau(-\frac{1}{7z})$ satisfies
that $\tau \tau^{\prime}=49$.
$$C_X/\text{PSL}(2, 7) \cong \text{Spec}\left(\left(\mathbb{C}
  [z_1, z_2, z_3]/(z_1^3 z_2+z_2^3 z_3+z_3^3 z_1)\right)^{\text{PSL}(2, 7)}\right),
  \eqno{(1.42)}$$
where $X$ is the Klein quartic curve $C_4$: $z_1^3 z_2+z_2^3 z_3+z_3^3 z_1=0$
in $\mathbb{P}^2(\mathbb{C})$.
$$C_4 \rightarrow C_4/\text{PSL}(2, 7) \cong \mathbb{P}^1.\eqno{(1.43)}$$
Galois resolvent of the modular equation $\Phi_7(\cdot, j)=0$ of level seven.

(3) $\text{SL}(2, 13)$: $g(X_0(13))=0$.

  In the modular equation $\Phi_{13}(X, Y)=0$ of level thirteen:
$$\aligned
  X=j(z) &=\frac{(\tau^2+5\tau+13)(\tau^4+247 \tau^3+3380 \tau^2
           +15379 \tau+28561)^3}{\tau^{13}}\\
         &=\frac{({\tau^{\prime}}^2+5 \tau^{\prime}+13)({\tau^{\prime}}^4
           +7 {\tau^{\prime}}^3+20 {\tau^{\prime}}^2+19 \tau^{\prime}+1)^3}
          {\tau^{\prime}},
\endaligned\eqno{(1.44)}$$
and
$$Y=j(13z)=\frac{(\tau^2+5 \tau+13)(\tau^4+7 \tau^3+20 \tau^2+19 \tau+1)^3}{\tau},
  \eqno{(1.45)}$$
where the Hauptmodul is given by
$$\tau=\left(\frac{\eta(z)}{\eta(13z)}\right)^2,\eqno{(1.46)}$$
and the Fricke involution $\tau^{\prime}=\tau(-\frac{1}{13 z})$ satisfies
that $\tau \tau^{\prime}=13$.
$$C_Y/\text{SL}(2, 13) \cong \text{Spec}\left(\left(\mathbb{C}[z_1,
  z_2, z_3, z_4, z_5, z_6]/I(Y)\right)^{\text{SL}(2, 13)}\right).\eqno{(1.47)}$$
$$Y \rightarrow Y/\text{SL}(2, 13) \cong \mathbb{P}^1.\eqno{(1.48)}$$
Galois resolvent of the modular equation $\Phi_{13}(\cdot, j)=0$ of level
thirteen.

  However, for the group $\text{PSL}(2, 11)$, $g(X_0(11))=1$, there
does not exist such structure as above. This leads us to study the
following ring of invariant polynomials:
$$\left(\mathbb{C}[z_1, z_2, z_3, z_4, z_5, z_6]/I(Y)\right)^{\text{SL}(2, 13)}$$
and the cone $C_Y/\text{SL}(2, 13)$ over the modular curve $X(13)$ which
gives rise to a new perspective on the theory of $E_6$, $E_7$ and $E_8$-singularities.
Roughly speaking, in his celebrated icosahedral book, Felix Klein showed that
the equations of $E_6$, $E_7$ and $E_8$-singularities come from the polyhedral
equations and hence have the modularity coming from the modular curves $X(3)$,
$X(4)$ and $X(5)$, respectively. In the present paper, we show that the equations
of $E_6$, $E_7$ and $E_8$ as well as $Q_{18}$ and $E_{20}$-singularities possess
a distinct modularity: they can be obtained uniformly from a quotient
$C_Y/\text{SL}(2, 13)$ over the modular curve $X(13)$, where $C_Y$ is a cone
over $Y$. In particular, the equations of $E_6$, $E_7$, $E_8$, $Q_{18}$ and
$E_{20}$-singularities possess infinitely many kinds of distinct modular
parameterizations. They form variations of $E_6$, $E_7$, $E_8$, $Q_{18}$ and
$E_{20}$-singularity structures over the modular curve $X(13)$, for which we
give their algebraic version, geometric version and $j$-function version. In
particular, they give the same singularity structures as the Brieskorn's
spheres but with different links. This shows that there exist algebraic
varieties defined over $\mathbb{Q}$ which have infinitely many distinct
modular parameterizations. As an application, this shows that $E_6$, $E_7$
and $E_8$-singularities are not necessarily Kleinian polyhedral singularities
and they can be constructed directly from $\text{SL}(2, 13)$ without using $E_6$,
$E_7$ and $E_8$ groups. In particular, this gives a negative answer to Arnol'd
and Brieskorn's questions about the mysterious relation between the icosahedron
and $E_8$, i.e., the $E_8$-singularity is not necessarily the Kleinian icosahedral
singularity, and it also gives a relation to Witten's question about the mysterious
relation between the $A$-$D$-$E$ groups and the $A$-$D$-$E$ singularities. Now, we
give the details.

  For the following equations:
$$\aligned
  &E_6: &\quad x^4+y^3+z^2=0, \quad &\text{quartic equation},\\
  &E_7: &\quad x^3 y+y^3+z^2=0, \quad &\text{quartic equation},\\
  &E_8: &\quad x^5+y^3+z^2=0, \quad &\text{quintic equation},\\
  &Q_{18}: &\quad x^8+y^3+x z^2=0, \quad &\text{octic equation},\\
  &E_{20}: &\quad x^{11}+y^3+z^2=0, \quad &\text{eleventh equation}.
\endaligned$$
The significance of the first three of the above equations goes back to
Hermite's solution of quintic equations by modular functions (see \cite{He})
and Schwarz's solution of Gauss hypergeometric differential equations by
means of algebraic functions (see \cite{Sch}). The key idea is attributed
to Galois. Combining ideas of Galois and Riemann, in his celebrated
icosahedral book as well as his monographs on elliptic modular functions
(see \cite{K}, \cite{KF1} and \cite{KF2}), Klein gave the modular
parameterizations of these equations by the tetrahedral equation (the
tetrahedral group), the octahedral equation (the octahedral group) and
the icosahedral equation (the icosahedral group), i.e., the modular
curves $X(3)$, $X(4)$ and $X(5)$, respectively.

  By a rational double point or a simple singularity we understand the
singularity of the quotient of $\mathbb{C}^2$ by the action of a finite
subgroup of $\text{SL}(2, \mathbb{C})$ (see \cite{Sl1}). Let $\Gamma$ be
a finite subgroup of $\text{SL}(2, \mathbb{C})$. Then $\Gamma$ is one of
the following: a cyclic group of order $\ell \geq 1$ $(A_{\ell})$, a
binary dihedral group of order $4(\ell-2)$, $\ell \geq 4$ $(D_{\ell})$,
the binary tetrahedral group $(E_6)$, the binary octahedral group $(E_7)$,
or the binary icosahedral group $(E_8)$. In 1874, Klein showed that the
ring of polynomials in two variables which are invariant under $\Gamma$
is generated by three elements $x$, $y$ and $z$, which satisfy the
following relation
$$\aligned
 &A_{\ell \geq 1} \quad  & x^{\ell+1}+y^2+z^2=0,\\
 &D_{\ell \geq 4} \quad  & x^{\ell-1}+x y^2+z^2=0,\\
 &E_6             \quad  & x^4+y^3+z^2=0,\\
 &E_7             \quad  & x^3 y+y^3+z^2=0,\\
 &E_8             \quad  & x^5+y^3+z^2=0.
\endaligned$$
These results of Klein on the invariant theory of the binary polyhedral
groups were a starting point for later developments. In the minimal
resolution of such a singularity an intersection configuration of the
components of the exceptional divisor appears which can be described in
a simple way by a Dynkin diagram of type $A_{\ell}$, $D_{\ell}$, $E_6$,
$E_7$ or $E_8$. Up to analytic isomorphism, these diagrams classify the
corresponding singularities. In other words, the $ADE$ singularities are
the Kleinian singularities, i.e., the quotient singularities of $\mathbb{C}^2$
by a finite subgroup of $\text{SL}(2, \mathbb{C})$. In particular, the
$E_8$-singularity is the icosahedral singularity $\mathbb{C}^2/\Gamma$,
where $\Gamma$ is the binary icosahedral group.

  In the present paper, we will show that the $E_6$, $E_7$ and
$E_8$-singularities can be obtained from a quotient $C_Y/\text{SL}(2, 13)$
over the modular curve $X(13)$, where $C_Y$ is a cone over $Y$. This gives
infinitely many kinds of distinct constructions of the $E_6$, $E_7$ and
$E_8$-singularities which are different from the polyhedral singularities,
i.e., the (binary) tetrahedral singularity, the (binary) octahedral singularity
and the (binary) icosahedral singularity. They form variations of the $E_6$,
$E_7$ and $E_8$-singularity structures over the modular curve $X(13)$, for
which we give their algebraic version, geometric version and $j$-function
version. Our construction is based on the invariant theory for the group
$\text{SL}(2, 13)$. We construct $G$-invariant homogeneous polynomials
$\Phi_{m, n}$ of degrees $d=4$, $8$, $10$, $12$, $14$, $16$, $18$, $20$,
$22$, $26$, $30$, $32$, $34$, $42$ and $44$. Furthermore, over the modular
curve $X(13)$, these invariants are modular forms. Hence, there are some
algebraic relations among these invariants. In particular, some of them
satisfy the equations of $E_6$, $E_7$ and $E_8$-singularities. Hence, we
obtain three homomorphisms from the ring corresponding to the equations
of $E_6$, $E_7$ and $E_8$-singularities to the ring of invariant polynomials
$$\left(\mathbb{C}[z_1, z_2, z_3, z_4, z_5, z_6]/I(Y) \right)^{\text{SL}(2, 13)}$$
over the modular curve $X(13)$.

\textbf{Theorem 1.21.} (Main Theorem 11) (Variation structures of the
$E_6$, $E_7$ and $E_8$-singularities over the modular curve $X$:
algebraic version) {\it The equations of $E_6$, $E_7$ and $E_8$-singularities
$$E_6: \quad \left(\frac{\Phi_{20}^3}{\Phi_{12}^4}+1728 \Phi_{12}\right)^2
  -\left(\frac{\Phi_{26}}{\Phi_{20}}\right)^4-4 \cdot 1728
  \left(\frac{\Phi_{20}}{\Phi_{12}}\right)^3=0,\eqno{(1.49)}$$
$$E_7: \quad \Phi_{12} \cdot \left(\frac{\Phi_{26}}{\Phi_{18}}\right)^3
  -\Phi_{18}^2-1728 \Phi_{12}^3=0,\eqno{(1.50)}$$
$$E_8: \quad \Phi_{20}^3-\left(\frac{\Phi_{12}^2 \Phi_{26}}{\Phi_{20}}\right)^2
  -1728 \Phi_{12}^5=0\eqno{(1.51)}$$
possess infinitely many kinds of distinct modular parametrizations
$($with the cardinality of the continuum in ZFC set theory$)$
$$(\Phi_{12}, \Phi_{18}, \Phi_{20}, \Phi_{26})=(\Phi_{12}^{\lambda},
   \Phi_{18}^{\mu}, \Phi_{20}^{\gamma}, \Phi_{26}^{\kappa})\eqno{(1.52)}$$
over the modular curve $X$ as follows$:$
$$\left\{\aligned
  \Phi_{12}^{\lambda} &=\lambda \Phi_{3, 0}+(1-\lambda) \Phi_{0, 2}
              \quad \text{mod $\mathfrak{a}_1$},\\
  \Phi_{18}^{\mu} &=\mu \Phi_{3, 1}+(1-\mu) \Phi_{0, 3}
              \quad \text{mod $\mathfrak{a}_2$},\\
  \Phi_{20}^{\gamma} &=\gamma \Phi_{5, 0}+(1-\gamma) \Phi_{2, 2}
              \quad \text{mod $\mathfrak{a}_3$},\\
  \Phi_{26}^{\kappa} &=\kappa \Phi_{2, 3}+(1-\kappa) \Phi_{5, 1}
             \quad \text{mod $\mathfrak{a}_4$},
\endaligned\right.\eqno{(1.53)}$$
where $\Phi_{12}$, $\Phi_{18}$, $\Phi_{20}$ and $\Phi_{30}$ are invariant
homogeneous polynomials of degree $12$, $18$, $20$ and $26$, respectively.
The ideals are given by
$$\left\{\aligned
  \mathfrak{a}_1 &=(\Phi_4, \Phi_8),\\
  \mathfrak{a}_2 &=(\Phi_4, \Phi_8, \Phi_{10}, \Phi_{14}),\\
  \mathfrak{a}_3 &=(\Phi_4, \Phi_8, \Phi_{10}, \Phi_{4, 0}, \Phi_{1, 2}),\\
  \mathfrak{a}_4 &=(\Phi_4, \Phi_8, \Phi_{10}, \Phi_{14}, \Phi_{4, 0},
                    \Phi_{1, 2}, \Phi_{1, 3}, \Phi_{4, 1}),
\endaligned\right.\eqno{(1.54)}$$
and the parameter space $\{ (\lambda, \mu, \gamma, \kappa) \} \cong \mathbb{C}^4$.
They form variations of the $E_6$, $E_7$ and $E_8$-singularity structures
over the modular curve $X$.}

  The relations between the link at $O$ for the cone $C_Y/G$ over $X(13)$
and the link at $O$ for $E_6$, $E_7$ and $E_8$-singularities are given as
follows:

\textbf{Theorem 1.22.} (Main Theorem 12) (Variation structures of the
$E_6$, $E_7$ and $E_8$-singularities over the modular curve $X$:
geometric version) {\it There are three morphisms from the cone $C_Y$
over $Y$ to the $E_6$, $E_7$ and $E_8$-singularities:
$$\aligned
  &f_6: C_Y/G \rightarrow\\
  &\text{Spec}\left(\mathbb{C}\{\Phi_{12}, \Phi_{20}\}
  [\Phi_{26}]/((\frac{\Phi_{20}^3}{\Phi_{12}^4}+1728 \Phi_{12})^2
  -(\frac{\Phi_{26}}{\Phi_{20}})^4-4 \cdot 1728
  (\frac{\Phi_{20}}{\Phi_{12}})^3)\right).
\endaligned\eqno{(1.55)}$$
$$f_7: C_Y/G \rightarrow \text{Spec}\left(\mathbb{C}\{\Phi_{18}\}
  [\Phi_{12}, \Phi_{26}]/((\Phi_{12} \cdot(\frac{\Phi_{26}}{\Phi_{18}})^3
  -\Phi_{18}^2-1728 \Phi_{12}^3)\right).\eqno{(1.56)}$$
$$f_8: C_Y/G \rightarrow \text{Spec}\left(\mathbb{C}\{\Phi_{20}\}
  [\Phi_{12}, \Phi_{26}]/(\Phi_{20}^3-(\frac{\Phi_{12}^2 \Phi_{26}}{\Phi_{20}})^2
  -1728 \Phi_{12}^5)\right).\eqno{(1.57)}$$
over the modular curve $X$. In particular, there are infinitely many
such quadruples $(\Phi_{12}, \Phi_{18}, \Phi_{20}, \Phi_{26})$ $=$
$(\Phi_{12}^{\lambda}, \Phi_{18}^{\mu}, \Phi_{20}^{\gamma}, \Phi_{26}^{\kappa})$
whose parameter space $\{ (\lambda, \mu, \gamma, \kappa) \}$
$\cong \mathbb{C}^4$. They form variations of the $E_6$, $E_7$ and
$E_8$-singularity structures over the modular curve $X$.}

  In particular, the links at $O$ for the $E_6$, $E_7$ and
$E_8$-singularities by a triple of polynomials are given as follows:

  The link at $O$ for the $E_6$-singularity by the triple
$(\Phi_{12}, \Phi_{20}, \Phi_{26})$:
$$\left\{\aligned
  &|\Phi_{12}|^2+|\Phi_{20}|^2+|\Phi_{26}|^2=1,\\
  &\left(\frac{\Phi_{20}^3}{\Phi_{12}^4}+1728 \Phi_{12}\right)^2
  -\left(\frac{\Phi_{26}}{\Phi_{20}}\right)^4-4 \cdot 1728
  \left(\frac{\Phi_{20}}{\Phi_{12}}\right)^3=0.
\endaligned\right.\eqno{(1.58)}$$

  The link at $O$ for the $E_7$-singularity by the triple
$(\Phi_{12}, \Phi_{18}, \Phi_{26})$:
$$\left\{\aligned
  &|\Phi_{12}|^2+|\Phi_{18}|^2+|\Phi_{26}|^2=1,\\
  &\Phi_{12} \cdot \left(\frac{\Phi_{26}}{\Phi_{18}}\right)^3
   -\Phi_{18}^2-1728 \Phi_{12}^3=0.
\endaligned\right.\eqno{(1.59)}$$

  The link at $O$ for the $E_8$-singularity by the triple
$(\Phi_{12}, \Phi_{20}, \Phi_{26})$:
$$\left\{\aligned
  &|\Phi_{12}|^2+|\Phi_{20}|^2+|\Phi_{26}|^2=1,\\
  &\Phi_{20}^3-\left(\frac{\Phi_{12}^2 \Phi_{26}}{\Phi_{20}}\right)^2
   -1728 \Phi_{12}^5=0.
\endaligned\right.\eqno{(1.60)}$$

  Theorem 1.21 implies the following:

\textbf{Theorem 1.23.} (Main Theorem 13) {\it There exist algebraic varieties
defined over $\mathbb{Q}$ which have infinitely many distinct modular
parameterizations.}

  Theorem 1.21 and Theorem 1.22 give a new perspective on the theory of
singularities: they show that there exist infinitely many kinds of
distinct constructions of the $E_6$, $E_7$ and $E_8$-singularities: one and
only one are given by the Kleinian singularities $\mathbb{C}^2/\Gamma$, where
$\Gamma$ is the binary tetrahedral group, the binary octahedral group and the
binary icosahedral group (see \cite{K}), i.e., the tetrahedral singularity, the
octahedral singularity and the icosahedral singularity, the other infinitely
many kinds of constructions are given from the quotient $C_Y/\text{SL}(2, 13)$
over the modular curve $X$. Hence, the equations of $E_6$, $E_7$ and
$E_8$-singularities possess infinitely many kinds of distinct modular
parametrizations.

  In his talk at ICM 1970 \cite{Br2}, Brieskorn showed how to construct
the singularity of type $ADE$ directly from the simple complex Lie group
of the same type. Namely, assume that $G$ is of type $ADE$, Brieskorn
proved a conjecture made by Grothendieck that the intersection of a
transversal slice to the sub-regular unipotent orbit with the unipotent
variety has a simple surface singularity of the same type as $G$. A
fuller treatment was given by Slodowy (see \cite{Br2} and \cite{Sl1}). A
clarification of the occurrence of the polyhedral groups in Brieskorn's
construction (see \cite{GrP} and \cite{GrP2}),  and thus a direct
relationship between the simple Lie groups and the finite subgroups of
$\text{SL}(2, \mathbb{C})$, was achieved by Kronheimer (see \cite{Kr1} and
\cite{Kr2}) using differential geometric methods. His construction starts
directly from the finite subgroups of $\text{SL}(2, \mathbb{C})$ and uses
hyper-K\"{a}hler quotient constructions. Kronheimer also gave an algebraic
approach using McKay correspondence. However, Brieskorn had still written
at the end of \cite{Br2}: ``Thus we see that there is a relation between
exotic spheres, the icosahedron and $E_8$. But I still do not understand
why the regular polyhedra come in.'' (see also \cite{Gr}, \cite{GrP},
\cite{GrP2} and \cite{Br3}). On the other hand, Arnol'd pointed out that
the theory of singularities is even linked (in a quite mysterious way) to
the classification of regular polyhedra in three-dimensional Euclidean
space (see \cite{Ar3}, p. 43). In his survey article on Platonic solids,
Kleinian singularities and Lie groups \cite{Sl2}, Slodowy found that the
objects of these different classifications are related to each other by
mathematical constructions. However, up to now, these constructions do
not explain why the different classifications should be related at all.

  From the viewpoint of string theory, according to Witten (see \cite{Wi1}
or \cite{H}): ``Mathematicians thought of them ($A$-$D$-$E$ singularities)
as somehow connected with $A$-$D$-$E$ groups $\ldots$. Dynkin diagrams are
all very well and good but they are a method of studying groups, not the
other way around.'' An extensive mathematical theory relates the $A$-$D$-$E$
singularity to the $A$-$D$-$E$ Dynkin diagram and many associated bits of
geometry and algebra. However, the role of the $A$-$D$-$E$ group in relation
to the singularity is elusive. In particular, The role of the $A$-$D$-$E$
group in relation to the $A$-$D$-$E$ singularity is somewhat mysterious
classically (see \cite{Wi2}).

  As a consequence, Theorem 1.21 and Theorem 1.22 show that the $E_6$, $E_7$
and $E_8$-singularities are not necessarily the Kleinian binary polyhedral
singularities, where the binary polyhedral groups are closely related to
the $E_6$, $E_7$ and $E_8$-groups. The group $\text{SL}(2, 13)$ can also
be related to the $E_6$, $E_7$ and $E_8$ in a unified way. This provides
a new perspective about the relation between the $A$-$D$-$E$ group and
the $A$-$D$-$E$ singularity (see \cite{Wi2}).

  In particular, the icosahedron does not necessarily appear in the triple
(exotic spheres, icosahedron, $E_8$) of Brieskorn \cite{Br2}. The group
$\text{SL}(2, 13)$ can take its place and there are infinitely many kinds
of the other triples (exotic spheres, $\text{SL}(2, 13)$, $E_8$). The links
of these infinitely many kinds of distinct constructions of the $E_8$-singularity:
$\mathbb{C}^2/\text{SL}(2, 5)$ and a variation of the $E_8$-singularity
structure
$$C_Y/G \rightarrow \text{Spec}\left(\mathbb{C}\{\Phi_{20}\}
  [\Phi_{12}, \Phi_{26}]/(\Phi_{20}^3-(\frac{\Phi_{12}^2 \Phi_{26}}{\Phi_{20}})^2
  -1728 \Phi_{12}^5)\right)$$
over the modular curve $X$ give the same singularity structure but with different
links. For $\text{SL}(2, 5)$, it leads to the Poincar\'{e} homology $3$-sphere,
which is a special case of the Brieskorn sphere $\Sigma (2, 3, p)$, i.e., the link
of the singularity $x^2+y^3+z^p=0$:
$$\Sigma (2, 3, p):
  \left\{\aligned
  |x|^2+|y|^2+|z|^2 &=1,\\
  x^2+y^3+z^p &=0,
  \endaligned\right. \quad (x, y, z \in \mathbb{C})$$
for $p=5$, whose higher dimensional lifting:
$$z_1^5+z_2^3+z_3^2+z_4^2+z_5^2=0, \quad \sum_{i=1}^{5} z_i \overline{z_i}=1,
  \quad z_i \in \mathbb{C} \quad (1 \leq i \leq 5)$$
gives the Milnor's standard generator of $\Theta_7$. On the other hand, for
$\text{SL}(2, 13)$, they are the links of $E_8$-singularity by the triple
$(\Phi_{12}, \Phi_{20}, \Phi_{26})$ given by (1.60). Hence, this gives
a negative answer to Arnol'd and Brieskorn's questions about the mysterious
relation between the icosahedron and $E_8$. Moreover, the relation between
Platonic solids, Kleinian singularities and Lie groups appearing in Slodowy's
survey \cite{Sl2} can be replaced by the relations between $\text{SL}(2, 13)$,
variations of the $E_6$, $E_7$ and $E_8$-singularity structures
$$\aligned
  &f_6: C_Y/G \rightarrow\\
  &\text{Spec}\left(\mathbb{C}\{\Phi_{12}, \Phi_{20}\}
  [\Phi_{26}]/((\frac{\Phi_{20}^3}{\Phi_{12}^4}+1728 \Phi_{12})^2
  -(\frac{\Phi_{26}}{\Phi_{20}})^4-4 \cdot 1728
  (\frac{\Phi_{20}}{\Phi_{12}})^3)\right),
\endaligned$$
$$f_7: C_Y/G \rightarrow \text{Spec}\left(\mathbb{C}\{\Phi_{18}\}
  [\Phi_{12}, \Phi_{26}]/((\Phi_{12} \cdot(\frac{\Phi_{26}}{\Phi_{18}})^3
  -\Phi_{18}^2-1728 \Phi_{12}^3)\right),$$
$$f_8: C_Y/G \rightarrow \text{Spec}\left(\mathbb{C}\{\Phi_{20}\}
  [\Phi_{12}, \Phi_{26}]/(\Phi_{20}^3-(\frac{\Phi_{12}^2 \Phi_{26}}{\Phi_{20}})^2
  -1728 \Phi_{12}^5)\right)$$
over the modular curve $X$ and $E_6$, $E_7$, $E_8$, respectively.

  Moreover, Theorem 1.21 and Theorem 1.22 can be extended to the following
two kinds of singularities:
$$\left\{\aligned
  &Q_{18}: & x^8+y^3+xz^2=0,\\
  &E_{20}: & x^{11}+y^3+z^2=0,
\endaligned\right.$$
where $Q_{18}$ and $E_{20}$ are two bimodal singularities in the
pyramids of $14$ exceptional singularities (see \cite{Ar2}, p.255).
Theorem 1.21, Theorem 1.22, as well as Theorem 1.24 and Theorem 1.25
show that five different algebraic surfaces, the equations of $E_6$,
$E_7$, $E_8$, $Q_{18}$ and $E_{20}$-singularities can be realized
from the same quotients $C_Y/\text{SL}(2, 13)$ over the modular
curve $X$ and have the same modular parametrizations.

\textbf{Theorem 1.24.} (Variations of $Q_{18}$ and $E_{20}$-singularity
structures over the modular curve $X$: algebraic version) {\it The equations
of $Q_{18}$ and $E_{20}$-singularities
$$\Phi_{32}^3-\Phi_{12} \left(\frac{\Phi_{12}^4 \Phi_{26}}{\Phi_{32}}\right)^2
  -1728 \Phi_{12}^8=0,\eqno{(1.61)}$$
$$\Phi_{44}^3-\left(\frac{\Phi_{12}^7 \Phi_{26}}{\Phi_{44}}\right)^2
  -1728 \Phi_{12}^{11}=0\eqno{(1.62)}$$
possess infinitely many kinds of distinct modular parametrizations
$($with the cardinality of the continuum in ZFC set theory$)$
$$(\Phi_{12}, \Phi_{26}, \Phi_{32}, \Phi_{44})=(\Phi_{12}^{\lambda},
   \Phi_{26}^{\mu}, \Phi_{32}^{\gamma}, \Phi_{44})\eqno{(1.63)}$$
over the modular curve $X$ as follows$:$
$$\left\{\aligned
  \Phi_{12}^{\lambda} &=\lambda \Phi_{3, 0}+(1-\lambda) \Phi_{0, 2}
              \quad \text{mod $\mathfrak{a}_1$},\\
  \Phi_{26}^{\mu} &=\mu \Phi_{2, 3}+(1-\mu) \Phi_{5, 1} \quad \text{mod $\mathfrak{a}_4$},\\
  \Phi_{32}^{\gamma} &=\gamma_1 \Phi_{8, 0}+\gamma_2 \Phi_{5, 2}+
                       (1-\gamma_1-\gamma_2) \Phi_{2, 4} \quad \text{mod $\mathfrak{a}_4$},\\
  \Phi_{44} &=\Phi_{11, 0} \quad \text{mod $\mathfrak{a}_5$},
\endaligned\right.\eqno{(1.64)}$$
where $\Phi_{12}$, $\Phi_{26}$, $\Phi_{32}$ and $\Phi_{44}$ are invariant
homogeneous polynomials of degree $12$, $26$, $32$ and $44$, respectively.
The ideals are given by
$$\left\{\aligned
  \mathfrak{a}_1 &=(\Phi_4, \Phi_8),\\
  \mathfrak{a}_4 &=(\Phi_4, \Phi_8, \Phi_{10}, \Phi_{14}, \Phi_{4, 0},
                    \Phi_{1, 2}, \Phi_{1, 3}, \Phi_{4, 1}),\\
  \mathfrak{a}_5 &=(\Phi_4, \Phi_8, \Phi_{10}, \Phi_{14}, \Phi_{4, 0},
                    \Phi_{1, 2}, \Phi_{1, 3}, \Phi_{4, 1}, \Phi_{1, 5},
                    \Phi_{4, 3}, \Phi_{7, 1}),
\endaligned\right.\eqno{(1.65)}$$
and the parameter space $\{ (\lambda, \mu, \gamma) \} \cong \mathbb{C}^4$.
They form variations of $Q_{18}$ and $E_{20}$-singularity structures over
the modular curve $X$.}

  The relations between the link at $O$ for the cone $C_Y/G$ over $X(13)$
and the link at $O$ for $Q_{18}$ and $E_{20}$-singularities are given as
follows:

\textbf{Theorem 1.25.} (Variations of $Q_{18}$ and $E_{20}$-singularity
structures over the modular curve $X$: geometric version) {\it There
are two morphisms from the cone $C_Y$ over $Y$ to the $Q_{18}$ and
$E_{20}$-singularities:
$$f_{18}: C_Y/G \rightarrow \text{Spec} \left(\mathbb{C}\{\Phi_{32}\}
          [\Phi_{12}, \Phi_{26}]/(\Phi_{32}^3-\Phi_{12} (\frac{\Phi_{12}^4
          \Phi_{26}}{\Phi_{32}})^2-1728 \Phi_{12}^8)\right)\eqno{(1.66)}$$
and
$$f_{20}: C_Y/G \rightarrow \text{Spec} \left(\mathbb{C}\{\Phi_{44}\}
          [\Phi_{12}, \Phi_{26}]/(\Phi_{44}^3-(\frac{\Phi_{12}^7
          \Phi_{26}}{\Phi_{44}})^2-1728 \Phi_{12}^{11})\right)\eqno{(1.67)}$$
over the modular curve $X$. In particular, there are infinitely many such
quadruples $(\Phi_{12}, \Phi_{26}, \Phi_{32}, \Phi_{44})$ $=$ $(\Phi_{12}^{\lambda},
\Phi_{26}^{\mu}, \Phi_{32}^{\gamma}, \Phi_{44})$ whose parameter space
$\{ (\lambda, \mu, \gamma) \}$ $\cong \mathbb{C}^4$. They form variations
of $Q_{18}$ and $E_{20}$-singularity structures over the modular curve $X$.}

  In particular, the links at $O$ for the $Q_{18}$ and $E_{20}$-singularities
by a triple of polynomials are given as follows:

  The link at $O$ for the $Q_{18}$-singularity by the triple
$(\Phi_{12}, \Phi_{26}, \Phi_{32})$:
$$\left\{\aligned
  &|\Phi_{12}|^2+|\Phi_{26}|^2+|\Phi_{32}|^2=1,\\
  &\Phi_{32}^3-\Phi_{12} \left(\frac{\Phi_{12}^4 \Phi_{26}}{\Phi_{32}}\right)^2
   -1728 \Phi_{12}^8=0.
\endaligned\right.\eqno{(1.68)}$$

  The link at $O$ for the $E_{20}$-singularity by the triple
$(\Phi_{12}, \Phi_{26}, \Phi_{44})$:
$$\left\{\aligned
  &|\Phi_{12}|^2+|\Phi_{26}|^2+|\Phi_{44}|^2=1,\\
  &\Phi_{44}^3-\left(\frac{\Phi_{12}^7 \Phi_{26}}{\Phi_{44}}\right)^2
   -1728 \Phi_{12}^{11}=0.
\endaligned\right.\eqno{(1.69)}$$

  In fact, Klein had noticed the similarity between the relation of
the equation $x^5+y^3+z^2=0$ to the icosahedral group $\text{PSL}(2, 5)$
and the relation of the equation $x^7+y^3+z^2=0$  to the group
$\text{PSL}(2, 7)$ (see \cite{K1}, \cite{K2}, \cite{KF1} and \cite{KF2}).
This is the starting point of the work of Dolgachev (see \cite{Do}) to
which Arnol'd was referring when he spoke about the wonderful coincidences
with Lobatchevsky triangles and automorphic functions (see \cite{Ar1}).
The normal form of Arnol'd for the quasi-homogeneous singularity $E_{12}$
in three variables is $x^7+y^3+z^2$, which can be realized as the quotient
conical singularity as follows (see \cite{BrPR} and \cite{Do}): The
canonical model $Y$ of the modular curve $X(7)$ in $\mathbb{CP}^2$ is
the Klein quartic given by the homogeneous equation
$z_1^3 z_2+z_2^3 z_3+z_3^3 z_1=0$. The finite group $\text{PSL}(2, 7)$
acts linearly on $\mathbb{C}^3$ and on $\mathbb{CP}^2$ leaving invariant
$Y \subset \mathbb{CP}^2$ and the cone $C_Y \subset \mathbb{C}^3$.
Calculations of invariants by Klein and Gordan imply:
$$\left[\mathbb{C}[z_1, z_2, z_3]/(z_1^3 z_2+z_2^3 z_3+z_3^3 z_1)
  \right]^{\text{PSL}(2, 7)} \cong \mathbb{C}[x, y, z]/(x^7+y^3+z^2).
  \eqno{(1.70)}$$
This algebraic result can be interpreted geometrically as follows:
The affine algebraic surface defined by the equation $x^7+y^3+z^2=0$
is the quotient of the cone $C_Y$ by the group $\text{PSL}(2, 7)$
over the modular curve $X(7)$, where $C_Y$ is the cone over $Y$.
Similarly, Klein also obtained the structure of the $\mathbb{C}$-algebra
of $\mathbb{C}[z_1, z_2]^{\text{SL}(2, 5)}$ of $\text{SL}(2, 5)$-invariant
polynomials on $\mathbb{C}^2$:
$$\mathbb{C}[z_1, z_2]^{\text{SL}(2, 5)} \cong \mathbb{C}[x, y, z]/
  (x^5+y^3+z^2).\eqno{(1.71)}$$
This algebraic result can also be interpreted geometrically as follows:
The affine algebraic surface defined by the equation $x^5+y^3+z^2=0$ is
the quotient of the cone $C_Y$ by the group $\text{SL}(2, 5)$ over the
modular curve $X(5)$, where $Y=\mathbb{CP}^1$ is the canonical model
of the modular curve $X(5)$ and $C_Y$ is a cone over $Y$. Therefore,
(1.49), (1.50), (1.51), (1.52), (1.55), (1.56), (1.57), (1.61), (1.62),
(1.63), (1.66), (1.67), (1.70) and (1.71) give a complete and unified
description for the relation between the $G$-invariant homogeneous
polynomials and the associated singularities corresponding to the genus
zero modular curves $X_0(N)$, where $N=5$, $7$, $13$ and $G=\text{SL}(2, 5)$,
$\text{PSL}(2, 7)$ and $\text{SL}(2, 13)$, respectively:
$$\begin{matrix}
  &C_Y/\text{SL}(2, 5) & C_Y/\text{PSL}(2, 7) & C_Y/\text{SL}(2, 13)\\
  &\downarrow          &\downarrow            & \downarrow\\
  &X(5)                &X(7)                  & X(13)\\
  &Y=\mathbb{CP}^1     &Y=\text{Klein quartic curve} &Y=Y, Y_2, Y_3\\
  &E_8\text{-singularity} &E_{12}\text{-singularity}
  &\text{$E_6$, $E_7$, $E_8$, $Q_{18}$}\\
  &                       &                          &\text{and $E_{20}$-singularities}
\end{matrix}\eqno{(1.72)}$$

  Finally, recall that there is a decomposition formula of the
elliptic modular function $j$ in terms of the icosahedral invariants
$f$, $H$ and $T$ of degrees $12$, $20$ and $30$ over the modular
curve $X(5)$ (see section two, in particular (2.8) for the details):
$$\aligned
 &j(z): j(z)-1728: 1\\
=&H(x_1(z), x_2(z))^3: -T(x_1(z), x_2(z))^2: f(x_1(z), x_2(z))^5,
\endaligned\eqno{(1.73)}$$
which was discovered by Klein (see \cite{K}, \cite{KF1} and \cite{KF2})
and later by Ramanujan (see \cite{Du}). In contrast with (1.73), we
have the following $j$-function version:

\textbf{Theorem 1.26.} (Variations of the structure of decomposition
formulas of the elliptic modular functions $j$ over the modular curve
$X$) {\it There are infinitely many kinds of distinct decomposition
formulas of the elliptic modular function $j$ in terms of the invariants
$\Phi_{12}$, $\Phi_{18}$, $\Phi_{20}$ and $\Phi_{26}$ over the modular
curve $X$:
$$E_7\text{-type}: \quad
  j(z): j(z)-1728: 1=\Phi_{12} \left(\frac{\Phi_{26}}{\Phi_{18}}\right)^3:
  \Phi_{18}^2: \Phi_{12}^3,\eqno{(1.74)}$$
$$E_8\text{-type}: \quad
  j(z): j(z)-1728: 1=\Phi_{20}^3: \left(\frac{\Phi_{26} \Phi_{12}^2}
  {\Phi_{20}}\right)^2: \Phi_{12}^5,\eqno{(1.75)}$$
where $(\Phi_{12}, \Phi_{18}, \Phi_{20}, \Phi_{26})=(\Phi_{12}^{\lambda},
\Phi_{18}^{\mu}, \Phi_{20}^{\gamma}, \Phi_{26}^{\kappa})$ are given
by $(1.38)$. They form variations of the structure of decomposition
formulas of the elliptic modular functions $j$ over the modular curve
$X$.}

  In fact, these infinitely many kinds of distinct decompositions (1.74)
and (1.75) have the form of $E_7$-type and $E_8$-type, respectively.
They have the different geometric interpretation with respect to (1.73):
one and only one (1.73) is over the modular curve $X(5)$, the other
infinitely many kinds of decompositions (which form variations of
the structure of decomposition formulas) are over the modular curve $X(13)$.
They also have the different algebraic interpretation: one and only one is
invariant under the group $\text{SL}(2, 5)$, the other infinitely many kinds
of decompositions (which form variations of the structure of decomposition
formulas) are invariant under the group $\text{SL}(2, 13)$.

  Finally, we have the following:

\textbf{Probelem 1.27.} Find an analogue of Theorem 1.10 for the other
modular curves $X(p)$, where $p \geq 17$ is a prime number.

\textbf{Probelem 1.28.} Arithmetical realizations by cohomology groups
defined over some algebraic number fields for higher dimensional modular
varieties or Shimura varieties, i.e., are there similar phenomena for
higher dimensional modular varieties or Shimura varieties?

  This paper consists of eleven sections. In section two, we revisit the
standard constructions of the $E_6$, $E_7$ and $E_8$-singularities as the
well-known Kleinian polyhedral singularities and their modularity coming
from the modular curves $X(3)$, $X(4)$ and $X(5)$, respectively. In section
three, we study the invariant theory and modular forms for $\text{SL}(2, 13)$.
In particular, we construct the quadratic invariants ($\mathbf{A}$-terms),
the cubic invariants ($\mathbf{D}$-terms), and the sextic invariants
($\mathbf{G}$-terms). In section four, we recall Klein's construction
of the locus $\mathcal{L}_p$ corresponding to the modular curves $X(p)$
by means of quartic systems. In particular, we revisit the beautiful
construction of the modular curve $X(11)$ by Klein's invariant Fano
three-fold. Then we continue the study of section three by constructing
the quartic invariants: $\mathbf{B}$-terms and $\mathbf{C}$-terms. In
particular, we prove that there is a one-to-one correspondence between the
$\mathbb{B}$-terms and Klein's quartic systems, i.e. $\Phi_{abcd}$-terms.
In section five, we construct an ideal $I(Y)$ generated by our quartic
invariants ($\mathbf{B}$-terms), and prove that the corresponding curve
$Y$ is isomorphic to the modular curve $X(13)$. Moreover, we show that
this curve $Y$ can be constructed from an invariant quartic Fano-fold
$\Phi_4(z_1, z_2, z_3, z_4, z_5, z_6)=0$, and the ideal $I(Y)$ is invariant
under the action of the group isomorphic to $\text{SL}(2, 13)$. In section
six, the construction of the invariant ideal $I(Y)$ leads to a $21$-dimensional
reducible representation of $\text{SL}(2, 13)$. We give its decomposition as
the direct sum of $1$, $7$ and $13$-dimensional representations, which gives
the corresponding geometric construction of the curve $Y$ and a non-standard
geometric realization of the degenerate principal series representation for
$\text{SL}(2, 13)$ as well as a geometric realization of the Steinberg
representation for $\text{SL}(2, 13)$. In section seven, we establish the
modularity of the curve $Y$ by means of an explicit construction by theta
constants of order $13$. It also gives an example of hyperbolic uniformization
of arithmetic type with respect to a congruence subgroup of
$\text{SL}(2, \mathbb{Z})$ as well as an example of the explicit uniformization
for algebraic space curves of higher genus. As a consequence, this gives a new
solution to Hilbert's 22nd problem, i.e., an explicit uniformization for
algebraic relations among five variables by means of automorphic functions.
On the other hand, it leads to $21$ modular equations, which greatly improve
the result of Ramanujan and Evans on modular equations of order 13. In
particular, we determine both a small set of generators for a polynomial
defining ideal and the minimum number of equations needed to define this
algebraic space curve $Y$ in $\mathbb{P}^5$. In section eight, we construct
a Galois covering $Y \rightarrow Y/\text{SL}(2, 13) \cong \mathbb{CP}^1$
and give the geometric realization of the Galois resolvent for the modular
equation of order $13$ corresponding to a Hauptmodul for $\Gamma_0(13)$.
In section nine, we continue the study of invariant theory for
$\text{SL}(2, 13)$ by the computation of some invariant polynomials. We
construct $G$-invariant homogeneous polynomials $\Phi_{m, n}$ of degrees
$d=4$, $8$, $10$, $12$, $14$, $16$, $18$, $20$, $22$, $26$, $30$ and $34$.
Furthermore, over the modular curve $X(13)$, these invariants are modular
forms. In section ten, we find that there are some algebraic relations among
these invariants. In particular, some of them satisfy the equations of $E_6$,
$E_7$ and $E_8$-singularities. Hence, we obtain three homomorphisms from the
ring corresponding to the equations of $E_6$, $E_7$ and $E_8$-singularities
to the ring of invariant polynomials
$\left(\mathbb{C}[z_1, z_2, z_3, z_4, z_5, z_6]/I(Y) \right)^{\text{SL}(2, 13)}$
over the modular curve $X(13)$. In particular, there are infinitely many
such homogeneous polynomials $(\Phi_{12}^{\lambda}, \Phi_{18}^{\mu},
\Phi_{20}^{\gamma}, \Phi_{26}^{\kappa})$,
which form variations of the $E_6$, $E_7$ and $E_8$-singularity structures
over the modular curve $X$. We give their algebraic version, geometric version
and $j$-function version. Therefore, we give a different construction of the
$E_6$, $E_7$ and $E_8$-singularities coming from a quotient $C_Y/\text{SL}(2, 13)$
over the modular curve $X$, where $C_Y$ is the cone over the algebraic curve
$Y$. In section eleven, we construct $G$-invariant homogeneous polynomials
$\Phi_{m, n}$ of degrees $d=32$, $42$ and $44$, and extend our work to the cases
of $Q_{18}$ and $E_{20}$-singularities and obtain variations of $Q_{18}$ and
$E_{20}$-singularity structures over the modular curve $X$.

\textbf{Acknowledgements}. The author would like to thank P. Deligne and
J.-P. Serre for their very detailed and helpful comments as well as their
patience.

\begin{center}
{\large\bf 2. Modularity for equations of $E_6$, $E_7$ and $E_8$-singularities
              coming from modular curves $X(3)$, $X(4)$ and $X(5)$}
\end{center}

  Let us recall some classical result on the relation between the
polyhedra in $\mathbb{R}^3$ and the $E_6$, $E_7$, $E_8$-singularities
(see \cite{Mc}). Starting with the polynomial invariants of the finite
subgroup of $\text{SL}(2, \mathbb{C})$, a surface is defined from the
single syzygy which relates the three polynomials in two variables.
This surface has a singularity at the origin; the singularity can be
resolved by constructing a smooth surface which is isomorphic to
the original one except for a set of component curves which form
the pre-image of the origin. The components form a Dynkin curve
and the matrix of their intersections is the negative of the Cartan
matrix for the appropriate Lie algebra. The Dynkin curve is the
dual of the Dynkin graph. For example, if $\Gamma$ is the binary
icosahedral group, the corresponding Dynkin curve is that of $E_8$,
and $\mathbb{C}^2/\Gamma \subset \mathbb{C}^3$ is the set of zeros
of the equation
$$x^5+y^3+z^2=0.\eqno{(2.1)}$$
The link of this $E_8$-singularity, the Poincar\'{e} homology
$3$-sphere (see \cite{KS}), has a higher dimensional lifting:
$$z_1^5+z_2^3+z_3^2+z_4^2+z_5^2=0, \quad \sum_{i=1}^{5} z_i \overline{z_i}=1,
  \quad z_i \in \mathbb{C} \quad (1 \leq i \leq 5),\eqno{(2.2)}$$
which is the Brieskorn description of one of Milnor's exotic
$7$-dimensional spheres. In fact, it is an exotic $7$-sphere
representing Milnor's standard generator of $\Theta_7$ (see
\cite{Br1}, \cite{Br2} and \cite{Hi}).

  In his celebrated book \cite{K}, Klein gave a parametric solution
of the above singularity (2.1) by homogeneous polynomials $T$, $H$,
$f$ in two variables of degrees $30$, $20$, $12$ with integral
coefficients, where
$$f=z_1 z_2 (z_1^{10}+11 z_1^5 z_2^5-z_2^{10}),$$
$$H=\frac{1}{121} \begin{vmatrix}
    \frac{\partial^2 f}{\partial z_1^2} &
    \frac{\partial^2 f}{\partial z_1 \partial z_2}\\
    \frac{\partial^2 f}{\partial z_2 \partial z_1} &
    \frac{\partial^2 f}{\partial z_2^2}
    \end{vmatrix}
  =-(z_1^{20}+z_2^{20})+228 (z_1^{15} z_2^5-z_1^5 z_2^{15})
   -494 z_1^{10} z_2^{10},$$
$$T=-\frac{1}{20} \begin{vmatrix}
    \frac{\partial f}{\partial z_1} &
    \frac{\partial f}{\partial z_2}\\
    \frac{\partial H}{\partial z_1} &
    \frac{\partial H}{\partial z_2}
    \end{vmatrix}
  =(z_1^{30}+z_2^{30})+522 (z_1^{25} z_2^5-z_1^5 z_2^{25})
   -10005 (z_1^{20} z_2^{10}+z_1^{10} z_2^{20}).$$
They satisfy the famous (binary) icosahedral equation
$$T^2+H^3=1728 f^5.\eqno{(2.3)}$$
In fact, $f$, $H$ and $T$ are invariant polynomials under the
action of the binary icosahedral group. The above equation (2.3)
is closely related to Hermite's celebrated work (see \cite{He})
on the resolution of the quintic equations by elliptic modular
functions of order five. Essentially the same relation had been
found a few years earlier by Schwarz (see \cite{Sch}), who
considered three polynomials $\varphi_{12}$, $\varphi_{20}$ and
$\varphi_{30}$ whose roots correspond to the vertices, the
midpoints of the faces and the midpoints of the edges of an
icosahedron inscribed in the Riemann sphere. He obtained the
identity $\varphi_{20}^3-1728 \varphi_{12}^5=\varphi_{30}^2$.
We see this identity as well as (2.3) as the defining relation
between three generators $f$, $H$ and $T$ of the ring of invariants
$\mathbb{C}[z_1, z_2]^{\Gamma}$ of the binary icosahedral group
$\Gamma$ acting on $\mathbb{C}^2$, and we identify this ring with
the ring of functions on the affine variety $\mathbb{C}^2/\Gamma$
embedded in $\mathbb{C}^3$ and given by such an equation (see
\cite{Br3}). Namely,
$$\mathbb{C}[z_1, z_2]^{\Gamma} \cong \mathbb{C}[f, H, T]/
  (T^2+H^3-1728 f^5).\eqno{(2.4)}$$
Thus we see that from the very beginning there was a close
relation between the $E_8$-singularity and the icosahedron.
Moreover, the icosahedral equation (2.3) can be interpreted
in terms of the modular curve $X(5)$ which was also known by
Klein (see \cite{KF1}, p. 631). Let $x_1(z)=\eta(z) a(z)$ and
$x_2(z)=\eta(z) b(z)$, where
$$a(z)=e^{-\frac{3 \pi i}{10}} \theta \begin{bmatrix}
       \frac{3}{5}\\ 1 \end{bmatrix}(0, 5z), \quad
  b(z)=e^{-\frac{\pi i}{10}} \theta \begin{bmatrix}
       \frac{1}{5}\\ 1 \end{bmatrix}(0, 5z)$$
are theta constants of order five and
$\eta(z):=q^{\frac{1}{24}} \prod_{n=1}^{\infty} (1-q^n)$ with
$q=e^{2 \pi i z}$ is the Dedekind eta function which are all
defined in the upper-half plane
$\mathbb{H}=\{ z \in \mathbb{C}: \text{Im}(z)>0 \}$. Then
$$\left\{\aligned
  f(x_1(z), x_2(z)) &=-\Delta(z),\\
  H(x_1(z), x_2(z)) &=-\eta(z)^8 \Delta(z)E_4(z),\\
  T(x_1(z), x_2(z)) &=\Delta(z)^2 E_6(z),
\endaligned\right.\eqno{(2.5)}$$
where
$$E_4(z):=\frac{1}{2} \sum_{m, n \in \mathbb{Z}, (m, n)=1}
          \frac{1}{(mz+n)^4}, \quad
  E_6(z):=\frac{1}{2} \sum_{m, n \in \mathbb{Z}, (m, n)=1}
          \frac{1}{(mz+n)^6}$$
are Eisenstein series of weight $4$ and $6$, and
$\Delta(z)=\eta(z)^{24}$ is the discriminant. The relations
$$j(z)=\frac{E_4(z)^3}{\Delta(z)}=\frac{H(x_1(z), x_2(z))^3}
        {f(x_1(z), x_2(z))^5},\eqno{(2.6)}$$
$$j(z)-1728=\frac{E_6(z)^2}{\Delta(z)}=-\frac{T(x_1(z), x_2(z))^2}
            {f(x_1(z), x_2(z))^5}\eqno{(2.7)}$$
give the icosahedral equation (2.3) in terms of theta constants
of order five. Hence, we have the following decomposition formula
of the elliptic modular function $j$ in terms of the icosahedral
invariants $f$, $H$ and $T$ over the modular curve $X(5)$:
$$\aligned
 &j(z): j(z)-1728: 1\\
=&H(x_1(z), x_2(z))^3: -T(x_1(z), x_2(z))^2: f(x_1(z), x_2(z))^5.
\endaligned\eqno{(2.8)}$$

  Similarly, if $\Gamma$ is the binary tetrahedral group, the
corresponding Dynkin curve is that of $E_6$, and
$\mathbb{C}^2/\Gamma \subset \mathbb{C}^3$ is the set of zeros
of the equation
$$x^4+y^3+z^2=0.\eqno{(2.9)}$$
In his book \cite{K}, Klein also gave a parametric solution
of the above singularity (2.9) by homogeneous polynomials $\chi$,
$W$, $t$ in two variables of degrees $12$, $8$, $6$ with integral
coefficients, where
$$\Phi=z_1^4+2 \sqrt{-3} z_1^2 z_2^2+z_2^4, \quad
  \Psi=z_1^4-2 \sqrt{-3} z_1^2 z_2^2+z_2^4,$$
and
$$\left\{\aligned
  t &=z_1 z_2 (z_1^4-z_2^4),\\
  W &=\Phi \Psi=z_1^8+14 z_1^4 z_2^4+z_2^8,\\
  \chi &=z_1^{12}-33 z_1^8 z_2^4-33 z_1^4 z_2^8+z_2^{12}.
\endaligned\right.$$
They satisfy the (binary) tetrahedral equation
$$\chi^2-W^3+108 t^4=0.\eqno{(2.10)}$$
In fact, $t$, $W$ and $\chi$ are invariant polynomials under the
action of the binary tetrahedral group. We see this identity as
well as (2.9) as the defining relation between three generators
$t$, $W$ and $\chi$ of the ring of invariants
$\mathbb{C}[z_1, z_2]^{\Gamma}$ of the binary tetrahedral group
$\Gamma$ acting on $\mathbb{C}^2$, and we identify this ring with
the ring of functions on the affine variety $\mathbb{C}^2/\Gamma$
embedded in $\mathbb{C}^3$ and given by such an equation. Namely,
$$\mathbb{C}[z_1, z_2]^{\Gamma} \cong \mathbb{C}[t, W, \chi]/
  (\chi^2-W^3+108 t^4).\eqno{(2.11)}$$

  In the end, if $\Gamma$ is the binary octahedral group, the
corresponding Dynkin curve is that of $E_7$, and
$\mathbb{C}^2/\Gamma \subset \mathbb{C}^3$ is the set of zeros
of the equation
$$x^3 y+y^3+z^2=0.\eqno{(2.12)}$$
In his book \cite{K}, Klein gave a parametric solution of the above
singularity (2.12) by homogeneous polynomials $f_1$, $f_2$, $f_3$ in
two variables of degrees $18$, $8$, $12$ with integral coefficients,
where
$$\left\{\aligned
  f_1 &=t \chi=z_1 z_2 (z_1^4-z_2^4)(z_1^{12}
        -33 z_1^8 z_2^4-33 z_1^4 z_2^8+z_2^{12}),\\
  f_2 &=W=z_1^8+14 z_1^4 z_2^4+z_2^8,\\
  f_3 &=t^2=z_1^2 z_2^2 (z_1^4-z_2^4)^2.
\endaligned\right.$$
They satisfy the (binary) octahedral equation
$$f_1^2=f_3 \cdot f_2^3-108 f_3^3.\eqno{(2.13)}$$
In fact, $f_1$, $f_2$ and $f_3$ are invariant polynomials under the
action of the binary octahedral group. We see this identity as
well as (2.12) as the defining relation between three generators
$f_1$, $f_2$ and $f_3$ of the ring of invariants
$\mathbb{C}[z_1, z_2]^{\Gamma}$ of the binary octahedral group
$\Gamma$ acting on $\mathbb{C}^2$, and we identify this ring with
the ring of functions on the affine variety $\mathbb{C}^2/\Gamma$
embedded in $\mathbb{C}^3$ and given by such an equation. Namely,
$$\mathbb{C}[z_1, z_2]^{\Gamma} \cong \mathbb{C}[f_1, f_2, f_3]/
  (f_1^2-f_3 f_2^3+108 f_3^3).\eqno{(2.14)}$$

  In their subsequent monographs on elliptic modular functions \cite{KF1},
Section 3, Chap. 4 and \cite{KF2}, Section 5, Chap. 4, Klein and Fricke
gave the modular parametrizations for the equations of $E_6$, $E_7$ and
$E_8$-singularities by modular curves $X(3)$, $X(4)$ and $X(5)$, respectively.
Let $\theta_0$ and $\theta_1$ be two modular forms of weight $\frac{1}{2}$
for the principal congruence subgroup $\Gamma(4)$ of level $4$ given by
$$\theta_0(z)=\sum_{x \in 2 \mathbb{Z}} q^{\frac{x^2}{4}}, \quad
  \theta_1(z)=\sum_{x \in 2 \mathbb{Z}+1} q^{\frac{x^2}{4}}.$$
Then
$$\left\{\aligned
  t(\theta_0(z), \theta_1(z))^4 &=16 \Delta(z),\\
  W(\theta_0(z), \theta_1(z)) &=E_4(z),\\
  \chi(\theta_0(z), \theta_1(z)) &=E_6(z).
\endaligned\right.$$
This gives a modular parametrization for the equation (2.10) as well as
(2.13) by the modular curve $X(4)$.

  Similarly, let $\theta_0$ and $\theta_1$ be two modular forms of weight
$1$ for the principal congruence subgroup $\Gamma(3)$ of level $3$ given by
$$\theta_0(z)=\sum_{(x, y) \in \mathbb{Z}} q^{x^2-xy+y^2}, \quad
  \theta_1(z)=q^{\frac{1}{3}} \sum_{(x, y) \in \mathbb{Z}} q^{x^2-xy+y^2+x-y}.$$
Put
$$\left\{\aligned
  F_1 &=z_1^4+8 z_1 z_2^3,\\
  F_2 &=4 (z_2^4-z_1^3 z_2),\\
  F_3 &=z_1^6-20 z_1^3 z_2^3-8 z_2^6.
\endaligned\right.$$
Then
$$\left\{\aligned
  F_1(\theta_0(z), \theta_1(z)) &=E_4(z),\\
  F_2(\theta_0(z), \theta_1(z))^3 &=-1728 \Delta(z),\\
  F_3(\theta_0(z), \theta_1(z)) &=E_6(z).
\endaligned\right.$$
This gives a modular parametrization for the equation
$$F_1^3+F_2^3=F_3^2\eqno{(2.15)}$$
by the modular curve $X(3)$. Moreover, let
$$f_1=F_1 F_2, \quad f_2=F_3, \quad f_3=F_1^3-F_2^3.$$
Then
$$f_2^4=(F_1^3+F_2^3)^2=(F_1^3-F_2^3)^2+4 F_1^3 F_2^3
 =f_3^2+4 f_1^3.\eqno{(2.16)}$$
Let
$$g_1=f_2 f_3, \quad g_2=f_1, \quad g_3=f_2^2.$$
Then
$$g_1^2=g_3^3-4 g_2^3 g_3.\eqno{(2.17)}$$
This implies that the equation (2.10) as well as (2.13) can be
parameterized (over $\mathbb{C}$) by the modular curve $X(3)$.

  Note that the equations of $E_6$ and $E_7$-singularities possess two
distinct modular parameterizations by modular curves $X(3)$ and $X(4)$,
respectively. However, the equations (2.10) and (2.16) are not equivalent
over $\mathbb{Q}$. On the other hand, the equation of $E_8$-singularity
possesses only one modular parametrization by the modular curve $X(5)$.

\begin{center}
{\large\bf 3. Invariant theory and modular forms for $\text{SL}(2, 13)$ I:
              the quadratic and cubic invariants}
\end{center}

  The representation of $\text{SL}(2, 13)$ which we will consider is
the unique six-dimensional irreducible complex representation for which
the eigenvalues of $\left(\begin{matrix} 1 & 1\\ 0 & 1 \end{matrix}\right)$
are the $\exp (\underline{a} . 2 \pi i/13)$ for $\underline{a}$ a
non-square mod $13$. We will give an explicit realization of this
representation. This explicit realization will play a major role for
giving a complete system of invariants associated to $\text{SL}(2, 13)$.
At first, we will study the six-dimensional representation of
the finite group $\text{SL}(2, 13)$ of order $2184$, which acts
on the five-dimensional projective space
$\mathbb{P}^5=\{ (z_1, z_2, z_3, z_4, z_5, z_6): z_i \in \mathbb{C}
 \quad (i=1, 2, 3, 4, 5, 6) \}$. This representation is defined
over the cyclotomic field $\mathbb{Q}(e^{\frac{2 \pi i}{13}})$.
Put
$$S=-\frac{1}{\sqrt{13}} \begin{pmatrix}
  \zeta^{12}-\zeta & \zeta^{10}-\zeta^3 & \zeta^4-\zeta^9 &
  \zeta^5-\zeta^8 & \zeta^2-\zeta^{11} & \zeta^6-\zeta^7\\
  \zeta^{10}-\zeta^3 & \zeta^4-\zeta^9 & \zeta^{12}-\zeta &
  \zeta^2-\zeta^{11} & \zeta^6-\zeta^7 & \zeta^5-\zeta^8\\
  \zeta^4-\zeta^9 & \zeta^{12}-\zeta & \zeta^{10}-\zeta^3 &
  \zeta^6-\zeta^7 & \zeta^5-\zeta^8 & \zeta^2-\zeta^{11}\\
  \zeta^5-\zeta^8 & \zeta^2-\zeta^{11} & \zeta^6-\zeta^7 &
  \zeta-\zeta^{12} & \zeta^3-\zeta^{10} & \zeta^9-\zeta^4\\
  \zeta^2-\zeta^{11} & \zeta^6-\zeta^7 & \zeta^5-\zeta^8 &
  \zeta^3-\zeta^{10} & \zeta^9-\zeta^4 & \zeta-\zeta^{12}\\
  \zeta^6-\zeta^7 & \zeta^5-\zeta^8 & \zeta^2-\zeta^{11} &
  \zeta^9-\zeta^4 & \zeta-\zeta^{12} & \zeta^3-\zeta^{10}
\end{pmatrix}\eqno{(3.1)}$$
and
$$T=\text{diag}(\zeta^7, \zeta^{11}, \zeta^8, \zeta^6,
    \zeta^2, \zeta^5),\eqno{(3.2)}$$
where $\zeta=\exp(2 \pi i/13)$. We have
$$S^2=-I, \quad T^{13}=(ST)^3=I.\eqno{(3.3)}$$
In \cite{Y1}, we put $P=S T^{-1} S$ and $Q=S T^3$. Then
$(Q^3 P^4)^3=-I$ (see \cite{Y1}, the proof of Theorem 3.1).
Let $G=\langle S, T \rangle$, then $G \cong \text{SL}(2, 13)$.

  Put $\theta_1=\zeta+\zeta^3+\zeta^9$,
$\theta_2=\zeta^2+\zeta^6+\zeta^5$,
$\theta_3=\zeta^4+\zeta^{12}+\zeta^{10}$,
and $\theta_4=\zeta^8+\zeta^{11}+\zeta^7$. We find that
$$\left\{\aligned
  &\theta_1+\theta_2+\theta_3+\theta_4=-1,\\
  &\theta_1 \theta_2+\theta_1 \theta_3+\theta_1 \theta_4+
   \theta_2 \theta_3+\theta_2 \theta_4+\theta_3 \theta_4=2,\\
  &\theta_1 \theta_2 \theta_3+\theta_1 \theta_2 \theta_4+
   \theta_1 \theta_3 \theta_4+\theta_2 \theta_3 \theta_4=4,\\
  &\theta_1 \theta_2 \theta_3 \theta_4=3.
\endaligned\right.$$
Hence, $\theta_1$, $\theta_2$, $\theta_3$ and $\theta_4$ satisfy
the quartic equation $z^4+z^3+2 z^2-4z+3=0$,
which can be decomposed as two quadratic equations
$$\left(z^2+\frac{1+\sqrt{13}}{2} z+\frac{5+\sqrt{13}}{2}\right)
  \left(z^2+\frac{1-\sqrt{13}}{2} z+\frac{5-\sqrt{13}}{2}\right)=0$$
over the real quadratic field $\mathbb{Q}(\sqrt{13})$. Therefore, the
four roots are given as follows:
$$\left\{\aligned
  \theta_1=\frac{1}{4} \left(-1+\sqrt{13}+\sqrt{-26+6 \sqrt{13}}\right),\\
  \theta_2=\frac{1}{4} \left(-1-\sqrt{13}+\sqrt{-26-6 \sqrt{13}}\right),\\
  \theta_3=\frac{1}{4} \left(-1+\sqrt{13}-\sqrt{-26+6 \sqrt{13}}\right),\\
  \theta_4=\frac{1}{4} \left(-1-\sqrt{13}-\sqrt{-26-6 \sqrt{13}}\right).
\endaligned\right.$$
Moreover, we find that
$$\left\{\aligned
  \theta_1+\theta_3+\theta_2+\theta_4 &=-1,\\
  \theta_1+\theta_3-\theta_2-\theta_4 &=\sqrt{13},\\
  \theta_1-\theta_3-\theta_2+\theta_4 &=-\sqrt{-13+2 \sqrt{13}},\\
  \theta_1-\theta_3+\theta_2-\theta_4 &=\sqrt{-13-2 \sqrt{13}}.
\endaligned\right.$$

  Let us study the action of $S T^{\nu}$ on $\mathbb{P}^5$, where
$\nu=0, 1, \ldots, 12$. Put
$$\alpha=\zeta+\zeta^{12}-\zeta^5-\zeta^8, \quad
   \beta=\zeta^3+\zeta^{10}-\zeta^2-\zeta^{11}, \quad
   \gamma=\zeta^9+\zeta^4-\zeta^6-\zeta^7.$$
We find that
$$\aligned
  &13 ST^{\nu}(z_1) \cdot ST^{\nu}(z_4)\\
=&\beta z_1 z_4+\gamma z_2 z_5+\alpha z_3 z_6+\\
 &+\gamma \zeta^{\nu} z_1^2+\alpha \zeta^{9 \nu} z_2^2+\beta \zeta^{3 \nu} z_3^2
  -\gamma \zeta^{12 \nu} z_4^2-\alpha \zeta^{4 \nu} z_5^2-\beta \zeta^{10 \nu}
  z_6^2+\\
 &+(\alpha-\beta) \zeta^{5 \nu} z_1 z_2+(\beta-\gamma) \zeta^{6 \nu} z_2 z_3
  +(\gamma-\alpha) \zeta^{2 \nu} z_1 z_3+\\
 &+(\beta-\alpha) \zeta^{8 \nu} z_4 z_5+(\gamma-\beta) \zeta^{7 \nu} z_5 z_6
  +(\alpha-\gamma) \zeta^{11 \nu} z_4 z_6+\\
 &-(\alpha+\beta) \zeta^{\nu} z_3 z_4-(\beta+\gamma) \zeta^{9 \nu} z_1 z_5
  -(\gamma+\alpha) \zeta^{3 \nu} z_2 z_6+\\
 &-(\alpha+\beta) \zeta^{12 \nu} z_1 z_6-(\beta+\gamma) \zeta^{4 \nu} z_2 z_4
  -(\gamma+\alpha) \zeta^{10 \nu} z_3 z_5.
\endaligned$$
$$\aligned
  &13 ST^{\nu}(z_2) \cdot ST^{\nu}(z_5)\\
=&\gamma z_1 z_4+\alpha z_2 z_5+\beta z_3 z_6+\\
 &+\alpha \zeta^{\nu} z_1^2+\beta \zeta^{9 \nu} z_2^2+\gamma \zeta^{3 \nu} z_3^2
  -\alpha \zeta^{12 \nu} z_4^2-\beta \zeta^{4 \nu} z_5^2-\gamma \zeta^{10 \nu}
  z_6^2+\\
 &+(\beta-\gamma) \zeta^{5 \nu} z_1 z_2+(\gamma-\alpha) \zeta^{6 \nu} z_2 z_3
  +(\alpha-\beta) \zeta^{2 \nu} z_1 z_3+\\
 &+(\gamma-\beta) \zeta^{8 \nu} z_4 z_5+(\alpha-\gamma) \zeta^{7 \nu} z_5 z_6
  +(\beta-\alpha) \zeta^{11 \nu} z_4 z_6+\\
 &-(\beta+\gamma) \zeta^{\nu} z_3 z_4-(\gamma+\alpha) \zeta^{9 \nu} z_1 z_5
  -(\alpha+\beta) \zeta^{3 \nu} z_2 z_6+\\
 &-(\beta+\gamma) \zeta^{12 \nu} z_1 z_6-(\gamma+\alpha) \zeta^{4 \nu} z_2 z_4
  -(\alpha+\beta) \zeta^{10 \nu} z_3 z_5.
\endaligned$$
$$\aligned
  &13 ST^{\nu}(z_3) \cdot ST^{\nu}(z_6)\\
=&\alpha z_1 z_4+\beta z_2 z_5+\gamma z_3 z_6+\\
 &+\beta \zeta^{\nu} z_1^2+\gamma \zeta^{9 \nu} z_2^2+\alpha \zeta^{3 \nu} z_3^2
  -\beta \zeta^{12 \nu} z_4^2-\gamma \zeta^{4 \nu} z_5^2-\alpha \zeta^{10 \nu}
  z_6^2+\\
 &+(\gamma-\alpha) \zeta^{5 \nu} z_1 z_2+(\alpha-\beta) \zeta^{6 \nu} z_2 z_3
  +(\beta-\gamma) \zeta^{2 \nu} z_1 z_3+\\
 &+(\alpha-\gamma) \zeta^{8 \nu} z_4 z_5+(\beta-\alpha) \zeta^{7 \nu} z_5 z_6
  +(\gamma-\beta) \zeta^{11 \nu} z_4 z_6+\\
 &-(\gamma+\alpha) \zeta^{\nu} z_3 z_4-(\alpha+\beta) \zeta^{9 \nu} z_1 z_5
  -(\beta+\gamma) \zeta^{3 \nu} z_2 z_6+\\
 &-(\gamma+\alpha) \zeta^{12 \nu} z_1 z_6-(\alpha+\beta) \zeta^{4 \nu} z_2 z_4
  -(\beta+\gamma) \zeta^{10 \nu} z_3 z_5.
\endaligned$$
Note that $\alpha+\beta+\gamma=\sqrt{13}$, we find that
$$\aligned
  &\sqrt{13} \left[ST^{\nu}(z_1) \cdot ST^{\nu}(z_4)+
   ST^{\nu}(z_2) \cdot ST^{\nu}(z_5)+ST^{\nu}(z_3) \cdot ST^{\nu}(z_6)\right]\\
 =&(z_1 z_4+z_2 z_5+z_3 z_6)+(\zeta^{\nu} z_1^2+\zeta^{9 \nu} z_2^2+\zeta^{3 \nu}
   z_3^2)-(\zeta^{12 \nu} z_4^2+\zeta^{4 \nu} z_5^2+\zeta^{10 \nu} z_6^2)+\\
  &-2(\zeta^{\nu} z_3 z_4+\zeta^{9 \nu} z_1 z_5+\zeta^{3 \nu} z_2 z_6)
   -2(\zeta^{12 \nu} z_1 z_6+\zeta^{4 \nu} z_2 z_4+\zeta^{10 \nu} z_3 z_5).
\endaligned$$
Let
$$\varphi_{\infty}(z_1, z_2, z_3, z_4, z_5, z_6)=\sqrt{13} (z_1 z_4+z_2 z_5+z_3 z_6)
  \eqno{(3.4)}$$
and
$$\varphi_{\nu}(z_1, z_2, z_3, z_4, z_5, z_6)=\varphi_{\infty}(ST^{\nu}(z_1, z_2,
                                              z_3, z_4, z_5, z_6))\eqno{(3.5)}$$
for $\nu=0, 1, \ldots, 12$. Then
$$\aligned
  \varphi_{\nu}
=&(z_1 z_4+z_2 z_5+z_3 z_6)+\zeta^{\nu} (z_1^2-2 z_3 z_4)+\zeta^{4 \nu}
  (-z_5^2-2 z_2 z_4)+\\
 &+\zeta^{9 \nu} (z_2^2-2 z_1 z_5)+\zeta^{3 \nu} (z_3^2-2 z_2 z_6)+
   \zeta^{12 \nu} (-z_4^2-2 z_1 z_6)+\\
 &+\zeta^{10 \nu} (-z_6^2-2 z_3 z_5).
\endaligned\eqno{(3.6)}$$
This leads us to define the following senary quadratic forms
(quadratic forms in six variables):
$$\left\{\aligned
  \mathbf{A}_0 &=z_1 z_4+z_2 z_5+z_3 z_6,\\
  \mathbf{A}_1 &=z_1^2-2 z_3 z_4,\\
  \mathbf{A}_2 &=-z_5^2-2 z_2 z_4,\\
  \mathbf{A}_3 &=z_2^2-2 z_1 z_5,\\
  \mathbf{A}_4 &=z_3^2-2 z_2 z_6,\\
  \mathbf{A}_5 &=-z_4^2-2 z_1 z_6,\\
  \mathbf{A}_6 &=-z_6^2-2 z_3 z_5.
\endaligned\right.\eqno{(3.7)}$$
Hence,
$$\sqrt{13} ST^{\nu}(\mathbf{A}_0)=\mathbf{A}_0+\zeta^{\nu} \mathbf{A}_1
  +\zeta^{4 \nu} \mathbf{A}_2+\zeta^{9 \nu} \mathbf{A}_3+\zeta^{3 \nu}
  \mathbf{A}_4+\zeta^{12 \nu} \mathbf{A}_5+\zeta^{10 \nu} \mathbf{A}_6.
  \eqno{(3.8)}$$
Let $H:=Q^5 P^2 \cdot P^2 Q^6 P^8 \cdot Q^5 P^2 \cdot P^3 Q$ where
$P=S T^{-1} S$ and $Q=S T^3$. Then (see \cite{Y2}, p.27)
$$H=\begin{pmatrix}
  0 &  0 &  0 & 0 & 0 & 1\\
  0 &  0 &  0 & 1 & 0 & 0\\
  0 &  0 &  0 & 0 & 1 & 0\\
  0 &  0 & -1 & 0 & 0 & 0\\
 -1 &  0 &  0 & 0 & 0 & 0\\
  0 & -1 &  0 & 0 & 0 & 0
\end{pmatrix}.\eqno{(3.9)}$$
Note that $H^6=1$ and $H^{-1} T H=-T^4$. Thus,
$\langle H, T \rangle \cong \mathbb{Z}_{13} \rtimes \mathbb{Z}_6$.
Hence, it is a maximal subgroup of order $78$ of $\text{PSL}(2, 13)$
with index $14$. We find that $\varphi_{\infty}^2$ is invariant under
the action of the maximal subgroup $\langle H, T \rangle$. Note that
$$\varphi_{\infty}=\sqrt{13} \mathbf{A}_0, \quad
  \varphi_{\nu}=\mathbf{A}_0+\zeta^{\nu} \mathbf{A}_1+
  \zeta^{4 \nu} \mathbf{A}_2+\zeta^{9 \nu} \mathbf{A}_3+
  \zeta^{3 \nu} \mathbf{A}_4+\zeta^{12 \nu} \mathbf{A}_5+
  \zeta^{10 \nu} \mathbf{A}_6$$
for $\nu=0, 1, \ldots, 12$. Let $w=\varphi^2$,
$w_{\infty}=\varphi_{\infty}^2$ and $w_{\nu}=\varphi_{\nu}^2$.
Then $w_{\infty}$, $w_{\nu}$ for $\nu=0, \ldots, 12$ form an
algebraic equation of degree fourteen, which is just the Jacobian
equation of degree fourteen (see \cite{K}, pp.161-162), whose
roots are these $w_{\nu}$ and $w_{\infty}$:
$$w^{14}+a_1 w^{13}+\cdots+a_{13} w+a_{14}=0.$$
In particular, the coefficients
$$-a_1=w_{\infty}+\sum_{\nu=0}^{12} w_{\nu}=
  26 ({\bf A}_0^2+{\bf A}_1 {\bf A}_5+{\bf A}_2 {\bf A}_3+{\bf A}_4 {\bf A}_6).\eqno{(3.10)}$$
This leads to an invariant quadric
$\Psi_2:={\bf A}_0^2+{\bf A}_1 {\bf A}_5+{\bf A}_2 {\bf A}_3+{\bf A}_4 {\bf A}_6
    =2 \Phi_4(z_1, z_2, z_3, z_4, z_5, z_6)$,
where
$$\begin{array}{rl}
  \Phi_4:=&(z_3 z_4^3+z_1 z_5^3+z_2 z_6^3)-(z_6 z_1^3+z_4 z_2^3+z_5 z_3^3)+\\
          &+3(z_1 z_2 z_4 z_5+z_2 z_3 z_5 z_6+z_3 z_1 z_6 z_4),
\end{array}\eqno{(3.11)}$$
Hence, the variety $\Phi_4=0$ is a quartic four-fold, which
is invariant under the action of the group $G$.

  On the other hand, we have
$$\aligned
  &-13 \sqrt{13} ST^{\nu}(z_1) \cdot ST^{\nu}(z_2) \cdot ST^{\nu}(z_3)\\
 =&-r_4 (\zeta^{8 \nu} z_1^3+\zeta^{7 \nu} z_2^3+\zeta^{11 \nu} z_3^3)
   -r_2 (\zeta^{5 \nu} z_4^3+\zeta^{6 \nu} z_5^3+\zeta^{2 \nu} z_6^3)\\
  &-r_3 (\zeta^{12 \nu} z_1^2 z_2+\zeta^{4 \nu} z_2^2 z_3+\zeta^{10 \nu} z_3^2 z_1)
   -r_1 (\zeta^{\nu} z_4^2 z_5+\zeta^{9 \nu} z_5^2 z_6+\zeta^{3 \nu} z_6^2 z_4)\\
  &+2 r_1 (\zeta^{3 \nu} z_1 z_2^2+\zeta^{\nu} z_2 z_3^2+\zeta^{9 \nu} z_3 z_1^2)
   -2 r_3 (\zeta^{10 \nu} z_4 z_5^2+\zeta^{12 \nu} z_5 z_6^2+\zeta^{4 \nu} z_6 z_4^2)\\
  &+2 r_4 (\zeta^{7 \nu} z_1^2 z_4+\zeta^{11 \nu} z_2^2 z_5+\zeta^{8 \nu} z_3^2 z_6)
   -2 r_2 (\zeta^{6 \nu} z_1 z_4^2+\zeta^{2 \nu} z_2 z_5^2+\zeta^{5 \nu} z_3 z_6^2)+\\
  &+r_1 (\zeta^{3 \nu} z_1^2 z_5+\zeta^{\nu} z_2^2 z_6+\zeta^{9 \nu} z_3^2 z_4)
   +r_3 (\zeta^{10 \nu} z_2 z_4^2+\zeta^{12 \nu} z_3 z_5^2+\zeta^{4 \nu} z_1 z_6^2)+\\
  &+r_2 (\zeta^{6 \nu} z_1^2 z_6+\zeta^{2 \nu} z_2^2 z_4+\zeta^{5 \nu} z_3^2 z_5)
   +r_4 (\zeta^{7 \nu} z_3 z_4^2+\zeta^{11 \nu} z_1 z_5^2+\zeta^{8 \nu} z_2 z_6^2)+\\
  &+r_0 z_1 z_2 z_3+r_{\infty} z_4 z_5 z_6+\\
  &-r_4 (\zeta^{11 \nu} z_1 z_2 z_4+\zeta^{8 \nu} z_2 z_3 z_5+\zeta^{7 \nu} z_1 z_3 z_6)+\\
  &+r_2 (\zeta^{2 \nu} z_1 z_4 z_5+\zeta^{5 \nu} z_2 z_5 z_6+\zeta^{6 \nu} z_3 z_4 z_6)+\\
  &-3 r_4 (\zeta^{7 \nu} z_1 z_2 z_5+\zeta^{11 \nu} z_2 z_3 z_6+\zeta^{8 \nu} z_1 z_3 z_4)+\\
  &+3 r_2 (\zeta^{6 \nu} z_2 z_4 z_5+\zeta^{2 \nu} z_3 z_5 z_6+\zeta^{5 \nu} z_1 z_4 z_6)+\\
  &-r_3 (\zeta^{10 \nu} z_1 z_2 z_6+\zeta^{4 \nu} z_1 z_3 z_5+\zeta^{12 \nu} z_2 z_3 z_4)+\\
  &+r_1 (\zeta^{3 \nu} z_3 z_4 z_5+\zeta^{9 \nu} z_2 z_4 z_6+\zeta^{\nu} z_1 z_5 z_6),
\endaligned$$
where
$$r_0=2(\theta_1-\theta_3)-3(\theta_2-\theta_4), \quad
  r_{\infty}=2(\theta_4-\theta_2)-3(\theta_1-\theta_3),$$
$$r_1=\sqrt{-13-2 \sqrt{13}}, \quad r_2=\sqrt{\frac{-13+3 \sqrt{13}}{2}},$$
$$r_3=\sqrt{-13+2 \sqrt{13}}, \quad r_4=\sqrt{\frac{-13-3 \sqrt{13}}{2}}.$$
This leads us to define the following senary cubic forms (cubic forms in
six variables):
$$\left\{\aligned
  \mathbf{D}_0 &=z_1 z_2 z_3,\\
  \mathbf{D}_1 &=2 z_2 z_3^2+z_2^2 z_6-z_4^2 z_5+z_1 z_5 z_6,\\
  \mathbf{D}_2 &=-z_6^3+z_2^2 z_4-2 z_2 z_5^2+z_1 z_4 z_5+3 z_3 z_5 z_6,\\
  \mathbf{D}_3 &=2 z_1 z_2^2+z_1^2 z_5-z_4 z_6^2+z_3 z_4 z_5,\\
  \mathbf{D}_4 &=-z_2^2 z_3+z_1 z_6^2-2 z_4^2 z_6-z_1 z_3 z_5,\\
  \mathbf{D}_5 &=-z_4^3+z_3^2 z_5-2 z_3 z_6^2+z_2 z_5 z_6+3 z_1 z_4 z_6,\\
  \mathbf{D}_6 &=-z_5^3+z_1^2 z_6-2 z_1 z_4^2+z_3 z_4 z_6+3 z_2 z_4 z_5,\\
  \mathbf{D}_7 &=-z_2^3+z_3 z_4^2-z_1 z_3 z_6-3 z_1 z_2 z_5+2 z_1^2 z_4,\\
  \mathbf{D}_8 &=-z_1^3+z_2 z_6^2-z_2 z_3 z_5-3 z_1 z_3 z_4+2 z_3^2 z_6,\\
  \mathbf{D}_9 &=2 z_1^2 z_3+z_3^2 z_4-z_5^2 z_6+z_2 z_4 z_6,\\
  \mathbf{D}_{10} &=-z_1 z_3^2+z_2 z_4^2-2 z_4 z_5^2-z_1 z_2 z_6,\\
  \mathbf{D}_{11} &=-z_3^3+z_1 z_5^2-z_1 z_2 z_4-3 z_2 z_3 z_6+2 z_2^2 z_5,\\
  \mathbf{D}_{12} &=-z_1^2 z_2+z_3 z_5^2-2 z_5 z_6^2-z_2 z_3 z_4,\\
  \mathbf{D}_{\infty}&=z_4 z_5 z_6.
\endaligned\right.\eqno{(3.12)}$$
Then
$$\aligned
  &-13 \sqrt{13} ST^{\nu}(\mathbf{D}_0)\\
 =&r_0 \mathbf{D}_0+r_1 \zeta^{\nu} \mathbf{D}_1+
   r_2 \zeta^{2 \nu} \mathbf{D}_2+
   r_1 \zeta^{3 \nu} \mathbf{D}_3+r_3 \zeta^{4 \nu} \mathbf{D}_4+\\
  &+r_2 \zeta^{5 \nu} \mathbf{D}_5+r_2 \zeta^{6 \nu} \mathbf{D}_6+
   r_4 \zeta^{7 \nu} \mathbf{D}_7+r_4 \zeta^{8 \nu} \mathbf{D}_8+\\
  &+r_1 \zeta^{9 \nu} \mathbf{D}_9+r_3 \zeta^{10 \nu} \mathbf{D}_{10}
   +r_4 \zeta^{11 \nu} \mathbf{D}_{11}+r_3 \zeta^{12 \nu} \mathbf{D}_{12}
   +r_{\infty} \mathbf{D}_{\infty}.
\endaligned$$
$$\aligned
  &-13 \sqrt{13} ST^{\nu}(\mathbf{D}_{\infty})\\
 =&r_{\infty} \mathbf{D}_0-r_3 \zeta^{\nu} \mathbf{D}_1-
   r_4 \zeta^{2 \nu} \mathbf{D}_2-r_3 \zeta^{3 \nu} \mathbf{D}_3+
   r_1 \zeta^{4 \nu} \mathbf{D}_4+\\
  &-r_4 \zeta^{5 \nu} \mathbf{D}_5-r_4 \zeta^{6 \nu} \mathbf{D}_6+
   r_2 \zeta^{7 \nu} \mathbf{D}_7+r_2 \zeta^{8 \nu} \mathbf{D}_8+\\
  &-r_3 \zeta^{9 \nu} \mathbf{D}_9+r_1 \zeta^{10 \nu} \mathbf{D}_{10}+
   r_2 \zeta^{11 \nu} \mathbf{D}_{11}+r_1 \zeta^{12 \nu} \mathbf{D}_{12}-
   r_0 \mathbf{D}_{\infty}.
\endaligned$$

  Let
$$\delta_{\infty}(z_1, z_2, z_3, z_4, z_5, z_6)
 =13^2 (z_1^2 z_2^2 z_3^2+z_4^2 z_5^2 z_6^2)\eqno{(3.13)}$$
and
$$\delta_{\nu}(z_1, z_2, z_3, z_4, z_5, z_6)
 =\delta_{\infty}(ST^{\nu}(z_1, z_2, z_3, z_4, z_5, z_6))\eqno{(3.14)}$$
for $\nu=0, 1, \ldots, 12$. Then
$$\delta_{\nu}=13^2 ST^{\nu}(\mathbf{G}_0)
 =-13 \mathbf{G}_0+\zeta^{\nu} \mathbf{G}_1+\zeta^{2 \nu} \mathbf{G}_2+
  \cdots+\zeta^{12 \nu} \mathbf{G}_{12},\eqno{(3.15)}$$
where the senary sextic forms (i.e., sextic forms in six
variables) are given as follows:
$$\left\{\aligned
  \mathbf{G}_0 =&\mathbf{D}_0^2+\mathbf{D}_{\infty}^2,\\
  \mathbf{G}_1 =&-\mathbf{D}_7^2+2 \mathbf{D}_0 \mathbf{D}_1+10 \mathbf{D}_{\infty}
                 \mathbf{D}_1+2 \mathbf{D}_2 \mathbf{D}_{12}+\\
                &-2 \mathbf{D}_3 \mathbf{D}_{11}-4 \mathbf{D}_4 \mathbf{D}_{10}-2
                 \mathbf{D}_9 \mathbf{D}_5,\\
  \mathbf{G}_2 =&-2 \mathbf{D}_1^2-4 \mathbf{D}_0 \mathbf{D}_2+6 \mathbf{D}_{\infty}
                 \mathbf{D}_2-2 \mathbf{D}_4 \mathbf{D}_{11}+\\
                &+2 \mathbf{D}_5 \mathbf{D}_{10}-2 \mathbf{D}_6 \mathbf{D}_9-2
                 \mathbf{D}_7 \mathbf{D}_8,\\
  \mathbf{G}_3 =&-\mathbf{D}_8^2+2 \mathbf{D}_0 \mathbf{D}_3+10 \mathbf{D}_{\infty}
                 \mathbf{D}_3+2 \mathbf{D}_6 \mathbf{D}_{10}+\\
                &-2 \mathbf{D}_9 \mathbf{D}_7-4 \mathbf{D}_{12} \mathbf{D}_4-2
                 \mathbf{D}_1 \mathbf{D}_2,\\
  \mathbf{G}_4 =&-\mathbf{D}_2^2+10 \mathbf{D}_0 \mathbf{D}_4-2 \mathbf{D}_{\infty}
                 \mathbf{D}_4+2 \mathbf{D}_5 \mathbf{D}_{12}+\\
                &-2 \mathbf{D}_9 \mathbf{D}_8-4 \mathbf{D}_1 \mathbf{D}_3-2
                 \mathbf{D}_{10} \mathbf{D}_7,\\
  \mathbf{G}_5 =&-2 \mathbf{D}_9^2-4 \mathbf{D}_0 \mathbf{D}_5+6 \mathbf{D}_{\infty}
                 \mathbf{D}_5-2 \mathbf{D}_{10} \mathbf{D}_8+\\
                &+2 \mathbf{D}_6 \mathbf{D}_{12}-2 \mathbf{D}_2 \mathbf{D}_3-2
                 \mathbf{D}_{11} \mathbf{D}_7,\\
  \mathbf{G}_6 =&-2 \mathbf{D}_3^2-4 \mathbf{D}_0 \mathbf{D}_6+6 \mathbf{D}_{\infty}
                 \mathbf{D}_6-2 \mathbf{D}_{12} \mathbf{}_7+\\
                &+2 \mathbf{D}_2 \mathbf{D}_4-2 \mathbf{D}_5 \mathbf{D}_1-2
                 \mathbf{D}_8 \mathbf{D}_{11},\\
  \mathbf{G}_7 =&-2 \mathbf{D}_{10}^2+6 \mathbf{D}_0 \mathbf{D}_7+4 \mathbf{D}_{\infty}
                 \mathbf{D}_7-2 \mathbf{D}_1 \mathbf{D}_6+\\
                &-2 \mathbf{D}_2 \mathbf{D}_5-2 \mathbf{D}_8 \mathbf{D}_{12}-2
                 \mathbf{D}_9 \mathbf{D}_{11},\\
  \mathbf{G}_8 =&-2 \mathbf{D}_4^2+6 \mathbf{D}_0 \mathbf{D}_8+4 \mathbf{D}_{\infty}
                 \mathbf{D}_8-2 \mathbf{D}_3 \mathbf{D}_5+\\
                &-2 \mathbf{D}_6 \mathbf{D}_2-2 \mathbf{D}_{11} \mathbf{D}_{10}-2
                 \mathbf{D}_1 \mathbf{D}_7,\\
  \mathbf{G}_9 =&-\mathbf{D}_{11}^2+2 \mathbf{D}_0 \mathbf{D}_9+10 \mathbf{D}_{\infty}
                 \mathbf{D}_9+2 \mathbf{D}_5 \mathbf{D}_4+\\
                &-2 \mathbf{D}_1 \mathbf{D}_8-4 \mathbf{D}_{10} \mathbf{D}_{12}-2
                 \mathbf{D}_3 \mathbf{D}_6,
\endaligned\right.\eqno{(3.16)}$$
$$\left\{\aligned
  \mathbf{G}_{10} =&-\mathbf{D}_5^2+10 \mathbf{D}_0 \mathbf{D}_{10}-2 \mathbf{D}_{\infty}
                    \mathbf{D}_{10}+2 \mathbf{D}_6 \mathbf{D}_4+\\
                   &-2 \mathbf{D}_3 \mathbf{D}_7-4 \mathbf{D}_9 \mathbf{D}_1-2
                    \mathbf{D}_{12} \mathbf{D}_{11},\\
  \mathbf{G}_{11} =&-2 \mathbf{D}_{12}^2+6 \mathbf{D}_0 \mathbf{D}_{11}+4 \mathbf{D}_{\infty}
                    \mathbf{D}_{11}-2 \mathbf{D}_9 \mathbf{D}_2+\\
                   &-2 \mathbf{D}_5 \mathbf{D}_6-2 \mathbf{D}_7 \mathbf{D}_4-2
                    \mathbf{D}_3 \mathbf{D}_8,\\
  \mathbf{G}_{12} =&-\mathbf{D}_6^2+10 \mathbf{D}_0 \mathbf{D}_{12}-2 \mathbf{D}_{\infty}
                    \mathbf{D}_{12}+2 \mathbf{D}_2 \mathbf{D}_{10}+\\
                   &-2 \mathbf{D}_1 \mathbf{D}_{11}-4 \mathbf{D}_3 \mathbf{D}_9-2
                    \mathbf{D}_4 \mathbf{D}_8.
\endaligned\right.\eqno{(3.16)}$$
We have that $\mathbf{G}_0$ is invariant under the action of
$\langle H, T \rangle$, a maximal subgroup of order $78$ of
$\text{PSL}(2, 13)$ with index $14$. Note that $\delta_{\infty}$,
$\delta_{\nu}$ for $\nu=0, \ldots, 12$ form an algebraic equation
of degree fourteen. However, we have
$$\Phi_6=\Phi_{0, 1}=\delta_{\infty}+\sum_{\nu=0}^{12} \delta_{\nu}=0.$$
Hence, it is not the Jacobian equation of degree fourteen.

\begin{center}
{\large\bf 4. Invariant theory and modular forms for $\text{SL}(2, 13)$ II:
              the quartic invariants and the modular curve $X(13)$}
\end{center}

  In this section, we will study the following:

\textbf{Problem 4.1.} {\it Give an explicit construction of the modular
curve of level $p$ from the invariant theory associated with
$\text{SL}(2, \mathbb{F}_p)$ and its representations using projective
algebraic geometry.}

  In \cite{K4}, Klein gave the following formula for the genus of $X(p)$:
$$g(X(p))=\frac{1}{24} (p+2) (p-3) (p-5).\eqno{(4.1)}$$
In particular,
$$g(X(7))=3, \quad g(X(11))=26, \quad g(X(13))=50.$$
Klein also obtained the following formula for the degree of $X(p)$ in the
projective coordinates $(z_1, \ldots, z_{(p-1)/2})$ (i.e., Klein's $z$-curve)
(see \cite{K4}):
$$\text{deg}(X(p))=\frac{1}{48} (p-3)(p^2-1).\eqno{(4.2)}$$
In particular,
$$\text{deg}(X(7))=4, \quad \text{deg}(X(11))=20, \quad \text{deg}(X(13))=35.$$

  The equations defining the universal elliptic curve with prime level
$p>3$ structure over the field $\mathbb{C}$ of complex numbers were
obtained by Klein in his paper on elliptic normal curves of order $p$
(see \cite{K5}, p.230, equation (77)). Explicitly, consider the space
$V$ of all complex-valued functions on the cyclic group
$\mathbb{Z}/p \mathbb{Z}$ of integers modulo $p$ and denote by
$\mathbb{P}(V)$ the projective space of lines in $V$. If $f$ is a
non-zero element of $V$, we will denote by $[f]$ the corresponding
element of $\mathbb{P}(V)$. Denote by $\mathcal{L}_p$ the locus of
all points $[f]$ of $\mathbb{P}(V)$ such that $f$ is an odd function
and such that
$$\aligned
 0 &=f(a+b) \cdot f(a-b) \cdot f(c+d) \cdot f(c-d)\\
   &+f(a+c) \cdot f(a-c) \cdot f(d+b) \cdot f(d-b)\\
   &+f(a+d) \cdot f(a-d) \cdot f(b+c) \cdot f(b-c)
\endaligned$$
for all $a$, $b$, $c$, $d$ in $\mathbb{Z}/p \mathbb{Z}$. Denote by
$\mathcal{E}_p$ the locus of all points $([z], [X])$ of $\mathcal{L}_p
\times \mathbb{P}(V)$ such that
$$\aligned
 0 &=z(a-b) \cdot z(c-d) \cdot X(a+b) \cdot X(c+d)\\
   &+z(a-c) \cdot z(d-b) \cdot X(a+c) \cdot X(d+b)\\
   &+z(a-d) \cdot z(b-c) \cdot X(a+d) \cdot X(b+c).
\endaligned$$
Note that $([z], [z])$ belongs to $\mathcal{E}_p$ for all $[z]$ in
$\mathcal{L}_p$. For generic $[z]$ in $\mathcal{L}_p$, the locus
$\mathcal{E}_p([z])$ of all points $([w], [X])$ of $\mathcal{E}_p$
with $[w]=[z]$ is an elliptic curve with $[z]$ as its identity
element. The locus $\mathcal{L}_p$ is the modular curve of level
$p$ and $\mathcal{E}_p$ is the universal elliptic curve with level
$p$ structure over $\mathcal{L}_p$. The points of order $p$ on the
elliptic curve $\mathcal{E}_n([z])$ correspond, via projection onto
the second factor of $\mathcal{L}_p \times \mathbb{P}(V)$, to the
orbit of $[z]$ under the action of the Heisenberg group of
$\mathbb{Z}/p \mathbb{Z}$ on $\mathbb{P}(V)$ (see \cite{Ad1}).

  In the modern terminology, denote by $L^2(\mathbb{F}_p)$ the
$p$-dimensional complex vector space of all square-integrable complex
valued functions on $\mathbb{F}_p$ with respect to counting measure,
that is, all functions from $\mathbb{F}_p$ to the complex numbers.
We can decompose $L^2(\mathbb{F}_p)$ as the direct sum of the space
$V^{+}$ of even functions and the space $V^{-}$ of odd functions.
The space $V^{-}$ has dimension $(p-1)/2$ and its associated projective
space $\mathbb{P}(V^{-})$ has dimension $(p-3)/2$. If $f$ is a nonzero
element of $V^{-}$, we will denote by $[f]$ the corresponding element
of $\mathbb{P}(V^{-})$, in keeping with the classical notation for
homogeneous coordinates (see \cite{Ad2}).

  Klein discovered the following general result:

\textbf{Theorem 4.2.} {\it The modular curve $X(p)$ is isomorphic to
the locus of all $[f]$ in $\mathbb{P}(V^{-})$ which for all $w$, $x$,
$y$, $z$ in $\mathbb{F}_p$ satisfy the identities
$$\aligned
 0 &=f(w+x) f(w-x) f(y+z) f(y-z)\\
   &+f(w+y) f(w-y) f(z+x) f(z-x)\\
   &+f(w+z) f(w-z) f(x+y) f(x-y).
\endaligned\eqno{(4.3)}$$}

  Thus, $X(p)$ is defined by a collection of quartics which we can write
down explicitly. As much as we may admire this theorem, it is natural to
feel somewhat daunted by it. For even though we know the equations, there
are awful lot of equations and it is not clear that they really do us any
good.

  Hence, Klein did not derive the properties of his loci from a direct
study of the equations defining them. Instead, he used explicit modular
forms to map the moduli spaces into these loci and investigated the
geometry of the images. In particular, for special values of $p$, notably
$p=5$, $7$, $11$, he examined the geometry with great care and in doing
so created a legacy of beautiful geometry, which even now is a source of
considerable inspiration. When $p=7$, the locus $\mathcal{L}_p$ is
isomorphic to the celebrated Klein quartic curve
$$x^3 y+y^3 z+z^3 x=0.\eqno{(4.4)}$$
In the case $p=11$, Klein proved the following remarkable theorem
regarding the modular curve $X(11)$ over the complex numbers.

\textbf{Theorem 4.3.} (see \cite{K3}) {\it The modular curve $X(11)$
is the singular locus of the Hessian of the cubic threefold
$$v^2 w+w^2 x+x^2 y+y^2 z+z^2 v=0,\eqno{(4.5)}$$
which is invariant under the action of a group isomorphic to
$\text{PSL}(2, 11)$. Namely, let $k$ be an algebraically closed field
of characteristic not equal to $11$. The modular curve $X(11)$ over
$k$ is isomorphic to the locus of all points $[v, w, x, y, z]$ in
$\mathbb{P}^4(k)$ such that the matrix
$$\left(\begin{matrix}
  w & v & 0 & 0 & z\\
  v & x & w & 0 & 0\\
  0 & w & y & x & 0\\
  0 & 0 & x & z & y\\
  z & 0 & 0 & y & v
\end{matrix}\right)\eqno{(4.6)}$$
has rank $3$.}

  There are $15$ distinct $4 \times 4$ minors of the matrix $(4.6)$.
They are obtained from the following three by successive applications
of the cyclic permutation $(vwxyz)$:
$$-v w^2 z-v x^3+vxyz+w^2 y^2-x y^3, \eqno{(4.7)}$$
$$v^2 x^2-v^2 yz+v y^3+wxyz,\eqno{(4.8)}$$
$$v^2 wz-vwy^2-x^2 yz,\eqno{(4.9)}$$
these being respectively the first three $4 \times 4$ minors along the
top row. If we apply the cyclic permutation $(zyxwv)$ to (4.7) and
subtract the result from (4.8), we obtain
$$v y^3+w^3 z+w x^3. \eqno{(4.10)}$$
Thus, the span of the $15$ distinct $4 \times 4$ minors of the matrix
(4.6) is the same as that of the polynomials (4.8), (4.9), (4.10)
and their transforms under successive applications of the cyclic
permutation $(vwxyz)$.

  If we substitute $y_1$, $y_9$, $y_4$, $y_3$, $y_5$ respectively for
$v$, $w$, $x$, $y$, $z$ then the polynomials (4.8), (4.9), (4.10)
become respectively the polynomials (2), (3), (4) on page 188 of Klein
(see \cite{K4}). On the other hand, as Klein remarks on pages 195-196
of Klein (see \cite{K4}), there are only $10$ distinct quartics (4.3)
defining the locus $\mathcal{L}_{11}$, namely (4.9) and (4.10) and
their images under successive applications of the cyclic permutation
$(vwxyz)$. Verification of this assertion is facilitated by the observation
that, denoting by $\Phi_{a, b, c, d}$ the locus defined by (4.3) for a
particular choice of $a$, $b$, $c$, $d$ in $\mathbb{Z}/p \mathbb{Z}$, we
have
$$\Phi_{b, a, c, d}=\Phi_{b, c, d, a}=-\Phi_{a, b, c, d}$$
and
$$\Phi_{-a, b, c, d}=\Phi_{a, b, c, d}.$$
Furthermore, if two of $a$, $b$, $c$, $d$ are equal then
$\Phi_{a, b, c, d}=0$. It follows that up to a sign, each of
the quartics appearing in (4.3) is equal to a quartic
$\Phi_{a, b, c, d}$ with
$$0 \leq a<b<c<d \leq \frac{1}{2}(p-1).$$
For $p=11$, this leads to $15$ quartics, but as Klein remarked,
only $10$ of them are actually distinct by the following observation
(see \cite{K4}, p. 196):
$$\aligned
 &v^2 x^2-v^2 yz+v y^3+wxyz\\
=&\frac{1}{w} [-y(v^2 zw-x^2 zy-y^2 vw)-x(y^2 xz-v^2 xw-w^2 yz)].
\endaligned$$

  Thus, the locus $\mathcal{L}_{11}$ contains the rank $3$ locus of
(4.6). To see that the two coincide, we note, as Klein does (see
\cite{K4}, p. 196), that
$$\aligned
 &-v w^2 z-v x^3+vxyz+w^2 y^2-x y^3\\
=&\frac{1}{w} [y(v^2 wz-vwy^2-x^2 yz)+x(y^2 xz-v^2 xw-w^2 yz)].
  \endaligned\eqno{(4.11)}$$
The left-hand side of (4.11) is just (4.7). The numerator of the
right-hand side of (4.11) consist of two parts, the first of which
is $y$ times (4.9) and the second of which is $x$ times a polynomial
obtained from (4.9) by the cube of the cyclic permutation $(vwxyz)$.
The denominator of the right-hand side of (4.11) is just $w$. And we
obtain five similar equations from (4.11) by successive applications
of the cyclic permutation $(vwxyz)$. Therefore, in the affine open set
defined by $vwxyz \neq 0$, the five quartics obtained from (4.7) by
successive application of the cyclic permutation $(vwxyz)$ lie in the
ideal generated by the ten quartics obtained via $(vwxyz)$ from (4.9)
and (4.10). On the other hand, it is easy to see that if $vwxyz=0$ then
the equation (4.11) imply that $[v, w, x, y, z]$ is a point with four of
its entries zero, that is, $[1, 0, 0, 0, 0]$ or one of its cyclic permutations.
These points are easily seen to lie on the rank $3$ locus of (4.6), so
the two loci have the same underlying set and therefore coincide as varieties
(see \cite{Ad1}).

  In other words, we can construct the modular curve $X(11)$ from the
cubic invariant threefold. However, there is nothing in Klein's general
theorem on $X(p)$ about any invariant hyper-surface. It just gives a
bunch of quartic equations that define $X(p)$. On the other hand, Adler
and Ramanan (see \cite{AR}) proved that whenever $p>3$ is a prime
congruent to $3$ modulo $8$, there is a unique cubic invariant for the
representation of $\text{SL}(2, p)$ on $V^{-}$ and we have the modular
curve $X(p)$. More generally, for any $p>3$ there is a unique $3$-tensor
$\Theta$ on $L^2(\mathbb{F}_p)$ invariant under the Weil representation
of $\text{SL}(2, p)$. Thus, with essentially no restriction on $p$,
we can ask: is there a way to construct the modular curve $X(p)$
geometrically from the invariant $3$-tensor $\Theta$? In \cite{Ad2},
Adler proved that the modular curve $X(p)$ can be constructed
geometrically from that $3$-tensor provided $p$ is a prime $\geq 11$
and $\neq 13$. More precisely, if $-1$ is a square modulo $p$, then
one can construct Klein's $\mathbf{A}$-curve of level $p$ from the
restriction of $\Theta$ to the even part of the Weil representation
of $\text{SL}(2, \mathbb{F}_p)$, while if $-1$ is not a square
modulo $p$ then one can construct Klein's $z$-curve from the
restriction of $\Theta$ to the odd part of the Weil representation
(see \cite{Ad2}). In fact, the case $p=13$ requires some additional
concern since the automorphism group of $\Theta$ on $V^{+}$ is the
complex Lie group $G_2$ in that case. A similar phenomenon occurs
in connection with the case $p=7$, the automorphism group of $\Theta$
on $V^{-}$ is $\text{SL}(3, \mathbb{C})$.

\textbf{Theorem 4.4.} (see \cite{Ad2}). {\it Suppose $p \geq 11$. Let
$\varepsilon$ equal the quadratic character of $-1$ in $\mathbb{F}_p$
and let $\eta$ equal the quadratic character of $-2$ in $\mathbb{F}_p$.
Let $n=(p+\varepsilon)/2$. Denote by $\Theta$ the unique $\eta$-symmetric
$3$-tensor on the $\varepsilon$ part of the Weil representation of
$\text{SL}(2, \mathbb{F}_p)$. Then the group of collineations which
preserve the $3$-tensor is isomorphic to $\text{PSL}(2, \mathbb{F}_p)
\cdot \text{Aut}(\mathbb{F}_p)$ unless $p=13$, in which case the group
in question is $G_2(\mathbb{C})$.}

\textbf{Theorem 4.5.} (see \cite{Ad2}). {\it Suppose $p$ is an odd prime which
is $\geq 11$ and $\neq 13$. Then the modular curve $X(p)$ may be constructed
geometrically from the $3$-tensor $\Theta$ in the set theoretic sense. More
precisely, if $\varepsilon=1$ then the $\mathbf{A}$-curve may be constructed
from the restriction $\Theta|V^{+}$ of $\Theta$ to $V^{+}$, while if
$\varepsilon=-1$, the $z$-curve may be constructed from the restriction
$\Theta|V^{-}$ of $\Theta$ to $V^{-}$.}

{\it Definition} 4.6. (see \cite{Ad2}) . Let $H$: $f=0$ be a hyper-surface
in a projective space $\mathbb{P}^d(\mathbb{C})$ of dimension $d$ and let $Z$
be a subvariety of $\mathbb{P}^d(\mathbb{C})$. We say that $Z$ can be constructed
geometrically from $H$ if the ideal defining $Z$ is generated by covariants
of $f$. We say that $Z$ can be constructed geometrically from $H$ in the set
theoretic sense if $Z$ is the set theoretic intersection of covariants of $f$.

  Note that the invariants of $\text{SL}(2, \mathbb{F}_p)$ separate orbits
of $\text{SL}(2, \mathbb{F}_p)$. It follows that the modular curve is the
set theoretic intersection of all of the invariant hyper-surfaces containing
it. This shows that the notion of geometric constructibility which Adler uses
is much too restrictive: if $X$ is a $G$-invariant hyper-surface, it requires
the locus to be an intersection of $G$-invariant hyper-surfaces. While this
may be true set theoretically, contemporary algebraic geometry requires us to
consider the locus from a scheme theoretic point of view and it is certainly
not reasonable to require the ideal defining $Z$ to be generated by invariants.
For example, if we take $Z$ to be the singular locus of $X$, the ideal defining
$Z$ will be generated by the first partial derivatives of the form $f$ defining
$X$. Even if $f$ is an invariant of $G$, the first partials of $f$ in general
will not be. Thus, passing to the singular locus of something geometrically
constructible is not geometrical according to the definition they used.

  Now, we give the geometry of the modular curve via the fundamental relation
(see \cite{AR}). For every element $t$ of $K(\delta)=\mathbb{Z}/p \mathbb{Z}$,
let $E_t$ denote the linear form on the space $V^{-}$ of odd functions
on $K(\delta)$ given by $E_t(h)=h(t)$. If $w$, $x$, $y$, $z$ are
any elements of $K(\delta)$, denote by $\Phi_{w, x, y, z}$ the
quartic form given by
$$\aligned
  \Phi_{w, x, y, z} &=E_{w+x} \cdot E_{w-x} \cdot E_{y+z} \cdot E_{y-z}+\\
                    &+E_{w+y} \cdot E_{w-y} \cdot E_{z+x} \cdot E_{z-x}+\\
                    &+E_{w+z} \cdot E_{w-z} \cdot E_{x+y} \cdot E_{x-y}.
\endaligned\eqno{(4.12)}$$
Then $\mathcal{L}$ is defined by the equations
$$\Phi_{w, x, y, z}=0$$
with $w$, $x$, $y$, $z$ in $\mathbb{Z}/p \mathbb{Z}$. Since $E_{-t}=-E_t$ for
every $t \in K(\delta)$, we have $\Phi_{-w, x, y, z}=\Phi_{w, x, y, z}$.
Furthermore, we have
$$\Phi_{x, w, y, z}=\Phi_{x, y, z, w}=-\Phi_{w, x, y, z}.\eqno{(4.13)}$$
Since the odd permutations
$$\left(\begin{matrix}
  1 & 2 & 3 & 4\\
  2 & 1 & 3 & 4
  \end{matrix}\right) \quad \text{and} \quad
  \left(\begin{matrix}
  1 & 2 & 3 & 4\\
  2 & 3 & 4 & 1
  \end{matrix}\right)$$
generate the group of all permutations on four objects, it follows that
if $(a, b, c, d)=\sigma(w, x, y, z)$ is a permutation of $(w, x, y, z)$
then
$$\Phi_{a, b, c, d}=(-1)^{\sigma} \Phi_{w, x, y, z},\eqno{(4.14)}$$
where $(-1)^{\sigma}$ denotes the sign of the permutation $\sigma$. It
follows from this and from the sentence preceding (4.13) that for all
choices of signs we have
$$\Phi_{\pm w, \pm x, \pm y, \pm z}=\Phi_{w, x, y, z}.\eqno{(4.15)}$$
It follows that each quartic $\Phi_{w, x, y, z}$ is equal, up to a sign,
to a quartic $\Phi_{a, b, c, d}$ with
$$0 \leq a < b < c < d \leq \frac{p-1}{2}.$$
The locus $\mathcal{L}$ is therefore defined by
$\left(\begin{matrix} m\\ 4 \end{matrix}\right)$ quartics, where $2m-1=p$.
In general, these quartics are not distinct.

  For example, when $p=13$, the cardinality is $35$ but there are only
$21$ distinct quartics, namely:

(1) $\Phi_{0123}=E_1^3 E_5-E_2^3 E_4+E_3^3 E_1$.

(2) $\Phi_{0124}=\Phi_{3456}=E_1^2 E_2 E_6-E_2^2 E_3 E_5+E_4^2 E_1 E_3$.

(3) $\Phi_{0125}=\Phi_{1234}=-E_1^2 E_3 E_6-E_2^2 E_4 E_6+E_5^2 E_1 E_3$.

(4) $\Phi_{0126}=\Phi_{2345}=-E_1^2 E_4 E_5+E_2^2 E_5 E_6+E_6^2 E_1 E_3$.

(5) $\Phi_{0134}=-E_1^3 E_6-E_3^3 E_5+E_4^3 E_2$.

(6) $\Phi_{0135}=\Phi_{2356}=-E_1^2 E_2 E_5-E_3^2 E_4 E_6+E_5^2 E_2 E_4$.

(7) $\Phi_{0136}=\Phi_{1245}=-E_1^3 E_3 E_4+E_3^2 E_5 E_6+E_6^2 E_2 E_4$.

(8) $\Phi_{0145}=-E_1^3 E_4-E_4^3 E_6+E_5^3 E_3$.

(9) $\Phi_{0146}=\Phi_{1256}=-E_1^2 E_2 E_3+E_4^2 E_5 E_6+E_6^2 E_3 E_5$.

(10) $\Phi_{0156}=-E_1^3 E_2+E_5^3 E_6+E_6^3 E_4$.

(11) $\Phi_{0234}=\Phi_{2456}=-E_2^2 E_1 E_6-E_3^2 E_2 E_6+E_4^2 E_1 E_5$.

(12) $\Phi_{0235}=-E_2^3 E_5+E_3^3 E_6+E_5^3 E_1$.

(13) $\Phi_{0236}=\Phi_{1346}=-E_2^2 E_3 E_4+E_3^2 E_4 E_5+E_6^2 E_1 E_5$.

(14) $\Phi_{0245}=\Phi_{1356}=-E_2^2 E_1 E_4+E_4^2 E_3 E_6+E_5^2 E_2 E_6$.

(15) $\Phi_{0246}=-E_2^3 E_3+E_4^3 E_5+E_6^3 E_2$.

(16) $\Phi_{0256}=-E_2^3 E_1+E_5^3 E_4-E_6^3 E_3$.

(17) $\Phi_{0345}=\Phi_{1236}=-E_3^2 E_1 E_4+E_4^2 E_2 E_5-E_5^2 E_1 E_6$.

(18) $\Phi_{0346}=-E_3^3 E_2+E_4^3 E_3-E_6^3 E_1$.

(19) $\Phi_{0356}=\Phi_{1246}=-E_3^2 E_1 E_2+E_5^2 E_3 E_4-E_6^2 E_2 E_5$.

(20) $\Phi_{0456}=\Phi_{1345}=-E_4^2 E_1 E_2+E_5^2 E_2 E_3-E_6^2 E_1 E_4$.

(21) $\Phi_{1235}=\Phi_{1456}=\Phi_{2346}=-E_1 E_2 E_3 E_5+E_2 E_3 E_4 E_6+E_1 E_4 E_5 E_6$.

  This proves the following:

\textbf{Theorem 4.7.} {\it The modular curve $X(13)$ is isomorphic to the
above $21$ quartic equations ($\Phi_{abcd}$-terms).}

  Now, we will give a connection between the above quartic $21$ equations
with the invariant theory for the group $G$. This is a crucial step in the
present paper. Note that
$$\mathbf{A}_0^2=(z_1^2 z_4^2+z_2^2 z_5^2+z_3^2 z_6^2)+
                 2 (z_1 z_2 z_4 z_5+z_2 z_3 z_5 z_6+z_3 z_1 z_6 z_4).\eqno{(4.16)}$$
Let $ST^{\nu}$ act on the above two quartic terms, respectively. We obtain
a decomposition of $\mathbf{A}$-terms into the sum of $\mathbf{B}$-terms
and $\mathbf{C}$-terms as follows:
$$\aligned
  &13 [ST^{\nu}(z_1)^2 \cdot ST^{\nu}(z_4)^2+ST^{\nu}(z_2)^2 \cdot ST^{\nu}(z_5)^2
   +ST^{\nu}(z_3)^2 \cdot ST^{\nu}(z_6)^2]\\
 =&\mathbf{C}_0+2 \zeta^{\nu} \mathbf{C}_1+2 \zeta^{3 \nu} \mathbf{C}_3
  +2\zeta^{9 \nu} \mathbf{C}_9+2 \zeta^{12 \nu} \mathbf{C}_{12}+2
   \zeta^{10 \nu} \mathbf{C}_{10}+2 \zeta^{4 \nu} \mathbf{C}_4+\\
  &+2 \zeta^{5 \nu} \mathbf{C}_5+2 \zeta^{2 \nu} \mathbf{C}_2+2 \zeta^{6 \nu}
   \mathbf{C}_6+2 \zeta^{8 \nu} \mathbf{C}_8+2 \zeta^{11 \nu} \mathbf{C}_{11}
   +2 \zeta^{7 \nu} \mathbf{C}_7,
\endaligned\eqno{(4.17)}$$
where
$$\left\{\aligned
  \mathbf{C}_0 &=-(z_1^2 z_4^2+z_2^2 z_5^2+z_3^2 z_6^2),\\
  \mathbf{C}_1 &=z_5^2 z_6^2+2 z_2 z_3 z_4 z_5+2 z_1^2 z_2 z_5-2 z_3 z_1 z_4^2
                -z_3^2 z_6 z_4+\\
               &+z_1 z_2^3+z_3 z_5^3+z_2 z_4 z_6^2+z_1^2 z_3 z_6,\\
  \mathbf{C}_3 &=z_4^2 z_5^2+2 z_1 z_2 z_6 z_4+2 z_3^2 z_1 z_4-2 z_2 z_3 z_6^2
                -z_2^2 z_5 z_6+\\
               &+z_3 z_1^3+z_2 z_4^3+z_1 z_6 z_5^2+z_3^2 z_2 z_5,\\
  \mathbf{C}_9 &=z_6^2 z_4^2+2 z_3 z_1 z_5 z_6+2 z_2^2 z_3 z_6-2 z_1 z_2 z_5^2
                -z_1^2 z_4 z_5+\\
               &+z_2 z_3^3+z_1 z_6^3+z_3 z_5 z_4^2+z_2^2 z_1 z_4,\\
  \mathbf{C}_{12} &=z_2^2 z_3^2+2 z_1 z_2 z_5 z_6-2 z_4^2 z_2 z_5-2 z_1^2 z_6 z_4
                   -z_6^2 z_3 z_1+\\
                  &-z_2^3 z_6-z_3^2 z_1 z_5+z_4 z_5^3-z_4^2 z_3 z_6,\\
  \mathbf{C}_{10} &=z_1^2 z_2^2+2 z_3 z_1 z_4 z_5-2 z_6^2 z_1 z_4-2 z_3^2 z_5 z_6
                   -z_5^2 z_2 z_3+\\
                  &-z_1^3 z_5-z_2^2 z_3 z_4+z_6 z_4^3-z_6^2 z_2 z_5,\\
  \mathbf{C}_4 &=z_3^2 z_1^2+2 z_2 z_3 z_6 z_4-2 z_5^2 z_3 z_6-2 z_2^2 z_4 z_5
                -z_4^2 z_1 z_2+\\
               &-z_3^3 z_4-z_1^2 z_2 z_6+z_5 z_6^3-z_5^2 z_1 z_4,\\
  \mathbf{C}_5 &=\frac{1}{2} z_2^4-2 z_4 z_5 z_6^2+2 z_3 z_4 z_5^2+2 z_2 z_3 z_4^2
                 -2 z_1^2 z_2 z_4+z_1^2 z_5^2,\\
  \mathbf{C}_2 &=\frac{1}{2} z_1^4-2 z_6 z_4 z_5^2+2 z_2 z_6 z_4^2+2 z_1 z_2 z_6^2
                 -2 z_3^2 z_1 z_6+z_3^2 z_4^2,\\
  \mathbf{C}_6 &=\frac{1}{2} z_3^4-2 z_5 z_6 z_4^2+2 z_1 z_5 z_6^2+2 z_3 z_1 z_5^2
                 -2 z_2^2 z_3 z_5+z_2^2 z_6^2,\\
  \mathbf{C}_8 &=\frac{1}{2} z_5^4-2 z_1 z_2 z_3^2+2 z_1^2 z_5 z_6+2 z_1 z_5 z_4^2
                 -2 z_2^2 z_1 z_6+z_2^2 z_4^2,\\
  \mathbf{C}_{11} &=\frac{1}{2} z_4^4-2 z_3 z_1 z_2^2+2 z_3^2 z_4 z_5+2 z_3 z_4 z_6^2
                    -2 z_1^2 z_3 z_5+z_1^2 z_6^2,\\
  \mathbf{C}_7 &=\frac{1}{2} z_6^4-2 z_2 z_3 z_1^2+2 z_2^2 z_6 z_4+2 z_2 z_6 z_5^2
                 -2 z_3^2 z_2 z_4+z_3^2 z_5^2.
\endaligned\right.\eqno{(4.18)}$$
On the other hand,
$$\aligned
  &13 [ST^{\nu}(z_1) \cdot ST^{\nu}(z_4) \cdot ST^{\nu}(z_2) \cdot ST^{\nu}(z_5)+\\
  &+ST^{\nu}(z_2) \cdot ST^{\nu}(z_5) \cdot ST^{\nu}(z_3) \cdot ST^{\nu}(z_6)+\\
  &+ST^{\nu}(z_3) \cdot ST^{\nu}(z_6) \cdot ST^{\nu}(z_1) \cdot ST^{\nu}(z_4)]\\
 =&\mathbf{B}_0+\zeta^{\nu} \mathbf{B}_1+\zeta^{3 \nu} \mathbf{B}_3
  +\zeta^{9 \nu} \mathbf{B}_9+\zeta^{12 \nu} \mathbf{B}_{12}+\zeta^{10 \nu}
   \mathbf{B}_{10}+\zeta^{4 \nu} \mathbf{B}_4+\\
  &+2 \zeta^{5 \nu} \mathbf{B}_5+2 \zeta^{2 \nu} \mathbf{B}_2+2 \zeta^{6 \nu}
   \mathbf{B}_6+2 \zeta^{8 \nu} \mathbf{B}_8+2 \zeta^{11 \nu} \mathbf{B}_{11}
   +2 \zeta^{7 \nu} \mathbf{B}_7,
\endaligned\eqno{(4.19)}$$
where
$$\left\{\aligned
  \mathbf{B}_0 &=5(z_1 z_2 z_4 z_5+z_2 z_3 z_5 z_6+z_3 z_1 z_6 z_4)
                +2(z_1 z_5^3+z_2 z_6^3+z_3 z_4^3)+\\
               &\quad -2(z_1^3 z_6+z_2^3 z_4+z_3^3 z_5),\\
  \mathbf{B}_1 &=z_1^3 z_4-z_1 z_2^3+z_3 z_5^3-z_3^2 z_6 z_4
                +z_2 z_4 z_6^2-z_1^2 z_2 z_5,\\
  \mathbf{B}_3 &=z_3^3 z_6-z_3 z_1^3+z_2 z_4^3-z_2^2 z_5 z_6
                +z_1 z_6 z_5^2-z_3^2 z_1 z_4,\\
  \mathbf{B}_9 &=z_2^3 z_5-z_2 z_3^3+z_1 z_6^3-z_1^2 z_4 z_5
                +z_3 z_5 z_4^2-z_2^2 z_3 z_6,\\
  \mathbf{B}_{12} &=-z_1 z_4^3-z_4 z_5^3-z_2^3 z_6-z_3^2 z_1 z_5
                   -z_3 z_1 z_6^2+z_4^2 z_2 z_5,\\
  \mathbf{B}_{10} &=-z_3 z_6^3-z_6 z_4^3-z_1^3 z_5-z_2^2 z_3 z_4
                   -z_2 z_3 z_5^2+z_6^2 z_1 z_4,\\
  \mathbf{B}_4 &=-z_2 z_5^3-z_5 z_6^3-z_3^3 z_4-z_1^2 z_2 z_6
                -z_1 z_2 z_4^2+z_5^2 z_3 z_6,\\
  \mathbf{B}_5 &=-z_2^2 z_1 z_5+z_4 z_5 z_6^2+z_2 z_3 z_4^2,\\
  \mathbf{B}_2 &=-z_1^2 z_3 z_4+z_6 z_4 z_5^2+z_1 z_2 z_6^2,\\
  \mathbf{B}_6 &=-z_3^2 z_2 z_6+z_5 z_6 z_4^2+z_3 z_1 z_5^2,\\
  \mathbf{B}_8 &=z_2 z_4 z_5^2+z_1 z_2 z_3^2+z_1^2 z_5 z_6,\\
  \mathbf{B}_{11} &=z_1 z_6 z_4^2+z_3 z_1 z_2^2+z_3^2 z_4 z_5,\\
  \mathbf{B}_7 &=z_3 z_5 z_6^2+z_2 z_3 z_1^2+z_2^2 z_6 z_4.
\endaligned\right.\eqno{(4.20)}$$
This leads us to define the following decompositions:
$$\mathbf{B}_0=5 \mathbf{B}_0^{(0)}+2 \mathbf{B}_0^{(1)}-2 \mathbf{B}_0^{(2)},
  \eqno{(4.21)}$$
where
$$\left\{\aligned
  \mathbf{B}_0^{(0)} &=z_1 z_2 z_4 z_5+z_2 z_3 z_5 z_6+z_3 z_1 z_6 z_4,\\
  \mathbf{B}_0^{(1)} &=z_1 z_5^3+z_2 z_6^3+z_3 z_4^3,\\
  \mathbf{B}_0^{(2)} &=z_1^3 z_6+z_2^3 z_4+z_3^3 z_5.
\endaligned\right.\eqno{(4.22)}$$
$$\mathbf{B}_1=\mathbf{B}_1^{(1)}+\mathbf{B}_1^{(2)},\eqno{(4.23)}$$
where
$$\left\{\aligned
  \mathbf{B}_1^{(1)} &=z_3 z_5^3+z_1^3 z_4-z_1 z_2^3,\\
  \mathbf{B}_1^{(2)} &=z_2 z_4 z_6^2-z_3^2 z_6 z_4-z_1^2 z_2 z_5.
\endaligned\right.\eqno{(4.24)}$$
$$\mathbf{B}_3=\mathbf{B}_3^{(1)}+\mathbf{B}_3^{(2)},\eqno{(4.25)}$$
where
$$\left\{\aligned
  \mathbf{B}_3^{(1)} &=z_2 z_4^3+z_3^3 z_6-z_3 z_1^3,\\
  \mathbf{B}_3^{(2)} &=z_1 z_6 z_5^2-z_2^2 z_5 z_6-z_3^2 z_1 z_4.
\endaligned\right.\eqno{(4.26)}$$
$$\mathbf{B}_9=\mathbf{B}_9^{(1)}+\mathbf{B}_9^{(2)},\eqno{(4.27)}$$
where
$$\left\{\aligned
  \mathbf{B}_9^{(1)} &=z_1 z_6^3+z_2^3 z_5-z_2 z_3^3,\\
  \mathbf{B}_9^{(2)} &=z_3 z_5 z_4^2-z_1^2 z_4 z_5-z_2^2 z_3 z_6.
\endaligned\right.\eqno{(4.28)}$$
$$\mathbf{B}_{12}=-\mathbf{B}_{12}^{(1)}+\mathbf{B}_{12}^{(2)},\eqno{(4.29)}$$
where
$$\left\{\aligned
  \mathbf{B}_{12}^{(1)} &=z_1 z_4^3+z_2^3 z_6+z_4 z_5^3,\\
  \mathbf{B}_{12}^{(2)} &=z_2 z_5 z_4^2-z_3^2 z_1 z_5-z_6^2 z_3 z_1.
\endaligned\right.\eqno{(4.30)}$$
$$\mathbf{B}_{10}=-\mathbf{B}_{10}^{(1)}+\mathbf{B}_{10}^{(2)},\eqno{(4.31)}$$
where
$$\left\{\aligned
  \mathbf{B}_{10}^{(1)} &=z_3 z_6^3+z_1^3 z_5+z_6 z_4^3,\\
  \mathbf{B}_{10}^{(2)} &=z_1 z_4 z_6^2-z_2^2 z_3 z_4-z_5^2 z_2 z_3.
\endaligned\right.\eqno{(4.32)}$$
$$\mathbf{B}_4=-\mathbf{B}_4^{(1)}+\mathbf{B}_4^{(2)},\eqno{(4.33)}$$
where
$$\left\{\aligned
  \mathbf{B}_4^{(1)} &=z_2 z_5^3+z_3^3 z_4+z_5 z_6^3,\\
  \mathbf{B}_4^{(2)} &=z_3 z_6 z_5^2-z_1^2 z_2 z_6-z_4^2 z_1 z_2.
\endaligned\right.\eqno{(4.34)}$$

  The significance of these quartic polynomials come from the following:

\textbf{Theorem 4.8.} {\it There is a one-to-one correspondence
between the above quartic equations $($$\Phi_{abcd}$-terms$)$
and the quartic polynomials $($$\mathbf{B}$-terms$)$.}

{\it Proof}. Under the following correspondence
$$z_1 \longleftrightarrow -e^{\frac{\pi i}{4}} E_1, \quad
  z_2 \longleftrightarrow -e^{\frac{\pi i}{4}} E_3, \quad
  z_3 \longleftrightarrow e^{\frac{\pi i}{4}} E_4,\eqno{(4.35)}$$
and
$$z_4 \longleftrightarrow e^{\frac{\pi i}{4}} E_5, \quad
  z_5 \longleftrightarrow e^{\frac{\pi i}{4}} E_2, \quad
  z_6 \longleftrightarrow e^{\frac{\pi i}{4}} E_6,\eqno{(4.36)}$$
we have

(1) $\mathbf{B}_0^{(0)} \longleftrightarrow \Phi_{1235}=\Phi_{1456}=\Phi_{2346}$.

(2) $\mathbf{B}_0^{(1)} \longleftrightarrow -\Phi_{0256}$.

(3) $\mathbf{B}_0^{(2)} \longleftrightarrow -\Phi_{0134}$.

(4) $\mathbf{B}_1^{(1)} \longleftrightarrow \Phi_{0123}$.

(5) $\mathbf{B}_1^{(2)} \longleftrightarrow \Phi_{0146}=\Phi_{1256}$.

(6) $\mathbf{B}_3^{(1)} \longleftrightarrow \Phi_{0145}$.

(7) $\mathbf{B}_3^{(2)} \longleftrightarrow -\Phi_{0234}=-\Phi_{2456}$.

(8) $\mathbf{B}_9^{(1)} \longleftrightarrow -\Phi_{0346}$.

(9) $\mathbf{B}_9^{(2)} \longleftrightarrow -\Phi_{0135}=-\Phi_{2356}$.

(10) $\mathbf{B}_{12}^{(1)} \longleftrightarrow \Phi_{0235}$.

(11) $\mathbf{B}_{12}^{(2)} \longleftrightarrow \Phi_{0456}=\Phi_{1345}$.

(12) $\mathbf{B}_{10}^{(1)} \longleftrightarrow -\Phi_{0156}$.

(13) $\mathbf{B}_{10}^{(2)} \longleftrightarrow \Phi_{0236}=\Phi_{1346}$.

(14) $\mathbf{B}_4^{(1)} \longleftrightarrow -\Phi_{0246}$.

(15) $\mathbf{B}_4^{(2)} \longleftrightarrow \Phi_{0125}=\Phi_{1234}$.

(16) $\mathbf{B}_5 \longleftrightarrow \Phi_{0356}=\Phi_{1246}$.

(17) $\mathbf{B}_2 \longleftrightarrow -\Phi_{0126}=-\Phi_{2345}$.

(18) $\mathbf{B}_6 \longleftrightarrow -\Phi_{0245}=-\Phi_{1356}$.

(19) $\mathbf{B}_8 \longleftrightarrow -\Phi_{0124}=-\Phi_{3456}$.

(20) $\mathbf{B}_{11} \longleftrightarrow -\Phi_{0345}=-\Phi_{1236}$.

(21) $\mathbf{B}_7 \longleftrightarrow -\Phi_{0136}=-\Phi_{1245}$.

This completes the proof of Theorem 4.9.

\noindent
$\qquad \qquad \qquad \qquad \qquad \qquad \qquad \qquad \qquad
 \qquad \qquad \qquad \qquad \qquad \qquad \qquad \boxed{}$

  Theorem 4.8 gives a connection between Klein's quartic system,
i.e., the locus for the modular curve $X(13)$ (algebraic geometry)
and our quartic invariants associated with $\text{SL}(2, 13)$
(invariant theory). This leads us to study Klein's quartic system
from the viewpoint of representation theory.

\begin{center}
{\large\bf 5. An invariant ideal defining the modular curve $X(13)$}
\end{center}

  Let $I=I(Y)$ be an ideal generated by the above twenty-one quartic
polynomials ($\mathbf{B}$-terms):
$$I=\langle \mathbf{B}_0^{(i)}, \mathbf{B}_1^{(j)}, \mathbf{B}_3^{(j)},
    \mathbf{B}_9^{(j)}, \mathbf{B}_{12}^{(j)}, \mathbf{B}_{10}^{(j)},
    \mathbf{B}_4^{(j)}, \mathbf{B}_5, \mathbf{B}_2, \mathbf{B}_6,
    \mathbf{B}_8, \mathbf{B}_{11}, \mathbf{B}_7 \rangle\eqno{(5.1)}$$
with $i=0. 1, 2$ and $j=1, 2$. The corresponding curve associated to
the ideal $I(Y)$ is denoted by $Y$.

\textbf{Theorem 5.1.} {\it The modular curve $X=X(13)$ is isomorphic
to the curve $Y$ in $\mathbb{CP}^5$.}

{\it Proof}. It is the consequence of Theorem 4.7 and Theorem 4.9.

\noindent
$\qquad \qquad \qquad \qquad \qquad \qquad \qquad \qquad \qquad
 \qquad \qquad \qquad \qquad \qquad \qquad \qquad \boxed{}$

\textbf{Theorem 5.2.} {\it The curve $Y$ can be constructed from the
invariant quartic Fano four-fold $\Phi_4=0$ in $\mathbb{CP}^5$, i.e.,
$$(z_3 z_4^3+z_1 z_5^3+z_2 z_6^3)-(z_6 z_1^3+z_4 z_2^3+z_5 z_3^3)
  +3(z_1 z_2 z_4 z_5+z_2 z_3 z_5 z_6+z_3 z_1 z_6 z_4)=0.\eqno{(5.2)}$$}

{\it Proof}. Note that the invariant quartic polynomial $\Phi_4$ has the
following decomposition:
$$\Phi_4=3 \mathbf{B}_0^{(0)}+\mathbf{B}_0^{(1)}-\mathbf{B}_0^{(2)}.\eqno{(5.3)}$$
Let $S$ act on the three terms $\mathbf{B}_0^{(0)}$, $\mathbf{B}_0^{(1)}$
and $\mathbf{B}_0^{(2)}$, respectively. This gives that
$$\aligned
   13 S(\mathbf{B}_0^{(0)})
 =&5 \mathbf{B}_0^{(0)}+2 \mathbf{B}_0^{(1)}-2 \mathbf{B}_0^{(2)}
   +2(\mathbf{B}_5+\mathbf{B}_2+\mathbf{B}_6+\mathbf{B}_8+\mathbf{B}_{11}+\mathbf{B}_7)+\\
  &+(\mathbf{B}_1^{(2)}+\mathbf{B}_3^{(2)}+\mathbf{B}_9^{(2)})
   +(\mathbf{B}_{12}^{(2)}+\mathbf{B}_{10}^{(2)}+\mathbf{B}_4^{(2)})+\\
  &+(\mathbf{B}_1^{(1)}+\mathbf{B}_3^{(1)}+\mathbf{B}_9^{(1)})
   -(\mathbf{B}_{12}^{(1)}+\mathbf{B}_{10}^{(1)}+\mathbf{B}_4^{(1)}).
\endaligned\eqno{(5.4)}$$
On the other hand,
$$\aligned
   13 S(\mathbf{B}_0^{(1)})
 =&12 \mathbf{B}_0^{(0)}+\frac{7+\sqrt{13}}{2} \mathbf{B}_0^{(1)}+
   \frac{-7+\sqrt{13}}{2} \mathbf{B}_0^{(2)}+\\
  &-3 (\mathbf{B}_5+\mathbf{B}_2+\mathbf{B}_6+\mathbf{B}_8+\mathbf{B}_{11}
   +\mathbf{B}_7)+\\
  &+\frac{-3-3 \sqrt{13}}{2} (\mathbf{B}_1^{(2)}+\mathbf{B}_3^{(2)}
   +\mathbf{B}_9^{(2)})+\\
  &+\frac{-3+3 \sqrt{13}}{2} (\mathbf{B}_{12}^{(2)}+\mathbf{B}_{10}^{(2)}
   +\mathbf{B}_4^{(2)})+\\
  &+\frac{-3+\sqrt{13}}{2} (\mathbf{B}_1^{(1)}+\mathbf{B}_3^{(1)}+\mathbf{B}_9^{(1)})+\\
  &+\frac{3+\sqrt{13}}{2} (\mathbf{B}_{12}^{(1)}+\mathbf{B}_{10}^{(1)}+\mathbf{B}_4^{(1)}).
\endaligned\eqno{(5.5)}$$
$$\aligned
   13 S(\mathbf{B}_0^{(2)})
 =&-12 \mathbf{B}_0^{(0)}+\frac{-7+\sqrt{13}}{2} \mathbf{B}_0^{(1)}+
   \frac{7+\sqrt{13}}{2} \mathbf{B}_0^{(2)}+\\
  &+3 (\mathbf{B}_5+\mathbf{B}_2+\mathbf{B}_6+\mathbf{B}_8+\mathbf{B}_{11}
   +\mathbf{B}_7)+\\
  &+\frac{3-3 \sqrt{13}}{2} (\mathbf{B}_1^{(2)}+\mathbf{B}_3^{(2)}
   +\mathbf{B}_9^{(2)})+\\
  &+\frac{3+3 \sqrt{13}}{2} (\mathbf{B}_{12}^{(2)}+\mathbf{B}_{10}^{(2)}
   +\mathbf{B}_4^{(2)})+\\
  &+\frac{3+\sqrt{13}}{2} (\mathbf{B}_1^{(1)}+\mathbf{B}_3^{(1)}+\mathbf{B}_9^{(1)})+\\
  &+\frac{-3+\sqrt{13}}{2} (\mathbf{B}_{12}^{(1)}+\mathbf{B}_{10}^{(1)}+\mathbf{B}_4^{(1)}).
\endaligned\eqno{(5.6)}$$
Combining with (4.19), this leads to all of the $\mathbf{B}$-terms (4.20), (4.22),
(4.24), (4.26), (4.28), (4.30), (4.32) and (4.34) which we have defined above.
Hence, the curve $Y$ can be constructed from the invariant quartic Fano four-fold
$\Phi_4=0$ by the action of the transformation $S$. Moreover,
$$\aligned
  &13 [S(\mathbf{B}_0^{(1)})-S(\mathbf{B}_0^{(2)})]\\
 =&24 \mathbf{B}_0^{(0)}+7 \mathbf{B}_0^{(1)}-7 \mathbf{B}_0^{(2)}
  -6(\mathbf{B}_5+\mathbf{B}_2+\mathbf{B}_6+\mathbf{B}_8+\mathbf{B}_{11}+\mathbf{B}_7)+\\
  &-3(\mathbf{B}_1^{(2)}+\mathbf{B}_3^{(2)}+\mathbf{B}_9^{(2)})
   -3(\mathbf{B}_{12}^{(2)}+\mathbf{B}_{10}^{(2)}+\mathbf{B}_4^{(2)})+\\
  &-3(\mathbf{B}_1^{(1)}+\mathbf{B}_3^{(1)}+\mathbf{B}_9^{(1)})
   +3(\mathbf{B}_{12}^{(1)}+\mathbf{B}_{10}^{(1)}+\mathbf{B}_4^{(1)}),
\endaligned$$
and
$$\aligned
  &13 S(3 \mathbf{B}_0^{(0)}+\mathbf{B}_0^{(1)}-\mathbf{B}_0^{(2)})\\
 =&15 \mathbf{B}_0^{(0)}+6 \mathbf{B}_0^{(1)}-6 \mathbf{B}_0^{(2)}
  +6(\mathbf{B}_5+\mathbf{B}_2+\mathbf{B}_6+\mathbf{B}_8+\mathbf{B}_{11}+\mathbf{B}_7)+\\
  &+24 \mathbf{B}_0^{(0)}+7 \mathbf{B}_0^{(1)}-7 \mathbf{B}_0^{(2)}
  -6(\mathbf{B}_5+\mathbf{B}_2+\mathbf{B}_6+\mathbf{B}_8+\mathbf{B}_{11}+\mathbf{B}_7)\\
 =&13 (3 \mathbf{B}_0^{(0)}+\mathbf{B}_0^{(1)}-\mathbf{B}_0^{(2)}),
\endaligned$$
 i.e.,
 $$S(3 \mathbf{B}_0^{(0)}+\mathbf{B}_0^{(1)}-\mathbf{B}_0^{(2)})
  =3 \mathbf{B}_0^{(0)}+\mathbf{B}_0^{(1)}-\mathbf{B}_0^{(2)}.$$
This proves again that $\Phi_4$ is invariant under the action of $S$.

\noindent
$\qquad \qquad \qquad \qquad \qquad \qquad \qquad \qquad \qquad
 \qquad \qquad \qquad \qquad \qquad \qquad \qquad \boxed{}$

\textbf{Theorem 5.3.} {\it The ideal $I(Y)$ is invariant under the action
of $G$, which gives a twenty-one dimensional representation of $G$.}

{\it Proof}. By the proof of Theorem 5.2, we have obtained the expression
of $S(\mathbf{B}_0^{(0)})$, $S(\mathbf{B}_0^{(1)})$ and $S(\mathbf{B}_0^{(2)})$.
The other $18$ terms can be divided into six triples:
$$(S(\mathbb{B}_1^{(1)}), S(\mathbb{B}_3^{(1)}), S(\mathbb{B}_9^{(1)})),$$
$$(S(\mathbb{B}_{12}^{(1)}), S(\mathbb{B}_{10}^{(1)}), S(\mathbb{B}_4^{(1)})),$$
$$(S(\mathbb{B}_1^{(2)}), S(\mathbb{B}_3^{(2)}), S(\mathbb{B}_9^{(1)})),$$
$$(S(\mathbb{B}_{12}^{(2)}), S(\mathbb{B}_{10}^{(2)}), S(\mathbb{B}_4^{(2)})),$$
$$(S(\mathbb{B}_5), S(\mathbb{B}_2), S(\mathbb{B}_6)),$$
$$(S(\mathbb{B}_8), S(\mathbb{B}_{11}), S(\mathbb{B}_7)).$$
Without loss of generality, we begin with the computation
of $S(\mathbf{B}_1^{(1)})$, $S(\mathbf{B}_{12}^{(1)})$, $S(\mathbf{B}_1^{(2)})$,
$S(\mathbf{B}_{12}^{(2)})$, $S(\mathbf{B}_5)$, and $S(\mathbf{B}_8)$. We have
$$\aligned
  &13 S(\mathbf{B}_1^{(1)})\\
 =&6 \mathbf{B}_0^{(0)}+\frac{-3+\sqrt{13}}{2} \mathbf{B}_0^{(1)}+
   \frac{3+\sqrt{13}}{2} \mathbf{B}_0^{(2)}+\\
  &+3 [(\zeta^9+\zeta^4+\zeta^6+\zeta^7+\zeta^2+\zeta^{11}) \mathbf{B}_5+\\
  &+(\zeta^3+\zeta^{10}+\zeta^2+\zeta^{11}+\zeta^5+\zeta^8) \mathbf{B}_2+\\
  &+(\zeta+\zeta^{12}+\zeta^5+\zeta^8+\zeta^6+\zeta^7) \mathbf{B}_6+\\
  &+(\zeta^9+\zeta^4+\zeta^6+\zeta^7+\zeta^3+\zeta^{10}) \mathbf{B}_8+\\
  &+(\zeta^3+\zeta^{10}+\zeta^2+\zeta^{11}+\zeta+\zeta^{12}) \mathbf{B}_{11}+\\
  &+(\zeta+\zeta^{12}+\zeta^5+\zeta^8+\zeta^9+\zeta^4) \mathbf{B}_7]+\\
  &+3 [(1+\zeta^9+\zeta^4) \mathbf{B}_1^{(2)}
   +(1+\zeta^3+\zeta^{10}) \mathbf{B}_3^{(2)}
   +(1+\zeta+\zeta^{12}) \mathbf{B}_9^{(2)}+\\
  &+(1+\zeta^6+\zeta^7) \mathbf{B}_{12}^{(2)}
   +(1+\zeta^2+\zeta^{11}) \mathbf{B}_{10}^{(2)}
   +(1+\zeta^5+\zeta^8) \mathbf{B}_4^{(2)}]+\\
  &+[1+(\zeta^9+\zeta^4)+2 (\zeta+\zeta^{12})+(\zeta^2+\zeta^{11})] \mathbf{B}_1^{(1)}+\\
  &+[1+(\zeta^3+\zeta^{10})+2 (\zeta^9+\zeta^4)+(\zeta^5+\zeta^8)] \mathbf{B}_3^{(1)}+\\
  &+[1+(\zeta+\zeta^{12})+2 (\zeta^3+\zeta^{10})+(\zeta^6+\zeta^7)] \mathbf{B}_9^{(1)}+\\
  &-[1+(\zeta^6+\zeta^7)+2 (\zeta^5+\zeta^8)+(\zeta^3+\zeta^{10})] \mathbf{B}_{12}^{(1)}+\\
  &-[1+(\zeta^2+\zeta^{11})+2 (\zeta^6+\zeta^7)+(\zeta+\zeta^{12})] \mathbf{B}_{10}^{(1)}+\\
  &-[1+(\zeta^5+\zeta^8)+2 (\zeta^2+\zeta^{11})+(\zeta^9+\zeta^4)] \mathbf{B}_4^{(1)}.
\endaligned\eqno{(5.7)}$$
$$\aligned
  &13 S(\mathbf{B}_{12}^{(1)})\\
 =&-6 \mathbf{B}_0^{(0)}+\frac{3+\sqrt{13}}{2} \mathbf{B}_0^{(1)}+
   \frac{-3+\sqrt{13}}{2} \mathbf{B}_0^{(2)}+\\
  &-3 [(\zeta^9+\zeta^4+\zeta^6+\zeta^7+\zeta^3+\zeta^{10}) \mathbf{B}_5+\\
  &+(\zeta^3+\zeta^{10}+\zeta^2+\zeta^{11}+\zeta+\zeta^{12}) \mathbf{B}_2+\\
  &+(\zeta+\zeta^{12}+\zeta^5+\zeta^8+\zeta^9+\zeta^4) \mathbf{B}_6+\\
  &+(\zeta^9+\zeta^4+\zeta^6+\zeta^7+\zeta^2+\zeta^{11}) \mathbf{B}_8+\\
  &+(\zeta^3+\zeta^{10}+\zeta^2+\zeta^{11}+\zeta^5+\zeta^8) \mathbf{B}_{11}+\\
  &+(\zeta+\zeta^{12}+\zeta^5+\zeta^8+\zeta^6+\zeta^7) \mathbf{B}_7]+\\
  &-3 [(1+\zeta^6+\zeta^7) \mathbf{B}_1^{(2)}
   +(1+\zeta^2+\zeta^{11}) \mathbf{B}_3^{(2)}
   +(1+\zeta^5+\zeta^8) \mathbf{B}_9^{(2)}+\\
  &+(1+\zeta^9+\zeta^4) \mathbf{B}_{12}^{(2)}
   +(1+\zeta^3+\zeta^{10}) \mathbf{B}_{10}^{(2)}
   +(1+\zeta+\zeta^{12}) \mathbf{B}_4^{(2)}]+\\
  &-[1+(\zeta^6+\zeta^7)+2 (\zeta^5+\zeta^8)+(\zeta^3+\zeta^{10})] \mathbf{B}_1^{(1)}+\\
  &-[1+(\zeta^2+\zeta^{11})+2 (\zeta^6+\zeta^7)+(\zeta+\zeta^{12})] \mathbf{B}_3^{(1)}+\\
  &-[1+(\zeta^5+\zeta^8)+2 (\zeta^2+\zeta^{11})+(\zeta^9+\zeta^4)] \mathbf{B}_9^{(1)}+\\
  &+[1+(\zeta^9+\zeta^4)+2 (\zeta+\zeta^{12})+(\zeta^2+\zeta^{11})] \mathbf{B}_{12}^{(1)}+\\
  &+[1+(\zeta^3+\zeta^{10})+2 (\zeta^9+\zeta^4)+(\zeta^5+\zeta^8)] \mathbf{B}_{10}^{(1)}+\\
  &+[1+(\zeta+\zeta^{12})+2 (\zeta^3+\zeta^{10})+(\zeta^6+\zeta^7)] \mathbf{B}_4^{(1)}.
\endaligned\eqno{(5.8)}$$
$$\aligned
  &13 S(\mathbf{B}_1^{(2)})\\
 =&2 \mathbf{B}_0^{(0)}+\frac{-1-\sqrt{13}}{2} \mathbf{B}_0^{(1)}+
   \frac{1-\sqrt{13}}{2} \mathbf{B}_0^{(2)}+\\
  &+(\zeta^9+\zeta^4+\zeta^6+\zeta^7+\zeta^2+\zeta^{11}) \mathbf{B}_5+\\
  &+(\zeta^3+\zeta^{10}+\zeta^2+\zeta^{11}+\zeta^5+\zeta^8) \mathbf{B}_2+\\
  &+(\zeta+\zeta^{12}+\zeta^5+\zeta^8+\zeta^6+\zeta^7) \mathbf{B}_6+\\
  &+(\zeta^9+\zeta^4+\zeta^6+\zeta^7+\zeta^3+\zeta^{10}) \mathbf{B}_8+\\
  &+(\zeta^3+\zeta^{10}+\zeta^2+\zeta^{11}+\zeta+\zeta^{12}) \mathbf{B}_{11}+\\
  &+(\zeta+\zeta^{12}+\zeta^5+\zeta^8+\zeta^9+\zeta^4) \mathbf{B}_7+\\
  &+[-1-(\zeta^9+\zeta^4)+2 (\zeta+\zeta^{12})+(\zeta^2+\zeta^{11})] \mathbf{B}_1^{(2)}+\\
  &+[-1-(\zeta^3+\zeta^{10})+2 (\zeta^9+\zeta^4)+(\zeta^5+\zeta^8)] \mathbf{B}_3^{(2)}+\\
  &+[-1-(\zeta+\zeta^{12})+2 (\zeta^3+\zeta^{10})+(\zeta^6+\zeta^7)] \mathbf{B}_9^{(2)}+\\
  &+[-1-(\zeta^6+\zeta^7)+2 (\zeta^5+\zeta^8)+(\zeta^3+\zeta^{10})] \mathbf{B}_{12}^{(2)}+\\
  &+[-1-(\zeta^2+\zeta^{11})+2 (\zeta^6+\zeta^7)+(\zeta+\zeta^{12})] \mathbf{B}_{10}^{(2)}+\\
  &+[-1-(\zeta^5+\zeta^8)+2 (\zeta^2+\zeta^{11})+(\zeta^9+\zeta^4)] \mathbf{B}_4^{(2)}+\\
  &+(1+\zeta^9+\zeta^4) \mathbf{B}_1^{(1)}+(1+\zeta^3+\zeta^{10}) \mathbf{B}_3^{(1)}
   +(1+\zeta+\zeta^{12}) \mathbf{B}_9^{(1)}+\\
  &-(1+\zeta^6+\zeta^7) \mathbf{B}_{12}^{(1)}-(1+\zeta^2+\zeta^{11}) \mathbf{B}_{10}^{(1)}
   -(1+\zeta^5+\zeta^8) \mathbf{B}_4^{(1)}.
\endaligned\eqno{(5.9)}$$
$$\aligned
  &13 S(\mathbf{B}_{12}^{(2)})\\
 =&2 \mathbf{B}_0^{(0)}+\frac{-1+\sqrt{13}}{2} \mathbf{B}_0^{(1)}+
   \frac{1+\sqrt{13}}{2} \mathbf{B}_0^{(2)}+\\
  &+(\zeta^9+\zeta^4+\zeta^6+\zeta^7+\zeta^3+\zeta^{10}) \mathbf{B}_5+\\
  &+(\zeta^3+\zeta^{10}+\zeta^2+\zeta^{11}+\zeta+\zeta^{12}) \mathbf{B}_2+\\
  &+(\zeta+\zeta^{12}+\zeta^5+\zeta^8+\zeta^9+\zeta^4) \mathbf{B}_6+\\
  &+(\zeta^9+\zeta^4+\zeta^6+\zeta^7+\zeta^2+\zeta^{11}) \mathbf{B}_8+\\
  &+(\zeta^3+\zeta^{10}+\zeta^2+\zeta^{11}+\zeta^5+\zeta^8) \mathbf{B}_{11}+\\
  &+(\zeta+\zeta^{12}+\zeta^5+\zeta^8+\zeta^6+\zeta^7) \mathbf{B}_7+\\
  &+[-1-(\zeta^6+\zeta^7)+2 (\zeta^5+\zeta^8)+(\zeta^3+\zeta^{10})] \mathbf{B}_1^{(2)}+\\
  &+[-1-(\zeta^2+\zeta^{11})+2 (\zeta^6+\zeta^7)+(\zeta+\zeta^{12})] \mathbf{B}_3^{(2)}+\\
  &+[-1-(\zeta^5+\zeta^8)+2 (\zeta^2+\zeta^{11})+(\zeta^9+\zeta^4)] \mathbf{B}_9^{(2)}+\\
  &+[-1-(\zeta^9+\zeta^4)+2 (\zeta+\zeta^{12})+(\zeta^2+\zeta^{11})] \mathbf{B}_{12}^{(2)}+\\
  &+[-1-(\zeta^3+\zeta^{10})+2 (\zeta^9+\zeta^4)+(\zeta^5+\zeta^8)] \mathbf{B}_{10}^{(2)}+\\
  &+[-1-(\zeta+\zeta^{12})+2 (\zeta^3+\zeta^{10})+(\zeta^6+\zeta^7)] \mathbf{B}_4^{(2)}+\\
  &+(1+\zeta^6+\zeta^7) \mathbf{B}_1^{(1)}+(1+\zeta^2+\zeta^{11}) \mathbf{B}_3^{(1)}
   +(1+\zeta^5+\zeta^8) \mathbf{B}_9^{(1)}+\\
  &-(1+\zeta^9+\zeta^4) \mathbf{B}_{12}^{(1)}-(1+\zeta^3+\zeta^{10}) \mathbf{B}_{10}^{(1)}
   -(1+\zeta+\zeta^{12}) \mathbf{B}_4^{(1)}.
\endaligned\eqno{(5.10)}$$
$$\aligned
  &13 S(\mathbf{B}_5)\\
 =&4 \mathbf{B}_0^{(0)}-\mathbf{B}_0^{(1)}+\mathbf{B}_0^{(2)}+\\
  &-[(\zeta^3+\zeta^{10})+2 (\zeta^2+\zeta^{11})+2 (\zeta^9+\zeta^4)
   +2 (\zeta+\zeta^{12})] \mathbf{B}_5+\\
  &-[(\zeta+\zeta^{12})+2 (\zeta^5+\zeta^8)+2 (\zeta^3+\zeta^{10})
   +2 (\zeta^9+\zeta^4)] \mathbf{B}_2+\\
  &-[(\zeta^9+\zeta^4)+2 (\zeta^6+\zeta^7)+2 (\zeta+\zeta^{12})
   +2 (\zeta^3+\zeta^{10})] \mathbf{B}_6+\\
  &-[(\zeta^2+\zeta^{11})+2 (\zeta^3+\zeta^{10})+2 (\zeta^6+\zeta^7)
   +2 (\zeta^5+\zeta^8)] \mathbf{B}_8+\\
  &-[(\zeta^5+\zeta^8)+2 (\zeta+\zeta^{12})+2 (\zeta^2+\zeta^{11})
   +2 (\zeta^6+\zeta^7)] \mathbf{B}_{11}+\\
  &-[(\zeta^6+\zeta^7)+2 (\zeta^9+\zeta^4)+2 (\zeta^5+\zeta^8)
   +2 (\zeta^2+\zeta^{11})] \mathbf{B}_7+\\
  &+(\zeta^9+\zeta^4+\zeta^6+\zeta^7+\zeta^2+\zeta^{11}) \mathbf{B}_1^{(2)}+\\
  &+(\zeta^3+\zeta^{10}+\zeta^2+\zeta^{11}+\zeta^5+\zeta^8) \mathbf{B}_3^{(2)}+\\
  &+(\zeta+\zeta^{12}+\zeta^5+\zeta^8+\zeta^6+\zeta^7) \mathbf{B}_9^{(2)}+\\
  &+(\zeta^9+\zeta^4+\zeta^6+\zeta^7+\zeta^3+\zeta^{10}) \mathbf{B}_{12}^{(2)}+\\
  &+(\zeta^3+\zeta^{10}+\zeta^2+\zeta^{11}+\zeta+\zeta^{12}) \mathbf{B}_{10}^{(2)}+\\
  &+(\zeta+\zeta^{12}+\zeta^5+\zeta^8+\zeta^9+\zeta^4) \mathbf{B}_4^{(2)}+\\
  &+(\zeta^9+\zeta^4+\zeta^6+\zeta^7+\zeta^2+\zeta^{11}) \mathbf{B}_1^{(1)}+\\
  &+(\zeta^3+\zeta^{10}+\zeta^2+\zeta^{11}+\zeta^5+\zeta^8) \mathbf{B}_3^{(1)}+\\
  &+(\zeta+\zeta^{12}+\zeta^5+\zeta^8+\zeta^6+\zeta^7) \mathbf{B}_9^{(1)}+\\
  &-(\zeta^9+\zeta^4+\zeta^6+\zeta^7+\zeta^3+\zeta^{10}) \mathbf{B}_{12}^{(1)}+\\
  &-(\zeta^3+\zeta^{10}+\zeta^2+\zeta^{11}+\zeta+\zeta^{12}) \mathbf{B}_{10}^{(1)}+\\
  &-(\zeta+\zeta^{12}+\zeta^5+\zeta^8+\zeta^9+\zeta^4) \mathbf{B}_4^{(1)}.
\endaligned\eqno{(5.11)}$$
$$\aligned
  &13 S(\mathbf{B}_8)\\
 =&4 \mathbf{B}_0^{(0)}-\mathbf{B}_0^{(1)}+\mathbf{B}_0^{(2)}+\\
  &-[(\zeta^2+\zeta^{11})+2 (\zeta^3+\zeta^{10})+2 (\zeta^6+\zeta^7)
   +2 (\zeta^5+\zeta^8)] \mathbf{B}_5+\\
  &-[(\zeta^5+\zeta^8)+2 (\zeta+\zeta^{12})+2 (\zeta^2+\zeta^{11})
   +2 (\zeta^6+\zeta^7)] \mathbf{B}_2+\\
  &-[(\zeta^6+\zeta^7)+2 (\zeta^9+\zeta^4)+2 (\zeta^5+\zeta^8)
   +2 (\zeta^2+\zeta^{11})] \mathbf{B}_6+\\
  &-[(\zeta^3+\zeta^{10})+2 (\zeta^2+\zeta^{11})+2 (\zeta^9+\zeta^4)
   +2 (\zeta+\zeta^{12})] \mathbf{B}_8+\\
  &-[(\zeta+\zeta^{12})+2 (\zeta^5+\zeta^8)+2 (\zeta^3+\zeta^{10})
   +2 (\zeta^9+\zeta^4)] \mathbf{B}_{11}+\\
  &-[(\zeta^9+\zeta^4)+2 (\zeta^6+\zeta^7)+2 (\zeta+\zeta^{12})
   +2 (\zeta^3+\zeta^{10})] \mathbf{B}_7+\\
  &+(\zeta^9+\zeta^4+\zeta^6+\zeta^7+\zeta^3+\zeta^{10}) \mathbf{B}_1^{(2)}+\\
  &+(\zeta^3+\zeta^{10}+\zeta^2+\zeta^{11}+\zeta+\zeta^{12}) \mathbf{B}_3^{(2)}+\\
  &+(\zeta+\zeta^{12}+\zeta^5+\zeta^8+\zeta^9+\zeta^4) \mathbf{B}_9^{(2)}+\\
  &+(\zeta^9+\zeta^4+\zeta^6+\zeta^7+\zeta^2+\zeta^{11}) \mathbf{B}_{12}^{(2)}+\\
  &+(\zeta^3+\zeta^{10}+\zeta^2+\zeta^{11}+\zeta^5+\zeta^8) \mathbf{B}_{10}^{(2)}+\\
  &+(\zeta+\zeta^{12}+\zeta^5+\zeta^8+\zeta^6+\zeta^7) \mathbf{B}_4^{(2)}+\\
  &+(\zeta^9+\zeta^4+\zeta^6+\zeta^7+\zeta^3+\zeta^{10}) \mathbf{B}_1^{(1)}+\\
  &+(\zeta^3+\zeta^{10}+\zeta^2+\zeta^{11}+\zeta+\zeta^{12}) \mathbf{B}_3^{(1)}+\\
  &+(\zeta+\zeta^{12}+\zeta^5+\zeta^8+\zeta^9+\zeta^4) \mathbf{B}_9^{(1)}+\\
  &-(\zeta^9+\zeta^4+\zeta^6+\zeta^7+\zeta^2+\zeta^{11}) \mathbf{B}_{12}^{(1)}+\\
  &-(\zeta^3+\zeta^{10}+\zeta^2+\zeta^{11}+\zeta^5+\zeta^8) \mathbf{B}_{10}^{(1)}+\\
  &-(\zeta+\zeta^{12}+\zeta^5+\zeta^8+\zeta^6+\zeta^7) \mathbf{B}_4^{(1)}.
\endaligned\eqno{(5.12)}$$

  The other terms: $S(\mathbf{B}_3^{(1)})$, $S(\mathbf{B}_{10}^{(1)})$,
$S(\mathbf{B}_3^{(2)})$, $S(\mathbf{B}_{10}^{(2)})$, $S(\mathbf{B}_2)$,
$S(\mathbf{B}_{11})$, $S(\mathbf{B}_9^{(1)})$, $S(\mathbf{B}_4^{(1)})$,
$S(\mathbf{B}_9^{(2)})$, $S(\mathbf{B}_4^{(2)})$, $S(\mathbf{B}_6)$, and
$S(\mathbf{B}_7)$ can be obtained by the permutation on the above six terms.
For each of the above triples, we can obtain the other two terms by the
following permutation
$$\zeta^9 \mapsto \zeta^3 \mapsto \zeta, \quad
  \zeta^4 \mapsto \zeta^{10} \mapsto \zeta^{12}, \quad
  \zeta^6 \mapsto \zeta^2 \mapsto \zeta^5, \quad
  \zeta^7 \mapsto \zeta^{11} \mapsto \zeta^8\eqno{(5.13)}$$
on the first term. For example, for the first triple
$(S(\mathbb{B}_1^{(1)}), S(\mathbb{B}_3^{(1)}), S(\mathbb{B}_9^{(1)}))$,
we can get the expression of $S(\mathbb{B}_3^{(1)})$ and $S(\mathbb{B}_9^{(1)})$
by the above permutation (5.13) on $S(\mathbb{B}_1^{(1)})$ as follows:
$$\aligned
  &13 S(\mathbf{B}_3^{(1)})\\
 =&6 \mathbf{B}_0^{(0)}+\frac{-3+\sqrt{13}}{2} \mathbf{B}_0^{(1)}+
   \frac{3+\sqrt{13}}{2} \mathbf{B}_0^{(2)}+\\
  &+3 [(\zeta^3+\zeta^{10}+\zeta^2+\zeta^{11}+\zeta^5+\zeta^8) \mathbf{B}_5+\\
  &+(\zeta+\zeta^{12}+\zeta^5+\zeta^8+\zeta^6+\zeta^7) \mathbf{B}_2+\\
  &+(\zeta^9+\zeta^4+\zeta^6+\zeta^7+\zeta^2+\zeta^{11}) \mathbf{B}_6+\\
  &+(\zeta^3+\zeta^{10}+\zeta^2+\zeta^{11}+\zeta+\zeta^{12}) \mathbf{B}_8+\\
  &+(\zeta+\zeta^{12}+\zeta^5+\zeta^8+\zeta^9+\zeta^4) \mathbf{B}_{11}+\\
  &+(\zeta^9+\zeta^4+\zeta^6+\zeta^7+\zeta^3+\zeta^{10}) \mathbf{B}_7]+\\
  &+3 [(1+\zeta^3+\zeta^{10}) \mathbf{B}_1^{(2)}
   +(1+\zeta+\zeta^{12}) \mathbf{B}_3^{(2)}
   +(1+\zeta^9+\zeta^4) \mathbf{B}_9^{(2)}+\\
  &+(1+\zeta^2+\zeta^{11}) \mathbf{B}_{12}^{(2)}
   +(1+\zeta^5+\zeta^8) \mathbf{B}_{10}^{(2)}
   +(1+\zeta^6+\zeta^7) \mathbf{B}_4^{(2)}]+\\
  &+[1+(\zeta^3+\zeta^{10})+2 (\zeta^9+\zeta^4)+(\zeta^5+\zeta^8)] \mathbf{B}_1^{(1)}+\\
  &+[1+(\zeta+\zeta^{12})+2 (\zeta^3+\zeta^{10})+(\zeta^6+\zeta^7)] \mathbf{B}_3^{(1)}+\\
  &+[1+(\zeta^9+\zeta^4)+2 (\zeta+\zeta^{12})+(\zeta^2+\zeta^{11})] \mathbf{B}_9^{(1)}+\\
  &-[1+(\zeta^2+\zeta^{11})+2 (\zeta^6+\zeta^7)+(\zeta+\zeta^{12})] \mathbf{B}_{12}^{(1)}+\\
  &-[1+(\zeta^5+\zeta^8)+2 (\zeta^2+\zeta^{11})+(\zeta^9+\zeta^4)] \mathbf{B}_{10}^{(1)}+\\
  &-[1+(\zeta^6+\zeta^7)+2 (\zeta^5+\zeta^8)+(\zeta^3+\zeta^{10})] \mathbf{B}_4^{(1)}.
\endaligned\eqno{(5.14)}$$
$$\aligned
  &13 S(\mathbf{B}_9^{(1)})\\
 =&6 \mathbf{B}_0^{(0)}+\frac{-3+\sqrt{13}}{2} \mathbf{B}_0^{(1)}+
   \frac{3+\sqrt{13}}{2} \mathbf{B}_0^{(2)}+\\
  &+3 [(\zeta+\zeta^{12}+\zeta^5+\zeta^8+\zeta^6+\zeta^7) \mathbf{B}_5+\\
  &+(\zeta^9+\zeta^4+\zeta^6+\zeta^7+\zeta^2+\zeta^{11}) \mathbf{B}_2+\\
  &+(\zeta^3+\zeta^{10}+\zeta^2+\zeta^{11}+\zeta^5+\zeta^8) \mathbf{B}_6+\\
  &+(\zeta+\zeta^{12}+\zeta^5+\zeta^8+\zeta^9+\zeta^4) \mathbf{B}_8+\\
  &+(\zeta^9+\zeta^4+\zeta^6+\zeta^7+\zeta^3+\zeta^{10}) \mathbf{B}_{11}+\\
  &+(\zeta^3+\zeta^{10}+\zeta^2+\zeta^{11}+\zeta+\zeta^{12}) \mathbf{B}_7]+\\
  &+3 [(1+\zeta+\zeta^{12}) \mathbf{B}_1^{(2)}
   +(1+\zeta^9+\zeta^4) \mathbf{B}_3^{(2)}
   +(1+\zeta^3+\zeta^{10}) \mathbf{B}_9^{(2)}+\\
  &+(1+\zeta^5+\zeta^8) \mathbf{B}_{12}^{(2)}
   +(1+\zeta^6+\zeta^7) \mathbf{B}_{10}^{(2)}
   +(1+\zeta^2+\zeta^{11}) \mathbf{B}_4^{(2)}]+\\
  &+[1+(\zeta+\zeta^{12})+2 (\zeta^3+\zeta^{10})+(\zeta^6+\zeta^7)] \mathbf{B}_1^{(1)}+\\
  &+[1+(\zeta^9+\zeta^4)+2 (\zeta+\zeta^{12})+(\zeta^2+\zeta^{11})] \mathbf{B}_3^{(1)}+\\
  &+[1+(\zeta^3+\zeta^{10})+2 (\zeta^9+\zeta^4)+(\zeta^5+\zeta^8)] \mathbf{B}_9^{(1)}+\\
  &-[1+(\zeta^5+\zeta^8)+2 (\zeta^2+\zeta^{11})+(\zeta^9+\zeta^4)] \mathbf{B}_{12}^{(1)}+\\
  &-[1+(\zeta^6+\zeta^7)+2 (\zeta^5+\zeta^8)+(\zeta^3+\zeta^{10})] \mathbf{B}_{10}^{(1)}+\\
  &-[1+(\zeta^2+\zeta^{11})+2 (\zeta^6+\zeta^7)+(\zeta+\zeta^{12})] \mathbf{B}_4^{(1)}.
\endaligned\eqno{(5.15)}$$
The other five triples can be obtained in the similar way.

  On the other hand, we have
$$T(\mathbf{B}_0^{(0)})=\mathbf{B}_0^{(0)}, \quad
  T(\mathbf{B}_0^{(1)})=\mathbf{B}_0^{(1)}, \quad
  T(\mathbf{B}_0^{(2)})=\mathbf{B}_0^{(2)}.$$
$$T(\mathbf{B}_1^{(1)})=\zeta \mathbf{B}_1^{(1)}, \quad
  T(\mathbf{B}_1^{(2)})=\zeta \mathbf{B}_1^{(2)}.$$
$$T(\mathbf{B}_3^{(1)})=\zeta^3 \mathbf{B}_3^{(1)}, \quad
  T(\mathbf{B}_3^{(2)})=\zeta^3 \mathbf{B}_3^{(2)}.$$
$$T(\mathbf{B}_9^{(1)})=\zeta^9 \mathbf{B}_9^{(1)}, \quad
  T(\mathbf{B}_9^{(2)})=\zeta^9 \mathbf{B}_9^{(2)}.$$
$$T(\mathbf{B}_{12}^{(1)})=\zeta^{12} \mathbf{B}_{12}^{(1)}, \quad
  T(\mathbf{B}_{12}^{(2)})=\zeta^{12} \mathbf{B}_{12}^{(2)}.$$
$$T(\mathbf{B}_{10}^{(1)})=\zeta^{10} \mathbf{B}_{10}^{(1)}, \quad
  T(\mathbf{B}_{10}^{(2)})=\zeta^{10} \mathbf{B}_{10}^{(2)}.$$
$$T(\mathbf{B}_4^{(1)})=\zeta^4 \mathbf{B}_4^{(1)}, \quad
  T(\mathbf{B}_4^{(2)})=\zeta^4 \mathbf{B}_4^{(2)}.$$
$$T(\mathbf{B}_5)=\zeta^5 \mathbf{B}_5, \quad
  T(\mathbf{B}_8)=\zeta^8 \mathbf{B}_8.$$
$$T(\mathbf{B}_2)=\zeta^2 \mathbf{B}_2, \quad
  T(\mathbf{B}_{11})=\zeta^{11} \mathbf{B}_{11}.$$
$$T(\mathbf{B}_6)=\zeta^6 \mathbf{B}_6, \quad
  T(\mathbf{B}_7)=\zeta^7 \mathbf{B}_7.$$

  Let $\widetilde{S}$ and $\widetilde{T}$ denote the corresponding matrices
of order $21$ corresponding to $S$ and $T$. Then they have the following trace:
$$\text{Tr}(\widetilde{S})=1, \quad \text{Tr}(\widetilde{T})=\frac{3+\sqrt{13}}{2}.
  \eqno{(5.16)}$$
This completes the proof of Theorem 5.3.

\noindent
$\qquad \qquad \qquad \qquad \qquad \qquad \qquad \qquad \qquad
 \qquad \qquad \qquad \qquad \qquad \qquad \qquad \boxed{}$

  This leads to the study of the following ring of invariant polynomials:
$$\left(\mathbb{C}[z_1, \ldots, z_6]/I(Y)\right)^{\text{SL}(2, 13)}\eqno{(5.17)}$$
and the cone over the corresponding curve $Y$:
$$C_Y/\text{SL}(2, 13) \cong \text{Spec}\left(\left(\mathbb{C}[z_1,
  \ldots, z_6]/I(Y)\right)^{\text{SL}(2, 13)}\right).\eqno{(5.18)}$$
Here, let us recall the cone over a projective variety: let $Y \subset
\mathbb{P}^n$ be a nonempty algebraic set over the field $k$, and let
$\theta: \mathbb{A}^{n+1} -\{ (0, \ldots, 0) \} \rightarrow \mathbb{P}^n$
be the map which sends the point with affine coordinates $(a_0, \ldots,
a_n)$ to the point with homogeneous coordinates $(a_0, \ldots, a_n)$.
We define the affine cone over $Y$ to be
$$C_Y=\theta^{-1}(Y) \cup \{ (0, \ldots, 0) \}.$$
Then $C_Y$ is an algebraic set in $\mathbb{A}^{n+1}$, whose ideal is
equal to $I(Y)$, considered as an ordinary ideal in $k[x_0, \ldots, x_n]$

\begin{center}
{\large\bf 6. The decomposition of a $21$-dimensional reducible representation
              of $\text{SL}(2, 13)$ and its geometric realization}
\end{center}

\begin{center}
{\bf 6.1. A $21$-dimensional reducible representation of $\text{SL}(2, 13)$
          and its decomposition: $\mathbf{21}=\mathbf{1} \oplus \mathbf{7}
          \oplus \mathbf{13}$}
\end{center}

  The irreducible representations of $\text{SL}(2, q)$ are well known
for a very long time (see \cite{Fr1}, \cite{J} and \cite{S2}) and are a
prototype example in many introductory courses on the subject (see
\cite{Bon}). For $G=\text{SL}(2, q)$, there are at most
$$1+1+\frac{q-3}{2}+2+\frac{q-1}{2}+2=q+4$$
conjugacy classes which can appear in the following formula on the
dimensions of ordinary representations:
$$|G|=1^2+q^2+\frac{q-3}{2} (q+1)^2+2 \cdot \left(\frac{q+1}{2}\right)^2
     +\frac{q-1}{2} (q-1)^2+2 \cdot \left(\frac{q-1}{2}\right)^2.\eqno{(6.1)}$$
In particular, Schur (see \cite{S2}) derived the character tables of
the groups $\text{SL}(2, q)$ for all values of $q$.

  In order to calculate the character table of $G$, algebraic methods
(in particular Harish-Chandra induction) give roughly half of the
irreducible characters (non-cuspidal characters with number $\frac{q+5}{2}$).
The other half (the cuspidal characters with number $\frac{q+3}{2}$)
can be obtained by the $\ell$-adic cohomology of the Drinfeld curve,
i.e., the Deligne-Lusztig induction (see \cite{Bon}).

  In particular, when $q=13$, $|G|=2184$. There are at most $q+4=17$
conjugacy classes and
$$1^2+13^2+5 \times 14^2+2 \times 7^2+6 \times 12^2+2 \times 6^2=2184$$
corresponding to the decomposition of conjugacy classes:
$$1+1+5+2+6+2=17.$$

  In general, let $G=G(\mathbb{F}_q)$ be the group of $\mathbb{F}_q$-rational
points of a reductive algebraic group defined over $\mathbb{F}_q$ where $q$
is a power of a prime number $p$. Let $\widehat{G}$ be the set of isomorphism
classes of irreducible representations of $G$. The Gel'fand-Harish-Chandra's
philosophy of cusp forms asserts that each member $\rho$ of $\widehat{G}$ can
be realized (see \cite{GH}), in some effective manner, insider the induction
$\text{Ind}_P^G(\sigma)$ from a cuspidal representation (i.e., one that is not
induced from a smaller parabolic subgroup) of a parabolic subgroup $P$ of $G$
$$\widehat{G}=\left\{\aligned
  &\rho < \text{Ind}_P^G(\sigma)\\
  &G \supset P-\text{parabolic}\\
  &\text{$\sigma$ cuspidal representation of $P$}
  \endaligned\right\}.$$
The above philosophy is a central one and leads to important developments in
representation theory of reductive groups over local and finite fields, in
particular the work of Deligne-Lusztig on the construction of representations
of finite reductive groups and Lusztig's classification of these representations.

  Let us see how this philosophy manifests itself in the cases when $G$ is
$\text{GL}_2(\mathbb{F}_q)$ of $\text{SL}_2(\mathbb{F}_q)$ (see \cite{GH}).
For simplicity, we ignore for the one dimensional representations. In both
cases around half of the irreducible representations are cuspidal, i.e.,
$P=G$, and form the so called discrete series, and half are the so called
principal series, i.e., induced from the Borel subgroup $P=B$ of upper
triangular matrices in $G$. Moreover, it is well known that the discrete
and principal series representations of $\text{SL}_2(\mathbb{F}_q)$ can be
obtained by restriction from the corresponding irreducible representations
of $\text{GL}_2(\mathbb{F}_q)$. These restriction are typically irreducible,
with only two exceptions: one discrete series representation of dimension
$q-1$ and one principal series representation of dimension $q+1$ split
into two pieces. Hence one has two discrete and two principal series
representations of dimensions $\frac{q-1}{2}$ and $\frac{q+1}{2}$,
respectively. These representations are called degenerate discrete series
and degenerate principal series. In general, For $G=\text{SL}_2(\mathbb{F}_q)$,
we have the trivial representation of dimension $1$, the Steinberg representation
of dimension $q$, induced representations $\text{Ind}_B^G \chi$ ($\chi$ a
character of the split torus) of dimension $q+1$ and discrete series (or
super-cuspidal) representations of dimension $q-1$. This discrete series
are associated to the characters of the non-split torus.
$$\begin{matrix}
 &G  & \text{discrete series ($P=G$)} \quad & \text{principal series ($P=B$)}\\
 & \text{GL}(2, \mathbb{F}_q) &\text{dim}(\rho)=q-1 \quad
 & \text{dim($\rho$)$=q+1$ or $q$}\\
 & \text{SL}(2, \mathbb{F}_q) &\text{dim}(\rho)=\left\{\aligned
                               &\text{$q-1$, or}\\
                               &\frac{q-1}{2}
                               \endaligned\right.
                              &\text{dim}(\rho)=\left\{\aligned
                               &\text{$q+1$ or $q$}\\
                               &\frac{q+1}{2}
                               \endaligned\right.
\end{matrix}$$

  In particular, when $q=13$ and $G=\text{SL}(2, 13)$. There are the
following representations $\rho$:

(1) discrete series
$$\text{dim}(\rho)=\left\{\aligned
  &12 \quad \text{(super-cuspidal representations) or}\\
  &6 \quad \text{(two degenerate discrete series representations)}
\endaligned\right.$$

(2) principal series
$$\text{dim}(\rho)=\left\{\aligned
  &14 \quad \text{(principal series representations) or}\\
  &13 \quad \text{(Steinberg representation) or}\\
  &7 \quad \text{(two degenerate principal series representations)}
\endaligned\right.$$

\textbf{Theorem 6.1}. {\it The $21$-dimensional representation is reducible,
which can be decomposed as the direct sum of a $1$-dimensional representation
$($the trivial representation$)$, a $7$-dimensional representation $($the
degenerate principal series representation) and a $13$-dimensional representation
$($the Steinberg representation):
$$\mathbf{21}=\mathbf{1} \oplus \mathbf{7} \oplus \mathbf{13}.\eqno{(6.2)}$$
More precisely, let $V$ be a complex vector space generated by the $21$ quartic
polynomials $\mathbf{B}_0^{(0)}$, $\mathbf{B}_0^{(1)}$, $\ldots$, $\mathbf{B}_7$.
Then $V$ has the following decomposition:
$$V=V_1 \oplus V_7 \oplus V_{13},\eqno{(6.3)}$$
where the $1$-dimensional subspace
$$V_1=\langle 3 \mathbf{B}_0^{(0)}+\mathbf{B}_0^{(1)}-\mathbf{B}_0^{(2)} \rangle,
      \eqno{(6.4)}$$
the $7$-dimensional subspace
$$\aligned
  V_7=\langle &\mathbf{B}_0^{(1)}+\mathbf{B}_0^{(2)},
               \mathbf{B}_1^{(1)}-3 \mathbf{B}_1^{(2)},
               \mathbf{B}_3^{(1)}-3 \mathbf{B}_3^{(2)},
               \mathbf{B}_9^{(1)}-3 \mathbf{B}_9^{(2)},\\
              &\mathbf{B}_{12}^{(1)}+3 \mathbf{B}_{12}^{(2)},
               \mathbf{B}_{10}^{(1)}+3 \mathbf{B}_{10}^{(2)},
               \mathbf{B}_4^{(1)}+3 \mathbf{B}_4^{(2)}
\rangle\endaligned\eqno{(6.5)}$$
and the $13$-dimensional subspace
$$\aligned
  V_{13}=\langle &4 \mathbf{B}_0^{(0)}-\mathbf{B}_0^{(1)}+\mathbf{B}_0^{(2)},
               \mathbf{B}_5, \mathbf{B}_2, \mathbf{B}_6, \mathbf{B}_8,
               \mathbf{B}_{11}, \mathbf{B}_7,\\
              &\mathbf{B}_1^{(1)}+\mathbf{B}_1^{(2)},
               \mathbf{B}_3^{(1)}+\mathbf{B}_3^{(2)},
               \mathbf{B}_9^{(1)}+\mathbf{B}_9^{(2)},\\
              &-\mathbf{B}_{12}^{(1)}+\mathbf{B}_{12}^{(2)},
               -\mathbf{B}_{10}^{(1)}+\mathbf{B}_{10}^{(2)},
               -\mathbf{B}_4^{(1)}+\mathbf{B}_4^{(2)} \rangle
\endaligned\eqno{(6.6)}$$
are stable under the action of $G \cong \text{SL}(2, 13)$.}

{\it Proof}. Our $21$-dimensional representation does not belong to the
discrete series or the principal series. By the proof of Theorem 5.2 and
Theorem 5.3, we have
$$S(3 \mathbf{B}_0^{(0)}+\mathbf{B}_0^{(1)}-\mathbf{B}_0^{(2)})
 =3 \mathbf{B}_0^{(0)}+\mathbf{B}_0^{(1)}-\mathbf{B}_0^{(2)}.\eqno{(6.7)}$$
$$T(3 \mathbf{B}_0^{(0)}+\mathbf{B}_0^{(1)}-\mathbf{B}_0^{(2)})
 =3 \mathbf{B}_0^{(0)}+\mathbf{B}_0^{(1)}-\mathbf{B}_0^{(2)}.\eqno{(6.8)}$$
Let $V$ be a complex vector space generated by the $21$ quartic polynomials
$\mathbf{B}_0^{(0)}$, $\mathbf{B}_0^{(1)}$, $\ldots$, $\mathbf{B}_7$. Put
$$V_1=\langle 3 \mathbf{B}_0^{(0)}+\mathbf{B}_0^{(1)}-\mathbf{B}_0^{(2)} \rangle.
  \eqno{(6.9)}$$
Then $V_1$ is a one-dimensional subspace which is stable under the action of
$G \cong \text{SL}(2, 13)$. Thus,
$$V=V_1 \oplus V_1^0,\eqno{(6.10)}$$
where $V_1^0$ is the complement of $V_1$ in $V$. This shows that this $21$-dimensional
representation is reducible. It can be decomposed as a direct sum of a
$1$-dimensional representation and a $20$-dimensional representation
$$\mathbf{21}=\mathbf{1} \oplus \mathbf{20},\eqno{(6.11)}$$
the latter one can be decomposed into the sum of irreducible representations
listed as above. More precisely, we have
$$\aligned
  &\sqrt{13} S(\mathbf{B}_0^{(1)}+\mathbf{B}_0^{(2)})\\
 =&(\mathbf{B}_0^{(1)}+\mathbf{B}_0^{(2)})
  +(\mathbf{B}_1^{(1)}-3 \mathbf{B}_1^{(2)})
  +(\mathbf{B}_3^{(1)}-3 \mathbf{B}_3^{(2)})
  +(\mathbf{B}_9^{(1)}-3 \mathbf{B}_9^{(2)})\\
 +&(\mathbf{B}_{12}^{(1)}+3 \mathbf{B}_{12}^{(2)})+
   (\mathbf{B}_{10}^{(1)}+3 \mathbf{B}_{10}^{(2)})+
   (\mathbf{B}_4^{(1)}+3 \mathbf{B}_4^{(2)}).
\endaligned\eqno{(6.12)}$$
$$T(\mathbf{B}_0^{(1)}+\mathbf{B}_0^{(2)})
 =\mathbf{B}_0^{(1)}+\mathbf{B}_0^{(2)}.\eqno{(6.13)}$$
Without loss of generality, we can only consider the action of $S$
and $T$ on $\mathbf{B}_1^{(1)}-3 \mathbf{B}_1^{(2)}$ and
$\mathbf{B}_{12}^{(1)}+3 \mathbf{B}_{12}^{(2)}$. The action of $S$
and $T$ on
$\mathbf{B}_3^{(1)}-3 \mathbf{B}_3^{(2)}$,
$\mathbf{B}_9^{(1)}-3 \mathbf{B}_9^{(2)}$,
$\mathbf{B}_{10}^{(1)}+3 \mathbf{B}_{10}^{(2)}$,
and $\mathbf{B}_4^{(1)}+3 \mathbf{B}_4^{(2)}$
can be done similarly. We have
$$\aligned
  &13 S(\mathbf{B}_1^{(1)}-3 \mathbf{B}_1^{(2)})\\
 =&2 \sqrt{13} (\mathbf{B}_0^{(1)}+\mathbf{B}_0^{(2)})\\
  &+[-2 (1+\zeta^9+\zeta^4)+2 (\zeta+\zeta^{12})+(\zeta^2+\zeta^{11})]
  (\mathbf{B}_1^{(1)}-3 \mathbf{B}_1^{(2)})\\
  &+[-2 (1+\zeta^3+\zeta^{10})+2 (\zeta^9+\zeta^4)+(\zeta^5+\zeta^8)]
  (\mathbf{B}_3^{(1)}-3 \mathbf{B}_3^{(2)})\\
  &+[-2 (1+\zeta+\zeta^{12})+2 (\zeta^3+\zeta^{10})+(\zeta^6+\zeta^7)]
  (\mathbf{B}_9^{(1)}-3 \mathbf{B}_9^{(2)})\\
  &+[2 (1+\zeta^6+\zeta^7)-2 (\zeta^5+\zeta^8)-(\zeta^3+\zeta^{10})]
  (\mathbf{B}_{12}^{(1)}+3 \mathbf{B}_{12}^{(2)})\\
  &+[2 (1+\zeta^2+\zeta^{11})-2 (\zeta^6+\zeta^7)-(\zeta+\zeta^{12})]
  (\mathbf{B}_{10}^{(1)}+3 \mathbf{B}_{10}^{(2)})\\
  &+[2 (1+\zeta^5+\zeta^8)-2 (\zeta^2+\zeta^{11})-(\zeta^9+\zeta^4)]
  (\mathbf{B}_4^{(1)}+3 \mathbf{B}_4^{(2)}).
\endaligned\eqno{(6.14)}$$
$$\aligned
  &13 S(\mathbf{B}_{12}^{(1)}+3 \mathbf{B}_{12}^{(2)})\\
 =&2 \sqrt{13} (\mathbf{B}_0^{(1)}+\mathbf{B}_0^{(2)})\\
  &+[2 (1+\zeta^6+\zeta^7)-2 (\zeta^5+\zeta^8)-(\zeta^3+\zeta^{10})]
  (\mathbf{B}_1^{(1)}-3 \mathbf{B}_1^{(2)})\\
  &+[2 (1+\zeta^2+\zeta^{11})-2 (\zeta^6+\zeta^7)-(\zeta+\zeta^{12})]
  (\mathbf{B}_3^{(1)}-3 \mathbf{B}_3^{(2)})\\
  &+[2 (1+\zeta^5+\zeta^8)-2 (\zeta^2+\zeta^{11})-(\zeta^9+\zeta^4)]
  (\mathbf{B}_9^{(1)}-3 \mathbf{B}_9^{(2)})\\
  &+[-2 (1+\zeta^9+\zeta^4)+2 (\zeta+\zeta^{12})+(\zeta^2+\zeta^{11})]
  (\mathbf{B}_{12}^{(1)}+3 \mathbf{B}_{12}^{(2)})\\
  &+[-2 (1+\zeta^3+\zeta^{10})+2 (\zeta^9+\zeta^4)+(\zeta^5+\zeta^8)]
  (\mathbf{B}_{10}^{(1)}+3 \mathbf{B}_{10}^{(2)})\\
  &+[-2 (1+\zeta+\zeta^{12})+2 (\zeta^3+\zeta^{10})+(\zeta^6+\zeta^7)]
  (\mathbf{B}_4^{(1)}+3 \mathbf{B}_4^{(2)}).
\endaligned\eqno{(6.15)}$$
$$T(\mathbf{B}_1^{(1)}-3 \mathbf{B}_1^{(2)})
 =\zeta (\mathbf{B}_1^{(1)}-3 \mathbf{B}_1^{(2)}).\eqno{(6.16)}$$
$$T(\mathbf{B}_{12}^{(1)}+3 \mathbf{B}_{12}^{(2)})
 =\zeta^{12} (\mathbf{B}_{12}^{(1)}+3 \mathbf{B}_{12}^{(2)}).\eqno{(6.17)}$$
Put
$$\aligned
  V_7=\langle &\mathbf{B}_0^{(1)}+\mathbf{B}_0^{(2)},
               \mathbf{B}_1^{(1)}-3 \mathbf{B}_1^{(2)},
               \mathbf{B}_3^{(1)}-3 \mathbf{B}_3^{(2)},
               \mathbf{B}_9^{(1)}-3 \mathbf{B}_9^{(2)},\\
              &\mathbf{B}_{12}^{(1)}+3 \mathbf{B}_{12}^{(2)},
               \mathbf{B}_{10}^{(1)}+3 \mathbf{B}_{10}^{(2)},
               \mathbf{B}_4^{(1)}+3 \mathbf{B}_4^{(2)}
\rangle.\endaligned\eqno{(6.18)}$$
Then $V_7$ is a seven-dimensional subspace which is stable under the
action of $G \cong \text{SL}(2, 13)$. Moreover,
$$\aligned
  &13 S(4 \mathbf{B}_0^{(0)}-\mathbf{B}_0^{(1)}+\mathbf{B}_0^{(2)})\\
 =&-(4 \mathbf{B}_0^{(0)}-\mathbf{B}_0^{(1)}+\mathbf{B}_0^{(2)})+14 \mathbf{B}_5
   +14 \mathbf{B}_2+14 \mathbf{B}_6+14 \mathbf{B}_8+14 \mathbf{B}_{11}+14 \mathbf{B}_7\\
  &+7 (\mathbf{B}_1^{(1)}+\mathbf{B}_1^{(2)})
   +7 (\mathbf{B}_3^{(1)}+\mathbf{B}_3^{(2)})
   +7 (\mathbf{B}_9^{(1)}+\mathbf{B}_9^{(2)})\\
  &+7 (-\mathbf{B}_{12}^{(1)}+\mathbf{B}_{12}^{(2)})
   +7 (-\mathbf{B}_{10}^{(1)}+\mathbf{B}_{10}^{(2)})
   +7 (-\mathbf{B}_4^{(1)}+\mathbf{B}_4^{(2)}).
\endaligned\eqno{(6.19)}$$
$$T(4 \mathbf{B}_0^{(0)}-\mathbf{B}_0^{(1)}+\mathbf{B}_0^{(2)})
 =4 \mathbf{B}_0^{(0)}-\mathbf{B}_0^{(1)}+\mathbf{B}_0^{(2)}.\eqno{(6.20)}$$
Without loss of generality, we can only consider the action of $S$
and $T$ on $\mathbf{B}_5$, $\mathbf{B}_8$, $\mathbf{B}_1^{(1)}+\mathbf{B}_1^{(2)}$,
$-\mathbf{B}_{12}^{(1)}+\mathbf{B}_{12}^{(2)}$. The action of $S$
and $T$ on $\mathbf{B}_2$, $\mathbf{B}_6$, $\mathbf{B}_{11}$, $\mathbf{B}_7$,
$\mathbf{B}_3^{(1)}+\mathbf{B}_3^{(2)}$, $\mathbf{B}_9^{(1)}+\mathbf{B}_9^{(2)}$,
$-\mathbf{B}_{10}^{(1)}+\mathbf{B}_{10}^{(2)}$ and $-\mathbf{B}_4^{(1)}+\mathbf{B}_4^{(2)}$
can be done similarly. We have
$$\aligned
  &13 S(\mathbf{B}_5)\\
 =&4 \mathbf{B}_0^{(0)}-\mathbf{B}_0^{(1)}+\mathbf{B}_0^{(2)}+\\
  &-[(\zeta^3+\zeta^{10})+2 (\zeta^2+\zeta^{11})+2 (\zeta^9+\zeta^4)
   +2 (\zeta+\zeta^{12})] \mathbf{B}_5+\\
  &-[(\zeta+\zeta^{12})+2 (\zeta^5+\zeta^8)+2 (\zeta^3+\zeta^{10})
   +2 (\zeta^9+\zeta^4)] \mathbf{B}_2+\\
  &-[(\zeta^9+\zeta^4)+2 (\zeta^6+\zeta^7)+2 (\zeta+\zeta^{12})
   +2 (\zeta^3+\zeta^{10})] \mathbf{B}_6+\\
  &-[(\zeta^2+\zeta^{11})+2 (\zeta^3+\zeta^{10})+2 (\zeta^6+\zeta^7)
   +2 (\zeta^5+\zeta^8)] \mathbf{B}_8+\\
  &-[(\zeta^5+\zeta^8)+2 (\zeta+\zeta^{12})+2 (\zeta^2+\zeta^{11})
   +2 (\zeta^6+\zeta^7)] \mathbf{B}_{11}+\\
  &-[(\zeta^6+\zeta^7)+2 (\zeta^9+\zeta^4)+2 (\zeta^5+\zeta^8)
   +2 (\zeta^2+\zeta^{11})] \mathbf{B}_7+\\
  &+(\zeta^9+\zeta^4+\zeta^6+\zeta^7+\zeta^2+\zeta^{11})
  (\mathbf{B}_1^{(1)}+\mathbf{B}_1^{(2)})+\\
  &+(\zeta^3+\zeta^{10}+\zeta^2+\zeta^{11}+\zeta^5+\zeta^8)
  (\mathbf{B}_3^{(1)}+\mathbf{B}_3^{(2)})+\\
  &+(\zeta+\zeta^{12}+\zeta^5+\zeta^8+\zeta^6+\zeta^7)
  (\mathbf{B}_9^{(1)}+\mathbf{B}_9^{(2)})+\\
  &+(\zeta^9+\zeta^4+\zeta^6+\zeta^7+\zeta^3+\zeta^{10})
  (-\mathbf{B}_{12}^{(1)}+\mathbf{B}_{12}^{(2)})+\\
  &+(\zeta^3+\zeta^{10}+\zeta^2+\zeta^{11}+\zeta+\zeta^{12})
  (-\mathbf{B}_{10}^{(1)}+\mathbf{B}_{10}^{(2)})+\\
  &+(\zeta+\zeta^{12}+\zeta^5+\zeta^8+\zeta^9+\zeta^4)
  (-\mathbf{B}_4^{(1)}+\mathbf{B}_4^{(2)}).
\endaligned\eqno{(6.21)}$$
$$\aligned
  &13 S(\mathbf{B}_8)\\
 =&4 \mathbf{B}_0^{(0)}-\mathbf{B}_0^{(1)}+\mathbf{B}_0^{(2)}+\\
  &-[(\zeta^2+\zeta^{11})+2 (\zeta^3+\zeta^{10})+2 (\zeta^6+\zeta^7)
   +2 (\zeta^5+\zeta^8)] \mathbf{B}_5+\\
  &-[(\zeta^5+\zeta^8)+2 (\zeta+\zeta^{12})+2 (\zeta^2+\zeta^{11})
   +2 (\zeta^6+\zeta^7)] \mathbf{B}_2+\\
  &-[(\zeta^6+\zeta^7)+2 (\zeta^9+\zeta^4)+2 (\zeta^5+\zeta^8)
   +2 (\zeta^2+\zeta^{11})] \mathbf{B}_6+\\
  &-[(\zeta^3+\zeta^{10})+2 (\zeta^2+\zeta^{11})+2 (\zeta^9+\zeta^4)
   +2 (\zeta+\zeta^{12})] \mathbf{B}_8+\\
  &-[(\zeta+\zeta^{12})+2 (\zeta^5+\zeta^8)+2 (\zeta^3+\zeta^{10})
   +2 (\zeta^9+\zeta^4)] \mathbf{B}_{11}+\\
  &-[(\zeta^9+\zeta^4)+2 (\zeta^6+\zeta^7)+2 (\zeta+\zeta^{12})
   +2 (\zeta^3+\zeta^{10})] \mathbf{B}_7+\\
  &+(\zeta^9+\zeta^4+\zeta^6+\zeta^7+\zeta^3+\zeta^{10})
   (\mathbf{B}_1^{(1)}+\mathbf{B}_1^{(2)})+\\
  &+(\zeta^3+\zeta^{10}+\zeta^2+\zeta^{11}+\zeta+\zeta^{12})
   (\mathbf{B}_3^{(1)}+\mathbf{B}_3^{(2)})+\\
  &+(\zeta+\zeta^{12}+\zeta^5+\zeta^8+\zeta^9+\zeta^4)
   (\mathbf{B}_9^{(1)}+\mathbf{B}_9^{(2)})+\\
  &+(\zeta^9+\zeta^4+\zeta^6+\zeta^7+\zeta^2+\zeta^{11})
   (-\mathbf{B}_{12}^{(1)}+\mathbf{B}_{12}^{(2)})+\\
  &+(\zeta^3+\zeta^{10}+\zeta^2+\zeta^{11}+\zeta^5+\zeta^8)
   (-\mathbf{B}_{10}^{(1)}+\mathbf{B}_{10}^{(2)})+\\
  &+(\zeta+\zeta^{12}+\zeta^5+\zeta^8+\zeta^6+\zeta^7)
   (-\mathbf{B}_4^{(1)}+\mathbf{B}_4^{(2)}).
\endaligned\eqno{(6.22)}$$
$$\aligned
  &13 S(\mathbf{B}_1^{(1)}+\mathbf{B}_1^{(2)})\\
 =&2 (4 \mathbf{B}_0^{(0)}-\mathbf{B}_0^{(1)}+\mathbf{B}_0^{(2)})+\\
  &+4 (\zeta^9+\zeta^4+\zeta^6+\zeta^7+\zeta^2+\zeta^{11}) \mathbf{B}_5+\\
  &+4 (\zeta^3+\zeta^{10}+\zeta^2+\zeta^{11}+\zeta^5+\zeta^8) \mathbf{B}_2+\\
  &+4 (\zeta+\zeta^{12}+\zeta^5+\zeta^8+\zeta^6+\zeta^7) \mathbf{B}_6+\\
  &+4 (\zeta^9+\zeta^4+\zeta^6+\zeta^7+\zeta^3+\zeta^{10}) \mathbf{B}_8+\\
  &+4 (\zeta^3+\zeta^{10}+\zeta^2+\zeta^{11}+\zeta+\zeta^{12}) \mathbf{B}_{11}+\\
  &+4 (\zeta+\zeta^{12}+\zeta^5+\zeta^8+\zeta^9+\zeta^4) \mathbf{B}_7]+\\
  &+[2 (1+\zeta^9+\zeta^4)+2 (\zeta+\zeta^{12})+(\zeta^2+\zeta^{11})
   (\mathbf{B}_1^{(1)}+\mathbf{B}_1^{(2)})+\\
  &+[2 (1+\zeta^3+\zeta^{10})+2 (\zeta^9+\zeta^4)+(\zeta^5+\zeta^8)]
   (\mathbf{B}_3^{(1)}+\mathbf{B}_3^{(2)})+\\
  &+[2 (1+\zeta+\zeta^{12})+2 (\zeta^3+\zeta^{10})+(\zeta^6+\zeta^7)]
   (\mathbf{B}_9^{(1)}+\mathbf{B}_9^{(2)})+\\
  &+[2 (1+\zeta^6+\zeta^7)+2 (\zeta^5+\zeta^8)+(\zeta^3+\zeta^{10})]
   (-\mathbf{B}_{12}^{(1)}+\mathbf{B}_{12}^{(2)})+\\
  &+[2 (1+\zeta^2+\zeta^{11})+2 (\zeta^6+\zeta^7)+(\zeta+\zeta^{12})]
   (-\mathbf{B}_{10}^{(1)}+\mathbf{B}_{10}^{(2)})+\\
  &+[2 (1+\zeta^5+\zeta^8)+2 (\zeta^2+\zeta^{11})+(\zeta^9+\zeta^4)]
   (-\mathbf{B}_4^{(1)}+\mathbf{B}_4^{(2)}).
\endaligned\eqno{(6.23)}$$
$$\aligned
  &13 S(-\mathbf{B}_{12}^{(1)}+\mathbf{B}_{12}^{(2)})\\
 =&2 (4 \mathbf{B}_0^{(0)}-\mathbf{B}_0^{(1)}+\mathbf{B}_0^{(2)})+\\
  &+4 (\zeta^9+\zeta^4+\zeta^6+\zeta^7+\zeta^3+\zeta^{10}) \mathbf{B}_5+\\
  &+4 (\zeta^3+\zeta^{10}+\zeta^2+\zeta^{11}+\zeta+\zeta^{12}) \mathbf{B}_2+\\
  &+4 (\zeta+\zeta^{12}+\zeta^5+\zeta^8+\zeta^9+\zeta^4) \mathbf{B}_6+\\
  &+4 (\zeta^9+\zeta^4+\zeta^6+\zeta^7+\zeta^2+\zeta^{11}) \mathbf{B}_8+\\
  &+4 (\zeta^3+\zeta^{10}+\zeta^2+\zeta^{11}+\zeta^5+\zeta^8) \mathbf{B}_{11}+\\
  &+4 (\zeta+\zeta^{12}+\zeta^5+\zeta^8+\zeta^6+\zeta^7) \mathbf{B}_7+\\
  &+[2 (1+\zeta^6+\zeta^7)+2 (\zeta^5+\zeta^8)+(\zeta^3+\zeta^{10})]
   (\mathbf{B}_1^{(1)}+\mathbf{B}_1^{(2)})+\\
  &+[2 (1+\zeta^2+\zeta^{11})+2 (\zeta^6+\zeta^7)+(\zeta+\zeta^{12})]
   (\mathbf{B}_3^{(1)}+\mathbf{B}_3^{(2)})+\\
  &+[2 (1+\zeta^5+\zeta^8)+2 (\zeta^2+\zeta^{11})+(\zeta^9+\zeta^4)]
   (\mathbf{B}_9^{(1)}+\mathbf{B}_9^{(2)})+\\
  &+[2 (1+\zeta^9+\zeta^4)+2 (\zeta+\zeta^{12})+(\zeta^2+\zeta^{11})]
   (-\mathbf{B}_{12}^{(1)}+\mathbf{B}_{12}^{(2)})+\\
  &+[2 (1+\zeta^3+\zeta^{10})+2 (\zeta^9+\zeta^4)+(\zeta^5+\zeta^8)]
   (-\mathbf{B}_{10}^{(1)}+\mathbf{B}_{10}^{(2)})+\\
  &+[2 (1+\zeta+\zeta^{12})+2 (\zeta^3+\zeta^{10})+(\zeta^6+\zeta^7)]
   (-\mathbf{B}_4^{(1)}+\mathbf{B}_4^{(2)}).
\endaligned\eqno{(6.24)}$$
$$T(\mathbf{B}_5)=\zeta^5 \mathbf{B}_5, \quad
  T(\mathbf{B}_8)=\zeta^8 \mathbf{B}_8.\eqno{(6.25)}$$
$$T(\mathbf{B}_1^{(1)}+\mathbf{B}_1^{(2)})
 =\zeta (\mathbf{B}_1^{(1)}+\mathbf{B}_1^{(2)}).\eqno{(6.26)}$$
$$T(-\mathbf{B}_{12}^{(1)}+\mathbf{B}_{12}^{(2)})
 =\zeta^{12} (-\mathbf{B}_{12}^{(1)}+\mathbf{B}_{12}^{(2)}).\eqno{(6.27)}$$
Let
$$\aligned
  V_{13}=\langle &4 \mathbf{B}_0^{(0)}-\mathbf{B}_0^{(1)}+\mathbf{B}_0^{(2)},
               \mathbf{B}_5, \mathbf{B}_2, \mathbf{B}_6, \mathbf{B}_8,
               \mathbf{B}_{11}, \mathbf{B}_7,\\
              &\mathbf{B}_1^{(1)}+\mathbf{B}_1^{(2)},
               \mathbf{B}_3^{(1)}+\mathbf{B}_3^{(2)},
               \mathbf{B}_9^{(1)}+\mathbf{B}_9^{(2)},\\
              &-\mathbf{B}_{12}^{(1)}+\mathbf{B}_{12}^{(2)},
               -\mathbf{B}_{10}^{(1)}+\mathbf{B}_{10}^{(2)},
               -\mathbf{B}_4^{(1)}+\mathbf{B}_4^{(2)} \rangle.
\endaligned\eqno{(6.28)}$$
Then $V_{13}$ is a thirteen-dimensional subspace which is stable under
the action of $G \cong \text{SL}(2, 13)$. It is easy to see that each
$v \in V$ can be represented by a linear combination of the basis of
$V_1 \oplus V_7 \oplus V_{13}$. This gives the following decomposition
of $V$:
$$V=V_1 \oplus V_7 \oplus V_{13},\eqno{(6.29)}$$
i.e.,
$$\mathbf{21}=\mathbf{1} \oplus \mathbf{7} \oplus \mathbf{13}.\eqno{(6.30)}$$
This completes the proof of Theorem 6.1.

\noindent
$\qquad \qquad \qquad \qquad \qquad \qquad \qquad \qquad \qquad
 \qquad \qquad \qquad \qquad \qquad \qquad \qquad \boxed{}$

\begin{center}
{\bf 6.2. A geometric construction of the curve $Y$ and the ring of
          invariant polynomials}
\end{center}

   Let $\mathfrak{a}_i$ be the ideals generated by the basis of $V_{6i-5}$ $(i=1, 2, 3)$
and $Y_i$ be the varieties corresponding to the ideals $\mathfrak{a}_i$. Theorem
6.1 implies that $\mathfrak{a}_i$ $(i=1, 2, 3)$ are $\text{SL}(2, 13)$-invariant
ideals and $Y_i$ $(i=1, 2, 3)$ are $\text{SL}(2, 13)$-invariant varieties. Let
$R=\mathbb{C}[z_1, z_2, z_3, z_4, z_5, z_6]$, we will study the three rings of
invariant polynomials $(R/\mathfrak{a}_i)^{\text{SL}(2, 13)}$ for $i=1, 2, 3$.
In particular, we construct the invariant homogeneous polynomials with minimal
degree in the invariant ideals $\mathfrak{a}_2$ and $\mathfrak{a}_3$, by which
we give a geometric construction of the curve $Y$.

\textbf{Theorem 6.2.} {\it The ideal $I(Y)$ can be decomposed as the sum of
ideals $\mathfrak{a}_i$ $(i=1, 2, 3)$. The curve $Y$ can be expressed as the
intersection of $Y_1$, $Y_2$ and $Y_3$:
$$I(Y)=\mathfrak{a}_1+\mathfrak{a}_2+\mathfrak{a}_3.\eqno{(6.31)}$$
$$Y=Y_1 \cap Y_2 \cap Y_3.\eqno{(6.32)}$$
In particular, the curve $Y_2$ is in a $\text{SL}(2, 13)$-invariant sextic
four-fold $W_2 \subset \mathbb{P}^5$, i.e., a Calabi-Yau four-fold, given by an
invariant homogeneous polynomial $\Psi_{1, 1}=0$, where
$$\aligned
  \Psi_{1, 1} =&2 \mathbf{A}_0 (\mathbf{B}_0^{(1)}+\mathbf{B}_0^{(2)})
                +\mathbf{A}_1 (\mathbf{B}_{12}^{(1)}+3 \mathbf{B}_{12}^{(2)})+\\
               &+\mathbf{A}_4 (\mathbf{B}_{10}^{(1)}+3 \mathbf{B}_{10}^{(2)})
                +\mathbf{A}_3 (\mathbf{B}_4^{(1)}+3 \mathbf{B}_4^{(2)})+\\
               &+\mathbf{A}_5 (\mathbf{B}_1^{(1)}-3 \mathbf{B}_1^{(2)})
                +\mathbf{A}_6 (\mathbf{B}_3^{(1)}-3 \mathbf{B}_3^{(2)})+\\
               &+\mathbf{A}_2 (\mathbf{B}_9^{(1)}-3 \mathbf{B}_9^{(2)})
                \in \mathfrak{a}_2.
\endaligned\eqno{(6.33)}$$
The curve $Y_3$ is in a $\text{SL}(2, 13)$-invariant octic four-fold $W_3
\subset \mathbb{P}^5$, i.e., a general type four-fold, given by an invariant
homogeneous polynomial $\Omega_{1, 1}=0$, where
$$\aligned
  &\Omega_{1, 1}=(4 \mathbf{B}_0^{(0)}-\mathbf{B}_0^{(1)}+\mathbf{B}_0^{(2)})
                 (6 \mathbf{A}_0^2-\mathbf{A}_1 \mathbf{A}_5-\mathbf{A}_2 \mathbf{A}_3
                  -\mathbf{A}_4 \mathbf{A}_6)+\\
                &+7 \mathbf{B}_5 (\mathbf{A}_2^2+2 \mathbf{A}_3 \mathbf{A}_5)
                 +7 \mathbf{B}_2 (\mathbf{A}_5^2+2 \mathbf{A}_1 \mathbf{A}_6)
                 +7 \mathbf{B}_6 (\mathbf{A}_6^2+2 \mathbf{A}_4 \mathbf{A}_2)+\\
                &+7 \mathbf{B}_8 (\mathbf{A}_3^2+2 \mathbf{A}_1 \mathbf{A}_2)+
                 +7 \mathbf{B}_{11} (\mathbf{A}_1^2+2 \mathbf{A}_4 \mathbf{A}_5)
                 +7 \mathbf{B}_7 (\mathbf{A}_4^2+2 \mathbf{A}_3 \mathbf{A}_6)+\\
                &+7 (\mathbf{B}_1^{(1)}+\mathbf{B}_1^{(2)}) (\mathbf{A}_0
                 \mathbf{A}_5+\mathbf{A}_3 \mathbf{A}_4)
                 +7 (\mathbf{B}_3^{(1)}+\mathbf{B}_3^{(2)}) (\mathbf{A}_0
                 \mathbf{A}_6+\mathbf{A}_1 \mathbf{A}_3)+\\
                &+7 (\mathbf{B}_9^{(1)}+\mathbf{B}_9^{(2)}) (\mathbf{A}_0
                 \mathbf{A}_2+\mathbf{A}_1 \mathbf{A}_4)
                 +7 (-\mathbf{B}_{12}^{(1)}+\mathbf{B}_{12}^{(2)}) (\mathbf{A}_0
                 \mathbf{A}_1+\mathbf{A}_2 \mathbf{A}_6)+\\
                &+7 (-\mathbf{B}_{10}^{(1)}+\mathbf{B}_{10}^{(2)}) (\mathbf{A}_0
                 \mathbf{A}_4+\mathbf{A}_2 \mathbf{A}_5)
                 +7 (-\mathbf{B}_4^{(1)}+\mathbf{B}_4^{(2)}) (\mathbf{A}_0
                 \mathbf{A}_3+\mathbf{A}_5 \mathbf{A}_6)\\
                &\in \mathfrak{a}_3.
\endaligned\eqno{(6.34)}$$}

{\it Proof}. (6.31) and (6.32) are obtained by the construction of the
ideals $\mathfrak{a}_1$, $\mathfrak{a}_2$ and $\mathfrak{a}_3$. In order
to prove (6.33), put
$$\phi_{\infty}=\sqrt{13} (\mathbf{B}_0^{(1)}+\mathbf{B}_0^{(2)}), \quad
  \phi_{\nu}=\phi_{\infty}(ST^{\nu}(z_1, z_2, z_3, z_4, z_5, z_6))\eqno{(6.35)}$$
for $\nu=0, 1, \cdots, 12$. By (6.12), we have
$$\aligned
  \phi_{\nu}=&(\mathbf{B}_0^{(1)}+\mathbf{B}_0^{(2)})
   +\zeta^{\nu} (\mathbf{B}_1^{(1)}-3 \mathbf{B}_1^{(2)})
   +\zeta^{3 \nu} (\mathbf{B}_3^{(1)}-3 \mathbf{B}_3^{(2)})\\
  &+\zeta^{9 \nu} (\mathbf{B}_9^{(1)}-3 \mathbf{B}_9^{(2)})
   +\zeta^{12 \nu} (\mathbf{B}_{12}^{(1)}+3 \mathbf{B}_{12}^{(2)})\\
  &+\zeta^{10 \nu} (\mathbf{B}_{10}^{(1)}+3 \mathbf{B}_{10}^{(2)})
   +\zeta^{4 \nu} (\mathbf{B}_4^{(1)}+3 \mathbf{B}_4^{(2)}).
\endaligned\eqno{(6.36)}$$
Together with (3.6) and (3.8), we have that $\varphi_0 \phi_0$ is invariant
under the action of $\langle H, T \rangle$. Put
$$\Psi_{m, n}=\varphi_0^m \phi_0^n+\varphi_1^m \phi_1^n+\cdots+\varphi_{12}^m
              \phi_{12}^n+\varphi_{\infty}^m \phi_{\infty}^n\eqno{(6.37)}$$
with $m+n$ is even. A straightforward computation gives (up to a constant
$13$) the expression of $\Psi_{1, 1}$ in (6.33). Similarly, put
$$\sigma_{\infty}=13 (4 \mathbf{B}_0^{(0)}-\mathbf{B}_0^{(1)}+\mathbf{B}_0^{(2)}),
  \quad \sigma_{\nu}=\sigma_{\infty}(ST^{\nu}(z_1, z_2, z_3, z_4, z_5, z_6))\eqno{(6.38)}$$
for $\nu=0, 1, \cdots, 12$. By (6.19), we have
$$\aligned
   \sigma_{\nu}
 =&-(4 \mathbf{B}_0^{(0)}-\mathbf{B}_0^{(1)}+\mathbf{B}_0^{(2)})\\
  &+14 \zeta^{5 \nu} \mathbf{B}_5
   +14 \zeta^{2 \nu} \mathbf{B}_2
   +14 \zeta^{6 \nu} \mathbf{B}_6
   +14 \zeta^{8 \nu} \mathbf{B}_8
   +14 \zeta^{11 \nu} \mathbf{B}_{11}
   +14 \zeta^{7 \nu} \mathbf{B}_7\\
  &+7 \zeta^{\nu} (\mathbf{B}_1^{(1)}+\mathbf{B}_1^{(2)})
   +7 \zeta^{3 \nu} (\mathbf{B}_3^{(1)}+\mathbf{B}_3^{(2)})
   +7 \zeta^{9 \nu} (\mathbf{B}_9^{(1)}+\mathbf{B}_9^{(2)})\\
  &+7 \zeta^{12 \nu} (-\mathbf{B}_{12}^{(1)}+\mathbf{B}_{12}^{(2)})
   +7 \zeta^{10 \nu} (-\mathbf{B}_{10}^{(1)}+\mathbf{B}_{10}^{(2)})
   +7 \zeta^{4 \nu} (-\mathbf{B}_4^{(1)}+\mathbf{B}_4^{(2)}).
\endaligned\eqno{(6.39)}$$
Note that $4 \mathbf{B}_0^{(0)}-\mathbf{B}_0^{(1)}+\mathbf{B}_0^{(2)}$
is invariant under the action of $\langle H, T \rangle$. Put
$$\Omega_{m, n}=w_0^m \sigma_0^n+w_1^m \sigma_1^n+\cdots+w_{12}^m
                \sigma_{12}^n+w_{\infty}^m \sigma_{\infty}^n.\eqno{(6.40)}$$
A straightforward computation gives (up to a constant $26$) the expression
of $\Omega_{1, 1}$ in (6.34).

\noindent
$\qquad \qquad \qquad \qquad \qquad \qquad \qquad \qquad \qquad
 \qquad \qquad \qquad \qquad \qquad \qquad \qquad \boxed{}$

\begin{center}
{\bf 6.3. A non-standard geometric realization of the degenerate principal
         series for $\text{SL}(2, 13)$ and a geometric realization of the
         Steinberg representation for $\text{SL}(2, 13)$}
\end{center}

  We have constructed four kinds of representations of $\text{SL}(2, 13)$
corresponding to the following bases:

(1) $(z_1, \ldots, z_6)$: discrete series of dimension six

(2) $(\mathbf{A}_0, \mathbf{A}_1, \ldots, \mathbf{A}_6)$: principal series
    of dimension seven

(3) $(\mathbf{D}_0, \mathbf{D}_1, \ldots, \mathbf{D}_{12}, \mathbf{D}_{\infty})$:
    principal series of dimension fourteen

(4) $(\mathbf{B}_0^{(0)}, \mathbf{B}_0^{(1)}, \mathbf{B}_0^{(2)}, \ldots, \mathbf{B}_7)$:
    reducible representation of dimension $21$, which can be decomposed as the
    direct sum of a one-dimensional representation, a seven-dimensional representation
    (principal series) and a thirteen-dimensional representation (principal series,
    Steinberg representation).

  Now, we prove that the above two $7$-dimensional representations are equivalent.
In \cite{Y1}, we obtain a seven-dimensional representation (degenerate principal
series)  of the group $\text{SL}(2, 13) \cong \langle \widetilde{S}, \widetilde{T}
\rangle$ which induces from the actions of $S$ and $T$ on the basis $(\mathbf{A}_0,
\mathbf{A}_1, \mathbf{A}_2, \mathbf{A}_3, \mathbf{A}_4, \mathbf{A}_5, \mathbf{A}_6)$.
Here,
$$\widetilde{S}=\frac{1}{\sqrt{13}} \left(\begin{matrix}
  1 & 1                  & 1               & 1
    & 1               & 1                  & 1\\
  2 & \zeta^2+\zeta^{11} & \zeta^9+\zeta^4 & \zeta^6+\zeta^7
    & \zeta^5+\zeta^8 & \zeta^3+\zeta^{10} & \zeta+\zeta^{12}\\
  2 & \zeta^9+\zeta^4 & \zeta^5+\zeta^8 & \zeta+\zeta^{12}
    & \zeta^3+\zeta^{10} & \zeta^6+\zeta^7 & \zeta^2+\zeta^{11}\\
  2 & \zeta^6+\zeta^7 & \zeta+\zeta^{12} & \zeta^5+\zeta^8
    & \zeta^2+\zeta^{11} & \zeta^9+\zeta^4 & \zeta^3+\zeta^{10}\\
  2 & \zeta^5+\zeta^8 & \zeta^3+\zeta^{10} & \zeta^2+\zeta^{11}
    & \zeta^6+\zeta^7 & \zeta+\zeta^{12} & \zeta^9+\zeta^4\\
  2 & \zeta^3+\zeta^{10} & \zeta^6+\zeta^7 & \zeta^9+\zeta^4
    & \zeta+\zeta^{12} & \zeta^2+\zeta^{11} & \zeta^5+\zeta^8\\
  2 & \zeta+\zeta^{12} & \zeta^2+\zeta^{11} & \zeta^3+\zeta^{10}
    & \zeta^9+\zeta^4 & \zeta^5+\zeta^8 & \zeta^6+\zeta^7
\end{matrix}\right),\eqno{(6.41)}$$
and
$$\widetilde{T}=\text{diag}(1, \zeta, \zeta^4, \zeta^9, \zeta^3,
                \zeta^{12}, \zeta^{10}).\eqno{(6.42)}$$
We have
$$\text{Tr}(\widetilde{S})=-1, \quad
  \text{Tr}(\widetilde{T})=\frac{1+\sqrt{13}}{2}.$$
On the other hand, for the other $7$-dimensional representation associated
with the ideal $\mathfrak{a}_2$, the trace of the action of $S$ and $T$ on
the basis corresponding to the ideal $\mathfrak{a}_2$ is:
$$\aligned
  \text{Tr}(S)=&\frac{1}{13}[\sqrt{13}
   +(-2(1+\zeta^9+\zeta^4)+2(\zeta+\zeta^{12})+(\zeta^2+\zeta^{11}))+\\
  &+(-2(1+\zeta+\zeta^{12})+2(\zeta^3+\zeta^{10})+(\zeta^6+\zeta^7))+\\
  &+(-2(1+\zeta^3+\zeta^{10})+2(\zeta^9+\zeta^4)+(\zeta^5+\zeta^8))+\\
  &+(-2(1+\zeta^9+\zeta^4)+2(\zeta+\zeta^{12})+(\zeta^2+\zeta^{11}))+\\
  &+(-2(1+\zeta+\zeta^{12})+2(\zeta^3+\zeta^{10})+(\zeta^6+\zeta^7))+\\
  &+(-2(1+\zeta^3+\zeta^{10})+2(\zeta^9+\zeta^4)+(\zeta^5+\zeta^8))]\\
 =&-1.
\endaligned$$
$$\text{Tr}(T)=1+\zeta+\zeta^3+\zeta^9+\zeta^{12}+\zeta^{10}+\zeta^4
              =\frac{1+\sqrt{13}}{2}.$$
This implies that the above two $7$-dimensional representations are
equivalent.

  By (3.8), we have
$$\sqrt{13} ST^{\nu}(\mathbf{A}_0)=\mathbf{A}_0+\zeta^{\nu} \mathbf{A}_1
  +\zeta^{3 \nu} \mathbf{A}_4+\zeta^{9 \nu} \mathbf{A}_3+\zeta^{12 \nu}
  \mathbf{A}_5+\zeta^{10 \nu} \mathbf{A}_6+\zeta^{4 \nu} \mathbf{A}_2.\eqno{(6.43)}$$
On the other hand, by (6.12), we have
$$\aligned
  &\sqrt{13} ST^{\nu}(\mathbf{B}_0^{(1)}+\mathbf{B}_0^{(2)})\\
 =&(\mathbf{B}_0^{(1)}+\mathbf{B}_0^{(2)})
  +\zeta^{\nu} (\mathbf{B}_1^{(1)}-3 \mathbf{B}_1^{(2)})
  +\zeta^{3 \nu} (\mathbf{B}_3^{(1)}-3 \mathbf{B}_3^{(2)})
  +\zeta^{9 \nu} (\mathbf{B}_9^{(1)}-3 \mathbf{B}_9^{(2)})\\
  +&\zeta^{12 \nu} (\mathbf{B}_{12}^{(1)}+3 \mathbf{B}_{12}^{(2)})+
   \zeta^{10 \nu} (\mathbf{B}_{10}^{(1)}+3 \mathbf{B}_{10}^{(2)})+
   \zeta^{4 \nu} (\mathbf{B}_4^{(1)}+3 \mathbf{B}_4^{(2)}).
\endaligned\eqno{(6.44)}$$
Note that both $\mathbf{A}_0^2$ and $(\mathbf{B}_0^{(1)}+\mathbf{B}_0^{(2)})^2$
are invariant under the action of the maximal subgroup $\langle H, T \rangle$.

  Let us recall a result which Jacobi had established as early as $1828$
in his ``Notices sur les fonctions elliptiques'' (see \cite{Ja1828} and
\cite{K}). Jacobi there considered, instead of the modular equation, the
so-called multiplier-equation, together with other equations equivalent
to it, and found that their $(n+1)$ roots are composed in a simple manner
of $\frac{n+1}{2}$ elements, with the help of merely numerical irrationalities.
Namely, if we denote these elements by $\mathbf{A}_0$, $\mathbf{A}_1$,
$\ldots$, $\mathbf{A}_{\frac{n-1}{2}}$, and further, for the roots
$z$ of the equation under consideration, apply the indices
employed by Galois, we have, with appropriate determination of the
square root occurring on the left-hand side:
$$\left\{\aligned
  \sqrt{z_{\infty}} &=\sqrt{(-1)^{\frac{n-1}{2}} \cdot n} \cdot \mathbf{A}_0,\\
  \sqrt{z_{\nu}} &=\mathbf{A}_0+\epsilon^{\nu} \mathbf{A}_1+\epsilon^{4 \nu} \mathbf{A}_2
  +\cdots+\epsilon^{(\frac{n-1}{2})^2 \nu} \mathbf{A}_{\frac{n-1}{2}}
\endaligned\right.\eqno{(6.45)}$$
for $\nu=0, 1, \cdots, n-1$ and $\epsilon=e^{\frac{2 \pi i}{n}}$.
Jacobi had himself emphasized the special significance of his
result by adding to his short communication: ``C'est un
th\'{e}or\`{e}me des plus importants dans la th\'{e}orie
alg\'{e}brique de la transformation et de la division des
fonctions elliptiques.''

  Now, (6.43) gives a standard geometric realization of (6.45), (6.44)
gives a non-standard geometric realization of (6.45). This prove the
following:

\textbf{Theorem 6.3.} {\it The degenerate principal series of
$\text{SL}(2, 13)$ has two distinct geometric realizations: The standard
realization is given by the basis $(\mathbf{A}_0, \mathbf{A}_1$, $\mathbf{A}_2$,
$\mathbf{A}_3$, $\mathbf{A}_4$, $\mathbf{A}_5$, $\mathbf{A}_6)$ corresponding
to the Klein $\mathbf{A}$-curve. The non-standard realization is given by
the basis $(\mathbf{B}_0^{(1)}+\mathbf{B}_0^{(2)}$, $\mathbf{B}_1^{(1)}-3
\mathbf{B}_1^{(2)}$, $\mathbf{B}_3^{(1)}-3 \mathbf{B}_3^{(2)}$,
$\mathbf{B}_9^{(1)}-3 \mathbf{B}_9^{(2)}$, $\mathbf{B}_{12}^{(1)}+3
\mathbf{B}_{12}^{(2)}$, $\mathbf{B}_{10}^{(1)}+3 \mathbf{B}_{10}^{(2)}$,
$\mathbf{B}_4^{(1)}+3 \mathbf{B}_4^{(2)})$ corresponding to the curve $Y_2$
lying over $Y$. In particular,
$$\sqrt{13} ST^{\nu}(\mathbf{A}_0)=\mathbf{A}_0+\zeta^{\nu} \mathbf{A}_1
  +\zeta^{3 \nu} \mathbf{A}_4+\zeta^{9 \nu} \mathbf{A}_3+\zeta^{12 \nu}
  \mathbf{A}_5+\zeta^{10 \nu} \mathbf{A}_6+\zeta^{4 \nu} \mathbf{A}_2.
  \eqno{(6.46)}$$
$$\aligned
  &\sqrt{13} ST^{\nu}(\mathbf{B}_0^{(1)}+\mathbf{B}_0^{(2)})\\
 =&(\mathbf{B}_0^{(1)}+\mathbf{B}_0^{(2)})
  +\zeta^{\nu} (\mathbf{B}_1^{(1)}-3 \mathbf{B}_1^{(2)})
  +\zeta^{3 \nu} (\mathbf{B}_3^{(1)}-3 \mathbf{B}_3^{(2)})
  +\zeta^{9 \nu} (\mathbf{B}_9^{(1)}-3 \mathbf{B}_9^{(2)})\\
  +&\zeta^{12 \nu} (\mathbf{B}_{12}^{(1)}+3 \mathbf{B}_{12}^{(2)})+
   \zeta^{10 \nu} (\mathbf{B}_{10}^{(1)}+3 \mathbf{B}_{10}^{(2)})+
   \zeta^{4 \nu} (\mathbf{B}_4^{(1)}+3 \mathbf{B}_4^{(2)}).
\endaligned\eqno{(6.47)}$$}

  Similarly, the curve $Y_3$ over $Y$ gives a geometric realization of the
Steinberg representation for $\text{SL}(2, 13)$. For the $13$-dimensional
representation associated with the ideal $\mathfrak{a}_3$, the trace of
the action of $S$ and $T$ on the basis corresponding to the ideal
$\mathfrak{a}_3$ is:
$$\aligned
  \text{Tr}(S)=&\frac{1}{13}[-1
   -5(\zeta+\zeta^{12}+\zeta^3+\zeta^{10}+\zeta^9+\zeta^4)\\
  &-2(\zeta^5+\zeta^8+\zeta^2+\zeta^{11}+\zeta^6+\zeta^7)+\\
  &-5(\zeta+\zeta^{12}+\zeta^3+\zeta^{10}+\zeta^9+\zeta^4)+\\
  &-2(\zeta^5+\zeta^8+\zeta^2+\zeta^{11}+\zeta^6+\zeta^7)+\\
  &+6+4(\zeta+\zeta^{12}+\zeta^3+\zeta^{10}+\zeta^9+\zeta^4)+\\
  &+(\zeta^5+\zeta^8+\zeta^2+\zeta^{11}+\zeta^6+\zeta^7)+\\
  &+6+4(\zeta+\zeta^{12}+\zeta^3+\zeta^{10}+\zeta^9+\zeta^4)+\\
  &+(\zeta^5+\zeta^8+\zeta^2+\zeta^{11}+\zeta^6+\zeta^7)\\
 =&1.
\endaligned$$
$$\text{Tr}(T)=1+\zeta^5+\zeta^2+\zeta^6+\zeta^8+\zeta^{11}+\zeta^7+
                \zeta+\zeta^3+\zeta^9+\zeta^{12}+\zeta^{10}+\zeta^4
              =0.$$
In particular, (6.39) gives that
$$\aligned
  &13S T^{\nu}(4 \mathbf{B}_0^{(0)}-\mathbf{B}_0^{(1)}+\mathbf{B}_0^{(2)})\\
 =&-(4 \mathbf{B}_0^{(0)}-\mathbf{B}_0^{(1)}+\mathbf{B}_0^{(2)})\\
  &+14 \zeta^{5 \nu} \mathbf{B}_5
   +14 \zeta^{2 \nu} \mathbf{B}_2
   +14 \zeta^{6 \nu} \mathbf{B}_6
   +14 \zeta^{8 \nu} \mathbf{B}_8
   +14 \zeta^{11 \nu} \mathbf{B}_{11}
   +14 \zeta^{7 \nu} \mathbf{B}_7\\
  &+7 \zeta^{\nu} (\mathbf{B}_1^{(1)}+\mathbf{B}_1^{(2)})
   +7 \zeta^{3 \nu} (\mathbf{B}_3^{(1)}+\mathbf{B}_3^{(2)})
   +7 \zeta^{9 \nu} (\mathbf{B}_9^{(1)}+\mathbf{B}_9^{(2)})\\
  &+7 \zeta^{12 \nu} (-\mathbf{B}_{12}^{(1)}+\mathbf{B}_{12}^{(2)})
   +7 \zeta^{10 \nu} (-\mathbf{B}_{10}^{(1)}+\mathbf{B}_{10}^{(2)})
   +7 \zeta^{4 \nu} (-\mathbf{B}_4^{(1)}+\mathbf{B}_4^{(2)}).
\endaligned\eqno{(6.48)}$$

\begin{center}
{\large\bf 7. Modularity for an invariant ideal and an explicit
              uniformization of algebraic space curves of higher genus}
\end{center}

  Recall that the theta functions with characteristic
$\begin{bmatrix} \epsilon\\ \epsilon^{\prime} \end{bmatrix} \in
\mathbb{R}^2$ is defined by the following series which converges
uniformly and absolutely on compact subsets of $\mathbb{C} \times
\mathbb{H}$ (see \cite{FK}, p. 73):
$$\theta \begin{bmatrix} \epsilon\\ \epsilon^{\prime} \end{bmatrix}(z, \tau)
 =\sum_{n \in \mathbb{Z}} \exp \left\{2 \pi i \left[\frac{1}{2}
  \left(n+\frac{\epsilon}{2}\right)^2 \tau+\left(n+\frac{\epsilon}{2}\right)
  \left(z+\frac{\epsilon^{\prime}}{2}\right)\right]\right\}.$$
The modified theta constants (see \cite{FK}, p. 215)
$\varphi_l(\tau):=\theta [\chi_l](0, k \tau)$,
where the characteristic $\chi_l=\begin{bmatrix} \frac{2l+1}{k}\\ 1
\end{bmatrix}$, $l=0, \ldots, \frac{k-3}{2}$, for odd $k$ and
$\chi_l=\begin{bmatrix} \frac{2l}{k}\\ 0 \end{bmatrix}$, $l=0,
\ldots, \frac{k}{2}$, for even $k$. Let $V(k)$ be the finite
dimensional vector space of holomorphic functions on $\mathbb{H}$
spanned by these modified theta constants $\varphi_l(\tau)$. We
have the following:

\textbf{Proposition 7.1.} (see \cite{FK}, p.236, Proposition 7.2 and
p.239, Proposition 7.9). {\it For each odd integer $k \geq 5$, the map
$$\Phi: \tau \mapsto (\varphi_0(\tau), \varphi_1(\tau), \ldots,
  \varphi_{\frac{k-5}{2}}(\tau), \varphi_{\frac{k-3}{2}}(\tau))$$
from $\mathbb{H} \cup \mathbb{Q} \cup \{ \infty \}$ to
$\mathbb{C}^{\frac{k-1}{2}}$, defines a holomorphic mapping from
$\overline{\mathbb{H}/\Gamma(k)}$ into $\mathbb{C} \mathbb{P}^{\frac{k-3}{2}}$.
If $k$ is prime, the distinguished punctures on the surface $\mathbb{H}/\Gamma(k)$
are sent (injectively) onto the coordinate vectors in projective space
$\mathbb{C} \mathbb{P}^{\frac{k-3}{2}}$. Further, in this case, the
map $\Phi$ restricted to the punctures is injective and of maximal rank
at each distinguished puncture. For all $k \in \mathbb{Z}^{+}$, the map
$\Phi$ has maximal rank everywhere.}

  In fact, the differential $d \Phi$ of the map $\Phi$ is shown to be
nonsingular everywhere, i.e., for all $x \in \mathbb{H} \cup \mathbb{Q}
\cup \{ \infty \}$, we can find functions $f$ and $g$ in $V(k)$ that
are regular at $x$ with $\frac{g}{f}$ having a simple zero at $x$ (see
\cite{FK}, p.239).

\noindent Remark. $\Phi(\overline{\mathbb{H}/\Gamma(k)})$ is a curve of
degree $\frac{(k^2-1)(k-3)}{48}$ in $\mathbb{CP}^{\frac{k-3}{2}}$ (see
\cite{FK}, p.237, Remark 7.4).

  In our case, the map
$\Phi: \tau \mapsto (\varphi_0(\tau), \varphi_1(\tau), \varphi_2(\tau),
 \varphi_3(\tau), \varphi_4(\tau), \varphi_5(\tau))$
gives a holomorphic mapping from the modular curve
$X(13)=\overline{\mathbb{H}/\Gamma(13)}$ into $\mathbb{C} \mathbb{P}^5$,
which corresponds to our six-dimensional representation, i.e., up to
the constants, $z_1, \ldots, z_6$ are just modular forms
$\varphi_0(\tau), \ldots, \varphi_5(\tau)$. Let
$$\left\{\aligned
  a_1(z) &:=e^{-\frac{11 \pi i}{26}} \theta
            \begin{bmatrix} \frac{11}{13}\\ 1 \end{bmatrix}(0, 13z)
           =q^{\frac{121}{104}} \sum_{n \in \mathbb{Z}} (-1)^n
            q^{\frac{1}{2}(13n^2+11n)},\\
  a_2(z) &:=e^{-\frac{7 \pi i}{26}} \theta
            \begin{bmatrix} \frac{7}{13}\\ 1 \end{bmatrix}(0, 13z)
           =q^{\frac{49}{104}} \sum_{n \in \mathbb{Z}} (-1)^n
            q^{\frac{1}{2}(13n^2+7n)},\\
  a_3(z) &:=e^{-\frac{5 \pi i}{26}} \theta
            \begin{bmatrix} \frac{5}{13}\\ 1 \end{bmatrix}(0, 13z)
           =q^{\frac{25}{104}} \sum_{n \in \mathbb{Z}} (-1)^n
            q^{\frac{1}{2}(13n^2+5n)},\\
  a_4(z) &:=-e^{-\frac{3 \pi i}{26}} \theta
            \begin{bmatrix} \frac{3}{13}\\ 1 \end{bmatrix}(0, 13z)
           =-q^{\frac{9}{104}} \sum_{n \in \mathbb{Z}} (-1)^n
            q^{\frac{1}{2}(13n^2+3n)},\\
  a_5(z) &:=e^{-\frac{9 \pi i}{26}} \theta
            \begin{bmatrix} \frac{9}{13}\\ 1 \end{bmatrix}(0, 13z)
           =q^{\frac{81}{104}} \sum_{n \in \mathbb{Z}} (-1)^n
            q^{\frac{1}{2}(13n^2+9n)},\\
  a_6(z) &:=e^{-\frac{\pi i}{26}} \theta
            \begin{bmatrix} \frac{1}{13}\\ 1 \end{bmatrix}(0, 13z)
           =q^{\frac{1}{104}} \sum_{n \in \mathbb{Z}} (-1)^n
            q^{\frac{1}{2}(13n^2+n)}
\endaligned\right.\eqno{(7.1)}$$
be the theta constants of order $13$ and
$$\mathbf{A}(z):=(a_1(z), a_2(z), a_3(z), a_4(z), a_5(z), a_6(z))^{T}.$$
The significance of our six dimensional representation of
$\text{SL}(2, 13)$ comes from the following:

\textbf{Proposition 7.2} (see \cite{Y2}, Proposition 2.5). {\it
If $z \in \mathbb{H}$, then the following relations hold:
$$\mathbf{A}(z+1)=e^{-\frac{3 \pi i}{4}} T \mathbf{A}(z), \quad
  \mathbf{A}\left(-\frac{1}{z}\right)=e^{\frac{\pi i}{4}} \sqrt{z}
  S \mathbf{A}(z),\eqno{(7.2)}$$
where $S$ and $T$ are given in (3.1) and (3.2), and
$0<\text{arg} \sqrt{z} \leq \pi/2$.}

  Note that it is (7.2) that gives an explicit realization of the
isomorphism between the unique sub-representation of parabolic
modular forms of weight $\frac{1}{2}$ on $X(13)$ and the six-dimensional
complex representation of $\text{SL}(2, 13)$ generated by $S$ and $T$.

  Recall that the principal congruence subgroup of level $13$ is the
normal subgroup $\Gamma(13)$ of $\Gamma=\text{SL}(2, \mathbb{Z})$
defined by the exact sequence
$1 \rightarrow \Gamma(13) \rightarrow \Gamma(1) \stackrel{f}{\rightarrow}
 \text{SL}(2, 13) \rightarrow 1$
where $f(\gamma) \equiv \gamma$ (mod $13$) for $\gamma \in \Gamma=\Gamma(1)$.
There is a representation $\rho: \Gamma \rightarrow \text{GL}(6, \mathbb{C})$
with kernel $\Gamma(13)$ defined as follows: if
$t=\begin{pmatrix} 1 & 1\\ 0 & 1 \end{pmatrix}$ and
$s=\begin{pmatrix} 0 & -1\\ 1 & 0 \end{pmatrix}$,
then $\rho(t)=T$ and $\rho(s)=S$. To see that such a representation exists,
note that $\text{PSL}(2, \mathbb{Z}):=\Gamma/\{ \pm I \}$ is defined by
the presentation $\langle s, t; s^2=(st)^3=1 \rangle$ satisfied by $s$ and
$t$ and we have proved that $S$ and $T$ satisfy these relations (in the
projective coordinates). Moreover, we have proved that $\text{PSL}(2, 13)$
is defined by the presentation $\langle S, T; S^2=T^{13}=(ST)^3=1 \rangle$.
Let $p=s t^{-1} s$ and $q=s t^3$. Then
$$h:=q^5 p^2 \cdot p^2 q^6 p^8 \cdot q^5 p^2 \cdot p^3 q
    =\begin{pmatrix}
     4, 428, 249 & -10, 547, 030\\
    -11, 594, 791 & 27, 616, 019
     \end{pmatrix}$$
satisfies that $\rho(h)=H$. The off-diagonal elements of the
matrix $h$, which corresponds to $H$, are congruent to $0$
mod $13$. The connection to $\Gamma_0(13)$ should be obvious.

  It is well-known that theta functions form a central object of study
both in complex analysis (studying Riemann surfaces) and in algebraic
geometry (studying abelian varieties). In particular, there exist some
quartic relations between theta constants, i.e., linear dependence
between products of four theta constants.

\textbf{Proposition 7.3.} (see \cite{FK}, p. 256). {\it Let
$j \in \mathbb{Z}$, $0<j \leq \frac{k-1}{2}$. For $i=1, 2, 3$,
choose odd integers $\alpha_i$ with $|\alpha_i| \leq k-2$,
$|\alpha_i| \neq k-2j$ and $\sum_{i=1}^{3} \alpha_i \equiv 0$}
(mod $k$). {\it The following quartic relation holds among theta
constants:
$$\aligned
 &\theta \begin{bmatrix} \frac{k-4 j}{k}\\ 1 \end{bmatrix} \prod_{i=1, 2, 3}
  \theta \begin{bmatrix} \frac{\alpha_i}{k}\\ 1 \end{bmatrix}+e^{\pi i \frac{k-4j}{k}}
  \theta \begin{bmatrix} \frac{k-2j}{k}\\ 1 \end{bmatrix} \prod_{i=1, 2, 3}
  \theta \begin{bmatrix} \frac{\alpha_i+2j}{k}\\ 1 \end{bmatrix}+\\
 &+e^{\pi i \frac{-k+2j}{k}} \theta \begin{bmatrix} \frac{k-2j}{k}\\ 1 \end{bmatrix}
  \prod_{i=1, 2, 3} \theta \begin{bmatrix} \frac{\alpha_i-2j}{k}\\ 1 \end{bmatrix}=0.
\endaligned\eqno{(7.3)}$$}

\textbf{Problem 7.4.} (see \cite{FK}, p. 256, Problem 10.2).  For a fixed odd
integer $k$, the above identity (7.3) is determined by the quadruple $(j; \alpha_1,
\alpha_2, \alpha_3)$. It is obvious that both $(j; \alpha_1, \alpha_2, \alpha_3)$
and $(j; -\alpha_1, -\alpha_2, -\alpha_3)$ always determine the same identity.
There are other, less clear, relations among the identities. It is of interest
to determine minimal generators for the ideal of identities and whether the
identities determine the curve $\Phi(\overline{\mathbb{H}/\Gamma(k)}) \subset
\mathbb{P}^{\frac{k-3}{2}}$.

  Now, we use the standard notation $(a)_{\infty}:=(a; q):=\prod_{k=0}^{\infty}
(1-aq^k)$ for $a \in \mathbb{C}^{\times}$ and $|q|<1$ so that $\eta(z)=q^{1/24}(q; q)$.
By the Jacobi triple product identity, we have
$$\theta \left[\begin{array}{c} \frac{l}{k}\\ 1 \end{array}\right](0, kz)
 =\exp \left(\frac{\pi i l}{2k}\right) q^{\frac{l^2}{8k}}
  (q^{\frac{k-l}{2}}; q^k)(q^{\frac{k+l}{2}}; q^k)(q^k; q^k),
  \eqno{(7.4)}$$
where $k$ is odd and $l=1, 3, 5, \cdots, k-2$. Hence,
$$\left\{\begin{array}{rl}
  a_1(z) &=q^{\frac{121}{104}} (q; q^{13})(q^{12}; q^{13})(q^{13}; q^{13}),\\
  a_2(z) &=q^{\frac{49}{104}} (q^3; q^{13})(q^{10}; q^{13})(q^{13}; q^{13}),\\
  a_3(z) &=q^{\frac{25}{104}} (q^9; q^{13})(q^4; q^{13})(q^{13}; q^{13}),\\
  a_4(z) &=-q^{\frac{9}{104}} (q^5; q^{13})(q^8; q^{13})(q^{13}; q^{13}),\\
  a_5(z) &=q^{\frac{81}{104}} (q^2; q^{13})(q^{11}; q^{13})(q^{13}; q^{13}),\\
  a_6(z) &=q^{\frac{1}{104}} (q^6; q^{13})(q^7; q^{13})(q^{13}; q^{13}).
\end{array}\right.\eqno{(7.5)}$$
Ramanujan's general theta function
$f(a, b):=\sum_{n=-\infty}^{\infty} a^{n(n+1)/2} b^{n(n-1)/2}$ and
$f(-q):=f(-q, -q^2)=\sum_{n=-\infty}^{\infty} (-1)^n q^{n(3n-1)/2}=(q; q)_{\infty}$.
In his notebooks (see \cite{R1}, p.326, \cite{R2}, p.244 and
\cite{B}, p.372), Ramanujan obtained his modular equations of
degree $13$.

\noindent{\bf Theorem 7.5.} {\it Define
$$\mu_1=\frac{f(-q^4, -q^9)}{q^{7/13} f(-q^2, -q^{11})},
  \mu_2=\frac{f(-q^6, -q^7)}{q^{6/13} f(-q^3, -q^{10})},
  \mu_3=\frac{f(-q^2, -q^{11})}{q^{5/13} f(-q, -q^{12})},$$
$$\mu_4=\frac{f(-q^5, -q^8)}{q^{2/13} f(-q^4, -q^9)},
  \mu_5=\frac{q^{5/13} f(-q^3, -q^{10})}{f(-q^5, -q^8)},
  \mu_6=\frac{q^{15/13} f(-q, -q^{12})}{f(-q^6, -q^7)}.$$
Then
$$1+\frac{f^2(-q)}{q f^2(-q^{13})}=\mu_1 \mu_2-\mu_3 \mu_5-
  \mu_4 \mu_6,\eqno{(7.6)}$$
$$-4-\frac{f^2(-q)}{q f^2(-q^{13})}=\frac{1}{\mu_1 \mu_2}-
  \frac{1}{\mu_3 \mu_5}-\frac{1}{\mu_4 \mu_6},\eqno{(7.7)}$$
$$3+\frac{f^2(-q)}{q f^2(-q^{13})}=\mu_2 \mu_3 \mu_4-\mu_1
  \mu_5 \mu_6,\eqno{(7.8)}$$
where $\mu_1 \mu_2 \mu_3 \mu_4 \mu_5 \mu_6=1$.}

  Let
$$G_m(z):=G_{m, p}(z):=(-1)^m q^{m(3m-p)/(2p^2)} \frac{f(-q^{2m/p},
          -q^{1-2m/p})}{f(-q^{m/p}, -q^{1-m/p})}$$
where $m$ is a positive integer, $p=13$ and $q=e^{2 \pi i z}$.
Then the above three formulas (7.6), (7.7) and (7.8) are
equivalent to the following (see \cite{Ev} or \cite{B},
p.373-375)
$$1+G_1(z) G_5(z)+G_2(z) G_3(z)+G_4(z) G_6(z)=-\frac{\eta^2(z/13)}{\eta^2(z)},
  \eqno{(7.9)}$$
$$\frac{1}{G_1(z) G_5(z)}+\frac{1}{G_2(z) G_3(z)}+\frac{1}{G_4(z) G_6(z)}
 =4+\frac{\eta^2(z/13)}{\eta^2(z)},\eqno{(7.10)}$$
$$G_1(z) G_3(z) G_4(z)-\frac{1}{G_1(z) G_3(z) G_4(z)}
 =3+\frac{\eta^2(z/13)}{\eta^2(z)},\eqno{(7.11)}$$
where $G_1(z) G_2(z) G_3(z) G_4(z) G_5(z) G_6(z)=-1$. Moreover,
there is the following beautiful formula (see \cite{Ev} or \cite{B},
p.375-376) in the spirit of (7.6), (7.7) and (7.8): for $t=q^{1/13}$,
$$\frac{1}{(t^2)_{\infty} (t^3)_{\infty} (t^{10})_{\infty} (t^{11})_{\infty}}
 +\frac{t}{(t^4)_{\infty} (t^6)_{\infty} (t^7)_{\infty} (t^9)_{\infty}}
 =\frac{1}{(t)_{\infty} (t^5)_{\infty} (t^8)_{\infty} (t^{12})_{\infty}},
  \eqno{(7.12)}$$
which is equivalent to, for $p=13$,
$G_1^{-1}(z) G_5^{-1}(z)+G_4(z) G_6(z)=1$.
In \cite{Y2}, we proved that it is equivalent to the following formula:
$$\frac{1}{a_1(z) a_4(z)}+\frac{1}{a_2(z) a_5(z)}+\frac{1}{a_3(z) a_6(z)}=0.
  \eqno{(7.13)}$$
Hence, it would be of interest to have similar elegant formulas as above.

  In fact, we show that there are twenty-one such formulas (i.e. modular
equations of order $13$) which generate the invariant ideal $I(Y)$ associated
with the modular curve $X=X(13)$.

\textbf{Theorem 7.6.} (Modularity for the ideal $I(Y)$) {\it The invariant
ideal $I(Y)$ can be parameterized by theta constants of order $13$, i.e.,
there are the following twenty-one modular equations of order $13$:
$$\left\{\aligned
  \mathbf{B}_0^{(i)}(a_1(z), \ldots, a_6(z)) &=0,\\
  \mathbf{B}_1^{(j)}(a_1(z), \ldots, a_6(z)) &=0,\\
  \mathbf{B}_3^{(j)}(a_1(z), \ldots, a_6(z)) &=0,\\
  \mathbf{B}_9^{(j)}(a_1(z), \ldots, a_6(z)) &=0,\\
  \mathbf{B}_{12}^{(j)}(a_1(z), \ldots, a_6(z)) &=0,\\
  \mathbf{B}_{10}^{(j)}(a_1(z), \ldots, a_6(z)) &=0,\\
  \mathbf{B}_4^{(j)}(a_1(z), \ldots, a_6(z)) &=0,\\
  \mathbf{B}_5(a_1(z), \ldots, a_6(z)) &=0,\\
  \mathbf{B}_2(a_1(z), \ldots, a_6(z)) &=0,\\
  \mathbf{B}_6(a_1(z), \ldots, a_6(z)) &=0,\\
  \mathbf{B}_8(a_1(z), \ldots, a_6(z)) &=0,\\
  \mathbf{B}_{11}(a_1(z), \ldots, a_6(z)) &=0,\\
  \mathbf{B}_7(a_1(z), \ldots, a_6(z)) &=0,
\endaligned\right.\eqno{(7.14)}$$
where $i=0, 1, 2$ and $j=1, 2$.}

{\it Numerical evidence by $q$-expansion}. Before obtaining a
rigorous proof, we give some computation by $q$-expansion to
provide some numerical evidence. We have
$$\left\{\aligned
  a_1(z) &=q^{\frac{121}{104}} (1-q-q^{12}+q^{15}+\cdots),\\
  a_2(z) &=q^{\frac{49}{104}} (1-q^3-q^{10}+q^{19}+\cdots),\\
  a_3(z) &=q^{\frac{25}{104}} (1-q^4-q^9+q^{21}+\cdots),\\
  a_4(z) &=-q^{\frac{9}{104}} (1-q^5-q^8+q^{23}+\cdots),\\
  a_5(z) &=q^{\frac{81}{104}} (1-q^2-q^{11}+q^{17}+\cdots),\\
  a_6(z) &=q^{\frac{1}{104}} (1-q^6-q^7+q^{25}+\cdots).
\endaligned\right.\eqno{(7.15)}$$
Hence,
$$\aligned
 &\mathbf{B}_0^{(0)}(a_1(z), \ldots, a_6(z))\\
=&q^{\frac{5}{2}} (-1) (1-q-q^{12}+q^{15}+\cdots)(1-q^3-q^{10}+q^{19}+\cdots)\\
 &\times (1-q^5-q^8+q^{23}+\cdots)(1-q^2-q^{11}+q^{17}+\cdots)\\
 &+q^{\frac{3}{2}} (1-q^3-q^{10}+q^{19}+\cdots)(1-q^4-q^9+q^{21}+\cdots)\\
 &\times (1-q^2-q^{11}+q^{17}+\cdots)(1-q^6-q^7+q^{25}+\cdots)\\
 &+q^{\frac{3}{2}} (-1) (1-q^4-q^9+q^{21}+\cdots)(1-q-q^{12}+q^{15}+\cdots)\\
 &\times (1-q^6-q^7+q^{25}+\cdots)(1-q^5-q^8+q^{23}+\cdots) \\
=&q^{\frac{1}{2}} [-q^2 (1-q-q^2+q^4+q^7-q^8-q^{10}+q^{11}+\\
 &-2 q^{13}+3 q^{14}+2 q^{15}-2 q^{17}-q^{19}-3 q^{20}+\cdots)\\
 &+q (1-q^2-q^3-q^4+q^5+q^8+q^{10}-2 q^{13}+\\
 &+3 q^{15}+2 q^{16}+2 q^{17}-3 q^{18}-q^{20}+\cdots)\\
 &-q (1-q-q^4+q^9+q^{10}+q^{11}-q^{12}-2 q^{13}+\\
 &+2 q^{14}+2 q^{17}-q^{18}+\cdots)]\\
=&q^{\frac{1}{2}}(0+O(q^{20})).
\endaligned$$
$$\aligned
 &\mathbf{B}_0^{(1)}(a_1(z), \ldots, a_6(z))\\
=&q^{\frac{7}{2}} (1-q-q^{12}+q^{15}+\cdots)(1-q^2-q^{11}+q^{17}+\cdots)^3\\
 &+q^{\frac{1}{2}} (1-q^3-q^{10}+q^{19}+\cdots)(1-q^6-q^7+q^{25}+\cdots)^3\\
 &-q^{\frac{1}{2}} (1-q^4-q^9+q^{21}+\cdots)(1-q^5-q^8+q^{23}+\cdots)^3\\
=&q^{\frac{1}{2}} [q^3 (1-q-3 q^2+3 q^3+3 q^4-3 q^5-q^6+q^7-3 q^{11}+2 q^{12}+\\
 &+6 q^{13}-3 q^{14}-2 q^{15}-2 q^{18}-3 q^{19}+6 q^{20}+\cdots)\\
 &+(1-q^3-3 q^6-3 q^7+3 q^9+2 q^{10}+3 q^{12}+6 q^{13}+3 q^{14}+\\
 &-3 q^{15}-3 q^{16}-q^{18}-2 q^{19}+\cdots)\\
 &-(1-q^4-3 q^5-3 q^8+2 q^9+3 q^{10}+3 q^{12}+6 q^{13}+\\
 &-q^{15}+3 q^{16}-3 q^{17}-3 q^{18}-2 q^{19}-3 q^{20}+\cdots)]\\
=&q^{\frac{1}{2}}(0+O(q^{20})).
\endaligned$$
$$\aligned
 &\mathbf{B}_0^{(2)}(a_1(z), \ldots, a_6(z))\\
=&q^{\frac{7}{2}} (1-q-q^{12}+q^{15}+\cdots)^3 (1-q^6-q^7+q^{25}+\cdots)\\
 &-q^{\frac{3}{2}} (1-q^3-q^{10}+q^{19}+\cdots)^3 (1-q^5-q^8+q^{23}+\cdots)\\
 &+q^{\frac{3}{2}} (1-q^4-q^9+q^{21}+\cdots)^3 (1-q^2-q^{11}+q^{17}+\cdots)\\
=&q^{\frac{1}{2}} [q^3 (1-3q+3 q^2-q^3-q^6+2 q^7-2 q^9+q^{10}-3 q^{12}+\\
 &+6 q^{13}-3 q^{14}+3 q^{15}-6 q^{16}+3 q^{17}+3 q^{18}-3 q^{19}+\cdots)\\
 &-q (1-3 q^3-q^5+3 q^6+2 q^8-q^9-3 q^{10}+6 q^{13}-2 q^{14}+\\
 &+3 q^{15}-3 q^{16}+q^{17}-3 q^{18}+3 q^{19}+\cdots)\\
 &+q (1-q^2-3 q^4+3 q^6+3 q^8-3 q^9-3 q^{10}+2 q^{11}-q^{12}+\\
 &+6 q^{13}+q^{14}-3 q^{15}-2 q^{17}+3 q^{18}+\cdots)]\\
=&q^{\frac{1}{2}}(0+O(q^{20})).
\endaligned$$
$$\aligned
 &\mathbf{B}_1^{(1)}(a_1(z), \ldots, a_6(z))\\
=&q^{\frac{67}{26}} (1-q^4-q^9+q^{21}+\cdots) (1-q^2-q^{11}+q^{17}+\cdots)^3\\
 &-q^{\frac{93}{26}} (1-q-q^{12}+q^{15}+\cdots)^3 (1-q^5-q^8+q^{23}+\cdots)\\
 &-q^{\frac{67}{26}} (1-q-q^{12}+q^{15}+\cdots) (1-q^3-q^{10}+q^{19}+\cdots)^3\\
=&q^{\frac{67}{26}} [(1-3 q^2+2 q^4+2 q^6-3 q^8-q^9+q^{10}+3 q^{13}\\
 &+q^{15}-3 q^{17}-3 q^{19}+\cdots)\\
 &-q (1-3 q+3 q^2-q^3-q^5+3 q^6-3 q^7+3 q^9-3 q^{10}+q^{11}\\
 &-3 q^{12}+6 q^{13}-3 q^{14}+3 q^{15}-6 q^{16}+6 q^{17}-6 q^{18}+3 q^{19}+\cdots)\\
 &-(1-q-3 q^3+3 q^4+3 q^6-3 q^7-q^9-2 q^{10}+3 q^{11}-q^{12}\\
 &+6 q^{13}-6 q^{14}+4 q^{15}-3 q^{16}+3 q^{17}-6 q^{18}+3 q^{19}+\cdots)]\\
=&q^{\frac{67}{26}}(0+O(q^{20})).
\endaligned$$
$$\aligned
 &\mathbf{B}_1^{(2)}(a_1(z), \ldots, a_6(z))\\
=&-q^{\frac{15}{26}} (1-q^3-q^{10}+q^{19}+\cdots)(1-q^5-q^8+q^{23}+\cdots)\\
 &\times (1-q^6-q^7+q^{25}+\cdots)^2\\
 &+q^{\frac{15}{26}} (1-q^4-q^9+q^{21}+\cdots)^2 (1-q^6-q^7+q^{25}+\cdots)\\
 &\times (1-q^5-q^8+q^{23}+\cdots)\\
 &-q^{\frac{93}{26}} (1-q-q^{12}+q^{15}+\cdots)^2 (1-q^3-q^{10}+q^{19}+\cdots)\\
 &\times (1-q^2-q^{11}+q^{17}+\cdots)\\
=&q^{\frac{15}{26}} [-(1-q^3-q^5-2 q^6-2 q^7+2 q^9+q^{10}+3 q^{11}\\
 &+3 q^{12}+2 q^{13}+q^{14}-2 q^{17}-3 q^{18}+\cdots)\\
 &+(1-2 q^4-q^5-q^6-q^7+2 q^{10}+3 q^{11}+3 q^{12}+\\
 &+q^{13}+2 q^{14}-q^{16}+2 q^{17}-3 q^{18}-3 q^{19}+\cdots)\\
 &-q^3 (1-2 q+q^3+q^4-2 q^6+q^7-q^{10}+q^{11}-q^{13}\\
 &+4 q^{14}-3 q^{16}-3 q^{17}+4 q^{19}+\cdots)]\\
=&q^{\frac{15}{26}}(0+O(q^{20})).
\endaligned$$
$$\aligned
 &\mathbf{B}_3^{(1)}(a_1(z), \ldots, a_6(z))\\
=&-q^{\frac{19}{26}} (1-q^3-q^{10}+q^{19}+\cdots) (1-q^5-q^8+q^{23}+\cdots)^3\\
 &+q^{\frac{19}{26}} (1-q^4-q^9+q^{21}+\cdots)^3 (1-q^6-q^7+q^{25}+\cdots)\\
 &-q^{\frac{97}{26}} (1-q^4-q^9+q^{21}+\cdots) (1-q-q^{12}+q^{15}+\cdots)^3\\
=&q^{\frac{19}{26}} [-(1-q^3-3 q^5+2 q^{10}+3 q^{11}+3 q^{13}+2 q^{15}-3 q^{16}
  +q^{18}-2 q^{19}+\cdots)\\
 &+(1-3 q^4-q^6-q^7+3 q^8-3 q^9+3 q^{10}+3 q^{11}\\
 &-q^{12}+6 q^{13}-3 q^{14}+3 q^{16}-3 q^{17}+4 q^{18}-5 q^{19}+\cdots)\\
 &-q^3 (1-3 q+3 q^2-q^3-q^4+3 q^5-3 q^6+q^7-q^9+3 q^{10}-3 q^{11}-2 q^{12}\\
 &+6 q^{13}-3 q^{14}+3 q^{15}-3 q^{16}-3 q^{17}+3 q^{18}-3 q^{19}+\cdots)]\\
=&q^{\frac{19}{26}}(0+O(q^{20})).
\endaligned$$
$$\aligned
 &\mathbf{B}_3^{(2)}(a_1(z), \ldots, a_6(z))\\
=&q^{\frac{71}{26}} (1-q-q^{12}+q^{15}+\cdots)(1-q^6-q^7+q^{25}+\cdots)\\
 &\times (1-q^2-q^{11}+q^{17}+\cdots)^2\\
 &-q^{\frac{45}{26}} (1-q^3-q^{10}+q^{19}+\cdots)^2 (1-q^2-q^{11}+q^{17}+\cdots)\\
 &\times (1-q^6-q^7+q^{25}+\cdots)\\
 &+q^{\frac{45}{26}} (1-q^4-q^9+q^{21}+\cdots)^2 (1-q-q^{12}+q^{15}+\cdots)\\
 &\times (1-q^5-q^8+q^{23}+\cdots)\\
=&q^{\frac{45}{26}} [q (1-q-2 q^2+2 q^3+q^4-q^5-q^6+3 q^8-3 q^{10}-2 q^{11}\\
 &+2 q^{12}+2 q^{13}+q^{15}-q^{16}+2 q^{17}-q^{18}-4 q^{19}+\cdots)\\
 &-(1-q^2-2 q^3+2 q^5-q^7+3 q^9-3 q^{11}-q^{12}+\\
 &+q^{13}+3 q^{14}-q^{15}+2 q^{16}+3 q^{17}-q^{18}-2 q^{19}+\cdots)\\
 &+(1-q-2 q^4+q^5+q^6+q^{12}-q^{13}+q^{14}-q^{15}\\
 &+q^{16}+4 q^{17}-3 q^{18}-q^{19}+\cdots)]\\
=&q^{\frac{71}{26}}(0+O(q^{20})).
\endaligned$$
$$\aligned
 &\mathbf{B}_9^{(1)}(a_1(z), \ldots, a_6(z))\\
=&q^{\frac{31}{26}} (1-q-q^{12}+q^{15}+\cdots) (1-q^6-q^7+q^{25}+\cdots)^3\\
 &+q^{\frac{57}{26}} (1-q^3-q^{10}+q^{19}+\cdots)^3 (1-q^2-q^{11}+q^{17}+\cdots)\\
 &-q^{\frac{31}{26}} (1-q^3-q^{10}+q^{19}+\cdots) (1-q^4-q^9+q^{21}+\cdots)^3\\
=&q^{\frac{31}{26}} [(1-q-3 q^6+3 q^8+2 q^{12}+3 q^{13}-3 q^{14}-2 q^{15}
  +2 q^{18}+q^{19}+\cdots)\\
 &+q (1-q^2-3 q^3+3 q^5+3 q^6-3 q^8-q^9-3 q^{10}+3 q^{12}\\
 &+6 q^{13}+3 q^{14}-6 q^{15}-3 q^{16}-2 q^{17}+3 q^{18}+3 q^{19}+\cdots)\\
 &-(1-q^3-3 q^4+3 q^7+3 q^8-3 q^9-q^{10}-3 q^{11}+2 q^{12}\\
 &+6 q^{13}+3 q^{14}+q^{15}-6 q^{16}-3 q^{17}+4 q^{19}+\cdots)]\\
=&q^{\frac{31}{26}}(0+O(q^{20})).
\endaligned$$
$$\aligned
 &\mathbf{B}_9^{(2)}(a_1(z), \ldots, a_6(z))\\
=&q^{\frac{31}{26}} (1-q^4-q^9+q^{21}+\cdots)(1-q^2-q^{11}+q^{17}+\cdots)\\
 &\times (1-q^5-q^8+q^{23}+\cdots)^2\\
 &+q^{\frac{83}{26}} (1-q-q^{12}+q^{15}+\cdots)^2 (1-q^5-q^8+q^{23}+\cdots)\\
 &\times (1-q^2-q^{11}+q^{17}+\cdots)\\
 &-q^{\frac{31}{26}} (1-q^3-q^{10}+q^{19}+\cdots)^2 (1-q^4-q^9+q^{21}+\cdots)\\
 &\times (1-q^6-q^7+q^{25}+\cdots)\\
=&q^{\frac{31}{26}} [(1-q^2-q^4-2 q^5+q^6+2 q^7-2 q^8+q^9+3 q^{10}-2 q^{11}\\
 &+q^{12}+2 q^{13}-q^{14}-q^{15}+2 q^{16}+q^{17}-q^{18}+q^{19}+\cdots)\\
 &+q^2 (1-2 q+2 q^3-q^4-q^5+2 q^6-3 q^8+3 q^9-3 q^{11}+q^{12}\\
 &+q^{13}+2 q^{14}-q^{16}-q^{17}-q^{18}+\cdots)\\
 &-(1-2 q^3-q^4+q^7+q^9+q^{11}+q^{12}-q^{13}\\
 &+4 q^{16}+q^{17}-2 q^{18}+\cdots)]\\
=&q^{\frac{31}{26}}(0+O(q^{20})).
\endaligned$$
$$\aligned
 &\mathbf{B}_{12}^{(1)}(a_1(z), \ldots, a_6(z))\\
=&-q^{\frac{37}{26}} (1-q-q^{12}+q^{15}+\cdots) (1-q^5-q^8+q^{23}+\cdots)^3\\
 &+q^{\frac{37}{26}} (1-q^3-q^{10}+q^{19}+\cdots)^3 (1-q^6-q^7+q^{25}+\cdots)\\
 &-q^{\frac{63}{26}} (1-q^5-q^8+q^{23}+\cdots) (1-q^2-q^{11}+q^{17}+\cdots)^3\\
=&q^{\frac{37}{26}} [-(1-q-3 q^5+3 q^6-3 q^8+3 q^9+3 q^{10}-3 q^{11}-q^{12}\\
 &+6 q^{13}-6 q^{14}+4 q^{16}-3 q^{18}+3 q^{19}+\cdots)\\
 &+(1-3 q^3+2 q^6-q^7+2 q^9-3 q^{12}+3 q^{13}+q^{15}+q^{16}+3 q^{17}-3 q^{19}+\cdots)\\
 &-q (1-3 q^2+3 q^4-q^5-q^6+3 q^7-q^8-3 q^9+3 q^{10}-2 q^{11}-3 q^{12}\\
 &+6 q^{13}+q^{14}-3 q^{15}+3 q^{16}+3 q^{17}-6 q^{18}-3 q^{19}+\cdots)]\\
=&q^{\frac{37}{26}}(0+O(q^{20})).
\endaligned$$
$$\aligned
 &\mathbf{B}_{12}^{(2)}(a_1(z), \ldots, a_6(z))\\
=&q^{\frac{37}{26}} (1-q^5-q^8+q^{23}+\cdots)^2 (1-q^3-q^{10}+q^{19}+\cdots)\\
 &\times (1-q^2-q^{11}+q^{17}+\cdots)\\
 &-q^{\frac{63}{26}} (1-q^4-q^9+q^{21}+\cdots)^2 (1-q-q^{12}+q^{15}+\cdots)\\
 &\times (1-q^2-q^{11}+q^{17}+\cdots)\\
 &-q^{\frac{37}{26}} (1-q-q^{12}+q^{15}+\cdots) (1-q^4-q^9+q^{21}+\cdots)\\
 &\times (1-q^6-q^7+q^{25}+\cdots)^2\\
=&q^{\frac{37}{26}} [(1-q^2-q^3-q^5+2 q^7+q^{11}-q^{13}+q^{14}+q^{15}\\
 &+q^{16}-q^{17}+3 q^{18}+\cdots)\\
 &-q (1-q-q^2+q^3-2 q^4+2 q^5+2 q^6-2 q^7+q^8-3 q^9+q^{10}+2 q^{11}\\
 &-2 q^{12}+2 q^{13}-q^{14}+q^{15}+2 q^{16}-2 q^{18}+\cdots)\\
 &-(1-q-q^4+q^5-2 q^6+2 q^8-q^9+3 q^{10}-2 q^{12}+q^{13}\\
 &-q^{14}+2 q^{15}-3 q^{17}+3 q^{18}+2 q^{19}+\cdots)]\\
=&q^{\frac{37}{26}}(0+O(q^{20})).
\endaligned$$
$$\aligned
 &\mathbf{B}_{10}^{(1)}(a_1(z), \ldots, a_6(z))\\
=&q^{\frac{7}{26}} (1-q^4-q^9+q^{21}+\cdots) (1-q^6-q^7+q^{25}+\cdots)^3\\
 &+q^{\frac{111}{26}} (1-q-q^{12}+q^{15}+\cdots)^3 (1-q^2-q^{11}+q^{17}+\cdots)\\
 &-q^{\frac{7}{26}} (1-q^6-q^7+q^{25}+\cdots) (1-q^5-q^8+q^{23}+\cdots)^3\\
=&q^{\frac{7}{26}} [(1-q^4-3 q^6-3 q^7-q^9+3 q^{10}+3 q^{11}+3 q^{12}\\
 &+6 q^{13}+3 q^{14}+3 q^{15}-6 q^{17}-4 q^{18}-3 q^{19}+\cdots)\\
 &+q^4 (1-3 q+2 q^2+2 q^3-3 q^4+q^5-q^{11}+3 q^{13}\\
 &+q^{14}-3 q^{15}-3 q^{16}+q^{17}+3 q^{18}+\cdots)\\
 &-(1-3 q^5-q^6-q^7-3 q^8+3 q^{10}+3 q^{11}+3 q^{12}\\
 &+6 q^{13}+3 q^{14}+2 q^{15}-3 q^{17}-3 q^{18}-6 q^{19}+\cdots)]\\
=&q^{\frac{7}{26}}(0+O(q^{20})).
\endaligned$$
$$\aligned
 &\mathbf{B}_{10}^{(2)}(a_1(z), \ldots, a_6(z))\\
=&-q^{\frac{33}{26}} (1-q^6-q^7+q^{25}+\cdots)^2 (1-q-q^{12}+q^{15}+\cdots)\\
 &\times (1-q^5-q^8+q^{23}+\cdots)\\
 &+q^{\frac{33}{26}} (1-q^3-q^{10}+q^{19}+\cdots)^2 (1-q^4-q^9+q^{21}+\cdots)\\
 &\times (1-q^5-q^8+q^{23}+\cdots)\\
 &-q^{\frac{59}{26}} (1-q^3-q^{10}+q^{19}+\cdots) (1-q^4-q^9+q^{21}+\cdots)\\
 &\times (1-q^2-q^{11}+q^{17}+\cdots)^2\\
=&q^{\frac{33}{26}} [-(1-q-q^5-q^6+q^8+q^9+2 q^{11}-q^{13}+q^{14}\\
 &-2 q^{16}+q^{18}+3 q^{19}+\cdots)\\
 &+(1-2 q^3-q^4-q^5+q^6+2 q^7+q^8-3 q^{10}+q^{11}\\
 &+q^{12}+2 q^{13}+2 q^{14}-3 q^{17}+q^{18}+2 q^{19}+\cdots)\\
 &-q (1-2 q^2-q^3+2 q^5+2 q^6-q^8-3 q^9-q^{10}+q^{11}\\
 &+3 q^{12}+q^{13}+2 q^{15}-3 q^{16}-q^{18}+\cdots)]\\
=&q^{\frac{33}{26}}(0+O(q^{20})).
\endaligned$$
$$\aligned
 &\mathbf{B}_4^{(1)}(a_1(z), \ldots, a_6(z))\\
=&q^{\frac{73}{26}} (1-q^3-q^{10}+q^{19}+\cdots) (1-q^2-q^{11}+q^{17}+\cdots)^3\\
 &-q^{\frac{21}{26}} (1-q^4-q^9+q^{21}+\cdots)^3 (1-q^5-q^8+q^{23}+\cdots)\\
 &+q^{\frac{21}{26}} (1-q^2-q^{11}+q^{17}+\cdots) (1-q^6-q^7+q^{25}+\cdots)^3\\
=&q^{\frac{21}{26}} [q^2 (1-3 q^2-q^3+3 q^4+3 q^5-q^6-3 q^7+q^9-q^{10}-3 q^{11}+3 q^{12}\\
 &+6 q^{13}-3 q^{15}-5 q^{16}+3 q^{17}+3 q^{18}-5 q^{19}+\cdots)\\
 &-(1-3 q^4-q^5+2 q^8+2 q^{12}+3 q^{13}+3 q^{14}-3 q^{16}+q^{17}-3 q^{18}+\cdots)\\
 &+(1-q^2-3 q^6-3 q^7+3 q^8+3 q^9-q^{11}+3 q^{12}+6 q^{13}\\
 &-6 q^{15}-3 q^{16}+4 q^{17}+2 q^{18}-3 q^{19}+\cdots)]\\
=&q^{\frac{21}{26}}(0+O(q^{20})).
\endaligned$$
$$\aligned
 &\mathbf{B}_4^{(2)}(a_1(z), \ldots, a_6(z))\\
=&q^{\frac{47}{26}} (1-q^2-q^{11}+q^{17}+\cdots)^2 (1-q^4-q^9+q^{21}+\cdots)\\
 &\times (1-q^6-q^7+q^{25}+\cdots)\\
 &-q^{\frac{73}{26}} (1-q-q^{12}+q^{15}+\cdots)^2 (1-q^3-q^{10}+q^{19}+\cdots)\\
 &\times (1-q^6-q^7+q^{25}+\cdots)\\
 &-q^{\frac{47}{26}} (1-q-q^{12}+q^{15}+\cdots) (1-q^3-q^{10}+q^{19}+\cdots)\\
 &\times (1-q^5-q^8+q^{23}+\cdots)^2\\
=&q^{\frac{47}{26}} [(1-2 q^2+q^6-q^7+q^8+q^9-2 q^{12}-q^{13}+q^{14}+4 q^{15}\\
 &+q^{16}-3 q^{19}+\cdots)\\
 &-q (1-2 q+q^2-q^3+2 q^4-q^5-q^6+q^7+q^8-2 q^{10}+q^{11}\\
 &-2 q^{12}+2 q^{13}+4 q^{15}-3 q^{16}-q^{17}-q^{18}+4 q^{19}+\cdots)\\
 &-(1-q-q^3+q^4-2 q^5+2 q^6+2 q^{11}-3 q^{12}+q^{13}\\
 &-q^{14}+4 q^{15}-3 q^{16}+3 q^{17}+q^{18}-2 q^{19}+\cdots)]\\
=&q^{\frac{47}{26}}(0+O(q^{20})).
\endaligned$$
$$\aligned
 &\mathbf{B}_5(a_1(z), \ldots, a_6(z))\\
=&-q^{\frac{75}{26}} (1-q^3-q^{10}+q^{19}+\cdots)^2 (1-q-q^{12}+q^{15}+\cdots)\\
 &\times (1-q^2-q^{11}+q^{17}+\cdots)\\
 &-q^{\frac{23}{26}} (1-q^5-q^8+q^{23}+\cdots) (1-q^2-q^{11}+q^{17}+\cdots)\\
 &\times (1-q^6-q^7+q^{25}+\cdots)^2\\
 &+q^{\frac{23}{26}} (1-q^3-q^{10}+q^{19}+\cdots) (1-q^4-q^9+q^{21}+\cdots)\\
 &\times (1-q^5-q^8+q^{23}+\cdots)^2\\
=&q^{\frac{23}{26}} [-q^2 (1-q-q^2-q^3+2 q^4+2 q^5-q^6-q^7-q^8+q^9\\
 &-2 q^{10}+q^{11}+2 q^{12}+q^{14}-q^{15}+2 q^{16}-3 q^{17}+\cdots)\\
 &-(1-q^2-q^5-2 q^6-q^7+q^8+2 q^9+q^{10}+\\
 &+q^{11}+3 q^{12}-2 q^{16}+q^{19}+\cdots)\\
 &+(1-q^3-q^4-2 q^5+q^7+q^9+2 q^{11}+q^{12}+\\
 &+q^{13}+2 q^{14}-q^{16}-q^{17}+2 q^{18}-2 q^{19}+\cdots)]\\
=&q^{\frac{23}{26}}(0+O(q^{20})).
\endaligned$$
$$\aligned
 &\mathbf{B}_2(a_1(z), \ldots, a_6(z))\\
=&q^{\frac{69}{26}} (1-q-q^{12}+q^{15}+\cdots)^2 (1-q^4-q^9+q^{21}+\cdots)\\
 &\times (1-q^5-q^8+q^{23}+\cdots)\\
 &-q^{\frac{43}{26}} (1-q^6-q^7+q^{25}+\cdots) (1-q^5-q^8+q^{23}+\cdots)\\
 &\times (1-q^2-q^{11}+q^{17}+\cdots)^2\\
 &+q^{\frac{43}{26}} (1-q-q^{12}+q^{15}+\cdots) (1-q^3-q^{10}+q^{19}+\cdots)\\
 &\times (1-q^6-q^7+q^{25}+\cdots)^2\\
=&q^{\frac{43}{26}} [q (1-2q+q^2-q^4+q^5+q^6-q^7-q^8+2 q^9-q^{10}-q^{12}\\
 &+2 q^{14}+q^{16}+q^{17}-4 q^{18}+\cdots)\\
 &-(1-2 q^2+q^4-q^5-q^6+q^7+q^8+q^9+q^{10}-2 q^{11}+\\
 &-q^{14}+2 q^{15}+q^{16}+2 q^{17}+q^{18}-q^{19}+\cdots)\\
 &+(1-q-q^3+q^4-2 q^6+2 q^8+2 q^9-q^{10}-q^{11}+\\
 &+q^{13}-q^{14}+q^{16}+q^{17}+3 q^{19}+\cdots)]\\
=&q^{\frac{43}{26}}(0+O(q^{20})).
\endaligned$$
$$\aligned
 &\mathbf{B}_6(a_1(z), \ldots, a_6(z))\\
=&-q^{\frac{25}{26}} (1-q^4-q^9+q^{21}+\cdots)^2 (1-q^3-q^{10}+q^{19}+\cdots)\\
 &\times (1-q^6-q^7+q^{25}+\cdots)\\
 &+q^{\frac{25}{26}} (1-q^2-q^{11}+q^{17}+\cdots) (1-q^6-q^7+q^{25}+\cdots)\\
 &\times (1-q^5-q^8+q^{23}+\cdots)^2\\
 &+q^{\frac{77}{26}} (1-q^4-q^9+q^{21}+\cdots) (1-q-q^{12}+q^{15}+\cdots)\\
 &\times (1-q^2-q^{11}+q^{17}+\cdots)^2\\
=&q^{\frac{25}{26}} [-(1-q^3-2 q^4-q^6+q^7+q^8-q^9+2 q^{10}+q^{11}+2 q^{12}\\
 &-q^{14}+q^{15}+q^{16}+2 q^{17}-q^{18}-q^{19}+\cdots)\\
 &+(1-q^2-2 q^5-q^6+q^7-q^8+q^9+3 q^{10}+q^{11}\\
 &+q^{12}-q^{17}+q^{18}+q^{19}+\cdots)\\
 &+q^2 (1-q-2 q^2+2 q^3+2 q^6-2 q^7-q^8+q^{10}-q^{12}+\\
 &+q^{13}+q^{14}+3 q^{15}-2 q^{16}-2 q^{17}+\cdots)]\\
=&q^{\frac{25}{26}}(0+O(q^{20})).
\endaligned$$
$$\aligned
 &\mathbf{B}_8(a_1(z), \ldots, a_6(z))\\
=&-q^{\frac{55}{26}} (1-q^3-q^{10}+q^{19}+\cdots) (1-q^5-q^8+q^{23}+\cdots)\\
 &\times (1-q^2-q^{11}+q^{17}+\cdots)^2\\
 &+q^{\frac{55}{26}} (1-q-q^{12}+q^{15}+\cdots) (1-q^3-q^{10}+q^{19}+\cdots)\\
 &\times (1-q^4-q^9+q^{21}+\cdots)^2\\
 &+q^{\frac{81}{26}} (1-q-q^{12}+q^{15}+\cdots)^2 (1-q^2-q^{11}+q^{17}+\cdots)\\
 &\times (1-q^6-q^7+q^{25}+\cdots)\\
=&q^{\frac{55}{26}} [-(1-2 q^2-q^3+q^4+q^5+q^7-q^9-q^{10}-q^{11}+2 q^{12}\\
 &+q^{14}+2 q^{15}-q^{18}+\cdots)\\
 &+(1-q-q^3-q^4+2 q^5+2 q^7-q^8-3 q^9+q^{10}+2 q^{12}\\
 &+2 q^{17}-q^{18}-q^{19}+\cdots)\\
 &+q (1-2 q+2 q^3-q^4-q^6+q^7+2 q^8-2 q^9-q^{10}+\\
 &+q^{13}+2 q^{14}-2 q^{16}+q^{18}+\cdots)]\\
=&q^{\frac{55}{26}}(0+O(q^{20})).
\endaligned$$
$$\aligned
 &\mathbf{B}_{11}(a_1(z), \ldots, a_6(z))\\
=&q^{\frac{35}{26}} (1-q-q^{12}+q^{15}+\cdots) (1-q^6-q^7+q^{25}+\cdots)\\
 &\times (1-q^5-q^8+q^{23}+\cdots)^2\\
 &+q^{\frac{61}{26}} (1-q^4-q^9+q^{21}+\cdots) (1-q-q^{12}+q^{15}+\cdots)\\
 &\times (1-q^3-q^{10}+q^{19}+\cdots)^2\\
 &-q^{\frac{35}{26}} (1-q^4-q^9+q^{21}+\cdots)^2 (1-q^5-q^8+q^{23}+\cdots)\\
 &\times (1-q^2-q^{11}+q^{17}+\cdots)\\
=&q^{\frac{35}{26}} [(1-q-2 q^5+q^6-q^8+2 q^9+q^{10}+q^{11}-q^{12}\\
 &+q^{15}-2 q^{16}+q^{17}+2 q^{18}-q^{19}+\cdots)\\
 &+q (1-q-2 q^3+q^4+q^5+q^6+q^7-2 q^8-q^9-2 q^{10}+3 q^{11}\\
 &+q^{12}+2 q^{16}-2 q^{17}-q^{18}+q^{19}+\cdots)\\
 &-(1-q^2-2 q^4-q^5+2 q^6+q^7-q^{11}+2 q^{12}+q^{13}\\
 &+q^{15}-2 q^{16}+3 q^{17}-2 q^{19}+\cdots)]\\
=&q^{\frac{35}{26}}(0+O(q^{20})).
\endaligned$$
$$\aligned
 &\mathbf{B}_7(a_1(z), \ldots, a_6(z))\\
=&q^{\frac{27}{26}} (1-q^4-q^9+q^{21}+\cdots) (1-q^2-q^{11}+q^{17}+\cdots)\\
 &\times (1-q^6-q^7+q^{25}+\cdots)^2\\
 &+q^{\frac{79}{26}} (1-q^3-q^{10}+q^{19}+\cdots) (1-q^4-q^9+q^{21}+\cdots)\\
 &\times (1-q-q^{12}+q^{15}+\cdots)^2\\
 &-q^{\frac{27}{26}} (1-q^3-q^{10}+q^{19}+\cdots)^2 (1-q^6-q^7+q^{25}+\cdots)\\
 &\times (1-q^5-q^8+q^{23}+\cdots)\\
=&q^{\frac{27}{26}} [(1-q^2-q^4-q^6-2 q^7+2 q^8+q^9+2 q^{10}+2 q^{11}-q^{12}\\
 &+q^{15}-q^{17}+2 q^{19}+\cdots)\\
 &+q^2 (1-2 q+q^2-q^3+q^4+q^5-q^6+q^7-2 q^8+q^{10}+q^{11}\\
 &-2 q^{12}+2 q^{14}+2 q^{15}-q^{16}-2 q^{17}-2 q^{18}+\cdots)\\
 &-(1-2 q^3-q^5-q^7+q^8+2 q^9+2 q^{11}+q^{13}-2 q^{14}\\
 &+q^{15}+2 q^{16}+q^{17}-q^{18}+\cdots)]\\
=&q^{\frac{27}{26}}(0+O(q^{20})).
\endaligned$$

{\it Proof}. For simplicity, we set
$\theta \begin{bmatrix} \frac{k}{13}\\ 1 \end{bmatrix}:=
  \theta \begin{bmatrix} \frac{k}{13}\\ 1 \end{bmatrix}(0, 13z)$.
In the proof, we will only use Proposition 7.3. Let $k=13$ and
$(j; \alpha_1, \alpha_2, \alpha_3)=(1; 1, 5, 7)$. Then
$13-4j=9$, $13-2j=11$ and
$$\aligned
 &\theta \begin{bmatrix} \frac{9}{13}\\ 1 \end{bmatrix}
  \theta \begin{bmatrix} \frac{1}{13}\\ 1 \end{bmatrix}
  \theta \begin{bmatrix} \frac{5}{13}\\ 1 \end{bmatrix}
  \theta \begin{bmatrix} \frac{7}{13}\\ 1 \end{bmatrix}
 +e^{\frac{9 \pi i}{13}}
  \theta \begin{bmatrix} \frac{11}{13}\\ 1 \end{bmatrix}
  \theta \begin{bmatrix} \frac{3}{13}\\ 1 \end{bmatrix}
  \theta \begin{bmatrix} \frac{7}{13}\\ 1 \end{bmatrix}
  \theta \begin{bmatrix} \frac{9}{13}\\ 1 \end{bmatrix}+\\
 &+e^{-\frac{11 \pi i}{13}}
  \theta \begin{bmatrix} \frac{11}{13}\\ 1 \end{bmatrix}
  \theta \begin{bmatrix} \frac{-1}{13}\\ 1 \end{bmatrix}
  \theta \begin{bmatrix} \frac{3}{13}\\ 1 \end{bmatrix}
  \theta \begin{bmatrix} \frac{5}{13}\\ 1 \end{bmatrix}
 =0.
\endaligned$$
Note that
$$\theta \begin{bmatrix} \frac{-1}{13}\\ 1 \end{bmatrix}
 =e^{-\frac{\pi i}{13}} \theta \begin{bmatrix} \frac{1}{13}\\ 1 \end{bmatrix}.$$
This proves that
$$\aligned
 &\mathbf{B}_0^{(0)}(a_1(z), \ldots, a_6(z))\\
=&a_1(z) a_2(z) a_4(z) a_5(z)+a_2(z) a_3(z) a_5(z) a_6(z)+
  a_3(z) a_1(z) a_6(z) a_4(z)=0,
\endaligned$$
which is just the formula (7.13).

  Let $k=13$ and $(j; \alpha_1, \alpha_2, \alpha_3)=(1; 9, 9, -5)$.
Then $13-4j=9$, $13-2j=11$ and
$$\theta \begin{bmatrix} \frac{9}{13}\\ 1 \end{bmatrix}^3
  \theta \begin{bmatrix} \frac{-5}{13}\\ 1 \end{bmatrix}
 +e^{\frac{9 \pi i}{13}}
  \theta \begin{bmatrix} \frac{11}{13}\\ 1 \end{bmatrix}^3
  \theta \begin{bmatrix} \frac{-3}{13}\\ 1 \end{bmatrix}
 +e^{-\frac{11 \pi i}{13}}
  \theta \begin{bmatrix} \frac{11}{13}\\ 1 \end{bmatrix}
  \theta \begin{bmatrix} \frac{7}{13}\\ 1 \end{bmatrix}^2
  \theta \begin{bmatrix} \frac{-7}{13}\\ 1 \end{bmatrix}
 =0.$$
Note that
$$\theta \begin{bmatrix} \frac{-5}{13}\\ 1 \end{bmatrix}
 =e^{-\frac{5 \pi i}{13}} \theta \begin{bmatrix} \frac{5}{13}\\ 1 \end{bmatrix}, \quad
  \theta \begin{bmatrix} \frac{-3}{13}\\ 1 \end{bmatrix}
 =e^{-\frac{3 \pi i}{13}} \theta \begin{bmatrix} \frac{3}{13}\\ 1 \end{bmatrix}, \quad
  \theta \begin{bmatrix} \frac{-7}{13}\\ 1 \end{bmatrix}
 =e^{-\frac{7 \pi i}{13}} \theta \begin{bmatrix} \frac{7}{13}\\ 1 \end{bmatrix}.$$
This proves that
$$\mathbf{B}_1^{(1)}(a_1(z), \ldots, a_6(z))
 =a_3(z) a_5(z)^3+a_1(z)^3 a_4(z)-a_1(z) a_2(z)^3=0.$$

  Let $k=13$ and $(j; \alpha_1, \alpha_2, \alpha_3)=(3; 5, 5, 3)$.
Then $13-4j=1$, $13-2j=7$ and
$$\aligned
 &\theta \begin{bmatrix} \frac{1}{13}\\ 1 \end{bmatrix}
  \theta \begin{bmatrix} \frac{5}{13}\\ 1 \end{bmatrix}^2
  \theta \begin{bmatrix} \frac{3}{13}\\ 1 \end{bmatrix}
 +e^{\frac{\pi i}{13}}
  \theta \begin{bmatrix} \frac{7}{13}\\ 1 \end{bmatrix}
  \theta \begin{bmatrix} \frac{11}{13}\\ 1 \end{bmatrix}^2
  \theta \begin{bmatrix} \frac{9}{13}\\ 1 \end{bmatrix}+\\
 &+e^{-\frac{7 \pi i}{13}}
  \theta \begin{bmatrix} \frac{7}{13}\\ 1 \end{bmatrix}
  \theta \begin{bmatrix} \frac{-1}{13}\\ 1 \end{bmatrix}^2
  \theta \begin{bmatrix} \frac{-3}{13}\\ 1 \end{bmatrix}
 =0.
\endaligned$$
Note that
$$\theta \begin{bmatrix} \frac{-1}{13}\\ 1 \end{bmatrix}
 =e^{-\frac{\pi i}{13}} \theta \begin{bmatrix} \frac{1}{13}\\ 1 \end{bmatrix}, \quad
  \theta \begin{bmatrix} \frac{-3}{13}\\ 1 \end{bmatrix}
 =e^{-\frac{3 \pi i}{13}} \theta \begin{bmatrix} \frac{3}{13}\\ 1 \end{bmatrix}.$$
This proves that
$$\aligned
 &\mathbf{B}_1^{(2)}(a_1(z), \ldots, a_6(z))\\
=&a_2(z) a_4(z) a_6(z)^2-a_3(z)^2 a_6(z) a_4(z)-
  a_1(z)^2 a_2(z) a_5(z)=0.
\endaligned$$

  Let $k=13$ and $(j; \alpha_1, \alpha_2, \alpha_3)=(4; 3, 3, 7)$.
Then $13-4j=-3$, $13-2j=5$ and
$$\aligned
 &\theta \begin{bmatrix} \frac{-3}{13}\\ 1 \end{bmatrix}
  \theta \begin{bmatrix} \frac{3}{13}\\ 1 \end{bmatrix}^2
  \theta \begin{bmatrix} \frac{7}{13}\\ 1 \end{bmatrix}
 +e^{-\frac{3 \pi i}{13}}
  \theta \begin{bmatrix} \frac{5}{13}\\ 1 \end{bmatrix}
  \theta \begin{bmatrix} \frac{11}{13}\\ 1 \end{bmatrix}^2
  \theta \begin{bmatrix} \frac{15}{13}\\ 1 \end{bmatrix}+\\
 &+e^{-\frac{5 \pi i}{13}}
  \theta \begin{bmatrix} \frac{5}{13}\\ 1 \end{bmatrix}
  \theta \begin{bmatrix} \frac{-5}{13}\\ 1 \end{bmatrix}^2
  \theta \begin{bmatrix} \frac{-1}{13}\\ 1 \end{bmatrix}
 =0.
\endaligned$$
Note that
$$\theta \begin{bmatrix} \frac{-1}{13}\\ 1 \end{bmatrix}
 =e^{-\frac{\pi i}{13}} \theta \begin{bmatrix} \frac{1}{13}\\ 1 \end{bmatrix}, \quad
  \theta \begin{bmatrix} \frac{-3}{13}\\ 1 \end{bmatrix}
 =e^{-\frac{3 \pi i}{13}} \theta \begin{bmatrix} \frac{3}{13}\\ 1 \end{bmatrix}, \quad
  \theta \begin{bmatrix} \frac{-5}{13}\\ 1 \end{bmatrix}
 =e^{-\frac{5 \pi i}{13}} \theta \begin{bmatrix} \frac{5}{13}\\ 1 \end{bmatrix}$$
and
$$\theta \begin{bmatrix} \frac{15}{13}\\ 1 \end{bmatrix}
 =e^{-\frac{11 \pi i}{13}} \theta \begin{bmatrix} \frac{11}{13}\\ 1 \end{bmatrix}.$$
This proves that
$$\mathbf{B}_3^{(1)}(a_1(z), \ldots, a_6(z))
 =a_2(z) a_4(z)^3+a_3(z)^3 a_6(z)-a_3(z) a_1(z)^3=0.$$

  Let $k=13$ and $(j; \alpha_1, \alpha_2, \alpha_3)=(1; 7, 7, -1)$.
Then $13-4j=9$, $13-2j=11$ and
$$\aligned
 &\theta \begin{bmatrix} \frac{9}{13}\\ 1 \end{bmatrix}
  \theta \begin{bmatrix} \frac{7}{13}\\ 1 \end{bmatrix}^2
  \theta \begin{bmatrix} \frac{-1}{13}\\ 1 \end{bmatrix}
 +e^{\frac{9 \pi i}{13}}
  \theta \begin{bmatrix} \frac{11}{13}\\ 1 \end{bmatrix}
  \theta \begin{bmatrix} \frac{9}{13}\\ 1 \end{bmatrix}^2
  \theta \begin{bmatrix} \frac{1}{13}\\ 1 \end{bmatrix}+\\
 &+e^{-\frac{11 \pi i}{13}}
  \theta \begin{bmatrix} \frac{11}{13}\\ 1 \end{bmatrix}
  \theta \begin{bmatrix} \frac{5}{13}\\ 1 \end{bmatrix}^2
  \theta \begin{bmatrix} \frac{-3}{13}\\ 1 \end{bmatrix}
 =0.
\endaligned$$
Note that
$$\theta \begin{bmatrix} \frac{-1}{13}\\ 1 \end{bmatrix}
 =e^{-\frac{\pi i}{13}} \theta \begin{bmatrix} \frac{1}{13}\\ 1 \end{bmatrix}, \quad
  \theta \begin{bmatrix} \frac{-3}{13}\\ 1 \end{bmatrix}
 =e^{-\frac{3 \pi i}{13}} \theta \begin{bmatrix} \frac{3}{13}\\ 1 \end{bmatrix}.$$
This proves that
$$\aligned
 &\mathbf{B}_3^{(2)}(a_1(z), \ldots, a_6(z))\\
=&a_1(z) a_6(z) a_5(z)^2-a_2(z)^2 a_5(z) a_6(z)-a_3(z)^2 a_1(z) a_4(z)=0.
\endaligned$$

  Let $k=13$ and $(j; \alpha_1, \alpha_2, \alpha_3)=(3; 1, 1, 11)$.
Then $13-4j=1$, $13-2j=7$ and
$$\theta \begin{bmatrix} \frac{1}{13}\\ 1 \end{bmatrix}^3
  \theta \begin{bmatrix} \frac{11}{13}\\ 1 \end{bmatrix}
 +e^{\frac{\pi i}{13}}
  \theta \begin{bmatrix} \frac{7}{13}\\ 1 \end{bmatrix}^3
  \theta \begin{bmatrix} \frac{17}{13}\\ 1 \end{bmatrix}
 +e^{-\frac{7 \pi i}{13}}
  \theta \begin{bmatrix} \frac{7}{13}\\ 1 \end{bmatrix}
  \theta \begin{bmatrix} \frac{-5}{13}\\ 1 \end{bmatrix}^2
  \theta \begin{bmatrix} \frac{5}{13}\\ 1 \end{bmatrix}
 =0.$$
Note that
$$\theta \begin{bmatrix} \frac{-5}{13}\\ 1 \end{bmatrix}
 =e^{-\frac{5 \pi i}{13}} \theta \begin{bmatrix} \frac{5}{13}\\ 1 \end{bmatrix}, \quad
  \theta \begin{bmatrix} \frac{17}{13}\\ 1 \end{bmatrix}
 =e^{-\frac{9 \pi i}{13}} \theta \begin{bmatrix} \frac{9}{13}\\ 1 \end{bmatrix}.$$
This proves that
$$\mathbf{B}_9^{(1)}(a_1(z), \ldots, a_6(z))
 =a_1(z) a_6(z)^3+a_2(z)^3 a_5(z)-a_2(z) a_3(z)^3=0.$$

  Let $k=13$ and $(j; \alpha_1, \alpha_2, \alpha_3)=(2; 7, 7, -1)$.
Then $13-4j=5$, $13-2j=9$ and
$$\aligned
 &\theta \begin{bmatrix} \frac{5}{13}\\ 1 \end{bmatrix}
  \theta \begin{bmatrix} \frac{7}{13}\\ 1 \end{bmatrix}^2
  \theta \begin{bmatrix} \frac{-1}{13}\\ 1 \end{bmatrix}
 +e^{\frac{5 \pi i}{13}}
  \theta \begin{bmatrix} \frac{9}{13}\\ 1 \end{bmatrix}
  \theta \begin{bmatrix} \frac{11}{13}\\ 1 \end{bmatrix}^2
  \theta \begin{bmatrix} \frac{3}{13}\\ 1 \end{bmatrix}+\\
 &+e^{-\frac{9 \pi i}{13}}
  \theta \begin{bmatrix} \frac{9}{13}\\ 1 \end{bmatrix}
  \theta \begin{bmatrix} \frac{3}{13}\\ 1 \end{bmatrix}^2
  \theta \begin{bmatrix} \frac{-5}{13}\\ 1 \end{bmatrix}
 =0.
\endaligned$$
Note that
$$\theta \begin{bmatrix} \frac{-1}{13}\\ 1 \end{bmatrix}
 =e^{-\frac{\pi i}{13}} \theta \begin{bmatrix} \frac{1}{13}\\ 1 \end{bmatrix}, \quad
  \theta \begin{bmatrix} \frac{-5}{13}\\ 1 \end{bmatrix}
 =e^{-\frac{5 \pi i}{13}} \theta \begin{bmatrix} \frac{5}{13}\\ 1 \end{bmatrix}.$$
This proves that
$$\aligned
 &\mathbf{B}_9^{(2)}(a_1(z), \ldots, a_6(z))\\
=&a_3(z) a_5(z) a_4(z)^2-a_1(z)^2 a_4(z) a_5(z)-a_2(z)^2 a_3(z) a_6(z)=0.
\endaligned$$

  Let $k=13$ and $(j; \alpha_1, \alpha_2, \alpha_3)=(5; 7, 7, -1)$.
Then $13-4j=-7$, $13-2j=3$ and
$$\aligned
 &\theta \begin{bmatrix} \frac{-7}{13}\\ 1 \end{bmatrix}
  \theta \begin{bmatrix} \frac{7}{13}\\ 1 \end{bmatrix}^2
  \theta \begin{bmatrix} \frac{-1}{13}\\ 1 \end{bmatrix}
 +e^{\frac{-7 \pi i}{13}}
  \theta \begin{bmatrix} \frac{3}{13}\\ 1 \end{bmatrix}
  \theta \begin{bmatrix} \frac{17}{13}\\ 1 \end{bmatrix}^2
  \theta \begin{bmatrix} \frac{9}{13}\\ 1 \end{bmatrix}+\\
 &+e^{-\frac{3 \pi i}{13}}
  \theta \begin{bmatrix} \frac{3}{13}\\ 1 \end{bmatrix}
  \theta \begin{bmatrix} \frac{-3}{13}\\ 1 \end{bmatrix}^2
  \theta \begin{bmatrix} \frac{-11}{13}\\ 1 \end{bmatrix}
 =0.
\endaligned$$
Note that
$$\theta \begin{bmatrix} \frac{-1}{13}\\ 1 \end{bmatrix}
 =e^{-\frac{\pi i}{13}} \theta \begin{bmatrix} \frac{1}{13}\\ 1 \end{bmatrix}, \quad
  \theta \begin{bmatrix} \frac{-3}{13}\\ 1 \end{bmatrix}
 =e^{-\frac{3 \pi i}{13}} \theta \begin{bmatrix} \frac{3}{13}\\ 1 \end{bmatrix}, \quad
  \theta \begin{bmatrix} \frac{-7}{13}\\ 1 \end{bmatrix}
 =e^{-\frac{7 \pi i}{13}} \theta \begin{bmatrix} \frac{7}{13}\\ 1 \end{bmatrix},$$
$$\theta \begin{bmatrix} \frac{-11}{13}\\ 1 \end{bmatrix}
 =e^{-\frac{11 \pi i}{13}} \theta \begin{bmatrix} \frac{11}{13}\\ 1 \end{bmatrix}, \quad
  \theta \begin{bmatrix} \frac{17}{13}\\ 1 \end{bmatrix}
 =e^{-\frac{9 \pi i}{13}} \theta \begin{bmatrix} \frac{9}{13}\\ 1 \end{bmatrix}.$$
This proves that
$$\mathbf{B}_{12}^{(1)}(a_1(z), \ldots, a_6(z))
 =a_1(z) a_4(z)^3+a_2(z)^3 a_6(z)+a_4(z) a_5(z)^3=0.$$

  Let $k=13$ and $(j; \alpha_1, \alpha_2, \alpha_3)=(2; 1, 1, 11)$.
Then $13-4j=5$, $13-2j=9$ and
$$\aligned
 &\theta \begin{bmatrix} \frac{5}{13}\\ 1 \end{bmatrix}
  \theta \begin{bmatrix} \frac{1}{13}\\ 1 \end{bmatrix}^2
  \theta \begin{bmatrix} \frac{11}{13}\\ 1 \end{bmatrix}
 +e^{\frac{5 \pi i}{13}}
  \theta \begin{bmatrix} \frac{9}{13}\\ 1 \end{bmatrix}
  \theta \begin{bmatrix} \frac{5}{13}\\ 1 \end{bmatrix}^2
  \theta \begin{bmatrix} \frac{15}{13}\\ 1 \end{bmatrix}+\\
 &+e^{-\frac{9 \pi i}{13}}
  \theta \begin{bmatrix} \frac{9}{13}\\ 1 \end{bmatrix}
  \theta \begin{bmatrix} \frac{-3}{13}\\ 1 \end{bmatrix}^2
  \theta \begin{bmatrix} \frac{7}{13}\\ 1 \end{bmatrix}
 =0.
\endaligned$$
Note that
$$\theta \begin{bmatrix} \frac{-3}{13}\\ 1 \end{bmatrix}
 =e^{-\frac{3 \pi i}{13}} \theta \begin{bmatrix} \frac{3}{13}\\ 1 \end{bmatrix}, \quad
  \theta \begin{bmatrix} \frac{15}{13}\\ 1 \end{bmatrix}
 =e^{-\frac{11 \pi i}{13}} \theta \begin{bmatrix} \frac{11}{13}\\ 1 \end{bmatrix}.$$
This proves that
$$\aligned
 &\mathbf{B}_{12}^{(2)}(a_1(z), \ldots, a_6(z))\\
=&a_2(z) a_5(z) a_4(z)^2-a_3(z)^2 a_1(z) a_5(z)-
  a_6(z)^2 a_1(z) a_3(z)=0.
\endaligned$$

  Let $k=13$ and $(j; \alpha_1, \alpha_2, \alpha_3)=(6; 11, 11, -9)$.
Then $13-4j=-11$, $13-2j=1$ and
$$\aligned
 &\theta \begin{bmatrix} \frac{-11}{13}\\ 1 \end{bmatrix}
  \theta \begin{bmatrix} \frac{11}{13}\\ 1 \end{bmatrix}^2
  \theta \begin{bmatrix} \frac{-9}{13}\\ 1 \end{bmatrix}
 +e^{-\frac{11 \pi i}{13}}
  \theta \begin{bmatrix} \frac{1}{13}\\ 1 \end{bmatrix}
  \theta \begin{bmatrix} \frac{23}{13}\\ 1 \end{bmatrix}^2
  \theta \begin{bmatrix} \frac{3}{13}\\ 1 \end{bmatrix}+\\
 &+e^{-\frac{\pi i}{13}}
  \theta \begin{bmatrix} \frac{1}{13}\\ 1 \end{bmatrix}
  \theta \begin{bmatrix} \frac{-1}{13}\\ 1 \end{bmatrix}^2
  \theta \begin{bmatrix} \frac{-21}{13}\\ 1 \end{bmatrix}
 =0.
\endaligned$$
Note that
$$\theta \begin{bmatrix} \frac{-1}{13}\\ 1 \end{bmatrix}
 =e^{-\frac{\pi i}{13}} \theta \begin{bmatrix} \frac{1}{13}\\ 1 \end{bmatrix}, \quad
  \theta \begin{bmatrix} \frac{-9}{13}\\ 1 \end{bmatrix}
 =e^{-\frac{9 \pi i}{13}} \theta \begin{bmatrix} \frac{9}{13}\\ 1 \end{bmatrix}, \quad
  \theta \begin{bmatrix} \frac{-11}{13}\\ 1 \end{bmatrix}
 =e^{-\frac{11 \pi i}{13}} \theta \begin{bmatrix} \frac{11}{13}\\ 1 \end{bmatrix},
$$
$$\theta \begin{bmatrix} \frac{-21}{13}\\ 1 \end{bmatrix}
 =\theta \begin{bmatrix} \frac{5}{13}\\ 1 \end{bmatrix}, \quad
  \theta \begin{bmatrix} \frac{23}{13}\\ 1 \end{bmatrix}
 =e^{-\frac{3 \pi i}{13}} \theta \begin{bmatrix} \frac{3}{13}\\ 1 \end{bmatrix}.$$
This proves that
$$\mathbf{B}_{10}^{(1)}(a_1(z), \ldots, a_6(z))
 =a_3(z) a_6(z)^3+a_1(z)^3 a_5(z)+a_6(z) a_4(z)^3=0.$$

  Let $k=13$ and $(j; \alpha_1, \alpha_2, \alpha_3)=(4; 1, 1, 11)$.
Then $13-4j=-3$, $13-2j=5$ and
$$\aligned
 &\theta \begin{bmatrix} \frac{-3}{13}\\ 1 \end{bmatrix}
  \theta \begin{bmatrix} \frac{1}{13}\\ 1 \end{bmatrix}^2
  \theta \begin{bmatrix} \frac{11}{13}\\ 1 \end{bmatrix}
 +e^{-\frac{3 \pi i}{13}}
  \theta \begin{bmatrix} \frac{5}{13}\\ 1 \end{bmatrix}
  \theta \begin{bmatrix} \frac{9}{13}\\ 1 \end{bmatrix}^2
  \theta \begin{bmatrix} \frac{19}{13}\\ 1 \end{bmatrix}+\\
 &+e^{-\frac{5 \pi i}{13}}
  \theta \begin{bmatrix} \frac{5}{13}\\ 1 \end{bmatrix}
  \theta \begin{bmatrix} \frac{-7}{13}\\ 1 \end{bmatrix}^2
  \theta \begin{bmatrix} \frac{3}{13}\\ 1 \end{bmatrix}
 =0.
\endaligned$$
Note that
$$\theta \begin{bmatrix} \frac{-3}{13}\\ 1 \end{bmatrix}
 =e^{-\frac{3 \pi i}{13}} \theta \begin{bmatrix} \frac{3}{13}\\ 1 \end{bmatrix}, \quad
  \theta \begin{bmatrix} \frac{19}{13}\\ 1 \end{bmatrix}
 =\theta \begin{bmatrix} \frac{-7}{13}\\ 1 \end{bmatrix}
 =e^{-\frac{7 \pi i}{13}} \theta \begin{bmatrix} \frac{7}{13}\\ 1 \end{bmatrix}.$$
This proves that
$$\aligned
 &\mathbf{B}_{10}^{(2)}(a_1(z), \ldots, a_6(z))\\
=&a_1(z) a_4(z) a_6(z)^2-a_2(z)^2 a_3(z) a_4(z)-a_5(z)^2 a_2(z) a_3(z)=0.
\endaligned$$

  Let $k=13$ and $(j; \alpha_1, \alpha_2, \alpha_3)=(2; 5, 5, 3)$.
Then $13-4j=5$, $13-2j=9$ and
$$\theta \begin{bmatrix} \frac{5}{13}\\ 1 \end{bmatrix}^3
  \theta \begin{bmatrix} \frac{3}{13}\\ 1 \end{bmatrix}
 +e^{\frac{5 \pi i}{13}}
  \theta \begin{bmatrix} \frac{9}{13}\\ 1 \end{bmatrix}^3
  \theta \begin{bmatrix} \frac{7}{13}\\ 1 \end{bmatrix}
 +e^{-\frac{9 \pi i}{13}}
  \theta \begin{bmatrix} \frac{9}{13}\\ 1 \end{bmatrix}
  \theta \begin{bmatrix} \frac{1}{13}\\ 1 \end{bmatrix}^2
  \theta \begin{bmatrix} \frac{-1}{13}\\ 1 \end{bmatrix}
=0.$$
Note that
$$\theta \begin{bmatrix} \frac{-1}{13}\\ 1 \end{bmatrix}
 =e^{-\frac{\pi i}{13}} \theta \begin{bmatrix} \frac{1}{13}\\ 1 \end{bmatrix}.$$
This proves that
$$\mathbf{B}_4^{(1)}(a_1(z), \ldots, a_6(z))
 =a_2(z) a_5(z)^3+a_3(z)^3 a_4(z)+a_5(z) a_6(z)^3=0.$$

  Let $k=13$ and $(j; \alpha_1, \alpha_2, \alpha_3)=(1; 5, -1, 9)$.
Then $13-4j=9$, $13-2j=11$ and
$$\aligned
 &\theta \begin{bmatrix} \frac{9}{13}\\ 1 \end{bmatrix}^2
  \theta \begin{bmatrix} \frac{5}{13}\\ 1 \end{bmatrix}
  \theta \begin{bmatrix} \frac{-1}{13}\\ 1 \end{bmatrix}
 +e^{\frac{9 \pi i}{13}}
  \theta \begin{bmatrix} \frac{11}{13}\\ 1 \end{bmatrix}^2
  \theta \begin{bmatrix} \frac{7}{13}\\ 1 \end{bmatrix}
  \theta \begin{bmatrix} \frac{1}{13}\\ 1 \end{bmatrix}+\\
 &+e^{-\frac{11 \pi i}{13}}
  \theta \begin{bmatrix} \frac{11}{13}\\ 1 \end{bmatrix}
  \theta \begin{bmatrix} \frac{3}{13}\\ 1 \end{bmatrix}
  \theta \begin{bmatrix} \frac{-3}{13}\\ 1 \end{bmatrix}
  \theta \begin{bmatrix} \frac{7}{13}\\ 1 \end{bmatrix}
 =0.
\endaligned$$
Note that
$$\theta \begin{bmatrix} \frac{-1}{13}\\ 1 \end{bmatrix}
 =e^{-\frac{\pi i}{13}} \theta \begin{bmatrix} \frac{1}{13}\\ 1 \end{bmatrix}, \quad
  \theta \begin{bmatrix} \frac{-3}{13}\\ 1 \end{bmatrix}
 =e^{-\frac{3 \pi i}{13}} \theta \begin{bmatrix} \frac{3}{13}\\ 1 \end{bmatrix}.$$
This proves that
$$\aligned
 &\mathbf{B}_4^{(2)}(a_1(z), \ldots, a_6(z))\\
=&a_3(z) a_6(z) a_5(z)^2-a_1(z)^2 a_2(z) a_6(z)-a_4(z)^2 a_1(z) a_2(z)=0.
\endaligned$$

  Let $k=13$ and $(j; \alpha_1, \alpha_2, \alpha_3)=(3; 3, 9, 1)$.
Then $13-4j=1$, $13-2j=7$ and
$$\aligned
 &\theta \begin{bmatrix} \frac{1}{13}\\ 1 \end{bmatrix}^2
  \theta \begin{bmatrix} \frac{3}{13}\\ 1 \end{bmatrix}
  \theta \begin{bmatrix} \frac{9}{13}\\ 1 \end{bmatrix}
 +e^{\frac{\pi i}{13}}
  \theta \begin{bmatrix} \frac{7}{13}\\ 1 \end{bmatrix}^2
  \theta \begin{bmatrix} \frac{9}{13}\\ 1 \end{bmatrix}
  \theta \begin{bmatrix} \frac{15}{13}\\ 1 \end{bmatrix}+\\
 &+e^{-\frac{7 \pi i}{13}}
  \theta \begin{bmatrix} \frac{7}{13}\\ 1 \end{bmatrix}
  \theta \begin{bmatrix} \frac{-3}{13}\\ 1 \end{bmatrix}
  \theta \begin{bmatrix} \frac{3}{13}\\ 1 \end{bmatrix}
  \theta \begin{bmatrix} \frac{-5}{13}\\ 1 \end{bmatrix}
 =0.
\endaligned$$
Note that
$$\theta \begin{bmatrix} \frac{-3}{13}\\ 1 \end{bmatrix}
 =e^{-\frac{3 \pi i}{13}} \theta \begin{bmatrix} \frac{3}{13}\\ 1 \end{bmatrix}, \quad
  \theta \begin{bmatrix} \frac{-5}{13}\\ 1 \end{bmatrix}
 =e^{-\frac{5 \pi i}{13}} \theta \begin{bmatrix} \frac{5}{13}\\ 1 \end{bmatrix}, \quad
  \theta \begin{bmatrix} \frac{15}{13}\\ 1 \end{bmatrix}
 =e^{-\frac{11 \pi i}{13}} \theta \begin{bmatrix} \frac{11}{13}\\ 1 \end{bmatrix}.$$
This proves that
$$\aligned
 &\mathbf{B}_5(a_1(z), \ldots, a_6(z))\\
=&-a_2(z)^2 a_1(z) a_5(z)+a_4(z) a_5(z) a_6(z)^2+a_2(z) a_3(z) a_4(z)^2=0.
\endaligned$$

  Let $k=13$ and $(j; \alpha_1, \alpha_2, \alpha_3)=(1; 1, 3, 9)$.
Then $13-4j=9$, $13-2j=11$ and
$$\aligned
 &\theta \begin{bmatrix} \frac{9}{13}\\ 1 \end{bmatrix}^2
  \theta \begin{bmatrix} \frac{1}{13}\\ 1 \end{bmatrix}
  \theta \begin{bmatrix} \frac{3}{13}\\ 1 \end{bmatrix}
 +e^{\frac{9 \pi i}{13}}
  \theta \begin{bmatrix} \frac{11}{13}\\ 1 \end{bmatrix}^2
  \theta \begin{bmatrix} \frac{3}{13}\\ 1 \end{bmatrix}
  \theta \begin{bmatrix} \frac{5}{13}\\ 1 \end{bmatrix}+\\
 &+e^{-\frac{11 \pi i}{13}}
  \theta \begin{bmatrix} \frac{11}{13}\\ 1 \end{bmatrix}
  \theta \begin{bmatrix} \frac{-1}{13}\\ 1 \end{bmatrix}
  \theta \begin{bmatrix} \frac{1}{13}\\ 1 \end{bmatrix}
  \theta \begin{bmatrix} \frac{7}{13}\\ 1 \end{bmatrix}
 =0.
\endaligned$$
Note that
$$\theta \begin{bmatrix} \frac{-1}{13}\\ 1 \end{bmatrix}
 =e^{-\frac{\pi i}{13}} \theta \begin{bmatrix} \frac{1}{13}\\ 1 \end{bmatrix}.$$
This proves that
$$\aligned
 &\mathbf{B}_2(a_1(z), \ldots, a_6(z))\\
=&-a_1(z)^2 a_3(z) a_4(z)+a_6(z) a_4(z) a_5(z)^2+a_1(z) a_2(z) a_6(z)^2=0.
\endaligned$$

  Let $k=13$ and $(j; \alpha_1, \alpha_2, \alpha_3)=(2; 5, 7, 1)$.
Then $13-4j=5$, $13-2j=9$ and
$$\aligned
 &\theta \begin{bmatrix} \frac{5}{13}\\ 1 \end{bmatrix}^2
  \theta \begin{bmatrix} \frac{7}{13}\\ 1 \end{bmatrix}
  \theta \begin{bmatrix} \frac{1}{13}\\ 1 \end{bmatrix}
 +e^{\frac{5 \pi i}{13}}
  \theta \begin{bmatrix} \frac{9}{13}\\ 1 \end{bmatrix}^2
  \theta \begin{bmatrix} \frac{11}{13}\\ 1 \end{bmatrix}
  \theta \begin{bmatrix} \frac{5}{13}\\ 1 \end{bmatrix}+\\
 &+e^{-\frac{9 \pi i}{13}}
  \theta \begin{bmatrix} \frac{9}{13}\\ 1 \end{bmatrix}
  \theta \begin{bmatrix} \frac{1}{13}\\ 1 \end{bmatrix}
  \theta \begin{bmatrix} \frac{3}{13}\\ 1 \end{bmatrix}
  \theta \begin{bmatrix} \frac{-3}{13}\\ 1 \end{bmatrix}
 =0.
\endaligned$$
Note that
$$\theta \begin{bmatrix} \frac{-3}{13}\\ 1 \end{bmatrix}
 =e^{-\frac{3 \pi i}{13}} \theta \begin{bmatrix} \frac{3}{13}\\ 1 \end{bmatrix}.$$
This proves that
$$\aligned
 &\mathbf{B}_6(a_1(z), \ldots, a_6(z))\\
=&-a_3(z)^2 a_2(z) a_6(z)+a_5(z) a_6(z) a_4(z)^2+a_3(z) a_1(z) a_5(z)^2=0.
\endaligned$$

  Let $k=13$ and $(j; \alpha_1, \alpha_2, \alpha_3)=(2; 11, 7, -5)$.
Then $13-4j=5$, $13-2j=9$ and
$$\aligned
 &\theta \begin{bmatrix} \frac{5}{13}\\ 1 \end{bmatrix}
  \theta \begin{bmatrix} \frac{11}{13}\\ 1 \end{bmatrix}
  \theta \begin{bmatrix} \frac{7}{13}\\ 1 \end{bmatrix}
  \theta \begin{bmatrix} \frac{-5}{13}\\ 1 \end{bmatrix}
 +e^{\frac{5 \pi i}{13}}
  \theta \begin{bmatrix} \frac{9}{13}\\ 1 \end{bmatrix}
  \theta \begin{bmatrix} \frac{15}{13}\\ 1 \end{bmatrix}
  \theta \begin{bmatrix} \frac{11}{13}\\ 1 \end{bmatrix}
  \theta \begin{bmatrix} \frac{-1}{13}\\ 1 \end{bmatrix}+\\
 &+e^{-\frac{9 \pi i}{13}}
  \theta \begin{bmatrix} \frac{9}{13}\\ 1 \end{bmatrix}
  \theta \begin{bmatrix} \frac{7}{13}\\ 1 \end{bmatrix}
  \theta \begin{bmatrix} \frac{3}{13}\\ 1 \end{bmatrix}
  \theta \begin{bmatrix} \frac{-9}{13}\\ 1 \end{bmatrix}
 =0.
\endaligned$$
Note that
$$\theta \begin{bmatrix} \frac{-1}{13}\\ 1 \end{bmatrix}
 =e^{-\frac{\pi i}{13}} \theta \begin{bmatrix} \frac{1}{13}\\ 1 \end{bmatrix}, \quad
  \theta \begin{bmatrix} \frac{-5}{13}\\ 1 \end{bmatrix}
 =e^{-\frac{5 \pi i}{13}} \theta \begin{bmatrix} \frac{5}{13}\\ 1 \end{bmatrix}, \quad
  \theta \begin{bmatrix} \frac{-9}{13}\\ 1 \end{bmatrix}
 =e^{-\frac{9 \pi i}{13}} \theta \begin{bmatrix} \frac{9}{13}\\ 1 \end{bmatrix},$$
$$\theta \begin{bmatrix} \frac{15}{13}\\ 1 \end{bmatrix}
 =e^{-\frac{11 \pi i}{13}} \theta \begin{bmatrix} \frac{11}{13}\\ 1 \end{bmatrix}.$$
This proves that
$$\aligned
 &\mathbf{B}_8(a_1(z), \ldots, a_6(z))\\
=&a_2(z) a_4(z) a_5(z)^2+a_1(z) a_2(z) a_3(z)^2+a_1(z)^2 a_5(z) a_6(z)=0.
\endaligned$$

  Let $k=13$ and $(j; \alpha_1, \alpha_2, \alpha_3)=(1; 5, 5, 3)$.
Then $13-4j=9$, $13-2j=11$ and
$$\aligned
 &\theta \begin{bmatrix} \frac{9}{13}\\ 1 \end{bmatrix}
  \theta \begin{bmatrix} \frac{5}{13}\\ 1 \end{bmatrix}^2
  \theta \begin{bmatrix} \frac{3}{13}\\ 1 \end{bmatrix}
 +e^{\frac{9 \pi i}{13}}
  \theta \begin{bmatrix} \frac{11}{13}\\ 1 \end{bmatrix}
  \theta \begin{bmatrix} \frac{7}{13}\\ 1 \end{bmatrix}^2
  \theta \begin{bmatrix} \frac{5}{13}\\ 1 \end{bmatrix}+\\
 &+e^{-\frac{11 \pi i}{13}}
  \theta \begin{bmatrix} \frac{11}{13}\\ 1 \end{bmatrix}
  \theta \begin{bmatrix} \frac{3}{13}\\ 1 \end{bmatrix}^2
  \theta \begin{bmatrix} \frac{1}{13}\\ 1 \end{bmatrix}
 =0.
\endaligned$$
This proves that
$$\aligned
 &\mathbf{B}_{11}(a_1(z), \ldots, a_6(z))\\
=&a_1(z) a_6(z) a_4(z)^2+a_3(z) a_1(z) a_2(z)^2+a_3(z)^2 a_4(z) a_5(z)=0.
\endaligned$$

  Let $k=13$ and $(j; \alpha_1, \alpha_2, \alpha_3)=(3; -1, 5, 9)$.
Then $13-4j=1$, $13-2j=7$ and
$$\aligned
 &\theta \begin{bmatrix} \frac{1}{13}\\ 1 \end{bmatrix}
  \theta \begin{bmatrix} \frac{-1}{13}\\ 1 \end{bmatrix}
  \theta \begin{bmatrix} \frac{5}{13}\\ 1 \end{bmatrix}
  \theta \begin{bmatrix} \frac{9}{13}\\ 1 \end{bmatrix}
 +e^{\frac{\pi i}{13}}
  \theta \begin{bmatrix} \frac{7}{13}\\ 1 \end{bmatrix}
  \theta \begin{bmatrix} \frac{5}{13}\\ 1 \end{bmatrix}
  \theta \begin{bmatrix} \frac{11}{13}\\ 1 \end{bmatrix}
  \theta \begin{bmatrix} \frac{15}{13}\\ 1 \end{bmatrix}+\\
 &+e^{-\frac{7 \pi i}{13}}
  \theta \begin{bmatrix} \frac{7}{13}\\ 1 \end{bmatrix}
  \theta \begin{bmatrix} \frac{-7}{13}\\ 1 \end{bmatrix}
  \theta \begin{bmatrix} \frac{-1}{13}\\ 1 \end{bmatrix}
  \theta \begin{bmatrix} \frac{3}{13}\\ 1 \end{bmatrix}
 =0.
\endaligned$$
Note that
$$\theta \begin{bmatrix} \frac{-1}{13}\\ 1 \end{bmatrix}
 =e^{-\frac{\pi i}{13}} \theta \begin{bmatrix} \frac{1}{13}\\ 1 \end{bmatrix}, \quad
  \theta \begin{bmatrix} \frac{-7}{13}\\ 1 \end{bmatrix}
 =e^{-\frac{7 \pi i}{13}} \theta \begin{bmatrix} \frac{7}{13}\\ 1 \end{bmatrix}, \quad
  \theta \begin{bmatrix} \frac{15}{13}\\ 1 \end{bmatrix}
 =e^{-\frac{11 \pi i}{13}} \theta \begin{bmatrix} \frac{11}{13}\\ 1 \end{bmatrix}.$$
This proves that
$$\aligned
 &\mathbf{B}_7(a_1(z), \ldots, a_6(z))\\
=&a_3(z) a_5(z) a_6(z)^2+a_2(z) a_3(z) a_1(z)^2+a_2(z)^2 a_6(z) a_4(z)=0.
\endaligned$$

  By $\mathbf{B}_6(a_1(z), \ldots, a_6(z))=0$, we have
$$a_3(z) a_1(z) a_5(z)^2=a_3(z)^2 a_2(z) a_6(z)-a_5(z) a_6(z) a_4(z)^2.$$
Hence,
$$a_1(z) a_5(z)^3=a_2(z) a_3(z) a_5(z) a_6(z)-\frac{1}{a_3(z)} a_6(z) a_4(z)^2 z_5(z)^2.$$
By $\mathbf{B}_2(a_1(z), \ldots, a_6(z))=0$, we have
$$a_1(z) a_2(z) a_6(z)^2=a_1(z)^2 a_3(z) a_4(z)-a_6(z) a_4(z) a_5(z)^2.$$
Hence,
$$a_2(z) a_6(z)^3=a_3(z) a_1(z) a_6(z) a_4(z)-\frac{1}{a_1(z)} a_4(z) a_5(z)^2 z_6(z)^2.$$
By $\mathbf{B}_5(a_1(z), \ldots, a_6(z))=0$, we have
$$a_2(z) a_3(z) a_4(z)^2=a_2(z)^2 a_1(z) a_5(z)-a_4(z) a_5(z) a_6(z)^2.$$
Hence,
$$a_3(z) a_4(z)^3=a_1(z) a_2(z) a_4(z) a_5(z)-\frac{1}{a_2(z)} a_5(z) a_6(z)^2 z_4(z)^2.$$
Thus,
$$\aligned
 &a_1(z) a_5(z)^3+a_2(z) a_6(z)^3+a_3(z) a_4(z)^3
 =\left[1-\frac{a_4(z) a_5(z) a_6(z)}{a_1(z) a_2(z) a_3(z)}\right]\\
 &\times \left[a_1(z) a_2(z) a_4(z) a_5(z)+a_2(z) a_3(z) a_5(z) a_6(z)
         +a_3(z) a_1(z) a_6(z) a_4(z)\right].
\endaligned$$
Therefore, the four equations
$$\left\{\aligned
  \mathbf{B}_5(a_1(z), \ldots, a_6(z)) &=0,\\
  \mathbf{B}_2(a_1(z), \ldots, a_6(z)) &=0,\\
  \mathbf{B}_6(a_1(z), \ldots, a_6(z)) &=0,\\
  \mathbf{B}_0^{(0)}(a_1(z), \ldots, a_6(z)) &=0.
\endaligned\right.$$
imply that $\mathbf{B}_0^{(1)}(a_1(z), \ldots, a_6(z))=0$.

  By $\mathbf{B}_8(a_1(z), \ldots, a_6(z))=0$, we have
$$a_1(z)^2 a_5(z) a_6(z)=-a_2(z) a_4(z) a_5(z)^2-a_1(z) a_2(z) a_3(z)^2.$$
Hence,
$$a_1(z)^3 a_6(z)=-a_1(z) a_2(z) a_4(z) a_5(z)-\frac{1}{a_5(z)} a_2(z) a_3(z)^2 a_1(z)^2.$$
By $\mathbf{B}_7(a_1(z), \ldots, a_6(z))=0$, we have
$$a_2(z)^2 a_6(z) a_4(z)=-a_3(z) a_5(z) a_6(z)^2-a_2(z) a_3(z) a_1(z)^2.$$
Hence,
$$a_2(z)^3 a_4(z)=-a_2(z) a_3(z) a_5(z) a_6(z)-\frac{1}{a_6(z)} a_3(z) a_1(z)^2 a_2(z)^2.$$
By $\mathbf{B}_{11}(a_1(z), \ldots, a_6(z))=0$, we have
$$a_3(z)^2 a_4(z) a_5(z)=-a_1(z) a_6(z) a_4(z)^2-a_3(z) a_1(z) a_2(z)^2.$$
Hence,
$$a_3(z)^3 a_5(z)=-a_3(z) a_1(z) a_6(z) a_4(z)-\frac{1}{a_4(z)} a_1(z) a_2(z)^2 a_3(z)^2.$$
Thus,
$$\aligned
 &a_1(z)^3 a_6(z)+a_2(z)^3 a_4(z)+a_3(z)^3 a_5(z)
 =-\left[1+\frac{a_1(z) a_2(z) a_3(z)}{a_4(z) a_5(z) a_6(z)}\right]\\
 &\times \left[a_1(z) a_2(z) a_4(z) a_5(z)+a_2(z) a_3(z) a_5(z) a_6(z)
         +a_3(z) a_1(z) a_6(z) a_4(z)\right].
\endaligned$$
Therefore, the four equations
$$\left\{\aligned
  \mathbf{B}_8(a_1(z), \ldots, a_6(z)) &=0,\\
  \mathbf{B}_{11}(a_1(z), \ldots, a_6(z)) &=0,\\
  \mathbf{B}_7(a_1(z), \ldots, a_6(z)) &=0, \\
  \mathbf{B}_0^{(0)}(a_1(z), \ldots, a_6(z)) &=0.
\endaligned\right.$$
imply that $\mathbf{B}_0^{(2)}(a_1(z), \ldots, a_6(z))=0$.
This completes the proof of Theorem 7.6.

\noindent
$\qquad \qquad \qquad \qquad \qquad \qquad \qquad \qquad \qquad
 \qquad \qquad \qquad \qquad \qquad \qquad \qquad \boxed{}$

\textbf{Corollary 7.7.} {\it The space curve $Y$ is modular, i.e., it
is parameterized by theta constants of order $13$. Two curves $Y_2$ and
$Y_3$ lying over $Y$ have modular components given by (7.14).}

  Corollary 7.7 provides an example of modularity for higher genus space
curve $(g \geq 4)$ by means of the explicit modular parametrization.

\textbf{Corollary 7.8.} {\it The space curve $Y$ admits a hyperbolic
uniformization of arithmetic type with respect to the congruence
subgroup $\Gamma(13)$.}

  Corollary 7.8 gives an example of hyperbolic uniformization of arithmetic
type with respect to a congruence subgroup of $\text{SL}(2, \mathbb{Z})$ as
well as an example of the explicit uniformization of algebraic space curves
of higher genus.

\textbf{Corollary 7.9.} {\it The $21$ quartic polynomials in six variables
given by (7.14) form a system of over-determined algebraic equations which
can be uniformized by theta constants of order $13$.}

  Corollary 7.9 gives a new solution to Hilbert's 22nd problem, i.e., an
explicit uniformization for algebraic relations among five variables by means
of automorphic functions.

\textbf{Theorem 7.10.} {\it Besides (7.13), there are twenty formulas
given by (7.14) in the spirits of (7.13).}

  Theorem 7.10 greatly improves the result of Ramanujan and Evans on
modular equations of order $13$.

   Now, we study the geometric properties of $Y$. We will prove that at
each point $x$ of $Y \subset \mathbb{P}^5$, after passing to one of the
affine charts $z_i \neq 0$ containing $x$, and making $z_i=1$ to have affine
coordinate, they are $4$ quartic polynomials $F_j$ in the $21$-dimensional
space such that the $dF_j$ are linearly independent at $x$.

\noindent (1) In the affine coordinates $z_6=1$, we have
$$\left\{\aligned
  \mathbf{B}_0^{(1)} &=z_1 z_5^3+z_2+z_3 z_4^3,\\
  \mathbf{B}_9^{(1)} &=z_1+z_2^3 z_5-z_2 z_3^3,\\
  \mathbf{B}_4^{(1)} &=z_2 z_5^3+z_3^3 z_4+z_5,\\
  \mathbf{B}_{10}^{(1)} &=z_3+z_1^3 z_5+z_4^3.
\endaligned\right.\eqno{(7.16)}$$
Hence,
$$\left\{\aligned
  d \mathbf{B}_0^{(1)} &=z_5^3 dz_1+dz_2+z_4^3 dz_3+3 z_3 z_4^2 dz_4
                        +3 z_1 z_5^2 dz_5 \neq 0,\\
  d \mathbf{B}_9^{(1)} &=dz_1+(3 z_2^2 z_5-z_3^3) dz_2-3 z_2 z_3^2 dz_3
                        +z_2^3 dz_5 \neq 0,\\
  d \mathbf{B}_4^{(1)} &=z_5^3 dz_2+3 z_3^2 z_4 dz_3+z_3^3 dz_4+(1+3 z_2 z_5^2)
                         dz_5 \neq 0,\\
  d \mathbf{B}_{10}^{(1)} &=3 z_1^2 z_5 dz_1+dz_3+3 z_4^2 dz_4+z_1^3 dz_5
                            \neq 0.
\endaligned\right.\eqno{(7.17)}$$
Suppose that there exist $a$, $b$, $c$ and $d$, such that
$$a \cdot d \mathbf{B}_0^{(1)}+b \cdot d \mathbf{B}_9^{(1)}+
  c \cdot d \mathbf{B}_{10}^{(1)}+d \cdot d \mathbf{B}_4^{(1)}=0$$
for each $(z_1, z_2, z_3, z_4, z_5) \in Y$. Then
$$\aligned
  &a [z_5^3 dz_1+dz_2+z_4^3 dz_3+3 z_3 z_4^2 dz_4+3 z_1 z_5^2 dz_5]+\\
 +&b [dz_1+(3 z_2^2 z_5-z_3^3) dz_2-3 z_2 z_3^2 dz_3+z_2^3 dz_5]+\\
 +&c [3 z_1^2 z_5 dz_1+dz_3+3 z_4^2 dz_4+z_1^3 dz_5]+\\
 +&d [z_5^3 dz_2+3 z_3^2 z_4 dz_3+z_3^3 dz_4+(1+3 z_2 z_5^2) dz_5]=0.
\endaligned$$
This gives that
$$\left\{\aligned
  a z_5^3+b+c \cdot 3 z_1^2 z_5 &=0,\\
  a+b (3 z_2^2 z_5-z_3^3)+ d z_5^3 &=0,\\
  a z_4^3+b \cdot (-3 z_2 z_3^2)+c+d \cdot 3 z_3^2 z_4 &=0,\\
  a \cdot 3 z_3 z_4^2+c \cdot 3 z_4^2+d \cdot z_3^3 &=0,\\
  a \cdot 3 z_1 z_5^2+b z_2^3+c z_1^3+d(1+3 z_2 z_5^2) &=0.
\endaligned\right.$$
Let $(z_1, z_2, z_3, z_4, z_5)=(0, 0, 0, 0, 0) \in Y$. We have
$b=0$, $a=0$, $c=0$ and $d=0$. This implies that
$d \mathbf{B}_0^{(1)}$, $d \mathbf{B}_9^{(1)}$, $d \mathbf{B}_4^{(1)}$
and $d \mathbf{B}_{10}^{(1)}$ are linearly independent.

\noindent (2) In the affine coordinates $z_5=1$, we have
$$\left\{\aligned
  \mathbf{B}_0^{(1)} &=z_1+z_2 z_6^3+z_3 z_4^3,\\
  \mathbf{B}_1^{(1)} &=z_3+z_1^3 z_4-z_1 z_2^3,\\
  \mathbf{B}_{12}^{(1)} &=z_1 z_4^3+z_2^3 z_6+z_4,\\
  \mathbf{B}_4^{(1)} &=z_2+z_3^3 z_4+z_6^3.
\endaligned\right.\eqno{(7.18)}$$
Hence,
$$\left\{\aligned
  d \mathbf{B}_0^{(1)} &=dz_1+z_6^3 dz_2+z_4^3 dz_3+3 z_3 z_4^2 dz_4
                        +3 z_2 z_6^2 dz_6 \neq 0,\\
  d \mathbf{B}_1^{(1)} &=(3 z_1^2 z_4-z_2^3) dz_1-3 z_1 z_2^2 dz_2+dz_3
                        +z_1^3 dz_4 \neq 0,\\
  d \mathbf{B}_{12}^{(1)} &=z_4^3 dz_1+3 z_2^2 z_6 dz_2+(3 z_1 z_4^2+1) dz_4
                           +z_2^3 dz_6 \neq 0,\\
  d \mathbf{B}_4^{(1)} &=dz_2+3 z_3^2 z_4 dz_3+z_3^3 dz_4+3 z_6^2 dz_6 \neq 0.
\endaligned\right.\eqno{(7.19)}$$
Similarly, we have that $d \mathbf{B}_0^{(1)}$, $d \mathbf{B}_1^{(1)}$,
$d \mathbf{B}_{12}^{(1)}$ and $d \mathbf{B}_4^{(1)}$ are linearly independent.

\noindent (3) In the affine coordinates $z_4=1$, we have
$$\left\{\aligned
  \mathbf{B}_0^{(1)} &=z_1 z_5^3+z_2 z_6^3+z_3,\\
  \mathbf{B}_3^{(1)} &=z_2+z_3^3 z_6-z_3 z_1^3,\\
  \mathbf{B}_{10}^{(1)} &=z_3 z_6^3+z_1^3 z_5+z_6,\\
  \mathbf{B}_{12}^{(1)} &=z_1+z_2^3 z_6+z_5^3.
\endaligned\right.\eqno{(7.20)}$$
Hence,
$$\left\{\aligned
  d \mathbf{B}_0^{(1)} &=z_5^3 dz_1+z_6^3 dz_2+dz_3+3 z_1 z_5^2 dz_5
                        +3 z_2 z_6^2 dz_6 \neq 0,\\
  d \mathbf{B}_3^{(1)} &=-3 z_3 z_1^2 dz_1+dz_2+(3 z_3^2 z_6-z_1^3) dz_3
                        +z_3^3 dz_6 \neq 0,\\
  d \mathbf{B}_{10}^{(1)} &=3 z_1^2 z_5 dz_1+z_6^3 dz_3+z_1^3 dz_5+(1+3 z_3
                            z_6^2) dz_6 \neq 0,\\
  d \mathbf{B}_{12}^{(1)} &=dz_1+3 z_2^2 z_6 dz_2+3 z_5^2 dz_5+z_2^3 dz_6
                            \neq 0.
\endaligned\right.\eqno{(7.21)}$$
Similarly, we have that $d \mathbf{B}_0^{(1)}$, $d \mathbf{B}_3^{(1)}$,
$d \mathbf{B}_{10}^{(1)}$ and $d \mathbf{B}_{12}^{(1)}$ are linearly
independent.

\noindent (4) In the affine coordinates $z_3=1$, we have
$$\left\{\aligned
  \mathbf{B}_0^{(2)} &=z_1^3 z_6+z_2^3 z_4+z_5,\\
  \mathbf{B}_9^{(1)} &=z_1 z_6^3+z_2^3 z_5-z_2,\\
  \mathbf{B}_4^{(1)} &=z_2 z_5^3+z_4+z_5 z_6^3,\\
  \mathbf{B}_3^{(1)} &=z_2 z_4^3+z_6-z_1^3.
\endaligned\right.\eqno{(7.22)}$$
Hence,
$$\left\{\aligned
  d \mathbf{B}_0^{(2)} &=3 z_1^2 z_6 dz_1+3 z_2^2 z_4 dz_2+z_2^3 dz_4+dz_5
                         +z_1^3 dz_6 \neq 0,\\
  d \mathbf{B}_9^{(1)} &=z_6^3 dz_1+(3 z_2^2 z_5-1) dz_2+z_2^3 dz_5+3 z_1
                         z_6^2 dz_6 \neq 0,\\
  d \mathbf{B}_4^{(1)} &=z_5^3 dz_2+dz_4+(3 z_2 z_5^2+z_6^3) dz_5+3 z_5
                         z_6^2 dz_6 \neq 0,\\
  d \mathbf{B}_3^{(1)} &=-3 z_1^2 dz_1+z_4^3 dz_2+3 z_2 z_4^2 dz_4+dz_6
                         \neq 0.
\endaligned\right.\eqno{(7.23)}$$
Similarly, we have that $d \mathbf{B}_0^{(2)}$, $d \mathbf{B}_9^{(1)}$,
$d \mathbf{B}_4^{(1)}$ and $d \mathbf{B}_3^{(1)}$ are linearly independent.

\noindent (5) In the affine coordinates $z_2=1$, we have
$$\left\{\aligned
  \mathbf{B}_0^{(2)} &=z_1^3 z_6+z_4+z_3^3 z_5,\\
  \mathbf{B}_1^{(1)} &=z_3 z_5^3+z_1^3 z_4-z_1,\\
  \mathbf{B}_{12}^{(1)} &=z_1 z_4^3+z_6+z_4 z_5^3,\\
  \mathbf{B}_9^{(1)} &=z_1 z_6^3+z_5-z_3^3.
\endaligned\right.\eqno{(7.24)}$$
Hence,
$$\left\{\aligned
  d \mathbf{B}_0^{(2)} &=3 z_1^2 z_6 dz_1+3 z_3^2 z_5 dz_3+dz_4+z_3^3 dz_5
                        +z_1^3 dz_6 \neq 0,\\
  d \mathbf{B}_1^{(1)} &=(3 z_1^2 z_4-1) dz_1+z_5^3 dz_3+z_1^3 dz_4+3 z_3 z_5^2 dz_5
                         \neq 0,\\
  d \mathbf{B}_{12}^{(1)} &=z_4^3 dz_1+(3 z_1 z_4^2+z_5^3) dz_4+3 z_4 z_5^2 dz_5+dz_6
                            \neq 0,\\
  d \mathbf{B}_9^{(1)} &=z_6^3 dz_1-3 z_3^2 dz_3+dz_5+3 z_1 z_6^2 dz_6
                         \neq 0.
\endaligned\right.\eqno{(7.25)}$$
Similarly, we have that $d \mathbf{B}_0^{(2)}$, $d \mathbf{B}_1^{(1)}$,
$d \mathbf{B}_{12}^{(1)}$ and $d \mathbf{B}_9^{(1)}$ are linearly independent.

\noindent (6) In the affine coordinates $z_1=1$, we have
$$\left\{\aligned
  \mathbf{B}_0^{(2)} &=z_6+z_2^3 z_4+z_3^3 z_5,\\
  \mathbf{B}_3^{(1)} &=z_2 z_4^3+z_3^3 z_6-z_3,\\
  \mathbf{B}_{10}^{(1)} &=z_3 z_6^3+z_5+z_6 z_4^3,\\
  \mathbf{B}_1^{(1)} &=z_3 z_5^3+z_4-z_2^3.
\endaligned\right.\eqno{(7.26)}$$
Hence,
$$\left\{\aligned
  d \mathbf{B}_0^{(2)} &=3 z_2^2 z_4 dz_2+3 z_3^2 z_5 dz_3+z_2^3 dz_4+z_3^3 dz_5
                        +dz_6 \neq 0,\\
  d \mathbf{B}_3^{(1)} &=z_4^3 dz_2+(3 z_3^2 z_6-1) dz_3+3 z_2 z_4^2 dz_4+z_3^3 dz_6
                         \neq 0,\\
  d \mathbf{B}_{10}^{(1)} &=z_6^3 dz_3+3 z_6 z_4^2 dz_4+dz_5+(3 z_3 z_6^2+z_4^3) dz_6
                            \neq 0,\\
  d \mathbf{B}_1^{(1)} &=-3 z_2^2 dz_2+z_5^3 dz_3+dz_4+3 z_3 z_5^2 dz_5
                         \neq 0.
\endaligned\right.\eqno{(7.27)}$$
Similarly, we have that $d \mathbf{B}_0^{(2)}$, $d \mathbf{B}_3^{(1)}$,
$d \mathbf{B}_{10}^{(1)}$ and $d \mathbf{B}_1^{(1)}$ are linearly independent.

\textbf{Theorem 7.11.} (Ramanan) (see \cite{AR}, p.59, Theorem (20.7)).
{\it The locus $\mathcal{L}$ can be defined by
$\left(\begin{matrix} \frac{p-1}{2}\\ 3 \end{matrix}\right)$ quartics.}

  For $p=13$, this number is $20$, which is smaller than the number $21$
of distinct (hence linearly independent) quartics $\Phi_{abcd}$ as above.
As Ramanan has remarked (see \cite{AR}, p.59), one can actually reduce
the number further, namely to $(p-1)/2$. However, this does not lead to
explicit equations. Now, we give an explicit construction for $p=13$:

\textbf{Theorem 7.12}. {\it The locus $\mathcal{L}$ of the modular curve
$X(13)$ can be defined by six quartics which can be given explicitly.}

{\it Proof}. By the above argument, we have that if $z_6 \neq 0$, then
$d \mathbf{B}_0^{(1)}$, $d \mathbf{B}_9^{(1)}$, $d \mathbf{B}_4^{(1)}$
and $d \mathbf{B}_{10}^{(1)}$ are linearly independent. If $z_3 \neq 0$,
then $d \mathbf{B}_0^{(2)}$, $d \mathbf{B}_9^{(1)}$, $d \mathbf{B}_4^{(1)}$
and $d \mathbf{B}_3^{(1)}$ are linearly independent. Thus, the six quartics
$$(\mathbf{B}_0^{(1)}, \mathbf{B}_0^{(2)}, \mathbf{B}_3^{(1)},
   \mathbf{B}_{10}^{(1)}, \mathbf{B}_9^{(1)}, \mathbf{B}_4^{(1)})\eqno{(7.28)}$$
are sufficient to define the locus $\mathcal{L}$.

  Similarly, if $z_5 \neq 0$, then $d \mathbf{B}_0^{(1)}$, $d \mathbf{B}_1^{(1)}$,
$d \mathbf{B}_{12}^{(1)}$ and $d \mathbf{B}_4^{(1)}$ are linearly independent.
If $z_2 \neq 0$, then $d \mathbf{B}_0^{(2)}$, $d \mathbf{B}_1^{(1)}$,
$d \mathbf{B}_{12}^{(1)}$ and $d \mathbf{B}_9^{(1)}$ are linearly independent.
Thus, the six quartics
$$(\mathbf{B}_0^{(1)}, \mathbf{B}_0^{(2)}, \mathbf{B}_9^{(1)},
   \mathbf{B}_4^{(1)}, \mathbf{B}_1^{(1)}, \mathbf{B}_{12}^{(1)})\eqno{(7.29)}$$
are sufficient to define the locus $\mathcal{L}$.

  If $z_4 \neq 0$, then $d \mathbf{B}_0^{(1)}$, $d \mathbf{B}_3^{(1)}$,
$d \mathbf{B}_{10}^{(1)}$ and $d \mathbf{B}_{12}^{(1)}$ are linearly independent.
If $z_1 \neq 0$, then $d \mathbf{B}_0^{(2)}$, $d \mathbf{B}_3^{(1)}$,
$d \mathbf{B}_{10}^{(1)}$ and $d \mathbf{B}_1^{(1)}$ are linearly independent.
Thus, the six quartics
$$(\mathbf{B}_0^{(1)}, \mathbf{B}_0^{(2)}, \mathbf{B}_1^{(1)},
   \mathbf{B}_{12}^{(1)}, \mathbf{B}_3^{(1)}, \mathbf{B}_{10}^{(1)})\eqno{(7.30)}$$
are sufficient to define the locus $\mathcal{L}$.

  Note that $(z_1, z_2, z_3, z_4, z_5, z_6) \neq (0, 0, 0, 0, 0, 0)$.
This completes the proof of Theorem 7.12.

\noindent
$\qquad \qquad \qquad \qquad \qquad \qquad \qquad \qquad \qquad
 \qquad \qquad \qquad \qquad \qquad \qquad \qquad \boxed{}$

  In fact, Theorem 7.12 also gives a solution for Problem 7.4 when $p=13$.

  Now, we give the computation for the degree of the curve $Y$. There are
several ways to look at the degree of the curve $Y$. From the viewpoint
of Riemann surfaces, by the remark of Proposition 7.1,
$\Phi(\overline{\mathbb{H}/\Gamma(k)})$ is a curve of degree
$\frac{(k^2-1)(k-3)}{48}$ in $\mathbb{CP}^{\frac{k-3}{2}}$.
On the other hand, we can look at the geometry of the modular curve via the
fundamental relation (see \cite{AR}). In the same notation as in (4.12),
for every element $t$ of $K(\delta)=\mathbb{Z}/p \mathbb{Z}$, let $E_t$
denote the linear form on the space $V^{-}$ of odd functions
on $K(\delta)$ given by $E_t(h)=h(t)$.

\textbf{Theorem 7.13}. (see \cite{AR}, p.90, Corollary (23.28) and p.56,
Theorem (19.17)) {\it Let $\delta=(p)$. The curve part $\mathcal{C}$ of
the locus $\mathcal{L}$ defined by
$$\aligned
  0 &=h(u+v) \cdot h(u-v) \cdot h(w+z) \cdot h(w-z)+\\
    &+h(u+w) \cdot h(u-w) \cdot h(z+v) \cdot h(z-v)+\\
    &+h(u+z) \cdot h(u-z) \cdot h(v+w) \cdot h(v-w)
\endaligned\eqno{(7.31)}$$
for every $[h]$ in $\mathcal{L}$ and all $u$, $v$, $w$, $z \in K(\delta)$
has degree
$$\frac{(p-3)(p-1)(p+1)}{48}=\left(\begin{matrix}
  \frac{p+1}{2}\\ 3 \end{matrix}\right).\eqno{(7.32)}$$
The locus defined by (7.31) is the same as $\mathcal{L}$.}

  Therefore, $\Phi(\overline{\mathbb{H}/\Gamma(p)})$ and the locus
$\mathcal{L}$ of the modular curve $X(p)$ have the same degree
$\frac{(p^2-1)(p-3)}{48}$, which agrees with Klein's $z$-curve (4.2).

\begin{center}
{\large\bf 8. Galois covering $Y \rightarrow Y/\text{SL}(2, 13) \cong
              \mathbb{CP}^1$, Galois resolvent for the modular equation
              of level $13$ and a Hauptmodul for $\Gamma_0(13)$}
\end{center}

  Let us recall some basic facts about modular equations (see \cite{dSG}).
Given an integer $N \geq 2$, we seek an equation linking $j(\tau)$ and
$j^{\prime}(\tau)=j(N \tau)$ for $\tau \in \mathbb{H}$. It is easy to
check that $j^{\prime}$ is left invariant by the group
$$\Gamma_0(N)=\left\{ \left(\begin{matrix}
              a & b\\
              c & d
             \end{matrix}\right) \in \text{SL}(2, \mathbb{Z}):
             c \equiv 0 (\text{mod $N$}) \right\},$$
which is, in fact, precisely the stabilizer of $j^{\prime}$.

  On the other hand, $j^{\prime}$ is meromorphic at the cusps of $\Gamma_0(N)$.
Indeed, by means of the action of $\Gamma(1)=\text{SL}(2, \mathbb{Z})$, one
reduces the situation to the cusp $\infty$ and to a function of the form
$j \circ \left(\begin{matrix} a & b\\ c & d \end{matrix}\right)$ with $a$,
$b$ and $d$ integers; for sufficiently large $k$ the product of the latter
with $q^{k/m}$ is bounded in a neighborhood of $q^{1/m}=0$. The extension
$K(\Gamma_0(N))/\mathbb{C}(j)$ being finite, this implies the existence of
an algebraic relation between $j$ and $j^{\prime}$. In order to exhibit such
a relation, one considers the transforms of $j^{\prime}$ by the elements of
$\Gamma(1)$, that is, the $j \circ \alpha$ with $\alpha$ ranging over the
orbit $O_N$ of the point
$$p_N=\Gamma(1) \left(\begin{matrix} N & 0\\ 0 & 1 \end{matrix}\right)
      \in \Gamma(1) \backslash \Delta_N,$$
under the action of $\Gamma(1)$ on the right; here $\Delta_N$ denotes the
set of integer matrices of determinant $N$. One can check that the stabilizer
of the point $p_N$ in $\Gamma(1)$ is $\Gamma_0(N)$, so that the orbit $O_N$
may be identified with the quotient $\Gamma_0(N) \backslash \Gamma(1)$.
Denote $d_N$ for the index of $\Gamma_0(N)$ in $\Gamma(1)$ and
$\alpha_k \in \Delta_N$ $(k=1, \ldots, d_N)$ for a system of representatives
of the orbit $O_N$. Then the coefficients of the polynomial
$\prod_{k=1}^{d_N} (X-j \circ \alpha_k)$ are invariant under $\Gamma(1)$,
holomorphic on $\mathbb{H}$ and meromorphic at the cusp $\infty$. We have
thus found a polynomial $\Phi_N \in \mathbb{C}[X, Y]$ of degree $d_N$ in $X$
such that
$$\Phi_N(j^{\prime}, j)=0.\eqno{(8.1)}$$
This is the modular equation associated with transformations of order $N$.
The stabilizer of $j \circ \alpha_k$ is conjugate to $\Gamma_0(N)$ (the
stabilizer of $j^{\prime}$), whence the subgroup fixing all the $j \circ
\alpha_k$ coincides with $\Gamma(N)=\cap_{\gamma \in \Gamma(1)} \gamma
\Gamma_0(N) \gamma^{-1}$. It follows that the splitting field of
$\Phi_N \in \mathbb{C}[j][X]$ is $K(\Gamma(N))$. Moreover $\Gamma(1)$
acts as a set of automorphisms of $K(\Gamma(N))$ in permuting transitively
the roots of this polynomial, which is therefore irreducible, whence, in
particular, $K(\Gamma_0(N))=\mathbb{C}(j, j^{\prime})$. When $N=p$ is a
prime, one sees that the matrices $\left(\begin{matrix} 1 & k\\ 0 & p
\end{matrix}\right)$ $(0 \leq k < p)$ and $\left(\begin{matrix} p & 0\\
0 & 1 \end{matrix}\right)$ form a system of representatives of
$O_p=\Gamma(1) \backslash \Delta_p$; the index of $\Gamma_0(p)$ is thus
$d_p=p+1$.

  An elementary calculation shows that $\Gamma_0(N)$ is normalized by the
matrix
$$\left(\begin{matrix}
  0 & N^{-1/2}\\
  -N^{1/2} & 0
  \end{matrix}\right)$$
which induces an involutary automorphism of the modular curve $X_0(N)$ and
its function field: this is the Fricke involution interchanging $j$ and
$j^{\prime}$. One infers from its existence that $\Phi_N \in \mathbb{C}[X, Y]$
is symmetric. Klein relies on this symmetry in his investigation of the modular
equation for $N=2$, $3$, $4$, $5$, $7$ and $13$. For these values of $N$ the
modular curve $X_0(N)$ is of genus zero and there exists $\tau \in K(\Gamma_0(N))$
such that $K(\Gamma_0(N))=\mathbb{C}(j, j^{\prime})=\mathbb{C}(\tau)$; one then
has $j=F(\tau)$ and $j^{\prime}=F(\tau^{\prime})$ with $F \in \mathbb{C}(Z)$, the
function $\tau^{\prime}$ being linked to $\tau$ by the Fricke involution. In each
of these cases Klein describes a fundamental region for the action of $\Gamma_0(N)$
on the half-plane, then deduces from ramification data an expression for $F$ and
gives the relation between $\tau$ and $\tau^{\prime}$. Note that for $N \in \{2,
3, 4, 5 \}$, the modular curve $X(N)$ is also of genus zero, with respective
automorphism groups (leaving the set of cusps globally fixed) the dihedral group,
the tetrahedral group $A_4$, the group $S_4$ of the cube and the octahedron, and
the group $A_5$ of the dodecahedron and the icosahedron. Over the two year period
1878-1879, Klein published a series of papers on modular equations, devoted
respectively to transformations of order $p=5$, $7$ and $11$. In each case he
constructs by geometric means a Galois resolvent, gives its roots explicitly-using
modular forms-and shows how to find the modular equation itself of degree $p+1$
as well as a resolvent of degree $p$ for each of these particular values of $p$.

  In fact, the complex function field of $X=X_0(p)$ with $p$ a prime
consists of the modular functions $f(z)$ for $\Gamma_0(p)$ which are
meromorphic on the extended upper half-plane. A function $f$ lies in
the rational function field $\mathbb{Q}(X)$ if and only if the Fourier
coefficients in its expansion at $\infty$: $f(z)=\sum a_n q^n$ are all
rational numbers. The field $\mathbb{Q}(X)$ is known to be generated over
$\mathbb{Q}$ by the classical $j$-functions
$$\left\{\aligned
  j &=j(z)=q^{-1}+744+196884 q+\cdots,\\
  j_p &=j \left(\frac{-1}{pz}\right)=j(pz)=q^{-p}+744+\cdots.
\endaligned\right.$$
A further element in the function field $\mathbb{Q}(X)=\mathbb{Q}(j,
j_p)$ is the modular unit $u=\Delta(z)/\Delta(pz)$ with divisor
$(p-1) \{ (0)-(\infty) \}$ where $\Delta(z)$ is the discriminant.
If $m=\text{gcd}(p-1, 12)$, then an $m$th root of $u$ lies in
$\mathbb{Q}(X)$. This function has the Fourier expansion
$$t=\root m \of{u}=q^{(1-p)/m} \prod_{n \geq 1} \left(\frac{1-q^n}
    {1-q^{np}}\right)^{24/m}=\left(\frac{\eta(z)}{\eta(pz)}\right)^{24/m}.$$
When $p-1$ divides $12$, so $m=p-1$, the function $t$ is a Hauptmodul
for the curve $X$ which has genus zero. It is well-known that the genus
of the modular curve $X$ for prime $p$ is zero if and only if
$p=2, 3, 5, 7, 13$. In his paper, Klein studied the modular equations
of orders $2$, $3$, $5$, $7$, $13$ with degrees $3$, $4$, $6$, $8$, $14$,
respectively. They are uniformized, i.e., parametrized, by so-called
Hauptmoduln (principal moduli). A Hauptmodul is a function $J_{\Gamma}$
that is a modular function for some subgroup of $\Gamma(1)=SL(2, \mathbb{Z})$,
with any other modular function expressible as a rational function of it. In this
case, $\Gamma=\Gamma_0(p)$.

  The curve $X_0(p)$ can be given as a plane curve by the modular polynomial
$\Phi_p(X, Y)$. These can quickly get very complicated. For $p=2$ we have
$$\aligned
  \Phi_2(X, Y)=&X^3+Y^3-X^2 Y^2+1488 (X^2 Y+X Y^2)-162000 (X^2+Y^2)\\
               &+40773375 XY+8748000000 (X+Y)-157464000000000,
\endaligned$$
where $X$, $Y$ are the $j$-invariants of the two elliptic curves involved.
It is not so easy to guess that this is a genus $0$ curve. Hence, it is much
better, for conceptual understanding, to parametrize the curve by a different
modular function, and then write $X$ and $Y$ in terms of the parameter. This
leads to the Hauptmodul
$$\tau=\left(\frac{\eta(z)}{\eta(2z)}\right)^{24}.$$
Here,
$$X=j(z)=\frac{(\tau+256)^3}{\tau^2}=\frac{(\tau^{\prime}+16)^3}{\tau^{\prime}},$$
$$Y=j(2z)=\frac{(\tau+16)^3}{\tau}.$$
The Fricke involution involution $\tau^{\prime}=\tau(-\frac{1}{2z})$ satisfies
that
$$\tau \tau^{\prime}=2^{12}.$$

  For $p=3$, the modular polynomial is given by
$$\aligned
  \Phi_3(X, Y)=&X^4+Y^4-X^3 Y^3+2232 (X^3 Y^2+X^2 Y^3)\\
               &-1069956 (X^3 Y+X Y^3)+36864000 (X^3+Y^3)\\
               &+2587918086 X^2 Y^2+8900222976000 (X^2 Y+X Y^2)\\
               &+452984832000000 (X^2+Y^2)-770845966336000000 XY\\
               &+1855425871872000000000(X+Y),
\endaligned$$
We have
$$X=j(z)=\frac{(\tau+27) (\tau+243)^3}{\tau^3}
        =\frac{({\tau^{\prime}}+27) (\tau^{\prime}+3)^3}{\tau^{\prime}},$$
$$Y=j(3z)=\frac{(\tau+27) (\tau+3)^3}{\tau},$$
where the Hauptmodul is given by
$$\tau=\left(\frac{\eta(z)}{\eta(3 z)}\right)^{12},$$
and the Fricke involution $\tau^{\prime}=\tau(-\frac{1}{3z})$ satisfies
that
$$\tau \tau^{\prime}=3^6.$$

  For $p=5$, the geometric model of $X(5)$ he uses is the regular icosahedron,
the resolvent of degree five being linked, as had been shown by Hermite,
to the general quintic equation. In fact, Klein shows that the morphism
$X(5) \rightarrow X(1)$ is isomorphic to that taking the quotient of the
unit sphere in $\mathbb{R}^3$ by the action of the symmetry group of the
regular icosahedron. More precisely, in the modular equation $\Phi_5(X, Y)=0$
of level five, we have
$$X=j(z)=-\frac{(\tau^2-250 \tau+3125)^3}{\tau^5}
      =-\frac{({\tau^{\prime}}^2-10 \tau^{\prime}+5)^3}{\tau^{\prime}},\eqno{(8.2)}$$
$$Y=j(5z)=-\frac{(\tau^2-10 \tau+5)^3}{\tau},\eqno{(8.3)}$$
where the Hauptmodul is given by
$$\tau=\left(\frac{\eta(z)}{\eta(5 z)}\right)^6,\eqno{(8.4)}$$
and the Fricke involution $\tau^{\prime}=\tau(-\frac{1}{5z})$ satisfies
that
$$\tau \tau^{\prime}=125.$$

  For $p=7$, Klein shows that the modular curve $X(7)$ is isomorphic to
the smooth plane quartic curve $C_4$ with equation $x^3 y+y^3 z+z^3 x=0$
in $\mathbb{CP}^2$, invariant under the action of a group $G$ isomorphic
to $\text{PSL}(2, 7)$, which is the automorphism group of $X(7)$. In this
projective model the natural morphism from $X(7)$ onto
$X(1) \simeq \mathbb{CP}^1$ is made concrete as the projection of $C_4$ on
$C_4/G$ (identified with $\mathbb{CP}^1$). This is a Galois covering whose
generic fibre is considered by Klein as the Galois resolvent of the modular
equation $\Phi_7(\cdot, j)=0$ of level $7$, which means that the function
field of $C_4$ is the splitting field of this modular equation over $\mathbb{C}(j)$.
Moreover, just as the sphere has a regular tiling induced from the faces of
an inscribed icosahedron, so also does the modular curve $X(7)$ admits a
regular tiling by triangles. This tiling is inherited combinatorially from
a tiling of $\mathbb{H}$ of type $(2, 3, \infty)$, and its triangles are
of type $(2, 3, 7)$. Here, a tiling of $\mathbb{H}$ by triangles is said
to be of type $(a, b, c)$ if it is realized by hyperbolic triangles
$(a, b, c)$, that is, with angles $(\frac{2 \pi}{a}, \frac{2 \pi}{b},
\frac{2 \pi}{c})$. More precisely, in the modular equation $\Phi_7(X, Y)=0$
of level seven, we have
$$\aligned
  X=j(z) &=\frac{(\tau^2+13 \tau+49) (\tau^2+245 \tau+2401)^3}{\tau^7}\\
         &=\frac{({\tau^{\prime}}^2+13 \tau^{\prime}+49)({\tau^{\prime}}^2
         +5 \tau^{\prime}+1)^3}{\tau^{\prime}},
\endaligned\eqno{(8.5)}$$
$$Y=j(7z)=\frac{(\tau^2+13 \tau+49)(\tau^2+5 \tau+1)^3}{\tau},\eqno{(8.6)}$$
where the Hauptmodul is given by
$$\tau=\left(\frac{\eta(z)}{\eta(7 z)}\right)^4,\eqno{(8.7)}$$
and the Fricke involution $\tau^{\prime}=\tau(-\frac{1}{7z})$ satisfies
that
$$\tau \tau^{\prime}=49.$$

  Now, we come to the group $\text{PSL}(2, 13)$ and the modular curves
$X(13)$. The function field $K(\Gamma(13))$ of the modular curve $X(13)$
is the splitting field of the polynomial $\Phi_{13} \in \mathbb{C}(j)[X]$
associated with transformations of order $13$. Let
$\mathbb{F}_{13}=\mathbb{Z}/13 \mathbb{Z}$ be the field of thirteen elements.
Since $\text{SL}(2, 13)$ is generated by
$\left(\begin{matrix} 1 & 1\\ 0 & 1 \end{matrix}\right)$ and
$\left(\begin{matrix} 1 & 0\\ 1 & 1 \end{matrix}\right)$, the reduction
morphism modulo $13$ from $\text{SL}(2, \mathbb{Z})$ to $\text{SL}(2, 13)$
is surjective, whence the exact sequence
$$1 \rightarrow \overline{\Gamma}(13) \rightarrow \text{PSL}(2, \mathbb{Z})
  \rightarrow \text{PSL}(2, 13) \rightarrow 1,$$
where $\overline{\Gamma}:=\Gamma/ \{ \pm I \}$. In particular, the quotient
$\overline{G}=\text{PSL}(2, \mathbb{Z})/\overline{\Gamma}(13)$ is isomorphic
to $\text{PSL}(2, 13)$, a simple group of order $1092$. The group $\overline{G}$
acts on $X(13)$ via automorphisms and $X(13)/\overline{G}$ can be identified
with $X(1)$. Thus, the fibres of the projection $X(13) \rightarrow X(1)$ are
the orbits of the action of $\overline{G}$ on $X(13)$. There are therefore
three singular fibres corresponding to the values $J=\infty$, $0$ and $1$
(recall that $J=j/1728$), that is, $j=\infty$, $0$ and $1728$, whose elements
are called $A$-points, $B$-points and $C$-points in Klein's terminology, with
stabilizers of orders $13$, $3$ and $2$, respectively. These fibres have
cardinality $84$, $364$ and $546$; all others have $1092$ elements. By the
Riemann-Hurwitz formula, the genus of $X(13)$ satisfies the relation
$$2-2g=2 \cdot 1092-12 \cdot 84-2 \cdot 364-546,\eqno{(8.8)}$$
whence $g=50$. Hence, the modular curve $X(13)$ admits a regular tiling by
triangles. This tiling is inherited combinatorially from a tiling of
$\mathbb{H}$ of type $(2, 3, \infty)$, and its triangles are of type
$(2, 3, 13)$.

  The following facts about the modular subgroups of level $13$ should
be noted. One has $\overline{\Gamma}(13)<\overline{\Gamma}_0(13)<\overline{\Gamma}(1)$,
the respective subgroup indices being $78=6 \cdot 13$ and $14=13+1$. The
quotient $\overline{\Gamma}(1)/\overline{\Gamma}(13)$ is of order
$1092=78 \cdot 14$, and $\overline{\Gamma}_0(13)/\overline{\Gamma}(13)$
is a subgroup of order $78$, which is isomorphic to a semidirect product
of $\mathbb{Z}_{13}$ by $\mathbb{Z}_6$. The respective quotients of
$\mathbb{H}$ by these three groups (compactified) are the modular curves
$X(13)$, $X_0(13)$, $X(1)$, and the coverings $X(13) \rightarrow X_0(13)
\rightarrow X(1)$ are respectively $78$-sheeted and $14$-sheeted. The curve
$X_0(13)$ like $X(1)$ being of genus zero.

  In the modular equation $\Phi_{13}(X, Y)=0$ of level thirteen, we have
$$\aligned
  X=j(z) &=\frac{(\tau^2+5\tau+13)(\tau^4+247 \tau^3+3380 \tau^2
         +15379 \tau+28561)^3}{\tau^{13}}\\
         &=\frac{({\tau^{\prime}}^2+5 \tau^{\prime}+13)({\tau^{\prime}}^4
         +7 {\tau^{\prime}}^3+20 {\tau^{\prime}}^2+19 \tau^{\prime}+1)^3}
         {\tau^{\prime}},
\endaligned\eqno{(8.9)}$$
and
$$Y=j(13z)=\frac{(\tau^2+5 \tau+13)(\tau^4+7 \tau^3+20 \tau^2+19 \tau+1)^3}{\tau},
  \eqno{(8.10)}$$
where the Hauptmodul is given by
$$\tau=\left(\frac{\eta(z)}{\eta(13z)}\right)^2,\eqno{(8.11)}$$
and the Fricke involution $\tau^{\prime}=\tau(-\frac{1}{13 z})$ satisfies
that
$$\tau \tau^{\prime}=13.$$
Note that the above two quartic polynomials are defined over the quadratic
field $\mathbb{Q}(\sqrt{13})$:
$$\aligned
 &\tau^4+247 \tau^3+3380 \tau^2+15379 \tau+28561\\
=&\left(\tau^2+\frac{247+65 \sqrt{13}}{2} \tau+\frac{1859+507
  \sqrt{13}}{2}\right)\\
 &\times \left(\tau^2+\frac{247-65 \sqrt{13}}{2}
  \tau+\frac{1859-507 \sqrt{13}}{2}\right),
\endaligned$$
$$\aligned
 &\tau^4+7 \tau^3+20 \tau^2+19 \tau+1\\
=&\left(\tau^2+\frac{7+\sqrt{13}}{2} \tau+\frac{11+3
  \sqrt{13}}{2}\right)\left(\tau^2+\frac{7-\sqrt{13}}{2} \tau+\frac{11-3 \sqrt{13}}{2}\right).
\endaligned$$

  Now, we have shown that the modular curve $X(13)$ is isomorphic to the
algebraic space curve $Y$ with the $21$ equations in $\mathbb{CP}^5$,
invariant under the action of a group $G$ isomorphic to $\text{SL}(2, 13)$,
which is the automorphism group of $X(13)$. In this projective model the
natural morphism from $X(13)$ onto $X(1) \simeq \mathbb{CP}^1$ is made
concrete as the projection of $Y$ onto $Y/G$ (identified with $\mathbb{CP}^1$).
This is a Galois covering whose generic fibre is considered as the Galois
resolvent of the modular equation $\Phi_{13}( \cdot, j)=0$ of level $13$,
which means that the function field $K(Y)$ of $Y$ is the splitting field of
this modular equation over $\mathbb{C}(j)$. This proves the following:

\textbf{Theorem 8.1.} (Galois covering, Galois resolvent and their geometric
realizations). {\it In the projective model, $X(13)$ is isomorphic to $Y$ in
$\mathbb{CP}^5$. The natural morphism from $X(13)$ onto $X(1) \simeq \mathbb{CP}^1$
is realized as the projection of $Y$ onto $Y/G$ $($identified with $\mathbb{CP}^1$$)$.
This is a Galois covering whose generic fibre is interpreted as the Galois
resolvent of the modular equation $\Phi_{13}( \cdot, j)=0$ of level $13$, i.e.,
the function field of $Y$ is the splitting field of this modular equation
over $\mathbb{C}(j)$.}

\begin{center}
{\large\bf 9. Invariant theory and modular forms for $\text{SL}(2, 13)$ III:
               some computation for invariant polynomials}
\end{center}

  In this section, in order to study the ring of invariant polynomials
in six variables, we need to give some computation for invariant polynomials
in this ring.

  Because $\mathbf{B}$-terms vanish in the modular curves $X(13)$ according
to Theorem 7.6 and the $\mathbf{A}$-terms decompose into the sum of
$\mathbf{B}$-terms and $\mathbf{C}$-terms, we can only use $\mathbf{A}$-terms
instead of $\mathbf{C}$-terms.

  Put $x_i(z)=\eta(z) a_i(z)$, $y_i(z)=\eta^3(z) a_i(z)$ and
$u_i(z)=\eta^9(z) a_i(z)$ $(1 \leq i \leq 6)$. Let
$$\left\{\aligned
  X(z) &=(x_1(z), \ldots, x_6(z))^{T},\\
  Y(z) &=(y_1(z), \ldots, y_6(z))^{T}.\\
  U(z) &=(u_1(z), \ldots, u_6(z))^{T}.
\endaligned\right.$$
Then
$$\left\{\aligned
  X(z) &=\eta(z) \mathbf{A}(z)\\
  Y(z) &=\eta^3(z) \mathbf{A}(z),\\
  U(z) &=\eta^9(z) \mathbf{A}(z).
\endaligned\right.$$
Recall that $\eta(z)$ satisfies the following transformation formulas
$\eta(z+1)=e^{\frac{\pi i}{12}} \eta(z)$ and
$\eta\left(-\frac{1}{z}\right)=e^{-\frac{\pi i}{4}} \sqrt{z} \eta(z)$.
By Proposition 7.2, we have
$$X(z+1)=e^{-\frac{2 \pi i}{3}} \rho(t) X(z), \quad
  X\left(-\frac{1}{z}\right)=z \rho(s) X(z),$$
$$Y(z+1)=e^{-\frac{\pi i}{2}} \rho(t) Y(z), \quad
  Y\left(-\frac{1}{z}\right)=e^{-\frac{\pi i}{2}} z^2 \rho(s) Y(z).$$
$$U(z+1)=\rho(t) U(z), \quad
  U\left(-\frac{1}{z}\right)=z^5 \rho(s) U(z).$$
Define $j(\gamma, z):=cz+d$ if $z \in \mathbb{H}$ and
$\gamma=\begin{pmatrix} a & b\\ c & d \end{pmatrix} \in \Gamma(1)$.
Hence,
$$\left\{\aligned
  X(\gamma(z)) &=u(\gamma) j(\gamma, z) \rho(\gamma) X(z),\\
  Y(\gamma(z)) &=v(\gamma) j(\gamma, z)^2 \rho(\gamma) Y(z),\\
  U(\gamma(z)) &=j(\gamma, z)^5 \rho(\gamma) U(z)
\endaligned\right.\eqno{(9.1)}$$
for $\gamma \in \Gamma(1)$, where $u(\gamma)=1, \omega$ or $\omega^2$
with $\omega=e^{\frac{2 \pi i}{3}}$ and $v(\gamma)=\pm 1$ or $\pm i$.
Since $\Gamma(13)=\text{ker}$ $\rho$, we have
$X(\gamma(z))=u(\gamma) j(\gamma, z) X(z)$,
$Y(\gamma(z))=v(\gamma) j(\gamma, z)^2 Y(z)$ and
$U(\gamma(z))=j(\gamma, z)^5 U(z)$ for $\gamma \in \Gamma(13)$.
This means that the functions $x_1(z)$, $\ldots$, $x_6(z)$ are
modular forms of weight one for $\Gamma(13)$ with the same
multiplier $u(\gamma)=1, \omega$ or $\omega^2$ and $y_1(z)$,
$\ldots$, $y_6(z)$ are modular forms of weight two for
$\Gamma(13)$ with the same multiplier $v(\gamma)=\pm 1 $
or $\pm i$.

  From now on, we will use the following abbreviation
$$\mathbf{A}_j=\mathbf{A}_j(a_1(z), \ldots, a_6(z)) \quad
  (0 \leq j \leq 6),$$
$$\mathbf{D}_j=\mathbf{D}_j(a_1(z), \ldots, a_6(z)) \quad
  (j=0,1, \ldots, 12, \infty)$$
and
$$\mathbf{G}_j=\mathbf{G}_j(a_1(z), \ldots, a_6(z)) \quad
  (0 \leq j \leq 12).$$
We have
$$\left\{\aligned
  \mathbf{A}_0 &=q^{\frac{1}{4}} (1+O(q)),\\
  \mathbf{A}_1 &=q^{\frac{17}{52}} (2+O(q)),\\
  \mathbf{A}_2 &=q^{\frac{29}{52}} (2+O(q)),\\
  \mathbf{A}_3 &=q^{\frac{49}{52}} (1+O(q)),\\
  \mathbf{A}_4 &=q^{\frac{25}{52}} (-1+O(q)),\\
  \mathbf{A}_5 &=q^{\frac{9}{52}} (-1+O(q)),\\
  \mathbf{A}_6 &=q^{\frac{1}{52}} (-1+O(q)),
\endaligned\right.$$
and
$$\left\{\aligned
  \mathbf{D}_0 &=q^{\frac{15}{8}} (1+O(q)),\\
  \mathbf{D}_{\infty} &=q^{\frac{7}{8}} (-1+O(q)),\\
  \mathbf{D}_1 &=q^{\frac{99}{104}} (2+O(q)),\\
  \mathbf{D}_2 &=q^{\frac{3}{104}} (-1+O(q)),\\
  \mathbf{D}_3 &=q^{\frac{11}{104}} (1+O(q)),\\
  \mathbf{D}_4 &=q^{\frac{19}{104}} (-2+O(q)),\\
  \mathbf{D}_5 &=q^{\frac{27}{104}} (-1+O(q)),
\endaligned\right. \quad \quad
  \left\{\aligned
  \mathbf{D}_6 &=q^{\frac{35}{104}} (-1+O(q)),\\
  \mathbf{D}_7 &=q^{\frac{43}{104}} (1+O(q)),\\
  \mathbf{D}_8 &=q^{\frac{51}{104}} (3+O(q)),\\
  \mathbf{D}_9 &=q^{\frac{59}{104}} (-2+O(q)),\\
  \mathbf{D}_{10} &=q^{\frac{67}{104}} (1+O(q)),\\
  \mathbf{D}_{11} &=q^{\frac{75}{104}} (-4+O(q)),\\
  \mathbf{D}_{12} &=q^{\frac{83}{104}} (-1+O(q)).
\endaligned\right.$$
Hence,
$$\left\{\aligned
  \mathbf{G}_0 &=q^{\frac{7}{4}} (1+O(q)),\\
  \mathbf{G}_1 &=q^{\frac{43}{52}} (13+O(q)),\\
  \mathbf{G}_2 &=q^{\frac{47}{52}} (-22+O(q)),\\
  \mathbf{G}_3 &=q^{\frac{51}{52}} (-21+O(q)),\\
  \mathbf{G}_4 &=q^{\frac{3}{52}} (-1+O(q)),\\
  \mathbf{G}_5 &=q^{\frac{7}{52}} (2+O(q)),\\
  \mathbf{G}_6 &=q^{\frac{11}{52}} (2+O(q)),
\endaligned\right. \quad \quad
  \left\{\aligned
  \mathbf{G}_7 &=q^{\frac{15}{52}} (-2+O(q)),\\
  \mathbf{G}_8 &=q^{\frac{19}{52}} (-8+O(q)),\\
  \mathbf{G}_9 &=q^{\frac{23}{52}} (6+O(q)),\\
  \mathbf{G}_{10} &=q^{\frac{27}{52}} (1+O(q)),\\
  \mathbf{G}_{11} &=q^{\frac{31}{52}} (-8+O(q)),\\
  \mathbf{G}_{12} &=q^{\frac{35}{52}} (17+O(q)).
\endaligned\right.$$
Note that
$$\aligned
  w_{\nu} &=(\mathbf{A}_0+\zeta^{\nu} \mathbf{A}_1+\zeta^{4 \nu} \mathbf{A}_2
           +\zeta^{9 \nu} \mathbf{A}_3+\zeta^{3 \nu} \mathbf{A}_4
           +\zeta^{12 \nu} \mathbf{A}_5+\zeta^{10 \nu} \mathbf{A}_6)^2\\
          &=\mathbf{A}_0^2+2 (\mathbf{A}_1 \mathbf{A}_5+\mathbf{A}_2 \mathbf{A}_3
           +\mathbf{A}_4 \mathbf{A}_6)+\\
          &+2 \zeta^{\nu} (\mathbf{A}_0 \mathbf{A}_1+\mathbf{A}_2 \mathbf{A}_6)
           +2 \zeta^{3 \nu} (\mathbf{A}_0 \mathbf{A}_4+\mathbf{A}_2 \mathbf{A}_5)+\\
          &+2 \zeta^{9 \nu} (\mathbf{A}_0 \mathbf{A}_3+\mathbf{A}_5 \mathbf{A}_6)
           +2 \zeta^{12 \nu} (\mathbf{A}_0 \mathbf{A}_5+\mathbf{A}_3 \mathbf{A}_4)+\\
          &+2 \zeta^{10 \nu} (\mathbf{A}_0 \mathbf{A}_6+\mathbf{A}_1 \mathbf{A}_3)
           +2 \zeta^{4 \nu} (\mathbf{A}_0 \mathbf{A}_2+\mathbf{A}_1 \mathbf{A}_4)+\\
          &+\zeta^{2 \nu} (\mathbf{A}_1^2+2 \mathbf{A}_4 \mathbf{A}_5)
           +\zeta^{5 \nu} (\mathbf{A}_3^2+2 \mathbf{A}_1 \mathbf{A}_2)+\\
          &+\zeta^{6 \nu} (\mathbf{A}_4^2+2 \mathbf{A}_3 \mathbf{A}_6)
           +\zeta^{11 \nu} (\mathbf{A}_5^2+2 \mathbf{A}_1 \mathbf{A}_6)+\\
          &+\zeta^{8 \nu} (\mathbf{A}_2^2+2 \mathbf{A}_3 \mathbf{A}_5)
           +\zeta^{7 \nu} (\mathbf{A}_6^2+2 \mathbf{A}_4 \mathbf{A}_2),
\endaligned$$
where
$$\mathbf{A}_0^2+2 (\mathbf{A}_1 \mathbf{A}_5+\mathbf{A}_2 \mathbf{A}_3
  +\mathbf{A}_4 \mathbf{A}_6)=q^{\frac{1}{2}} (-1+O(q)),$$
$$\left\{\aligned
  \mathbf{A}_0 \mathbf{A}_1+\mathbf{A}_2 \mathbf{A}_6 &=q^{\frac{41}{26}} (-3+O(q)),\\
  \mathbf{A}_0 \mathbf{A}_4+\mathbf{A}_2 \mathbf{A}_5 &=q^{\frac{19}{26}} (-3+O(q)),\\
  \mathbf{A}_0 \mathbf{A}_3+\mathbf{A}_5 \mathbf{A}_6 &=q^{\frac{5}{26}} (1+O(q)),\\
  \mathbf{A}_0 \mathbf{A}_5+\mathbf{A}_3 \mathbf{A}_4 &=q^{\frac{11}{26}} (-1+O(q)),\\
  \mathbf{A}_0 \mathbf{A}_6+\mathbf{A}_1 \mathbf{A}_3 &=q^{\frac{7}{26}} (-1+O(q)),\\
  \mathbf{A}_0 \mathbf{A}_2+\mathbf{A}_1 \mathbf{A}_4 &=q^{\frac{47}{26}} (-1+O(q)),
\endaligned\right.$$
and
$$\left\{\aligned
  \mathbf{A}_1^2+2 \mathbf{A}_4 \mathbf{A}_5 &=q^{\frac{17}{26}} (6+O(q)),\\
  \mathbf{A}_3^2+2 \mathbf{A}_1 \mathbf{A}_2 &=q^{\frac{23}{26}} (8+O(q)),\\
  \mathbf{A}_4^2+2 \mathbf{A}_3 \mathbf{A}_6 &=q^{\frac{25}{26}} (-1+O(q)),\\
  \mathbf{A}_5^2+2 \mathbf{A}_1 \mathbf{A}_6 &=q^{\frac{9}{26}} (-3+O(q)),\\
  \mathbf{A}_2^2+2 \mathbf{A}_3 \mathbf{A}_5 &=q^{\frac{29}{26}} (2+O(q)),\\
  \mathbf{A}_6^2+2 \mathbf{A}_4 \mathbf{A}_2 &=q^{\frac{1}{26}} (1+O(q)).
\endaligned\right.$$

  In the introduction, the invariant homogeneous polynomials $\Phi_{m, n}$
are defined. Now we give the normalization for the following four families
of polynomials of degrees $d=12$, $16$, $20$ and $30$. For $d=12$, there are
two such invariant homogeneous polynomials $\Phi_{3, 0}$ and $\Phi_{0, 2}$:
$$\left\{\aligned
  \Phi_{3, 0} &:=-\frac{1}{13 \cdot 30} \left(\sum_{\nu=0}^{12}
                   w_{\nu}^3+w_{\infty}^3\right),\\
  \Phi_{0, 2} &:=-\frac{1}{13 \cdot 52} \left(\sum_{\nu=0}^{12}
                 \delta_{\nu}^2+\delta_{\infty}^2\right),
\endaligned\right.\eqno{(9.2)}$$
For $d=16$, there are two such invariant homogeneous polynomials
$\Phi_{4, 0}$ and $\Phi_{1, 2}$ which need not to be normalized.
For $d=20$, there are two such invariant homogeneous polynomials
$\Phi_{5, 0}$ and $\Phi_{2, 2}$:
$$\left\{\aligned
  \Phi_{5, 0} &:=\frac{1}{13 \cdot 25} \left(\sum_{\nu=0}^{12} w_{\nu}^5
                 +w_{\infty}^5\right),\\
  \Phi_{2, 2} &:=\frac{1}{13 \cdot 26} \left(\sum_{\nu=0}^{12} w_{\nu}^2
                 \delta_{\nu}^2+w_{\infty}^2 \delta_{\infty}^2\right)
\endaligned\right.\eqno{(9.3)}$$
For $d=30$, there are three such invariant homogeneous polynomials
$\Phi_{0, 5}$, $\Phi_{3, 3}$ and $\Phi_{6, 1}$:
$$\left\{\aligned
  \Phi_{0, 5} &:=-\frac{1}{13 \cdot 1315} \left(\sum_{\nu=0}^{12}
                 \delta_{\nu}^5+\delta_{\infty}^5\right),\\
  \Phi_{3, 3} &:=-\frac{1}{13 \cdot 27} \left(\sum_{\nu=0}^{12}
                 w_{\nu}^3 \delta_{\nu}^3+w_{\infty}^3 \delta_{\infty}^3\right),\\
  \Phi_{6, 1} &:=-\frac{1}{13 \cdot 285} \left(\sum_{\nu=0}^{12}
                 w_{\nu}^6 \delta_{\nu}+w_{\infty}^6 \delta_{\infty}\right),
\endaligned\right.\eqno{(9.4)}$$

\textbf{Theorem 9.1.} {\it The $G$-invariant homogeneous polynomials
$\Phi_{m, n}$ of degrees $d=4$, $8$, $10$, $12$, $14$, $16$, $20$ and
$30$ in $x_1(z)$, $\ldots$, $x_6(z)$ can be identified with modular
forms as follows$:$
$$\left\{\aligned
  \Phi_{4}(x_1(z), \ldots, x_6(z)) &=0,\\
  \Phi_{8}(x_1(z), \ldots, x_6(z)) &=0,\\
  \Phi_{10}(x_1(z), \ldots, x_6(z)) &=0,\\
  \Phi_{3, 0}(x_1(z), \ldots, x_6(z)) &=\Delta(z),\\
  \Phi_{0, 2}(x_1(z), \ldots, x_6(z)) &=\Delta(z),\\
  \Phi_{14}(x_1(z), \ldots, x_6(z)) &=0,\\
  \Phi_{4, 0}(x_1(z), \ldots, x_6(z)) &=0,\\
  \Phi_{1, 2}(x_1(z), \ldots, x_6(z)) &=0,\\
  \Phi_{5, 0}(x_1(z), \ldots, x_6(z)) &=\eta(z)^8 \Delta(z) E_4(z),\\
  \Phi_{2, 2}(x_1(z), \ldots, x_6(z)) &=\eta(z)^8 \Delta(z) E_4(z),\\
  \Phi_{0, 5}(x_1(z), \ldots, x_6(z)) &=\Delta(z)^2 E_6(z),\\
  \Phi_{3, 3}(x_1(z), \ldots, x_6(z)) &=\Delta(z)^2 E_6(z),\\
  \Phi_{6, 1}(x_1(z), \ldots, x_6(z)) &=\Delta(z)^2 E_6(z).\\
\endaligned\right.\eqno{(9.5)}$$}

{\it Proof}. We divide the proof into four parts (see also \cite{Y3}).
The first part is the calculation of $\Phi_{5, 0}$ and $\Phi_{3, 0}$.
Up to a constant, $\Phi_{5, 0}=w_0^5+w_1^5+\cdots+w_{12}^5+w_{\infty}^5$.
As a polynomial in six variables, $\Phi_{5, 0}(z_1, z_2, z_3, z_4, z_5, z_6)$
is a $G$-invariant polynomial. Moreover, for $\gamma \in \Gamma(1)$,
$$\aligned
  &\Phi_{5, 0}(Y(\gamma(z))^{T})=\Phi_{5, 0}(v(\gamma)
   j(\gamma, z)^2 (\rho(\gamma) Y(z))^{T})\\
 =&v(\gamma)^{20} j(\gamma, z)^{40} \Phi_{5, 0}((\rho(\gamma) Y(z))^{T})
 =j(\gamma, z)^{40} \Phi_{5, 0}((\rho(\gamma) Y(z))^{T}).
\endaligned$$
Note that $\rho(\gamma) \in \langle \rho(s), \rho(t) \rangle=G$ and
$\Phi_{5, 0}$ is a $G$-invariant polynomial, we have
$$\Phi_{5, 0}(Y(\gamma(z))^{T})=j(\gamma, z)^{40} \Phi_{5, 0}(Y(z)^{T}),
  \quad \text{for $\gamma \in \Gamma(1)$}.$$
This implies that $\Phi_{5, 0}(y_1(z), \ldots, y_6(z))$ is a modular form
of weight $40$ for the full modular group $\Gamma(1)$. Moreover, we will
show that it is a cusp form. In fact,
$$\aligned
  &\Phi_{5, 0}(a_1(z), \ldots, a_6(z))=13^5 q^{\frac{5}{2}} (1+O(q))^5+\\
  &+\sum_{\nu=0}^{12} [q^{\frac{1}{2}} (-1+O(q))+\\
  &+2 \zeta^{\nu} q^{\frac{41}{26}} (-3+O(q))+2 \zeta^{3 \nu}
    q^{\frac{19}{26}} (-3+O(q))
   +2 \zeta^{9 \nu} q^{\frac{5}{26}} (1+O(q))+\\
  &+2 \zeta^{12 \nu} q^{\frac{11}{26}} (-1+O(q))+2 \zeta^{10 \nu}
   q^{\frac{7}{26}} (-1+O(q))
   +2 \zeta^{4 \nu} q^{\frac{47}{26}} (-1+O(q))+\\
  &+\zeta^{2 \nu} q^{\frac{17}{26}} (6+O(q))+\zeta^{5 \nu}
   q^{\frac{23}{26}} (8+O(q))
   +\zeta^{6 \nu} q^{\frac{25}{26}} (-1+O(q))+\\
  &+\zeta^{11 \nu} q^{\frac{9}{26}} (-3+O(q))+\zeta^{8 \nu}
   q^{\frac{29}{26}} (2+O(q))
   +\zeta^{7 \nu} q^{\frac{1}{26}} (1+O(q))]^5.
\endaligned$$
We will calculate the $q^{\frac{1}{2}}$-term which is the lowest
degree. For the partition $13=4 \cdot 1+9$, the corresponding
term is
$$\begin{pmatrix} 5\\ 4, 1 \end{pmatrix} (\zeta^{7 \nu} q^{\frac{1}{26}})^4
  \cdot (-3) \zeta^{11 \nu} q^{\frac{9}{26}}=-15 q^{\frac{1}{2}}.$$
For the partition $13=3 \cdot 1+2 \cdot 5$, the corresponding term
is
$$\begin{pmatrix} 5\\ 3, 2 \end{pmatrix} (\zeta^{7 \nu} q^{\frac{1}{26}})^3
  \cdot (2 \zeta^{9 \nu} q^{\frac{5}{26}})^2=40 q^{\frac{1}{2}}.$$
Hence, for $\Phi_{5, 0}(y_1(z), \ldots, y_6(z))$ which is a modular
form for $\Gamma(1)$ with weight $40$, the lowest degree term is given by
$$(-15+40) q^{\frac{1}{2}} \cdot q^{\frac{3}{24} \cdot 20}=25 q^3.$$
Thus,
$$\Phi_{5, 0}(y_1(z), \ldots, y_6(z))=q^3 (13 \cdot 25+O(q)).$$
The leading term of $\Phi_{5, 0}(y_1(z), \ldots, y_6(z))$ together
with its weight $40$ suffice to identify this modular form with
$\Phi_{5, 0}(y_1(z), \ldots, y_6(z))=13 \cdot 25 \Delta(z)^3 E_4(z)$.
Consequently,
$$\Phi_{5, 0}(x_1(z), \ldots, x_6(z))=13 \cdot 25 \Delta(z)^3
  E_4(z)/\eta(z)^{40}=13 \cdot 25 \eta(z)^8 \Delta(z) E_4(z).$$

  Up to a constant, $\Phi_{3, 0}=w_0^3+w_1^3+\cdots+w_{12}^3+w_{\infty}^3$,
The calculation of $\Phi_{3, 0}$ is similar as that of $\Phi_{5, 0}$.
We find that
$$\Phi_{3, 0}(x_1(z), \ldots, x_6(z))=-13 \cdot 30 \Delta(z).$$

  The second part is the calculation of $\Phi_4$, $\Phi_8$ and
$\Phi_{4, 0}$. The calculation of $\Phi_4$ has been done in \cite{Y2},
Theorem 3.1. We will give the calculation of $\Phi_{4, 0}$. Up to a
constant, $\Phi_{4, 0}=w_0^4+w_1^4+\cdots+w_{12}^4+w_{\infty}^4$.
Similar as the above calculation for $\Phi_{5, 0}$, we find that
$\Phi_{4, 0}(y_1(z), \ldots, y_6(z))$ is a modular form of weight
$32$ for the full modular group $\Gamma(1)$. Moreover, we will
show that it is a cusp form. In fact,
$$\aligned
  &\Phi_{4, 0}(a_1(z), \ldots, a_6(z))=13^4 q^2 (1+O(q))^4+\\
  &+\sum_{\nu=0}^{12} [q^{\frac{1}{2}} (-1+O(q))+\\
  &+2 \zeta^{\nu} q^{\frac{41}{26}} (-3+O(q))+2 \zeta^{3 \nu} q^{\frac{19}{26}} (-3+O(q))
   +2 \zeta^{9 \nu} q^{\frac{5}{26}} (1+O(q))+\\
  &+2 \zeta^{12 \nu} q^{\frac{11}{26}} (-1+O(q))+2 \zeta^{10 \nu} q^{\frac{7}{26}} (-1+O(q))
   +2 \zeta^{4 \nu} q^{\frac{47}{26}} (-1+O(q))+\\
  &+\zeta^{2 \nu} q^{\frac{17}{26}} (6+O(q))+\zeta^{5 \nu} q^{\frac{23}{26}} (8+O(q))
   +\zeta^{6 \nu} q^{\frac{25}{26}} (-1+O(q))+\\
  &+\zeta^{11 \nu} q^{\frac{9}{26}} (-3+O(q))+\zeta^{8 \nu} q^{\frac{29}{26}} (2+O(q))
   +\zeta^{7 \nu} q^{\frac{1}{26}} (1+O(q))]^4.
\endaligned$$
We will calculate the $q$-term which is the lowest degree. For example,
consider the partition $26=3 \cdot 1+23$, the corresponding term is
$$\begin{pmatrix} 4\\ 3, 1 \end{pmatrix} (\zeta^{7 \nu} q^{\frac{1}{26}})^3
  \cdot 8 \zeta^{5 \nu} q^{\frac{23}{26}}=32 q.$$
For the other partitions, the calculation is similar. In conclusion,
we find that the coefficients of the $q$-term is an integer. Hence,
for $\Phi_{4, 0}(y_1(z), \ldots, y_6(z))$ which is a modular form for
$\Gamma(1)$ with weight $32$, the lowest degree term is given by
$$\text{some integer} \cdot q \cdot q^{\frac{3}{24} \cdot 16}
 =\text{some integer} \cdot q^3.$$
This implies that $\Phi_{4, 0}(y_1(z), \ldots, y_6(z))$ has a factor of
$\Delta(z)^3$, which is a cusp form of weight $36$. Therefore,
$\Phi_{4, 0}(y_1(z), \ldots, y_6(z))=0$. The calculation of $\Phi_8$ is
similar as that of $\Phi_{4, 0}$.

   The third part is the calculation of $\Phi_{0, 2}$ and $\Phi_{0, 5}$.
Up to a constant,
$\Phi_{0, 2}=\delta_0^2+\delta_1^2+\cdots+\delta_{12}^2+\delta_{\infty}^2$.
As a polynomial in six variables, $\Phi_{0, 2}(z_1, z_2, z_3, z_4, z_5, z_6)$
is a $G$-invariant polynomial. Moreover, for $\gamma \in \Gamma(1)$,
$$\aligned
  &\Phi_{0, 2}(X(\gamma(z))^{T})
 =\Phi_{0, 2}(u(\gamma) j(\gamma, z) (\rho(\gamma) X(z))^{T})\\
 =&u(\gamma)^{12} j(\gamma, z)^{12} \Phi_{0, 2}((\rho(\gamma) X(z))^{T})
 =j(\gamma, z)^{12} \Phi_{0, 2}((\rho(\gamma) X(z))^{T}).
\endaligned$$
Note that $\rho(\gamma) \in \langle \rho(s), \rho(t) \rangle=G$
and $\Phi_{0, 2}$ is a $G$-invariant polynomial, we have
$$\Phi_{0, 2}(X(\gamma(z))^{T})=j(\gamma, z)^{12}
  \Phi_{0, 2}(X(z)^{T}), \quad \text{for $\gamma \in
  \Gamma(1)$}.$$
This implies that $\Phi_{0, 2}(x_1(z), \ldots, x_6(z))$ is a modular
form of weight $12$ for the full modular group $\Gamma(1)$. Moreover,
we will show that it is a cusp form. In fact,
$$\aligned
  &\Phi_{0, 2}(a_1(z), \ldots, a_6(z))=13^4 q^{\frac{7}{2}} (1+O(q))^2+\\
  &+\sum_{\nu=0}^{12} [-13 q^{\frac{7}{4}} (1+O(q))+\\
  &+\zeta^{\nu} q^{\frac{43}{52}} (13+O(q))+\zeta^{2 \nu}
   q^{\frac{47}{52}} (-22+O(q))
   +\zeta^{3 \nu} q^{\frac{51}{52}} (-21+O(q))+\\
  &+\zeta^{4 \nu} q^{\frac{3}{52}} (-1+O(q))+\zeta^{5 \nu}
   q^{\frac{7}{52}} (2+O(q))
   +\zeta^{6 \nu} q^{\frac{11}{52}} (2+O(q))+\\
  &+\zeta^{7 \nu} q^{\frac{15}{52}} (-2+O(q))+\zeta^{8 \nu}
   q^{\frac{19}{52}} (-8+O(q))
   +\zeta^{9 \nu} q^{\frac{23}{52}} (6+O(q))+\\
  &+\zeta^{10 \nu} q^{\frac{27}{52}} (1+O(q))+\zeta^{11 \nu}
   q^{\frac{31}{52}} (-8+O(q))
   +\zeta^{12 \nu} q^{\frac{35}{52}} (17+O(q))]^2.
\endaligned$$
We will calculate the $q^{\frac{1}{2}}$-term which is the lowest
degree. For the partition $26=3+23$, the corresponding term is
$$\begin{pmatrix} 2\\ 1, 1 \end{pmatrix} \zeta^{4 \nu} q^{\frac{3}{52}}
  \cdot (-1) \cdot \zeta^{9 \nu} q^{\frac{23}{52}} \cdot 6=-12 q^{\frac{1}{2}}.$$
For the partition $26=7+19$, the corresponding term is
$$\begin{pmatrix} 2\\ 1, 1 \end{pmatrix} \zeta^{5 \nu} q^{\frac{7}{52}}
  \cdot 2 \cdot \zeta^{8 \nu} q^{\frac{19}{52}} \cdot (-8)=-32 q^{\frac{1}{2}}.$$
For the partition $26=11+15$, the corresponding term is
$$\begin{pmatrix} 2\\ 1, 1 \end{pmatrix} \zeta^{6 \nu} q^{\frac{11}{52}}
  \cdot 2 \cdot \zeta^{7 \nu} q^{\frac{15}{52}} \cdot (-2)=-8 q^{\frac{1}{2}}.$$
Hence, for $\Phi_{0, 2}(x_1(z), \ldots, x_6(z))$ which is a
modular form for $\Gamma(1)$ with weight $12$, the lowest degree
term is given by $(-12-32-8) q^{\frac{1}{2}} \cdot q^{\frac{12}{24}}=-52 q$.
Thus,
$$\Phi_{0, 2}(x_1(z), \ldots, x_6(z))=q (-13 \cdot 52+O(q)).$$
The leading term of $\Phi_{0, 2}(x_1(z), \ldots, x_6(z))$ together
with its weight $12$ suffice to identify this modular form with
$$\Phi_{0, 2}(x_1(z), \ldots, x_6(z))=-13 \cdot 52 \Delta(z).$$

  Up to a constant,
$\Phi_{0, 5}=\delta_0^5+\delta_1^5+\cdots+\delta_{12}^5+\delta_{\infty}^5$.
As a polynomial in six variables, $\Phi_{0, 5}(z_1, z_2, z_3, z_4, z_5, z_6)$
is a $G$-invariant polynomial. Similarly as above, we can show that
$\Phi_{0, 5}(x_1(z), \ldots, x_6(z))$ is a modular form of weight $30$ for
the full modular group $\Gamma(1)$. Moreover, we will show that it is a
cusp form. In fact,
$$\aligned
  &\Phi_{0, 5}(a_1(z), \ldots, a_6(z))=13^{10} q^{\frac{35}{4}} (1+O(q))^5+\\
  &+\sum_{\nu=0}^{12} [-13 q^{\frac{7}{4}} (1+O(q))+\\
  &+\zeta^{\nu} q^{\frac{43}{52}} (13+O(q))+\zeta^{2 \nu} q^{\frac{47}{52}} (-22+O(q))
   +\zeta^{3 \nu} q^{\frac{51}{52}} (-21+O(q))+\\
  &+\zeta^{4 \nu} q^{\frac{3}{52}} (-1+O(q))+\zeta^{5 \nu} q^{\frac{7}{52}} (2+O(q))
   +\zeta^{6 \nu} q^{\frac{11}{52}} (2+O(q))+\\
  &+\zeta^{7 \nu} q^{\frac{15}{52}} (-2+O(q))+\zeta^{8 \nu} q^{\frac{19}{52}} (-8+O(q))
   +\zeta^{9 \nu} q^{\frac{23}{52}} (6+O(q))+\\
  &+\zeta^{10 \nu} q^{\frac{27}{52}} (1+O(q))+\zeta^{11 \nu} q^{\frac{31}{52}} (-8+O(q))
   +\zeta^{12 \nu} q^{\frac{35}{52}} (17+O(q))]^5.
\endaligned$$
We will calculate the $q^{\frac{3}{4}}$-term which is the lowest degree.
(1) For the partition $39=4 \cdot 3+27$, the corresponding term is
$$\begin{pmatrix} 5\\ 4, 1 \end{pmatrix} (\zeta^{4 \nu}
  q^{\frac{3}{52}} \cdot (-1))^4 \cdot \zeta^{10 \nu}
  q^{\frac{27}{52}}=5 q^{\frac{3}{4}}.$$
(2) For the partition $39=3 \cdot 3+7+23$, the corresponding term is
$$\begin{pmatrix} 5\\ 3, 1, 1 \end{pmatrix} (\zeta^{4 \nu}
  q^{\frac{3}{52}} \cdot (-1))^3 \cdot \zeta^{5 \nu} q^{\frac{7}{52}}
  \cdot 2 \cdot \zeta^{9 \nu} q^{\frac{23}{52}} \cdot 6=-240 q^{\frac{3}{4}}.$$
(3) For the partition $39=3 \cdot 3+11+19$, the corresponding term is
$$\begin{pmatrix} 5\\ 3, 1, 1 \end{pmatrix} (\zeta^{4 \nu}
  q^{\frac{3}{52}} \cdot (-1))^3 \cdot
  \zeta^{6 \nu} q^{\frac{11}{52}} \cdot 2 \cdot \zeta^{8 \nu}
  q^{\frac{19}{52}} \cdot (-8)=320 q^{\frac{3}{4}}.$$
(4) For the partition $39=3 \cdot 3+2 \cdot 15$, the corresponding term is
$$\begin{pmatrix} 5\\ 3, 2 \end{pmatrix} (\zeta^{4 \nu}
  q^{\frac{3}{52}} \cdot (-1))^3 \cdot
  (\zeta^{7 \nu} q^{\frac{15}{52}} \cdot (-2))^2=-40 q^{\frac{3}{4}}.$$
(5) For the partition $39=2 \cdot 3+3 \cdot 11$, the corresponding term is
$$\begin{pmatrix} 5\\ 2, 3 \end{pmatrix} (\zeta^{4 \nu}
  q^{\frac{3}{52}} \cdot (-1))^2 \cdot
  (\zeta^{6 \nu} q^{\frac{11}{52}} \cdot 2)^3=80 q^{\frac{3}{4}}.$$
(6) For the partition $39=2 \cdot 3+2 \cdot 7+19$, the corresponding term is
$$\begin{pmatrix} 5\\ 2, 2, 1 \end{pmatrix} (\zeta^{4 \nu}
  q^{\frac{3}{52}} \cdot (-1))^2 \cdot
  (\zeta^{5 \nu} q^{\frac{7}{52}} \cdot 2)^2 \cdot \zeta^{8 \nu}
  q^{\frac{19}{52}} \cdot (-8)=-960 q^{\frac{3}{4}}.$$
(7) For the partition $39=2 \cdot 3+7+11+15$, the corresponding term is
$$\begin{pmatrix} 5\\ 2, 1, 1, 1 \end{pmatrix} (\zeta^{4 \nu}
  q^{\frac{3}{52}} \cdot (-1))^2 \cdot
  \zeta^{5 \nu} q^{\frac{7}{52}} \cdot 2 \cdot \zeta^{6 \nu}
  q^{\frac{11}{52}} \cdot 2 \cdot \zeta^{7 \nu}
  q^{\frac{15}{52}} \cdot (-2)=-480 q^{\frac{3}{4}}.$$
(8) For the partition $39=1 \cdot 3+3 \cdot 7+15$, the corresponding
term is
$$\begin{pmatrix} 5\\ 1, 3, 1 \end{pmatrix} \zeta^{4 \nu}
  q^{\frac{3}{52}} \cdot (-1) \cdot
  (\zeta^{5 \nu} q^{\frac{7}{52}} \cdot 2)^3 \cdot \zeta^{7 \nu}
  q^{\frac{15}{52}} \cdot (-2)=320 q^{\frac{3}{4}}.$$
(9) For the partition $39=1 \cdot 3+2 \cdot 7+2 \cdot 11$, the
corresponding term is
$$\begin{pmatrix} 5\\ 1, 2, 2 \end{pmatrix} \zeta^{4 \nu}
  q^{\frac{3}{52}} \cdot (-1) \cdot
  (\zeta^{5 \nu} q^{\frac{7}{52}} \cdot 2)^2 \cdot (\zeta^{6 \nu}
  q^{\frac{11}{52}} \cdot 2)^2=-480 q^{\frac{3}{4}}.$$
(10) For the partition $39=4 \cdot 7+11$, the corresponding term is
$$\begin{pmatrix} 5\\ 4, 1 \end{pmatrix} (\zeta^{5 \nu}
  q^{\frac{7}{52}} \cdot 2)^4 \cdot
  \zeta^{6 \nu} q^{\frac{11}{52}} \cdot 2=160 q^{\frac{3}{4}}.$$
Hence, for $\Phi_{0, 5}(x_1(z), \ldots, x_6(z))$ which is a modular
form for $\Gamma(1)$ with weight $30$, the lowest degree term is
given by
$$(5-240+320-40+80-960-480+320-480+160) q^{\frac{3}{4}} \cdot
  q^{\frac{30}{24}}=-1315 q^2.$$
Thus,
$$\Phi_{0, 5}(x_1(z), \ldots, x_6(z))=q^2 (-13 \cdot 1315+O(q)).$$
The leading term of $\Phi_{0, 5}(x_1(z), \ldots, x_6(z))$ together
with its weight $30$ suffice to identify this modular form with
$$\Phi_{0, 5}(x_1(z), \ldots, x_6(z))=-13 \cdot 1315 \Delta(z)^2 E_6(z).$$

  The last part is the calculation of $\Phi_{10}$, $\Phi_{14}$, $\Phi_{1, 2}$,
$\Phi_{2, 2}$, $\Phi_{3, 3}$ and $\Phi_{6, 1}$. For
$\Phi_{10}=w_0 \delta_0+w_1 \delta_1+\cdots+w_{12} \delta_{12}+w_{\infty} \delta_{\infty}$.
As a polynomial in six variables, $\Phi_{10}(z_1, z_2, z_3, z_4, z_5, z_6)$
is a $G$-invariant polynomial. Moreover, for $\gamma \in \Gamma(1)$,
$$\Phi_{10}(U(\gamma(z))^{T})=\Phi_{10}(j(\gamma, z)^5 (\rho(\gamma) U(z))^{T})
 =j(\gamma, z)^{50} \Phi_{10}((\rho(\gamma) U(z))^{T}).$$
Note that $\rho(\gamma) \in \langle \rho(s), \rho(t) \rangle=G$ and
$\Phi_{10}$ is a $G$-invariant polynomial, we have
$$\Phi_{10}(U(\gamma(z))^{T})=j(\gamma, z)^{50} \Phi_{10}(U(z)^{T}),
  \quad \text{for $\gamma \in \Gamma(1)$}.$$
This implies that $\Phi_{10}(u_1(z), \ldots, u_6(z))$ is a modular form
of weight $50$ for the full modular group $\Gamma(1)$. Moreover, we will
show that it is a cusp form. In fact,
$$\aligned
  &\Phi_{10}(a_1(z), \ldots, a_6(z))=13^3 q^{\frac{9}{4}} (1+O(q))^2+\\
  &+\sum_{\nu=0}^{12} [q^{\frac{1}{2}} (-1+O(q))+\\
  &+2 \zeta^{\nu} q^{\frac{41}{26}} (-3+O(q))+2 \zeta^{3 \nu}
    q^{\frac{19}{26}} (-3+O(q))
   +2 \zeta^{9 \nu} q^{\frac{5}{26}} (1+O(q))+\\
  &+2 \zeta^{12 \nu} q^{\frac{11}{26}} (-1+O(q))+2 \zeta^{10 \nu}
   q^{\frac{7}{26}} (-1+O(q))
   +2 \zeta^{4 \nu} q^{\frac{47}{26}} (-1+O(q))+\\
  &+\zeta^{2 \nu} q^{\frac{17}{26}} (6+O(q))+\zeta^{5 \nu}
   q^{\frac{23}{26}} (8+O(q))
   +\zeta^{6 \nu} q^{\frac{25}{26}} (-1+O(q))+\\
  &+\zeta^{11 \nu} q^{\frac{9}{26}} (-3+O(q))+\zeta^{8 \nu}
   q^{\frac{29}{26}} (2+O(q))
   +\zeta^{7 \nu} q^{\frac{1}{26}} (1+O(q))]\\
  &\times [-13 q^{\frac{7}{4}} (1+O(q))+\\
  &+\zeta^{\nu} q^{\frac{43}{52}} (13+O(q))+\zeta^{2 \nu}
   q^{\frac{47}{52}} (-22+O(q))
   +\zeta^{3 \nu} q^{\frac{51}{52}} (-21+O(q))+\\
  &+\zeta^{4 \nu} q^{\frac{3}{52}} (-1+O(q))+\zeta^{5 \nu}
   q^{\frac{7}{52}} (2+O(q))
   +\zeta^{6 \nu} q^{\frac{11}{52}} (2+O(q))+\\
  &+\zeta^{7 \nu} q^{\frac{15}{52}} (-2+O(q))+\zeta^{8 \nu}
   q^{\frac{19}{52}} (-8+O(q))
   +\zeta^{9 \nu} q^{\frac{23}{52}} (6+O(q))+\\
  &+\zeta^{10 \nu} q^{\frac{27}{52}} (1+O(q))+\zeta^{11 \nu}
   q^{\frac{31}{52}} (-8+O(q))
   +\zeta^{12 \nu} q^{\frac{35}{52}} (17+O(q))]
\endaligned$$
We will calculate the $q^{\frac{1}{4}}$-term which is the lowest
degree:
$$2 q^{\frac{5}{26}} \cdot q^{\frac{3}{52}} \cdot (-1)+
  q^{\frac{1}{26}} \cdot 1 \cdot q^{\frac{11}{52}} \cdot 2=0.$$
Hence, for $\Phi_{10}(u_1(z), \ldots, u_6(z))$ which is a modular form for
$\Gamma(1)$ with weight $50$, the lowest degree term is given by
$$\text{some integer} \cdot q^{\frac{5}{4}} \cdot q^{\frac{9}{24} \cdot 10}
 =\text{some integer} \cdot q^5.$$
This implies that $\Phi_{10}(u_1(z), \ldots, u_6(z))$ has a factor of
$\Delta(z)^5$, which is a cusp form of weight $60$. Therefore,
$\Phi_{10}(u_1(z), \ldots, u_6(z))=0$. Consequently,
$\Phi_{10}(x_1(z), \ldots, x_6(z))=0$. The calculation of $\Phi_{14}$
and $\Phi_{1, 2}$ is similar as that of $\Phi_{10}$.

  For
$\Phi_{2, 2}=w_0^2 \delta_0^2+w_1^2 \delta_1^2+\cdots+w_{12}^2 \delta_{12}^2+
w_{\infty}^2 \delta_{\infty}^2$.
As a polynomial in six variables, $\Phi_{2, 2}(z_1, z_2, z_3, z_4, z_5, z_6)$
is a $G$-invariant polynomial. Moreover, for $\gamma \in \Gamma(1)$,
$$\Phi_{2, 2}(U(\gamma(z))^{T})=\Phi_{2, 2}(j(\gamma, z)^5 (\rho(\gamma) U(z))^{T})
 =j(\gamma, z)^{100} \Phi_{2, 2}((\rho(\gamma) U(z))^{T}).$$
Note that $\rho(\gamma) \in \langle \rho(s), \rho(t) \rangle=G$ and
$\Phi_{2, 2}$ is a $G$-invariant polynomial, we have
$$\Phi_{2, 2}(U(\gamma(z))^{T})=j(\gamma, z)^{100} \Phi_{2, 2}(U(z)^{T}),
  \quad \text{for $\gamma \in \Gamma(1)$}.$$
This implies that $\Phi_{2, 2}(u_1(z), \ldots, u_6(z))$ is a modular form
of weight $100$ for the full modular group $\Gamma(1)$. Moreover, we will
show that it is a cusp form. In fact,
$$\aligned
  &\Phi_{2, 2}(a_1(z), \ldots, a_6(z))=13^6 q^{\frac{9}{2}} (1+O(q))^2+\\
  &+\sum_{\nu=0}^{12} [q^{\frac{1}{2}} (-1+O(q))+\\
  &+2 \zeta^{\nu} q^{\frac{41}{26}} (-3+O(q))+2 \zeta^{3 \nu}
    q^{\frac{19}{26}} (-3+O(q))
   +2 \zeta^{9 \nu} q^{\frac{5}{26}} (1+O(q))+\\
  &+2 \zeta^{12 \nu} q^{\frac{11}{26}} (-1+O(q))+2 \zeta^{10 \nu}
   q^{\frac{7}{26}} (-1+O(q))
   +2 \zeta^{4 \nu} q^{\frac{47}{26}} (-1+O(q))+\\
  &+\zeta^{2 \nu} q^{\frac{17}{26}} (6+O(q))+\zeta^{5 \nu}
   q^{\frac{23}{26}} (8+O(q))
   +\zeta^{6 \nu} q^{\frac{25}{26}} (-1+O(q))+\\
  &+\zeta^{11 \nu} q^{\frac{9}{26}} (-3+O(q))+\zeta^{8 \nu}
   q^{\frac{29}{26}} (2+O(q))
   +\zeta^{7 \nu} q^{\frac{1}{26}} (1+O(q))]^2\\
  &\times [-13 q^{\frac{7}{4}} (1+O(q))+\\
  &+\zeta^{\nu} q^{\frac{43}{52}} (13+O(q))+\zeta^{2 \nu}
   q^{\frac{47}{52}} (-22+O(q))
   +\zeta^{3 \nu} q^{\frac{51}{52}} (-21+O(q))+\\
  &+\zeta^{4 \nu} q^{\frac{3}{52}} (-1+O(q))+\zeta^{5 \nu}
   q^{\frac{7}{52}} (2+O(q))
   +\zeta^{6 \nu} q^{\frac{11}{52}} (2+O(q))+\\
  &+\zeta^{7 \nu} q^{\frac{15}{52}} (-2+O(q))+\zeta^{8 \nu}
   q^{\frac{19}{52}} (-8+O(q))
   +\zeta^{9 \nu} q^{\frac{23}{52}} (6+O(q))+\\
  &+\zeta^{10 \nu} q^{\frac{27}{52}} (1+O(q))+\zeta^{11 \nu}
   q^{\frac{31}{52}} (-8+O(q))
   +\zeta^{12 \nu} q^{\frac{35}{52}} (17+O(q))]^2
\endaligned$$
We will calculate the $q^{\frac{1}{2}}$-term which is the lowest
degree, there are five such terms:

\noindent (1)
$$\aligned
 &(\zeta^{7 \nu} q^{\frac{1}{26}} \cdot 1)^2 \times
  [(\zeta^{6 \nu} q^{\frac{11}{52}} \cdot 2)^2+
  2 \cdot \zeta^{5 \nu} q^{\frac{7}{52}} \cdot 2 \cdot
  \zeta^{7 \nu} q^{\frac{15}{52}} \cdot (-2)+\\
 &+2 \cdot \zeta^{4 \nu} q^{\frac{3}{52}} \cdot (-1) \cdot
  \zeta^{8 \nu} q^{\frac{19}{52}} \cdot (-8)]\\
=&12 q^{\frac{1}{2}}.
\endaligned$$
\noindent (2)
$$2 \cdot 2 \zeta^{9 \nu} q^{\frac{5}{26}} \cdot 1 \cdot \zeta^{7 \nu}
  q^{\frac{1}{26}} \cdot 1
  \times [(\zeta^{5 \nu} q^{\frac{7}{52}} \cdot 2)^2+
  2 \zeta^{4 \nu} q^{\frac{3}{52}} \cdot (-1) \cdot \zeta^{6 \nu}
  q^{\frac{11}{52}} \cdot 2]=0.$$
\noindent (3)
$$2 \cdot 2 \zeta^{10 \nu} q^{\frac{7}{26}} \cdot (-1) \cdot
  \zeta^{7 \nu} q^{\frac{1}{26}} \cdot 1 \times 2 \zeta^{4 \nu}
  q^{\frac{3}{52}} \cdot (-1) \cdot \zeta^{5 \nu} q^{\frac{7}{52}} \cdot 2
 =16 q^{\frac{1}{2}}.$$
\noindent (4)
$$2 \cdot \zeta^{11 \nu} q^{\frac{9}{26}} \cdot (-3) \cdot \zeta^{7 \nu}
  q^{\frac{1}{26}} \cdot 1 \times [\zeta^{4 \nu} q^{\frac{3}{52}} \cdot (-1)]^2
 =-6 q^{\frac{1}{2}}.$$
\noindent (5)
$$(2 \zeta^{9 \nu} q^{\frac{5}{26}} \cdot 1)^2 \times
  (\zeta^{4 \nu} q^{\frac{3}{52}} \cdot (-1))^2=4 q^{\frac{1}{2}}.$$
Hence, for $\Phi_{2, 2}(u_1(z), \ldots, u_6(z))$ which is a modular
form for $\Gamma(1)$ with weight $100$, the lowest degree term is given by
$$(12+0+16-6+4) q^{\frac{1}{2}} \cdot q^{\frac{9}{24} \cdot 20}=26 q^8.$$
Thus,
$$\Phi_{2, 2}(u_1(z), \ldots, u_6(z))=q^8 (13 \cdot 26+O(q)).$$
The leading term of $\Phi_{2, 2}(u_1(z), \ldots, u_6(z))$ together
with its weight $100$ suffice to identify this modular form with
$\Phi_{2, 2}(u_1(z), \ldots, u_6(z))=13 \cdot 26 \Delta(z)^8 E_4(z)$.
Consequently,
$$\Phi_{2, 2}(x_1(z), \ldots, x_6(z))=13 \cdot 26 \Delta(z)^8
  E_4(z)/\eta(z)^{160}=13 \cdot 26 \eta(z)^8 \Delta(z) E_4(z).$$
The calculation of $\Phi_{3, 3}$ and $\Phi_{6, 1}$ is similar as
that of $\Phi_{2, 2}$. We have
$$\Phi_{3, 3}(x_1(z), \ldots, x_6(z))=-13 \cdot 27 \Delta(z)^2 E_6(z).$$
$$\Phi_{6, 1}(x_1(z), \ldots, x_6(z))=-13 \cdot 285 \Delta(z)^2 E_6(z).$$
After normalization, this completes the proof of Theorem 9.1.

\noindent
$\qquad \qquad \qquad \qquad \qquad \qquad \qquad \qquad \qquad
 \qquad \qquad \qquad \qquad \qquad \qquad \qquad \boxed{}$

  Now we give the normalization for the following four families of
invariant homogeneous polynomials of degrees $d=18$, $22$, $26$ and
$34$. For $d=18$, there are two such invariant homogeneous polynomials
$\Phi_{3, 1}$ and $\Phi_{0, 3}$:
$$\left\{\aligned
  \Phi_{3, 1} &:=\frac{1}{13 \cdot 2} \left(\sum_{\nu=0}^{12}
  w_{\nu}^3 \delta_{\nu} +w_{\infty}^3 \delta_{\infty}\right),\\
  \Phi_{0, 3} &:=\frac{1}{13 \cdot 6} \left(\sum_{\nu=0}^{12}
                 \delta_{\nu}^3+\delta_{\infty}^3\right),
\endaligned\right.\eqno{(9.6)}$$
For $d=22$, there are two such invariant homogeneous polynomials
$\Phi_{1, 3}$ and $\Phi_{4, 1}$. For $d=26$, there are two such
invariant homogeneous polynomials $\Phi_{5, 1}$ and $\Phi_{2, 3}$:
$$\left\{\aligned
  \Phi_{5, 1} &:=-\frac{1}{13} \left(\sum_{\nu=0}^{12}
  w_{\nu}^5 \delta_{\nu}+w_{\infty}^5 \delta_{\infty}\right),\\
  \Phi_{2, 3} &:=-\frac{1}{13} \left(\sum_{\nu=0}^{12}
  w_{\nu}^2 \delta_{\nu}^3+w_{\infty}^2 \delta_{\infty}^3\right),
\endaligned\right.\eqno{(9.7)}$$
For $d=34$, there are three such invariant homogeneous polynomials
$\Phi_{1, 5}$, $\Phi_{4, 3}$ and $\Phi_{7, 1}$.

\textbf{Theorem 9.2.} {\it The $G$-invariant homogeneous polynomials
$\Phi_{m, n}$ of degrees $d=18$, $22$, $26$ and $34$ in $x_1(z)$,
$\ldots$, $x_6(z)$ can be identified with modular forms as follows$:$
$$\left\{\aligned
  \Phi_{3, 1}(x_1(z), \ldots, x_6(z)) &=\Delta(z) E_6(z),\\
  \Phi_{0, 3}(x_1(z), \ldots, x_6(z)) &=\Delta(z) E_6(z),\\
  \Phi_{1, 3}(x_1(z), \ldots, x_6(z)) &=0,\\
  \Phi_{4, 1}(x_1(z), \ldots, x_6(z)) &=0,\\
  \Phi_{5, 1}(x_1(z), \ldots, x_6(z)) &=\eta(z)^8 \Delta(z) E_4(z) E_6(z),\\
  \Phi_{2, 3}(x_1(z), \ldots, x_6(z)) &=\eta(z)^8 \Delta(z) E_4(z) E_6(z),\\
  \Phi_{1, 5}(x_1(z), \ldots, x_6(z)) &=0,\\
  \Phi_{4, 3}(x_1(z), \ldots, x_6(z)) &=0,\\
  \Phi_{7, 1}(x_1(z), \ldots, x_6(z)) &=0.\\
\endaligned\right.\eqno{(9.8)}$$}

{\it Proof}. The proof is similar as that of Theorem 9.1.

\noindent
$\qquad \qquad \qquad \qquad \qquad \qquad \qquad \qquad \qquad
 \qquad \qquad \qquad \qquad \qquad \qquad \qquad \boxed{}$

  As a consequence of Theorem 3.3 and Theorem 3.4, we find an image
of the modular curve $X=X(13)$ in $\mathbb{CP}^5$ given by a complete
intersection curve as follows:

\textbf{Theorem 9.3.} {\it There is a morphism
$$\Phi: X \rightarrow Z \subset \mathbb{CP}^5$$
with $\Phi(z)=(x_1(z), \ldots, x_6(z))$, where $Z$ is a complete intersection
curve corresponding to the ideal
$$Z=(\Phi_4, \Phi_8, \Phi_{10}, \Phi_{14}).\eqno{(9.9)}$$
Moreover, there is a morphism
$\Psi: X \to W \subset \mathbb{CP}^5$, where $W$ is the scheme
corresponding to the ideal
$$J=J(W)=(\Phi_4, \Phi_8, \Phi_{10}, \Phi_{14}, \Phi_{4, 0}, \Phi_{1, 2},
         \Phi_{1, 3}, \Phi_{4, 1}, \Phi_{1, 5}, \Phi_{4, 3}, \Phi_{7, 1}).
         \eqno{(9.10)}$$}

{\it Proof}. Theorem 9.1 and Theorem 9.2 imply that
$$\left\{\aligned
  \Phi_{4}(x_1(z), \ldots, x_6(z)) &=0,\\
  \Phi_{8}(x_1(z), \ldots, x_6(z)) &=0,\\
  \Phi_{10}(x_1(z), \ldots, x_6(z)) &=0,\\
  \Phi_{14}(x_1(z), \ldots, x_6(z)) &=0,\\
  \Phi_{4, 0}(x_1(z), \ldots, x_6(z)) &=0,\\
  \Phi_{1, 2}(x_1(z), \ldots, x_6(z)) &=0,\\
  \Phi_{1, 3}(x_1(z), \ldots, x_6(z)) &=0,\\
  \Phi_{4, 1}(x_1(z), \ldots, x_6(z)) &=0,\\
  \Phi_{1, 5}(x_1(z), \ldots, x_6(z)) &=0,\\
  \Phi_{4, 3}(x_1(z), \ldots, x_6(z)) &=0,\\
  \Phi_{7, 1}(x_1(z), \ldots, x_6(z)) &=0.
\endaligned\right.$$
This complete the proof of Theorem 9.3.

\noindent
$\qquad \qquad \qquad \qquad \qquad \qquad \qquad \qquad \qquad
 \qquad \qquad \qquad \qquad \qquad \qquad \qquad \boxed{}$

\begin{center}
{\large\bf 10. Modularity for equations of $E_6$, $E_7$ and $E_8$-singularities
              coming from $C_Y/\text{SL}(2, 13)$ and variations of $E_6$, $E_7$
              and $E_8$-singularity structures over $X(13)$}
\end{center}

  In this section, we will study the ring of invariant polynomials
$$\left(\mathbb{C}[z_1, z_2, z_3, z_4, z_5, z_6]/I(Y)\right)^{\text{SL}(2, 13)}.$$
By Theorem 6.2, there are three invariant cones $C_Y/\text{SL}(2, 13)$,
$C_{Y_2}/\text{SL}(2, 13)$ and $C_{Y_3}/\text{SL}(2, 13)$ over the curve
$Y$. This leads to a new perspective on the theory of singularities. In
particular, it gives a different construction of the $E_6$, $E_7$ and
$E_8$-singularities from a quotient $C_Y/G$ over the modular curve $X$.
Moreover, they have a simultaneous modular parametrizations in a unified
way.

\textbf{Theorem 10.1.} (Variation structures of the $E_6$, $E_7$ and
$E_8$-singularities over the modular curve $X$: algebraic version)
{\it The equations of $E_6$, $E_7$ and $E_8$-singularities
$$E_6: \quad \left(\frac{\Phi_{20}^3}{\Phi_{12}^4}+1728 \Phi_{12}\right)^2
  -\left(\frac{\Phi_{26}}{\Phi_{20}}\right)^4-4 \cdot 1728
  \left(\frac{\Phi_{20}}{\Phi_{12}}\right)^3=0,\eqno{(10.1)}$$
$$E_7: \quad \Phi_{12} \cdot \left(\frac{\Phi_{26}}{\Phi_{18}}\right)^3
  -\Phi_{18}^2-1728 \Phi_{12}^3=0,\eqno{(10.2)}$$
$$E_8: \quad \Phi_{20}^3-\left(\frac{\Phi_{12}^2 \Phi_{26}}{\Phi_{20}}\right)^2
  -1728 \Phi_{12}^5=0\eqno{(10.3)}$$
possess infinitely many kinds of distinct modular parametrizations
$($with the cardinality of the continuum in ZFC set theory$)$
$$(\Phi_{12}, \Phi_{18}, \Phi_{20}, \Phi_{26})=(\Phi_{12}^{\lambda},
   \Phi_{18}^{\mu}, \Phi_{20}^{\gamma}, \Phi_{26}^{\kappa})\eqno{(10.4)}$$
over the modular curve $X$ as follows$:$
$$\left\{\aligned
  \Phi_{12}^{\lambda} &=\lambda \Phi_{3, 0}+(1-\lambda) \Phi_{0, 2}
              \quad \text{mod $\mathfrak{a}_1$},\\
  \Phi_{18}^{\mu} &=\mu \Phi_{3, 1}+(1-\mu) \Phi_{0, 3}
              \quad \text{mod $\mathfrak{a}_2$},\\
  \Phi_{20}^{\gamma} &=\gamma \Phi_{5, 0}+(1-\gamma) \Phi_{2, 2}
              \quad \text{mod $\mathfrak{a}_3$},\\
  \Phi_{26}^{\kappa} &=\kappa \Phi_{2, 3}+(1-\kappa) \Phi_{5, 1}
             \quad \text{mod $\mathfrak{a}_4$},
\endaligned\right.\eqno{(10.5)}$$
where $\Phi_{12}$, $\Phi_{18}$, $\Phi_{20}$ and $\Phi_{30}$ are invariant
homogeneous polynomials of degree $12$, $18$, $20$ and $26$, respectively.
The ideals are given by
$$\left\{\aligned
  \mathfrak{a}_1 &=(\Phi_4, \Phi_8),\\
  \mathfrak{a}_2 &=(\Phi_4, \Phi_8, \Phi_{10}, \Phi_{14}),\\
  \mathfrak{a}_3 &=(\Phi_4, \Phi_8, \Phi_{10}, \Phi_{4, 0}, \Phi_{1, 2}),\\
  \mathfrak{a}_4 &=(\Phi_4, \Phi_8, \Phi_{10}, \Phi_{14}, \Phi_{4, 0},
                    \Phi_{1, 2}, \Phi_{1, 3}, \Phi_{4, 1}),
\endaligned\right.\eqno{(10.6)}$$
and the parameter space $\{ (\lambda, \mu, \gamma, \kappa) \} \cong \mathbb{C}^4$.
They form variations of the $E_6$, $E_7$ and $E_8$-singularity structures
over the modular curve $X$.}

{\it Proof}. By Theorem 9.3, the ideals $\mathfrak{a}_1$, $\mathfrak{a}_2$,
$\mathfrak{a}_3$ and $\mathfrak{a}_4$ are zero ideals over the modular curve
$X$. On the other hand, by Theorem 9.1 and Theorem 9.2, for degree $d=12$,
the invariant homogeneous polynomials $\Phi_{3, 0}$ and $\Phi_{0, 2}$
form a two-dimensional complex vector space, and
$\Phi_{3, 0}=\Phi_{0, 2}=\Delta(z)$ over the modular curve $X$ (after
normalization). Hence,
$$\Phi_{12}^{\lambda}=\Delta(z)\eqno{(10.7)}$$
over the modular curve $X$. For degree $d=18$, the invariant homogeneous
polynomials $\Phi_{3, 1}$ and $\Phi_{0, 3}$ form a two-dimensional complex
vector space, and $\Phi_{3, 1}=\Phi_{0, 3}=\Delta(z) E_6(z)$ over the modular
curve $X$ (after normalization). Hence,
$$\Phi_{18}^{\mu}=\Delta(z) E_6(z)\eqno{(10.8)}$$
over the modular curve $X$. For degree $d=20$, the invariant homogeneous
polynomials $\Phi_{5, 0}$ and $\Phi_{2, 2}$ form a two-dimensional complex
vector space, and $\Phi_{5, 0}=\Phi_{2, 2}=\eta(z)^8 \Delta(z) E_4(z)$ over
the modular curve $X$ (after normalization). Hence,
$$\Phi_{20}^{\gamma}=\eta(z)^8 \Delta(z) E_4(z)\eqno{(10.9)}$$
over the modular curve $X$. Finally, for degree $d=26$, the invariant homogeneous
polynomials $\Phi_{5, 1}$ and $\Phi_{2, 3}$ form a two-dimensional complex vector
space, and $\Phi_{5, 1}=\Phi_{2, 3}=\eta(z)^8 \Delta(z) E_4(z) E_6(z)$ over the
modular curve $X$ (after normalization). Hence,
$$\Phi_{26}^{\kappa}=\eta(z)^8 \Delta(z) E_4(z) E_6(z)\eqno{(10.10)}$$
over the modular curve $X$.

  Combining with (10.7), (10.8), (10.9) and (10.10), a straightforward
calculation shows that the equations of $E_6$, $E_7$ and $E_8$-singularities
can be parametrized by the invariant polynomials $\Phi_{12}$, $\Phi_{18}$,
$\Phi_{20}$ and $\Phi_{26}$:
$$\quad \left(\frac{\Phi_{20}^3}{\Phi_{12}^4}+1728 \Phi_{12}\right)^2
  -\left(\frac{\Phi_{26}}{\Phi_{20}}\right)^4-4 \cdot 1728
  \left(\frac{\Phi_{20}}{\Phi_{12}}\right)^3=0,$$
$$\quad \Phi_{12} \cdot \left(\frac{\Phi_{26}}{\Phi_{18}}\right)^3
  -\Phi_{18}^2-1728 \Phi_{12}^3=0,$$
$$\Phi_{20}^3-\left(\frac{\Phi_{12}^2 \Phi_{26}}{\Phi_{20}}\right)^2
  -1728 \Phi_{12}^5=0$$
over the modular curve $X$.

\noindent
$\qquad \qquad \qquad \qquad \qquad \qquad \qquad \qquad \qquad
 \qquad \qquad \qquad \qquad \qquad \qquad \qquad \boxed{}$

\textbf{Theorem 10.2.} (Variation structures of the $E_6$, $E_7$ and
$E_8$-singularities over the modular curve $X$: geometric version)
{\it There are three morphisms from the cone $C_Y$ over $Y$ to the
$E_6$, $E_7$ and $E_8$-singularities:
$$\aligned
  &f_6: C_Y/G \rightarrow\\
  &\text{Spec}\left(\mathbb{C}\{\Phi_{12}, \Phi_{20}\}
  [\Phi_{26}]/((\frac{\Phi_{20}^3}{\Phi_{12}^4}+1728 \Phi_{12})^2
  -(\frac{\Phi_{26}}{\Phi_{20}})^4-4 \cdot 1728
  (\frac{\Phi_{20}}{\Phi_{12}})^3)\right).
\endaligned\eqno{(10.11)}$$
$$f_7: C_Y/G \rightarrow \text{Spec}\left(\mathbb{C}\{\Phi_{18}\}
  [\Phi_{12}, \Phi_{26}]/((\Phi_{12} \cdot(\frac{\Phi_{26}}{\Phi_{18}})^3
  -\Phi_{18}^2-1728 \Phi_{12}^3)\right).\eqno{(10.12)}$$
$$f_8: C_Y/G \rightarrow \text{Spec}\left(\mathbb{C}\{\Phi_{20}\}
  [\Phi_{12}, \Phi_{26}]/(\Phi_{20}^3-(\frac{\Phi_{12}^2 \Phi_{26}}{\Phi_{20}})^2
  -1728 \Phi_{12}^5)\right).\eqno{(10.13)}$$
over the modular curve $X$. In particular, there are infinitely many
such quadruples $(\Phi_{12}, \Phi_{18}, \Phi_{20}, \Phi_{26})$ $=$
$(\Phi_{12}^{\lambda}, \Phi_{18}^{\mu}, \Phi_{20}^{\gamma}, \Phi_{26}^{\kappa})$
whose parameter space $\{ (\lambda, \mu, \gamma, \kappa) \}$
$\cong \mathbb{C}^4$. They form variations of the $E_6$, $E_7$ and
$E_8$-singularity structures over the modular curve $X$.}

{\it Proof}. By Theorem 9.1, Theorem 9.2, Theorem 9.3 and Theorem 10.1,
we have the following three ring homomorphisms
$$\aligned
 &\mathbb{C}\{\Phi_{12}, \Phi_{20}\}[\Phi_{26}]/
  ((\frac{\Phi_{20}^3}{\Phi_{12}^4}+1728 \Phi_{12})^2
  -(\frac{\Phi_{26}}{\Phi_{20}})^4-4 \cdot 1728
  (\frac{\Phi_{20}}{\Phi_{12}})^3)\\
 &\rightarrow \left(\mathbb{C}[z_1, z_2, z_3, z_4, z_5, z_6]/I(Y)
  \right)^{\text{SL}(2, 13)}
\endaligned$$
$$\aligned
 &\mathbb{C}\{\Phi_{18}\}[\Phi_{12}, \Phi_{26}]/
  (\Phi_{12} \cdot(\frac{\Phi_{26}}{\Phi_{18}})^3-\Phi_{18}^2-1728 \Phi_{12}^3)\\
 &\rightarrow \left(\mathbb{C}[z_1, z_2, z_3, z_4, z_5, z_6]/I(Y)
  \right)^{\text{SL}(2, 13)}
\endaligned$$
$$\aligned
 &\mathbb{C}\{\Phi_{20}\}[\Phi_{12}, \Phi_{26}]/
  (\Phi_{20}^3-(\frac{\Phi_{12}^2 \Phi_{26}}{\Phi_{20}})^2-1728 \Phi_{12}^5)\\
 &\rightarrow \left(\mathbb{C}[z_1, z_2, z_3, z_4, z_5, z_6]/I(Y)
  \right)^{\text{SL}(2, 13)}
\endaligned$$
over the modular curve $X$, which induce three morphisms of schemes
$$\aligned
  &f_6: C_Y/G \rightarrow\\
  &\text{Spec}\left(\mathbb{C}\{\Phi_{12}, \Phi_{20}\}
  [\Phi_{26}]/((\frac{\Phi_{20}^3}{\Phi_{12}^4}+1728 \Phi_{12})^2
  -(\frac{\Phi_{26}}{\Phi_{20}})^4-4 \cdot 1728
  (\frac{\Phi_{20}}{\Phi_{12}})^3)\right),
\endaligned$$
$$f_7: C_Y/G \rightarrow \text{Spec}\left(\mathbb{C}\{\Phi_{18}\}
  [\Phi_{12}, \Phi_{26}]/((\Phi_{12} \cdot(\frac{\Phi_{26}}{\Phi_{18}})^3
  -\Phi_{18}^2-1728 \Phi_{12}^3)\right),$$
$$f_8: C_Y/G \rightarrow \text{Spec}\left(\mathbb{C}\{\Phi_{20}\}
  [\Phi_{12}, \Phi_{26}]/(\Phi_{20}^3-(\frac{\Phi_{12}^2 \Phi_{26}}{\Phi_{20}})^2
  -1728 \Phi_{12}^5)\right)$$
over $X$.

  At the standard cusp $z=i \infty$, $X(13)$ has local coordinate
$q^{\frac{1}{13}}$. Near a cusp, for this to extend at the cusp,
$x_i(z)$ has to be replaced by $q^{-\frac{2}{39}} x_i(z)$
$(=q^{-\frac{1}{104}} a_i(z))$ for $1 \leq i \leq 6$. Then near
the cusp $z=i \infty$, we have the following:

\noindent (1) For $\Phi_{12}(x_1(z), \ldots, x_6(z))$, the
leading term is given by
$$q \cdot (q^{-\frac{2}{39}})^{12}=q^{\frac{5}{13}}.$$

\noindent (2) For $\Phi_{18}(x_1(z), \ldots, x_6(z))$, the
leading term is given by
$$q \cdot (q^{-\frac{2}{39}})^{18}=q^{\frac{1}{13}}.$$

\noindent (3) For $\Phi_{20}(x_1(z), \ldots, x_6(z))$, the
leading term is given by
$$q^{\frac{4}{3}} \cdot (q^{-\frac{2}{39}})^{20}=q^{\frac{4}{13}}.$$

\noindent (4) For $\Phi_{26}(x_1(z), \ldots, x_6(z))$, the
leading term is given by
$$q^{\frac{4}{3}} \cdot (q^{-\frac{2}{39}})^{26}=1.$$

  Hence, this shows that $\Phi_{26}$ gives functions not vanishing
on the rays of the cusp $z=i \infty$, while $\Phi_{12}$, $\Phi_{18}$
and $\Phi_{20}$ vanish exactly on the rays of the cusp $z=i \infty$.
Moreover, the triple $(\Phi_{12}, \Phi_{20}, \Phi_{26})$ of invariant
polynomials defines two maps from the cone $C_Y/G$ over $X(13)$ to the
equations $x^4+y^3+z^2=0$ and $x^5+y^3+z^2=0$. On the other hand, the
triple $(\Phi_{12}, \Phi_{18}, \Phi_{26})$ defines a map from the cone
$C_Y/G$ over $X(13)$ to the equation $x^3 y+y^3+z^2=0$. They map to
$(0, 0, 1)$ (after normalization) the generatrix of this cone
corresponding to the cusps. Hence, they give relations between the link
at $O$ for the cone $C_Y/G$ over $X(13)$ and the links at $O$ for $E_6$,
$E_7$ and $E_8$-singularities, which are given as follows:

  The link at $O$ for the $E_6$-singularity by the triple $(\Phi_{12},
\Phi_{20}, \Phi_{26})$:
$$\left\{\aligned
  &|\Phi_{12}|^2+|\Phi_{20}|^2+|\Phi_{26}|^2=1,\\
  &\left(\frac{\Phi_{20}^3}{\Phi_{12}^4}+1728 \Phi_{12}\right)^2
  -\left(\frac{\Phi_{26}}{\Phi_{20}}\right)^4-4 \cdot 1728
  \left(\frac{\Phi_{20}}{\Phi_{12}}\right)^3=0.
\endaligned\right.$$

  The link at $O$ for the $E_7$-singularity by the triple $(\Phi_{12},
\Phi_{18}, \Phi_{26})$:
$$\left\{\aligned
  &|\Phi_{12}|^2+|\Phi_{18}|^2+|\Phi_{26}|^2=1,\\
  &\Phi_{12} \cdot \left(\frac{\Phi_{26}}{\Phi_{18}}\right)^3
   -\Phi_{18}^2-1728 \Phi_{12}^3=0.
\endaligned\right.$$

  The link at $O$ for the $E_8$-singularity by the triple $(\Phi_{12},
\Phi_{20}, \Phi_{26})$:
$$\left\{\aligned
  &|\Phi_{12}|^2+|\Phi_{20}|^2+|\Phi_{26}|^2=1,\\
  &\Phi_{20}^3-\left(\frac{\Phi_{12}^2 \Phi_{26}}{\Phi_{20}}\right)^2
   -1728 \Phi_{12}^5=0.
\endaligned\right.$$
This completes the proof of Theorem 10.2.

\noindent
$\qquad \qquad \qquad \qquad \qquad \qquad \qquad \qquad \qquad
 \qquad \qquad \qquad \qquad \qquad \qquad \qquad \boxed{}$

\textbf{Theorem 10.3.} (Variations of the structure of decomposition
formulas of the elliptic modular functions $j$ over the modular curve
$X$) {\it There are infinitely many kinds of distinct decomposition
formulas of the elliptic modular function $j$ in terms of the invariants
$\Phi_{12}$, $\Phi_{18}$, $\Phi_{20}$ and $\Phi_{26}$ over the modular
curve $X$:
$$E_7\text{-type}: \quad
  j(z): j(z)-1728: 1=\Phi_{12} \left(\frac{\Phi_{26}}{\Phi_{18}}\right)^3:
  \Phi_{18}^2: \Phi_{12}^3,\eqno{(10.14)}$$
$$E_8\text{-type}: \quad
  j(z): j(z)-1728: 1=\Phi_{20}^3: \left(\frac{\Phi_{26} \Phi_{12}^2}
  {\Phi_{20}}\right)^2: \Phi_{12}^5,\eqno{(10.15)}$$
where $(\Phi_{12}, \Phi_{18}, \Phi_{20}, \Phi_{26})=(\Phi_{12}^{\lambda},
\Phi_{18}^{\mu}, \Phi_{20}^{\gamma}, \Phi_{26}^{\kappa})$ are given
by $(10.5)$. They form variations of the structure of decomposition
formulas of the elliptic modular functions $j$ over the modular curve
$X$.}

{\it Proof}. By Theorem 10.1, we have the following relations over the
modular curve $X$:
$$j(z)=\frac{E_4(z)^3}{\Delta(z)}=\frac{\Phi_{26}^3}{\Phi_{18}^3 \Phi_{12}^2}, \quad
  j(z)-1728=\frac{E_6(z)^2}{\Delta(z)}=\frac{\Phi_{18}^2}{\Phi_{12}^3},
  \eqno{(10.16)}$$
or
$$j(z)=\frac{E_4(z)^3}{\Delta(z)}=\frac{\Phi_{20}^3}{\Phi_{12}^5}, \quad
  j(z)-1728=\frac{E_6(z)^2}{\Delta(z)}=\frac{\Phi_{26}^2}{\Phi_{20}^2 \Phi_{12}}.
  \eqno{(10.17)}$$
Hence, we obtain an infinitely many kinds of distinct decomposition formulas
of the elliptic modular function $j$ in terms of the invariants $\Phi_{12}$,
$\Phi_{18}$, $\Phi_{20}$ and $\Phi_{26}$ over the modular curve $X$:
$$j(z): j(z)-1728: 1=\Phi_{12} \left(\frac{\Phi_{26}}{\Phi_{18}}\right)^3:
  \Phi_{18}^2: \Phi_{12}^3,$$
or
$$j(z): j(z)-1728: 1=\Phi_{20}^3: \left(\frac{\Phi_{26} \Phi_{12}^2}
  {\Phi_{20}}\right)^2: \Phi_{12}^5,$$
where $(\Phi_{12}, \Phi_{18}, \Phi_{20}, \Phi_{26})=(\Phi_{12}^{\lambda},
\Phi_{18}^{\mu}, \Phi_{20}^{\gamma}, \Phi_{26}^{\kappa})$ are given by
(10.5). This completes the proof of Theorem 10.3.

\noindent
$\qquad \qquad \qquad \qquad \qquad \qquad \qquad \qquad \qquad
 \qquad \qquad \qquad \qquad \qquad \qquad \qquad \boxed{}$

  Theorem 10.3 shows that there are infinitely many kinds of distinct
decomposition formulas of the elliptic modular function $j$ in terms of
the invariant polynomials over the modular curves, one and only one is
given by $\text{SL}(2, 5)$ corresponding to the modular curve $X(5)$,
the other infinitely many kinds of distinct decomposition formulas
(which form variations of the structure of decomposition formulas) are
given by $\text{SL}(2, 13)$ corresponding to the modular curve $X(13)$.

  In the end, let us recall some facts about exotic spheres (see \cite{Hi}).
A $k$-dimensional compact oriented differentiable manifold is called a
$k$-sphere if it is homeomorphic to the $k$-dimensional standard sphere.
A $k$-sphere not diffeomorphic to the standard $k$-sphere is said to be
exotic. The first exotic sphere was discovered by Milnor in 1956 (see
\cite{Mi}). Two $k$-spheres are called equivalent if there exists an
orientation preserving diffeomorphism between them. The equivalence classes
of $k$-spheres constitute for $k \geq 5$ a finite abelian group $\Theta_k$
under the connected sum operation. $\Theta_k$ contains the subgroup $b P_{k+1}$
of those $k$-spheres which bound a parallelizable manifold. $b P_{4m}$ ($m \geq 2$)
is cyclic of order $2^{2m-2}(2^{2m-1}-1)$ numerator $(4 B_m/m)$, where $B_m$
is the $m$-th Bernoulli number. Let $g_m$ be the Milnor generator of $b P_{4m}$.
If a $(4m-1)$-sphere $\Sigma$ bounds a parallelizable manifold $B$ of dimension
$4m$, then the signature $\tau(B)$ of the intersection form of $B$ is divisible
by $8$ and $\Sigma=\frac{\tau(B)}{8} g_m$. For $m=2$ we have
$b P_8=\Theta_7=\mathbb{Z} /28 \mathbb{Z}$. All these results are due to
Milnor-Kervaire (see \cite{KM}). In particular,
$$\sum_{i=0}^{2m} z_i \overline{z_i}=1, \quad
  z_0^3+z_1^{6k-1}+z_2^2+\cdots+z_{2m}^2=0$$
is a $(4m-1)$-sphere embedded in $S^{4m+1} \subset \mathbb{C}^{2n+1}$ which
represents the element $(-1)^m k \cdot g_m \in b P_{4m}$. For $m=2$ and
$k=1, 2, \cdots, 28$ we get the $28$ classes of $7$-spheres. Theorem 10.1 and
Theorem 10.2 show that the infinitely many kinds of distinct constructions of
the $E_8$-singularity: $\mathbb{C}^2/\text{SL}(2, 5)$ and a variation of the
$E_8$-singularity structure
$$C_Y/\text{SL}(2, 13) \rightarrow \text{Spec}\left(\mathbb{C}\{\Phi_{20}\}
  [\Phi_{12}, \Phi_{26}]/(\Phi_{20}^3-(\frac{\Phi_{12}^2 \Phi_{26}}{\Phi_{20}})^2
  -1728 \Phi_{12}^5)\right)$$
over the modular curve $X$ give the same singularity structure but with
different links.

\begin{center}
{\large\bf 11. Modularity for equations of $Q_{18}$ and $E_{20}$-singularities
              coming from $C_Y/\text{SL}(2, 13)$ and variations of $Q_{18}$
              and $E_{20}$-singularity structures over $X(13)$}
\end{center}

  In this section, we will extend our work on variations of the $E_6$,
$E_7$ and $E_8$-singularity structures over the modular curve $X$ to the
cases of $Q_{18}$ and $E_{20}$-singularities and obtain variations of
$Q_{18}$ and $E_{20}$-singularity structures over the modular curve $X$.

  Now we give the normalization for the following three families of
invariant homogeneous polynomials of degrees $d=32$, $42$ and $44$.
For $d=32$, there are three such invariant homogeneous polynomials
$\Phi_{8, 0}$, $\Phi_{5, 2}$ and $\Phi_{2, 4}$:
$$\left\{\aligned
  \Phi_{8, 0} &:=-\frac{1}{13 \cdot 1840} \left(\sum_{\nu=0}^{12} w_{\nu}^8
                 +w_{\infty}^8\right),\\
  \Phi_{5, 2} &:=-\frac{1}{13 \cdot 2064} \left(\sum_{\nu=0}^{12} w_{\nu}^5
                 \delta_{\nu}^2+w_{\infty}^5 \delta_{\infty}^2\right),\\
  \Phi_{2, 4} &:=-\frac{1}{13 \cdot 680} \left(\sum_{\nu=0}^{12} w_{\nu}^2
                 \delta_{\nu}^4+w_{\infty}^2 \delta_{\infty}^4\right).
\endaligned\right.\eqno{(11.1)}$$
For $d=42$, there are four such invariant homogeneous polynomials
$\Phi_{0, 7}$, $\Phi_{3, 5}$, $\Phi_{6, 3}$ and $\Phi_{9, 1}$:
$$\left\{\aligned
  \Phi_{0, 7} &:=\frac{1}{13 \cdot 226842} \left(\sum_{\nu=0}^{12}
                 \delta_{\nu}^7+\delta_{\infty}^7\right),\\
  \Phi_{3, 5} &:=\frac{1}{13 \cdot 634} \left(\sum_{\nu=0}^{12} w_{\nu}^3
                 \delta_{\nu}^5+w_{\infty}^3 \delta_{\infty}^5\right),\\
  \Phi_{6, 3} &:=\frac{1}{13 \cdot 10656} \left(\sum_{\nu=0}^{12} w_{\nu}^6
                 \delta_{\nu}^3+w_{\infty}^6 \delta_{\infty}^3\right),\\
  \Phi_{9, 1} &:=\frac{1}{13 \cdot 39134} \left(\sum_{\nu=0}^{12} w_{\nu}^9
                 \delta_{\nu}+w_{\infty}^9 \delta_{\infty}\right).
\endaligned\right.\eqno{(11.2)}$$
For $d=44$, there are four such invariant homogeneous polynomials
$\Phi_{11, 0}$, $\Phi_{8, 2}$, $\Phi_{5, 4}$ and $\Phi_{2, 6}$. Let
$$\Phi_{11, 0}:=\frac{1}{13 \cdot 146905} \left(\sum_{\nu=0}^{12} w_{\nu}^{11}
                +w_{\infty}^{11}\right).\eqno{(11.3)}$$

\textbf{Theorem 11.1.} {\it The $G$-invariant homogeneous polynomials
$\Phi_{m, n}$ of degrees $d=32$, $42$ and $44$ in $x_1(z)$, $\ldots$,
$x_6(z)$ can be identified with modular forms as follows$:$
$$\left\{\aligned
  \Phi_{8, 0}(x_1(z), \ldots, x_6(z)) &=\eta(z)^8 \Delta(z)^2 E_4(z),\\
  \Phi_{5, 2}(x_1(z), \ldots, x_6(z)) &=\eta(z)^8 \Delta(z)^2 E_4(z),\\
  \Phi_{2, 4}(x_1(z), \ldots, x_6(z)) &=\eta(z)^8 \Delta(z)^2 E_4(z),\\
  \Phi_{0, 7}(x_1(z), \ldots, x_6(z)) &=\Delta(z)^3 E_6(z),\\
  \Phi_{3, 5}(x_1(z), \ldots, x_6(z)) &=\Delta(z)^3 E_6(z),\\
  \Phi_{6, 3}(x_1(z), \ldots, x_6(z)) &=\Delta(z)^3 E_6(z),\\
  \Phi_{9, 1}(x_1(z), \ldots, x_6(z)) &=\Delta(z)^3 E_6(z),\\
  \Phi_{11, 0}(x_1(z), \ldots, x_6(z)) &=\eta(z)^8 \Delta(z)^3 E_4(z),\\
  \Phi_{8, 2}(x_1(z), \ldots, x_6(z)) &\in \eta(z)^8 \Delta(z)^2
  (\mathbb{C} E_4(z)^4 \oplus \mathbb{C} E_4(z) E_6(z)^2),\\
  \Phi_{5, 4}(x_1(z), \ldots, x_6(z)) &\in \eta(z)^8 \Delta(z)^2
  (\mathbb{C} E_4(z)^4 \oplus \mathbb{C} E_4(z) E_6(z)^2),\\
  \Phi_{2, 6}(x_1(z), \ldots, x_6(z)) &\in \eta(z)^8 \Delta(z)^2
  (\mathbb{C} E_4(z)^4 \oplus \mathbb{C} E_4(z) E_6(z)^2).
\endaligned\right.\eqno{(11.4)}$$}

{\it Proof}. The proof is similar as that of Theorem 9.1. For the
calculation of $\Phi_{8, 0}$, $\Phi_{0, 7}$ and $\Phi_{11, 0}$,
see also \cite{Y3}, section 4.

\noindent
$\qquad \qquad \qquad \qquad \qquad \qquad \qquad \qquad \qquad
 \qquad \qquad \qquad \qquad \qquad \qquad \qquad \boxed{}$

\textbf{Theorem 11.2.} (Variations of $Q_{18}$ and $E_{20}$-singularity
structures over the modular curve $X$: algebraic version) {\it The equations
of $Q_{18}$ and $E_{20}$-singularities
$$\Phi_{32}^3-\Phi_{12} \left(\frac{\Phi_{12}^4 \Phi_{26}}{\Phi_{32}}\right)^2
  -1728 \Phi_{12}^8=0, \quad
  \Phi_{44}^3-\left(\frac{\Phi_{12}^7 \Phi_{26}}{\Phi_{44}}\right)^2
  -1728 \Phi_{12}^{11}=0\eqno{(11.5)}$$
possess infinitely many kinds of distinct modular parametrizations
$($with the cardinality of the continuum in ZFC set theory$)$
$$(\Phi_{12}, \Phi_{26}, \Phi_{32}, \Phi_{44})=(\Phi_{12}^{\lambda},
   \Phi_{26}^{\mu}, \Phi_{32}^{\gamma}, \Phi_{44})$$
over the modular curve $X$ as follows$:$
$$\left\{\aligned
  \Phi_{12}^{\lambda} &=\lambda \Phi_{3, 0}+(1-\lambda) \Phi_{0, 2}
              \quad \text{mod $\mathfrak{a}_1$},\\
  \Phi_{26}^{\mu} &=\mu \Phi_{2, 3}+(1-\mu) \Phi_{5, 1} \quad \text{mod $\mathfrak{a}_4$},\\
  \Phi_{32}^{\gamma} &=\gamma_1 \Phi_{8, 0}+\gamma_2 \Phi_{5, 2}+
                       (1-\gamma_1-\gamma_2) \Phi_{2, 4} \quad \text{mod $\mathfrak{a}_4$},\\
  \Phi_{44} &=\Phi_{11, 0} \quad \text{mod $\mathfrak{a}_5$},
\endaligned\right.\eqno{(11.6)}$$
where $\Phi_{12}$, $\Phi_{26}$, $\Phi_{32}$ and $\Phi_{44}$ are invariant
homogeneous polynomials of degree $12$, $26$, $32$ and $44$, respectively.
The ideals are given by
$$\left\{\aligned
  \mathfrak{a}_1 &=(\Phi_4, \Phi_8),\\
  \mathfrak{a}_4 &=(\Phi_4, \Phi_8, \Phi_{10}, \Phi_{14}, \Phi_{4, 0},
                    \Phi_{1, 2}, \Phi_{1, 3}, \Phi_{4, 1}),\\
  \mathfrak{a}_5 &=(\Phi_4, \Phi_8, \Phi_{10}, \Phi_{14}, \Phi_{4, 0},
                    \Phi_{1, 2}, \Phi_{1, 3}, \Phi_{4, 1}, \Phi_{1, 5},
                    \Phi_{4, 3}, \Phi_{7, 1}),
\endaligned\right.\eqno{(11.7)}$$
and the parameter space $\{ (\lambda, \mu, \gamma) \} \cong \mathbb{C}^4$.
They form variations of $Q_{18}$ and $E_{20}$-singularity structures over
the modular curve $X$.}

{\it Proof}. Similar as in the proof of Theorem 10.1, the ideals
$\mathfrak{a}_1$, $\mathfrak{a}_4$, and $\mathfrak{a}_5$ are zero
ideals over the modular curve $X$, and we have proved
$$\Phi_{12}^{\lambda}=\Delta(z), \quad
  \Phi_{26}^{\mu}=\eta(z)^8 \Delta(z) E_4(z) E_6(z)\eqno{(11.8)}$$
over the modular curve $X$. For degree $d=32$, the invariant homogeneous
polynomials $\Phi_{8, 0}$, $\Phi_{5, 2}$ and $\Phi_{2, 4}$ form a
three-dimensional complex vector space, and
$\Phi_{8, 0}=\Phi_{5, 2}=\Phi_{2, 4}=\eta(z)^8 \Delta(z)^2 E_4(z)$
over the modular curve $X$ (after normalization). Hence,
$$\Phi_{32}^{\gamma}=\eta(z)^8 \Delta(z)^2 E_4(z)\eqno{(11.9)}$$
over the modular curve $X$.

  Combining with (11.8) and (11.9), a straightforward calculation
shows that the equations of $Q_{18}$ and $E_{20}$-singularities
can be parametrized by the invariant polynomials $\Phi_{12}$,
$\Phi_{26}$, $\Phi_{32}$ and $\Phi_{44}$:
$$\Phi_{32}^3-\Phi_{12} \cdot \left(\frac{\Phi_{12}^4 \Phi_{26}}
 {\Phi_{32}}\right)^2-1728 \Phi_{12}^8=0,$$
$$\Phi_{44}^3-\Phi_{12}^4 \cdot \left(\frac{\Phi_{12}^5 \Phi_{26}}
 {\Phi_{44}}\right)^2-1728 \Phi_{12}^{11}=0$$
over the modular curve $X$.

\noindent
$\qquad \qquad \qquad \qquad \qquad \qquad \qquad \qquad \qquad
 \qquad \qquad \qquad \qquad \qquad \qquad \qquad \boxed{}$

\textbf{Theorem 11.3.} (Variations of $Q_{18}$ and $E_{20}$-singularity
structures over the modular curve $X$: geometric version) {\it There
are two morphisms from the cone $C_Y$ over $Y$ to the $Q_{18}$ and
$E_{20}$-singularities:
$$f_{18}: C_Y/G \rightarrow \text{Spec} \left(\mathbb{C}\{\Phi_{32}\}
          [\Phi_{12}, \Phi_{26}]/(\Phi_{32}^3-\Phi_{12} (\frac{\Phi_{12}^4
          \Phi_{26}}{\Phi_{32}})^2-1728 \Phi_{12}^8)\right)\eqno{(11.10)}$$
and
$$f_{20}: C_Y/G \rightarrow \text{Spec} \left(\mathbb{C}\{\Phi_{44}\}
          [\Phi_{12}, \Phi_{26}]/(\Phi_{44}^3-(\frac{\Phi_{12}^7
          \Phi_{26}}{\Phi_{44}})^2-1728 \Phi_{12}^{11})\right)\eqno{(11.11)}$$
over the modular curve $X$. In particular, there are infinitely many such
quadruples $(\Phi_{12}, \Phi_{26}, \Phi_{32}, \Phi_{44})$ $=$ $(\Phi_{12}^{\lambda},
\Phi_{26}^{\mu}, \Phi_{32}^{\gamma}, \Phi_{44})$ whose parameter space
$\{ (\lambda, \mu, \gamma) \}$ $\cong \mathbb{C}^4$. They form variations
of $Q_{18}$ and $E_{20}$-singularity structures over the modular curve $X$.}

{\it Proof}. By Theorem 9.1, Theorem 9.2, Theorem 9.3 and Theorem 11.2, we
have the following two ring homomorphisms
$$\aligned
 &\mathbb{C}\{\Phi_{32}\}[\Phi_{12}, \Phi_{26}]/(\Phi_{32}^3-\Phi_{12}
  (\frac{\Phi_{12}^4 \Phi_{26}}{\Phi_{32}})^2-1728 \Phi_{12}^8)\\
 &\rightarrow \left(\mathbb{C}[z_1, z_2, z_3, z_4, z_5, z_6]/I(Y)
  \right)^{\text{SL}(2, 13)}
\endaligned$$
$$\aligned
 &\mathbb{C}\{\Phi_{44}\}[\Phi_{12}, \Phi_{26}]/(\Phi_{44}^3-
  (\frac{\Phi_{12}^7 \Phi_{26}}{\Phi_{44}})^2-1728 \Phi_{12}^{11}))\\
 &\rightarrow \left(\mathbb{C}[z_1, z_2, z_3, z_4, z_5, z_6]/I(Y)
  \right)^{\text{SL}(2, 13)}
\endaligned$$
over the modular curve $X$, which induce two morphisms of schemes
$$f_{18}: C_Y/G \rightarrow \text{Spec} \left(\mathbb{C}\{\Phi_{32}\}
          [\Phi_{12}, \Phi_{26}]/(\Phi_{32}^3-\Phi_{12} (\frac{\Phi_{12}^4
          \Phi_{26}}{\Phi_{32}})^2-1728 \Phi_{12}^8)\right)$$
and
$$f_{20}: C_Y/G \rightarrow \text{Spec} \left(\mathbb{C}\{\Phi_{44}\}
          [\Phi_{12}, \Phi_{26}]/(\Phi_{44}^3-(\frac{\Phi_{12}^7
          \Phi_{26}}{\Phi_{44}})^2-1728 \Phi_{12}^{11})\right)$$
over $X$.

  At the standard cusp $z=i \infty$, $X(13)$ has local coordinate
$q^{\frac{1}{13}}$. Near a cusp, for this to extend at the cusp,
$x_i(z)$ has to be replaced by $q^{-\frac{2}{39}} x_i(z)$
$(=q^{-\frac{1}{104}} a_i(z))$ for $1 \leq i \leq 6$. Then near
the cusp $z=i \infty$, we have the following:

\noindent (1) For $\Phi_{12}(x_1(z), \ldots, x_6(z))$, the
leading term is given by
$$q \cdot (q^{-\frac{2}{39}})^{12}=q^{\frac{5}{13}}.$$

\noindent (2) For $\Phi_{26}(x_1(z), \ldots, x_6(z))$, the
leading term is given by
$$q^{\frac{4}{3}} \cdot (q^{-\frac{2}{39}})^{26}=1.$$

\noindent (3) For $\Phi_{32}(x_1(z), \ldots, x_6(z))$, the
leading term is given by
$$q^{\frac{7}{3}} \cdot (q^{-\frac{2}{39}})^{32}=q^{\frac{9}{13}}.$$

\noindent (4) For $\Phi_{44}(x_1(z), \ldots, x_6(z))$, the
leading term is given by
$$q^{\frac{10}{3}} \cdot (q^{-\frac{2}{39}})^{44}=q^{\frac{14}{13}}.$$

  Hence, this shows that $\Phi_{26}$ gives functions not vanishing
on the rays of the cusp $z=i \infty$, while $\Phi_{12}$, $\Phi_{32}$
and $\Phi_{44}$ vanish exactly on the rays of the cusp $z=i \infty$.
Moreover, the triple $(\Phi_{12}, \Phi_{26}, \Phi_{32})$ of invariant
polynomials defines a map from the cone $C_Y/G$ over $X(13)$ to the
equation $x^8+y^3+x z^2=0$. On the other hand, the triple $(\Phi_{12},
\Phi_{26}, \Phi_{44})$ defines a map from the cone $C_Y/G$ over $X(13)$
to the equation $x^{11}+y^3+z^2=0$. They map to $(0, 1, 0)$ (after
normalization) the generatrix of this cone corresponding to the cusps.
Hence, they give relations between the link at $O$ for the cone $C_Y/G$
over $X(13)$ and the link at $O$ for $Q_{18}$ and $E_{20}$-singularities,
which are as follows:

  The link at $O$ for the $Q_{18}$-singularity by the triple $(\Phi_{12},
\Phi_{26}, \Phi_{32})$:
$$\left\{\aligned
  &|\Phi_{12}|^2+|\Phi_{26}|^2+|\Phi_{32}|^2=1,\\
  &\Phi_{32}^3-\Phi_{12} \left(\frac{\Phi_{12}^4 \Phi_{26}}{\Phi_{32}}\right)^2
   -1728 \Phi_{12}^8=0.
\endaligned\right.$$

  The link at $O$ for the $E_{20}$-singularity by the trile $(\Phi_{12},
\Phi_{26}, \Phi_{44})$:
$$\left\{\aligned
  &|\Phi_{12}|^2+|\Phi_{26}|^2+|\Phi_{44}|^2=1,\\
  &\Phi_{44}^3-\left(\frac{\Phi_{12}^7 \Phi_{26}}{\Phi_{44}}\right)^2
   -1728 \Phi_{12}^{11}=0.
\endaligned\right.$$
This completes the proof of Theorem 11.3.

\noindent
$\qquad \qquad \qquad \qquad \qquad \qquad \qquad \qquad \qquad
 \qquad \qquad \qquad \qquad \qquad \qquad \qquad \boxed{}$

\vskip 2.0 cm

\noindent{Department of Mathematics, Peking University}

\noindent{Beijing 100871, P. R. China}

\noindent{\it E-mail address}: yanglei$@$math.pku.edu.cn

\end{document}